\documentclass[12pt]{report}
\usepackage{german}
\usepackage{makeidx}
\usepackage{amssymb}
\usepackage{stmaryrd}
\usepackage{wasysym}
\usepackage{fontenc}
\makeindex
\author{Andreas Johann Raab\\
Luisenstrasse 60, 80798 M"unchen,\\
Federal Republic of Germany\\
E-mail:andreas@andreasjohannraab.de}
\title{Konzepte der abstrakten Ergodentheorie\\
Zweiter Teil:\\
Sensitive Cantor-Systeme}
\begin{document}
\maketitle
\newpage
{\Large {\bf The Subject}}\newline
\newline 
In for a penny, in for a pound:
Our conception of a generalized ergodic theory exceeded the
generality of general topology in the first part \cite{rabe} and it shall exceed general topology in the second part also.
The theorems $3.3^{1}$ and $3.8^{1}$ of first part of the conceptions of abstract ergodic theory gave 
some emphasis to 
certain forms of generalized attractors
within the pluralism of terms of 
generalized attractors, which we have found in the second chapter of the first part: At the beginning of the second chapter 
we constructed the set $\mbox{{\bf @}}(\xi,\mathcal{A})$ of attractors according to $(2.1)^{1}$ and $(2.2)^{1}$
for a given set $\mathcal{A}$ of subsets of the phase-space $Y=\bigcup \mathcal{A}$ and
for the evolution
\begin{displaymath}
\begin{array}{c}
\xi:Y\times \mathbb{R}\to Y\ ,\\
(x,t)\mapsto \xi(x,t):=\xi^{t}(x)\ ,
\end{array}
\end{displaymath}
which is determined by a given one-dimensional bundle of flows $\{\xi^{t}\}_{t\in\mathbb{R}}$:
Those attractors of the set $\mbox{{\bf @}}(\xi,\overline{\mathcal{A}})$
were involved in the theorems $3.3^{1}$, where
$$\overline{\mathcal{A}}:=\{\mathbf{cl}_{\mathcal{A}}(X):X\subset Y\}$$
is the set of sets, which are closed for the
operator $\mathbf{cl}_{\mathcal{A}}$ satisfying
$$\mathbf{cl}_{\mathcal{A}}(Z):=\bigcap\{X\in\mathcal{A}^{c}\setminus\{\emptyset\}:X\supset Z\}$$
for any subset $Z\subset Y$.
\newline
Now, in this second part of the generalized ergodic theory, 
we shall investigate the constitution
of the term of sensitivity of a given flow. We 
invite the reader to become acquainted with some 
alleged aporia in the first chapter: 
We try to feel the problem arising, when we try to
extend the term of metric sensitivity of an attractor
over the grade of generality of metric spaces. 
But, is metric sensitivity not enough?
We also explain,
how adapted metrization of phase-spaces of 
sensitive flows can be confronted with huge 
diffculties.\newline 
In the second chapter we first give some contrivance, when we discuss intinsic sensivity of single trajectories.
When we consider the possibility of generalizing sensitivity we come across the idea, that 
sensitive attractors are homogenenous.
In the third chapter we find, how to handle the  
problem of generalizing sensitivity. Generalizing sensitivity 
turns out to be a supplement of the main
result of the first part of the generalized ergodic theory, which can be identified 
with the assertion, that the set of subsets
$$[[\xi]]_{\mathcal{A}}:=\{\mathbf{cl}_{\mathcal{A}}(\xi(x,\mathbb{R})):x\in Y\}$$
is a partition of the phase-space $Y$,
if 
for any real number $t\in\mathbb{R}$
the implication
$$P\in \mathcal{A}\setminus\{\emptyset\}\Rightarrow \exists\ Q\in \mathcal{A}\setminus\{\emptyset\}:Q\subset(\xi^{t})^{-1}(P)$$
is true. So that $\xi^{t}$ is in the set of Cantor-continuous functions $\mathcal{C}_{+}(\mathcal{A})$ fulfilling 
the commutation
\begin{equation}\label{cinsym}
[\mathbf{cl}_{\mathcal{A}},\xi^{t}]=\emptyset\ ,
\end{equation}
where $[\mathbf{cl}_{\mathcal{A}},\xi^{t}]$ is the 
commutator we introduced in the first part of the generalized ergodic theory. Thus
we articulated the Cantor-continuosity of flows as a form of 
cinematic symmetry, which fixes the fact, that $[[\xi]]_{\mathcal{A}}$ forms a partition of the phase-space $Y$.
\newline
Our supplementation shall concern generalized sensitivity, 
which is based on the set of subsets $\mathcal{A}$ of the phase-space $Y=\bigcup\mathcal{A}$, since
this general form of sensitivity is our principal theme. It turns out, that
the violation of the cinematic symmetry  (\ref{cinsym}) can be identified
with the sought general form of sensitivity. 
If this general form of sensitivity appears within attractors, it
defines general chaos. Thus there
exist sources of chaos. If the special form of general ultrasensitivity is given, the generation 
of general ultrachaos can be located in one certain point of the phase-space.
\newline
As chaos is excluded by Cantor-continuosity we search for some other condition concerning the temporal continuousity of the flow,
which does not exclude sensitivity.
Whereat the temporal continuousity of the bundle of flows $\{\xi^{t}\}_{t\in\mathbb{R}}$ is
a generalized form of
continuousity of the functions $\xi(x,\mbox{id})$ for $x\in  Y$, which is referable to the set of subsets $\mathcal{A}$.
In the last chapter
we
have to initiate the general form of temporal continuousity, which is converse comanency. Converse comanency does exclude neither
general sensitivity nor ultrasensitivity. This general form of continuosity 
arranges a certain regularity within ultrachaos, which is generated by ultrasensitivity.  
\newline

\newpage
\tableofcontents

\chapter{Die Schwierigkeit der Abstraktion des Sensitivit"atsbegriffes}
{\small Versetzen wir uns in die Anf"ange der 
Diskussion der Sensitivit"at zur"uck!
Vor dem Moment des breiteren Aufbruches erregt Stephen Smale\index{Smale, Stephen} mit der Pr"asentation seiner Umst"ulpung 
der Sph"are im Jahre 1958
Aufsehen, die er kurz darauf in der Arbeit \cite{maal}
vorstellt.
Der Moment des breiteren Aufbruches ist in der Arbeit von David Ruelle und Floris Takens \cite{rute}
oder in der Arbeit M. H$\acute{e}$nons u. C. Heiles' \cite{hehi}
dokumentiert,
der 
sich in den weiteren
Arbeiten Smales \cite{male}, H$\acute{e}$nons \cite{hehe},  
R"osslers \cite{ross}, Lozis \cite{lozi}, Ruelles \cite{rule} oder \cite{ecke}
und anderer fortsetzt, welche zumeist
ber"uhmt gewordene seltsame Attraktoren vorstellen, die nach
deren jeweiligen Entdeckern benannt wurden. Nach David Ruelle ist unseres 
Wissens kein seltsamer Attraktor benannt, indess hat David Ruelle\index{Ruelle, David} den Begriff
des seltsamen Attraktors gepr"agt:\newline
Wir bestimmen als ein wesentliches Motiv der Konstitution des 
Sensitivit"atsbegriffes das {\em numerische Trauma}, dass sich jede noch so scharfe Lokalisierung
eines Zustandes innerhalb mancher Attraktoren aufl"ost. Die \glqq noch so scharfe Lokalisierung\grqq
ist dabei wesentlich
{\em quantitativ} gedacht. Sie ist also von vorne herein eine Lokalisierung, die kaum anders als auf der Grundlage
einer generalisierten
Metrisierung des Zustandsraumes verfassbar ist.
Der Sensitivit"atsbegriff soll gerade dieses 
sensitive Entwicklungsverhalten eines in einer noch so \glqq {\em kleinen} Umgebung\grqq lokalisierten Phasenvolumens objektivieren; und exakt jene
Attraktoren, innerhalb derer dieses Entwicklungsverhalten sich 
aufl"osender, noch so scharfer Lokalisierung auftritt, sollen dann die 
sensitiven Attraktoren sein.\newline
Wenn der Sensitivit"atsbegriff dieses
numerische Trauma formuliert, so besteht dessen numerischer Charakter nur dann,
wenn eine quantitative Lokalisierung m"oglich ist,
wenn also eine Distanzfunktion auf dem Zustandsraum angelegt ist.
Daher br"auchte es uns nicht zu verwundern, wenn sich abzeichnete, dass der Sensitivit"atsbegriff wesentlich
auf Distanzfunktionen basiert.
Das ist es, was uns nicht gerade verheisst, dass der metrische Sensitivit"atsbegriff
weiter generalisierbar ist, als es 
der Begriff der Metrik\index{Metrik} 
ist. 
Dabei ist der Sensitivit"atsbegriff,
wie wir bald sehen,
mit dem Begriff der Metrik kogeneralisierbar.\index{Kogeneralisierbarkeit des Sensitivit\"atsbegriffes}
Das Streben nach einer Fassung des Sensitivit"atsbegriffes f"ur beliebige 
Smale-Systeme oder sogar f"ur beliebige 
Cantor-Systeme, das
Verlangen nach einer Fassung des Sensitivit"atsbegriffes,
die "uber 
den Allgemeinheitsgrad
generalisierter 
Metriken 
hinaus geht, mag als fernab vom Schuss erscheinen.
Wir k"onnen doch 
mit dem Allgemeinheitsgrad des Begriffes metrischer Sensitivit"at 
zufrieden sein, weil wir mit demselben alle die Sensitivit"at betreffenden Fragen formulieren 
k"onnen, die beispielsweise
mit der in \cite{oran} dargestellten Ergodentheorie der klassischen Mechanik
verbunden sind, welche
in normierten, mithin in metrischen Phasenr"aumen darstellbar ist.
Indess legt uns unser Thema darauf fest, die 
Generalisierbarkeit des Sensitivit"atsbegriffes zu untersuchen und 
dessen Generalisierung
so weit wie nur
m"oglich zu treiben. Unsere erste Frage ist daher,
ob
die Grenze der Generalisierbarkeit des Sensitivit"atsbegriffes
nicht innerhalb, aber auch nicht ausserhalb derjeniger
Entsprechung verl"auft, welche die
Grenze der Generalisierbarkeit des Begriffes der Metrik 
zieht.
Kann es sein, dass der Sensitivit"atsbegriff -- im Unterschied zum Begriff des Attraktors --
die volle topologische und die sogar dar"uber hinausgehende Allgemeinheit nicht erreicht?}
\section{Metrische Sensitivit"at}
{\bf 1.1.1 Sensitive und seltsame Attraktoren sind zu unterscheiden}\newline\newline
Ganz an den Eingang des ersten Teiles der Konzepte der 
abstrakten Ergodentheorie \cite{rabe} stellten wir heraus,
dass die Ergodentheorie nicht \glqq in der Mathematik zur Welt gekommen\grqq sei und wir 
erl"auterten anschliessend deren Herkunft aus dem 
deterministisch-indeterministischen Dualismus.\index{deterministisch-indeterministischer Dualismus} 
Wir schlossen diese Erl"auterung mit dem Versprechen, dass wir 
etwas nicht vergessen wollen, 
n"amlich
die 
dualistische Herkunft der Ergodentheorie, die jenseits der Mathematik im physikalischen, im {\em Zeitlichen}, ja, 
im 
erlebbaren
deterministisch-indeterministischen Dualismus liegt; woraufhin wir, ganz innerhalb der Mathematik verbleibend, 
den Attraktorbegriff untersuchten und verallgemeinerten; und zwar auf eine Weise, 
der vorgehalten werden mag, dass sie im Zuge der Abstraktion 
die dualistische Herkunft der Ergodentheorie aus den Augen 
verliert. Ja, in der Abhandlung des elementaren Quasiergodensatzes \cite{raab}, da
sei die Bindung zur Herkunft der Ergodentheorie aus der
Spannung des deterministisch-indeterministischen Dualismus noch erkennbar, nicht mehr aber 
in der abstrakten Ergodentheorie \cite{rabe}. 
Es stimmt, es kann der Eindruck entstehen, dass wir unser Versprechen gaben, die Herkunft
der Ergodentheorie aus dem deterministisch-indeterministischen Dualismus nicht
vergessen zu wollen -- und unser Versprechen anschliessend brachen.   
\newline
Wer diesen Eindruck hat, beobachte jedoch auch dies:
Dasjenige Kontinuum,
das uns als der uns bislang einzig objektivierte Schatten
der physikalischen Zeit zuhanden ist, wenn wir die 
physikalische Zeit aus ihrem relativistischen Zusammenhang herausbrechen, ist der Zahlenstrahl $\mathbb{R}$.
Dessen logische Rekonstruktion ist errungen und ein Produkt nicht nur Cantors Obsession.
Wir weisen darauf hin, dass wir selbst in dem 
als recht abstrakt und weitab vom Schuss in irrelevante Entr"ucktheit abschiebbaren 
ersten Teil der Konzepte der 
abstrakten Ergodentheorie \cite{rabe} 
nicht etwa nur die Attraktoren 
einzelner Autobolismen untersuchten. Stattdessen liessen
wir
ein besonderes Gewicht 
der Untersuchung 
des Phasenflusses $\{\xi^{t}\}_{t\in\mathbb{R}}$ 
als einer Familie 
von Autobolismen zukommen,
welche --
der Zahlenstrahl indiziert. 
\newline
Dass wir den kontinuumsbezogenen 
Phasenfluss $\{\xi^{t}\}_{t\in\mathbb{R}}=\{\xi(\mbox{id},t)\}_{t\in\mathbb{R}}$
und die durch denselben
festgelegten Evolutionen hier wie
im ersten Teil der Konzepte der 
abstrakten Ergodentheorie bevorzugen,
ist der ferne Widerhall unseres Interesses an der physikalischen Zeit, deren 
Spannung uns als die Spannung des deterministisch-indeterministischen Dualismus erscheint.
Wo wir die Sensitivit"at von Fl"ussen ber"uhren, die
wir nun im zweiten Teil als Hauptthema 
er"ortern, ist
die Diskussion des kontinuumsbezogenen 
Phasenflusses $\{\xi^{t}\}_{t\in\mathbb{R}}$
ist mit folgender Missverst"andlichkeit
verbunden: Es ist die 
gem"ass der Definition 3.3 verfasste allgemeine Sensitivit"at
einzelner Fl"usse $\xi^{t}$, welche 
f"ur 
die allgemeine
Sensitivit"at der Flussfunktion $\xi=\xi^{\mathbf{P}_{2}}\circ\mathbf{P}_{1}$  
notwendig und hinreichend ist, die 
nur behauptet, dass es einen
Fluss $\xi^{t}$ gibt, der allgemein sensitiv ist.
Wie wir gleich sehen werden, ist 
das mit der metrischen Sensitivit"at anders.
Die metrische Sensitivit"at
einzelner Fl"usse $\xi^{t}$ gem"ass (\ref{seltb}) ist zwar hinreichend f"ur 
die metrische
Sensitivit"at einer metrischen Flussfunktion $\xi=\xi^{\mathbf{P}_{2}}\circ\mathbf{P}_{1}$
gem"ass (\ref{seltc}). Die Notwendigkeit metrischer Sensitivit"at
einzelner Fl"usse $\xi^{t}$ gem"ass (\ref{seltb}) daf"ur,
dass die metrische
Sensitivit"at der metrischen Flussfunktion $\xi$ vorliegt, ist aber fraglich.\newline
Die Attraktoren einzelner Fl"usse $\xi^{t}$
k"onnen nur dann seltsam sein, wenn $\xi^{t}$ allgemein sensitiv ist. Die allgemein sensitiven Attraktoren der Flussfunktion $\xi=\xi^{\mathbf{P}_{2}}\circ\mathbf{P}_{1}$
enthalten diese gegebenenfalls seltsamen Attraktoren einzelner sensitiver Fl"usse $\xi^{t}$, deren fraktale Beschaffenheit im 
jeweiligen
sensitiven Attraktor der Flussfunktion auf die Weise aufgehoben sei kann,
wie der geschw"arzte Schaustellter vor schwarzem Hintergrund verschwindet.
\newline
Es sei $h$ ein einzelner Autobolismus eines metrischen Raumes $(X,d)$, sodass 
dessen Wertemenge
$\mathbf{P}_{2}h=h^{-1}(X)=X= \mathbf{P}_{1}h$
mit dessen Definitionsmenge "ubereinstimmt.
$h$ sei bez"uglich der von der Metrik $d$ induzierten Topologie $\mathbf{T}(d)$ stetig, wobei
f"ur alle Teilmengen $Y\subset X$ deren jeweiliger Durchmesser 
$d(Y):=\sup\{d(a,b):a,b\in Y\}$ sei.
Wenn f"ur $x\in X$ und eine Umgebung ${\rm U}$ des Umgebungssystemes $\mathbf{T}(d)_{\{x\}}$ des Punktes $x$ die Aussage
\begin{equation}\label{selta}
\exists\ \varepsilon(x,{\rm U})\in\mathbb{R}^{+}: \forall\ y\in{\rm U}\
\exists\ n\in\mathbb{N}:\ 
d(h^{n}(x),h^{n}(y))\geq\varepsilon(x,{\rm U})
\end{equation}
wahr
ist, ist der Limes Superior f"ur alle $y\in{\rm U}\in\mathbf{T}(d)_{\{x\}}$
gem"ass
$$\lim_{n\to \infty}\sup d(h^{n}(x),h^{n}(y))\geq\varepsilon(x,{\rm U})$$
von unten beschr"ankt. Der G"ultigkeit von (\ref{selta}) steht 
die Stetigkeit des Autobolismus $h$ nicht im Wege: Denken wir an eine stetige,
lokale Expansion $q$ im Ursprung der zentrierten abgeschlossenen Einheitskugel, die
aber an dessen Berandung durch Kontraktion kompensiert ist.   
Indess, es kann dabei sein, dass bei der Kontraktion der Umgebungen um $x$, d.h, 
f"ur $d({\rm U})\to 0$ die Zahl $\varepsilon(x,{\rm U})$ ebenfalls verschwindet. Es kann 
(\ref{selta}) gelten 
und dabei
\begin{equation}\label{aselt}
\lim_{d({\rm U})\to 0}\sup\lim_{n\to \infty}\sup d(h^{n}({\rm U}))=0
\end{equation}
sein.\footnote{Wobei $\lim_{n\to \infty}\sup d(h^{n}({\rm U}))$ hier
als Menge jeweiliger Limites Superiores zu lesen ist: So, wie
mit dem Grenzwert $\lim_{x\to\infty}1/x=0$ eine Aussage "uber 
die Menge der Nullfolgen gemacht ist, behauptet (\ref{aselt}), dass f"ur alle Folgen von Umgebungen 
$\{{\rm U}_{j}\}_{j\in\mathbb{N}}$, deren Glieder 
Elemente des Umgebungssystemes $\mathbf{T}(d)_{\{x\}}$ sind,
die Implikation
$$\lim_{j\to \infty}d({\rm U}_{j})=0\Rightarrow
\lim_{j\to \infty}\lim_{n\to \infty}\sup d(h^{n}({\rm U}_{j}))=0$$
wahr ist.}
Die G"ultigkeit 
der Aussage (\ref{selta}) schliesst demnach weder aus noch impliziert sie, dass
der Fall vorliegt, exakt welchen wir unter der 
punktweisen Sensitivit"at eines Autobolismus $f$ eines metrischen Raumes $(X,d)$
im Punkt $x$ verstehen:
Genau dann n"amlich, wenn f"ur einen Punkt $x$ die
Aussage 
\begin{equation}\label{seltb}
\begin{array}{c}
\exists\ \varepsilon(x)\in\mathbb{R}^{+}: \forall\ ({\rm U},n)\in\mathbf{T}(d)_{\{x\}}\times\mathbb{N}\\ 
\exists\ (y,b)\in {\rm U}\times\{m\in\mathbb{N}:m>n\}:\\ 
d(f^{b}(x),f^{b}(y))\geq\varepsilon(x)
\end{array}
\end{equation}
gilt, nennen wir den Autobolismus $f$ im Punkt $x$ metrisch sensitiv.\index{punktweise metrische Sensitivit\"at eines metrischen Autobolismus}
Dass die punktweise metrische Sensitivit"at eines metrischen Autobolismus nicht zu dessen
beidseitiger metrischer Sensitivit"at "aquivalent ist, f"uhrt uns dies rasch vor Augen,
dass 
die gerade bedachte, stetige,
lokale Expansion $q$ im Ursprung der Ebene, die
aber an der Berandung der zentrierten abgeschlossenen Einheitskugel kontrahiert,   
metrisch sensitiv ist. $q$ ist aber nicht beidseitig metrisch sensitiv.
Genau dann, wenn f"ur einen Autobolismus $f$ und einen Punkt $x\in\mathbf{P}_{1}f$
\begin{equation}\label{thseltb}
\begin{array}{c}
\exists\ \varepsilon(x)\in\mathbb{R}^{+}: \forall\ ({\rm U},n)\in\mathbf{T}(d)_{\{x\}}\times\mathbb{N}\\ 
\exists\ (y,b)\in {\rm U}\times(\mathbb{Z}\setminus\{z\in\mathbb{Z}:|z|<n\}:\\ 
d(f^{b}(x),f^{b}(y))\geq\varepsilon(x)
\end{array}
\end{equation}
gilt, nennen wir den Autobolismus $f$ im Punkt $x$ beidseitig metrisch sensitiv.\index{beidseitige punktweise metrische Sensitivit\"at eines metrischen Autobolismus}
Es gibt dann 
eine positive Konstante $\Delta_{f}(x)\in\mathbb{R}^{+}$, f"ur die f"ur alle ${\rm U}\in\mathbf{T}(d)_{\{x\}}$ die
Ungleichung
\begin{equation}
d(f^{n}({\rm U}))\geq\Delta_{f}(x)
\end{equation}
gilt, wenn $n$ hinreichend gross gew"ahlt ist und wenn $f$ metrisch sensitiv ist.
Die punktweise Sensitivit"at 
des Autobolismus $f$
ist ein notwendiges Kriterium daf"ur, dass
er seltsame Attraktoren hat.
Die globale Version der punktweisen Sensitivit"at des Autobolismus $f$, dass dieser in allen Punkten
sensitiv ist, gilt\index{Sensitivit\"at eines metrischen Autobolismus}
"ublicherweise nicht
als dessen metrische Sensitivit"at. 
Die metrische
Sensitivit"at des Autobolismus $f$ liegt vor, wenn 
in seiner Definitionsmenge punktweise Sensitivit"at auftritt.
Es liegt daher diktiertermassen nahe, die punktweise
Sensitivit"at eines Phasenflusses $\{\xi^{t}\}_{t\in\mathbb{R}}=\{\xi(\mbox{id},t)\}_{t\in\mathbb{R}}$
und der durch diesen bestimmten Flussfunktion $\xi$ 
im Zustand $x\in\mathbf{P}_{2}\xi$ 
als die Gegebenheit festzulegen, dass es eine positive reelle Zahl $\varepsilon(x)$ gibt,
f"ur welche 
die Aussage
\begin{equation}\label{seltc}
\begin{array}{c}
\forall\ ({\rm U},t_{\star})\in\mathbf{T}(d)_{\{x\}}\times\mathbb{R}^{+} \ \exists\ (y,t)\in {\rm U}\times]-\infty,t_{\star}[\cup ]t_{\star},\infty[]:\\ 
d(\xi^{t}(x),\xi^{t}(y))\geq\varepsilon(x)
\end{array}
\end{equation}
wahr ist. Folgen wir diesem nat"urlichen Diktat, so riskieren wir, dass Konfusion entsteht:
Und zwar, weil 
sensitive Autobolismen $f$ Attraktoren haben k"onnen, die seltsam sind.
Die punktweise
Sensitivit"at einer Flussfunktion $\xi$ im Zustand $x$ behauptet dabei aber nicht unmittelbar,
dass auch nur einer der Autobolismen $\xi^{t}$ der Phasenflussgruppe $(\{\xi^{t}:t\in\mathbb{R}\},\circ)$
in $x$
punktweise
Sensitivit"at zeige. Wenn wir f"ur eine Flussfunktion 
Attraktoren formulieren, so folgen wir bei der Benennung 
derjeniger Attraktoren $\chi$ einer Flussfunktion, f"ur welche in einem  
Zustand $x\in\chi$ punktweise Sensitivit"at
gem"ass (\ref{seltc}) vorliegt,
der nat"urlichen Benennungslogik, wenn wir exakt solche Attraktoren, die in diesem Sinn
sensitives Verhalten zeigen, als 
die sensitiven Attraktoren der jeweiligen
{\em Flussfunktion} bezeichnen.
"Uber Sachverhalte streitet hier niemand, "uber Benennungen d"urften wir gelassen 
verschiedener 
Auffassung sein, wobei wir mit der hier geltenden Bezeichnungsweise 
nicht nur auf der Seite der nat"urlichen Benennungslogik sind. Um die erscheinungsbildliche Differenz der
sensitiven Attraktoren jeweiliger Flussfunktionen und der 
sensitiven Attraktoren
einzelner Autobolismen
hervorzuheben, verweisen
wir auf den Artikel David Ruelles \cite{rufe}, dessen Hauptanliegen
es ist, m"oglichst deutlich zu machen, was
ein seltsamer Attraktor ist, n"amlich ein sensitiver Attraktor 
eines einzelnen 
Autobolismus.\newline 
Wir bleiben deshalb auf der Seite der nat"urlichen Benennungslogik, weil 
dieselbe in diesem Fall mit dem 
Permanenzprinzip\index{Permanenzprinzip} einhergeht:
Dem Sachverhalt, 
dass 
unsere Definition der Sensitivit"at einer Flussfunktion
die Sensitivit"at eines einzelnen Autobolismus 
umfasst, steht n"amlich einzig 
folgendes Hindernis im Wege: Im ersten Teil der Konzepte der
abstrakten Ergodentheorie 
legten wir in $(1.51)^{1}$ fest,
dass das, was wir als eine Flussfunktion $\Phi$ 
gelten lassen,
so beschaffen sein soll, dass
die Funktion $\Phi(x,\mbox{id})$
f"ur alle Zust"ande $x\in\mathbf{P}_{2}\Phi$
invertierbar sei.
Nennen wir exakt alle Funktionen $\Theta$, die
bis auf dieses Kriterium der Invertierbarkeit
der Funktion $\Theta(x,\mbox{id})$ f"ur alle Zust"ande $x\in\mathbf{P}_{2}\Theta$
alle Merkmale einer Flussfunktion haben,
Flussfunktionen im weiteren Sinn,\index{Flussfunktion im weiteren Sinn} so k"onnen 
wir genausogut auch f"ur jede 
Flussfunktion im weiteren Sinn
festlegen, 
dass $\Theta$ 
genau dann punktweise bzw. globale 
Sensitivit"at habe, wenn 
die f"ur $\Theta$ statt f"ur 
$\Phi$ formulierte Aussage
(\ref{seltc}) f"ur alle 
$x\in\mathbf{P}_{2}\Theta$ wahr ist.
Genauso k"onnen wir auch den Begriff des freien 
Attraktors einer Flussfunktion, den die gem"ass $(2.1)^{1}$
und $(2.2)^{1}$ verfassten Mengensysteme $\mbox{{\bf @}}(\xi,\mathcal{A})$
f"ur eine jeweilige Flussfunktion $\xi$ und eine Zustandsraum"uberdeckung $\mathcal{A}$  
objektivieren, auch f"ur 
Flussfunktionen im weiteren Sinn verfassen. Der 
verallgemeinerte insensitve Ergodensatz\index{verallgemeinerter insensitver Ergodensatz} $(3.3)^{1}$ gilt offenbar ebenfalls f"ur
Flussfunktionen im weiteren Sinn. 
Zwischen den Flussfunktionen und der Flussfunktionen im weiteren Sinn
liegen die Entwicklungen, die wir aus darstellungspraktischen Gr"unden eigens herausheben:
\newline\newline
{\bf Festlegung 1.1: Entwicklungen und finit-additive Entwicklungen}\newline
{\em Wir bezeichnen jede Flussfunktion $\xi$ im weiteren Sinn genau dann als eine Entwicklung,
\index{Entwicklung}
wenn f"ur alle $t\in\mathbb{R}$ die Inversionen 
der Autbolismen} $\xi^{t}=\xi(\mbox{id},t)$
{\em existieren.\newline 
Da $\xi^{t}$ f"ur alle $t\in\mathbb{R}$ existiert, ist dann das Paar
\begin{equation}\label{seccltc}
(\mathbb{R},+_{\xi})
\end{equation}
durch die Festlegung, dass f"ur alle $a,b\in \mathbb{R}$ 
\begin{equation}\label{scelctc}
\xi^{\ a+_{\xi}b}:=\xi^{a}\circ\xi^{b}
\end{equation}
gelte,
bestimmt. $(\mathbb{R},+_{\xi})$ ist eine Gruppe. \newline
Genau dann bezeichnen wir jede Entwicklung $\xi$ als finit-additiv,\index{finit-additive Entwicklung}
wenn 
f"ur alle $a\in \mathbb{R}$ die Invarianz 
\begin{equation}\label{scelctc}
\{X+_{\xi}\{a\}:X\in[\mathbb{R}]\}=[\mathbb{R}]
\end{equation}
der beschr"ankten Teilmengen des Zahlenstrahles gilt, 
wobei $[\mathbb{R}]\subset 2^{\mathbb{R}}$
all deren Menge 
bezeichne.
Exakt jede strukturierte Flussfunktion $(\xi,\mathcal{A})$ im weiteren Sinn nennen wir dabei 
eine strukturierte
bzw. eine strukturierte finit-additive Entwicklung,\index{strukturierte Entwicklung}
\index{strukturierte finit-additive Entwicklung} deren erste Komponente eine  
Entwicklung bzw. eine finit-additive Entwicklung ist.}\newline 
\newline 
Diese einfache Festlegung bedarf keines Komentares.
Es sei sogleich die Betrachtung $f$ des metrischen 
des Autobolismus $f$ fortgesetzt und
$\underline{f}:X\times\mathbb{R}\to X$ die Funktion, deren 
jeweilige Werte f"ur alle $(x,t)\in \mathbf{P}_{1}\underline{f}$
\begin{equation}\label{seltd}
\underline{f}(x,t):=f^{\mbox{int}(t)}(x)
\end{equation}
sind, wobei
$$\mbox{int}(a):=\max\{n\in\mathbb{N}:n\leq a\}$$
den Wert der rechtseitig stetigen Abrundung bezeichne,
dann ist
$\underline{f}$ eine Flussfunktion 
im weiteren Sinn, die genau in denjenigen
Punkten des Zustandsraumes $X$ 
sensitiv ist, in denen $f$ punktweise 
Sensitivit"at zeigt.
Der Begriff der punktweisen Sensitivit"at gem"ass (\ref{seltc}) 
f"ur Flussfunktionen im weiteren Sinn
umfasst also den Begriff der punktweisen Sensitivit"at
jeweiliger Autobolismen.
\newline
Die Funktion
$\underline{f}(x,\mbox{id})$ ist f"ur alle Zust"ande $x\in \mathbf{P}_{2}f$
exakt auf den Intervallen 
$[j,j+1[$ f"ur $j\in \mathbb{Z}$ konstant. In dem Sinn, dass
zu jeder Teilmenge des Zahlenstrahles $J\subset\mathbb{R}$, auf der die
Restriktion
$(\underline{f}(x,\mbox{id})|J)^{-1}$
eine Funktion
ist, genau ein Partitionselement $J^{+}\in \{\mathbb{Z}+\{q\}:q\in]0,1[\}$ existiert, f"ur das
$J^{+}\supset J$ gilt,
ist die Funktion
$\underline{f}(x,\mbox{id})$
h"ochstens 
auf den Mengen
des Mengensystemes
\begin{equation}\label{brum}
\{\mathbb{Z}+\{q\}:q\in[0,1[\}\in\mathbf{part}(\mathbb{R})
\end{equation}
invertierbar. Wir
k"onnen auch so sagen: Die maximalen Bereiche
der Invertierbarkeit der Funktion $\underline{f}(x,\mbox{id})$
sind in Elementen der Partition (\ref{brum}) enthalten.\newline
Wir sehen auch, dass es eine l"osbare Aufgabe sein
kann, die Phasenflussgruppen 
$(f^{\mathbb{Z}},\circ)$, die einzelne
Autobolismen generieren, in 
Phasenflussgruppen $(F^{\mathbb{R}},\circ)$
stetig einzubetten; in der Form also, dass
die zu $(F^{\mathbb{R}},\circ)$ "aquivalente 
Flussfunktion bez"uglich einer vorgegebenen
Zustandsraumtopologie stetig ist. Die Untersuchung dieser Fragestellung
gehen wir in dieser Abhandlung nicht an, ebensowenig, wie wir die
Herausforderung annehmen,
das Verh"altnis zu untersuchen, das zwischen den fraktalen\index{Fraktal} 
sensitiven Attraktoren $a\in \mbox{{\bf @}}(\underline{\xi^{t}},\mathcal{A})$ der einzelnen Fl"usse $\xi^{t}$
einer Phasenflussgruppe $(\{\xi^{t}:t\in\mathbb{R}\},\circ)$
und den Attraktoren $\chi\in \mbox{{\bf @}}(\xi,\mathcal{A})$ der 
Flussfunktionen $\xi$ besteht. Wie wir leicht verifizieren, k"onnen 
wir das im ersten Teil der Konzepte der 
abstrakten Ergodentheorie
festgelegte 
Mengensystem der Attraktoren $\mbox{{\bf @}}(\tilde{\xi},\tilde{\mathcal{A}})$
genausogut f"ur strukturierte Flussfunktionen im weiteren Sinn 
$(\tilde{\xi},\tilde{\mathcal{A}})$ wie f"ur strukturierte Flussfunktionen 
verfassen.\index{Attraktoren strukturierter Flussfunktionen im weiteren Sinn}
F"ur die fraktalen sensitiven Attraktoren $a\in \mbox{{\bf @}}(\underline{\xi^{t}},\mathcal{A})$
gilt offensichtlich 
$$\mathbf{card}(\mbox{{\bf @}}(\xi,\mathcal{A})_{a})=1\ .$$
Die fraktalen sensitiven Attraktoren einzelner Fl"usse sind in den sensitiven Attraktoren
sie einbettender Flussfunktionen 
enthalten. 
Indess,
wenn wir diese Herausforderung ann"ahmen,
w"urden wir unser Versprechen, die Herkunft
der Ergodentheorie aus dem deterministisch-indeterministischen Dualismus nicht
vergessen zu wollen, keineswegs brechen,
das eigentlich das Interesse 
an der integrativen Sichtweise auf den
Determinismus ist, die sich nicht nur f"ur den einzelnen
Autobolismus begeistert. Wir w"urden dann aber 
die fraktalen sensitiven Attraktoren einzelner
Gruppenelemente 
der Phasenflussgruppe 
einer korrespondierenden Flussfunktion
innerhalb eines sensitiven Attraktors
derselben erforschen und die fraktalen sensitiven Attraktoren
in einen integralen Kontext stellen.
\newline\newline
{\bf 1.1.2 Sensitivit"at und Komanenz bez"uglich nicht-negativer Funktionen}
\newline\newline
Exakt die 
(\ref{seltc}) gem"asse
Auffassung 
von Sensitivit"at 
liegt der Einf"uhrung des Aufl"osungsfeldes einer Wellenfunktion
in dem Trakat "uber den elementaren Quasiergodensatz \cite{raab} 
zugrunde, die ebenda lediglich nicht in der Allgemeinheit metrischer R"aume verfasst ist, sondern
die da
f"ur normierte R"aume formuliert ist: 
Das Aufl"osungsfeld einer Wellenfunktion 
$\Phi\in B^{A\times\mathbb{R}}$ wurde
in dem Fall eingef"uhrt, dass die Mengen $A\subset X$ und $B\subset Y$ in normierten R"aumen $(X,||\mbox{id}||_{X})$ und $(Y,||\mbox{id}||_{Y})$ liegen. 
Das Aufl"osungsfeld einer Wellenfunktion
ist die Abbildung
\begin{equation}\label{aalzbb}
\begin{array}{c}
\Delta_{\Phi}:A\to\mathbb{R}^{+}\cup\{-\infty\}\ ,\\
x\mapsto\Delta_{\Phi}(x):=\sup M_{\Phi}(x)\ ,
\end{array}
\end{equation}
welche die Sensitivit"at der Wellenfunktion $\Phi$ quantitativ bewertet;
f"ur jede reelle Wellenfunktion $\Phi$ formuliert dann die 
Menge
\begin{equation}\label{aalzb}
\begin{array}{c}
M_{\Phi}(x):=
\Bigl\{\delta\in\mathbb{R}^{+}:\ \forall\ \varepsilon\in\mathbb{R}^{+}
\ \exists\ y\in \{z\in A:||y-x||_{X}<\varepsilon\}\\ \forall t_{\star}\in\mathbb{R}^{+}\ \exists\ 
t\in ]-\infty,t_{\star}[\cup ]t_{\star},\infty[\\
||\Phi(x,t)-\Phi(y,t)||_{Y}\ >\ \delta\Bigr\}
\end{array}
\end{equation}
die Sensitivit"at der Wellenfunktion $\Phi$ in dem Zustand $x\in A$
folgendermassen: Genau dann, wenn die
Menge $M_{\Phi}(x)$ nicht leer ist, nennen wir die Wellenfunktion $\Xi$
im Zustand $x\in A$
sensitiv.\index{Sensitivit\"at einer reellen Wellenfunktion}
Genau dann, wenn die Wellenfunktion $\Xi$ im Zustand 
$x$ nicht sensitiv ist, ist ihr Aufl"osungsfeld in $x$ singul"ar,
n"amlich 
das mit $-\infty$ identifizierte Supremum der leeren Menge.
Die Sensitivit"at einer reellen Wellenfunktion f"uhren wir also punkt- oder zustandsweise
ein. 
Dieser Sensitivit"atsbegriff, der f"ur 
auf normierten R"aumen erkl"arte Wellenfunktionen mit in
einem normierten Raum liegender Wertemenge
formuliert ist, umfasst damit
auch eine der Sensitivit"atsformen einer jeweiligen Flussfunktion $\Xi$
in einem jeweiligen Zustand eines normierten Raumes; und gerade diese 
Sensitivit"atsform einer jeweiligen Flussfunktion binnen eines normierten Raumes 
bezeichnen wir als deren 
Sensitivit"at im jeweiligen Zustand.\index{Sensitivit\"at einer Flussfunktion}  
\newline
Wir k"onnen $M_{\Phi}$ als mengenwertige Funktion 
auffassen, die 
auf dem ersten kartesischen Faktor 
$\mathbf{P}_{1}\mathbf{P}_{1}\Phi=A$
der Definitionsmenge $\mathbf{P}_{1}\Phi=A\times\mathbb{R}$
der Wellenfunktion $\Phi$ 
definiert ist.
Dabei ist $M_{\Phi}(x)$ 
f"ur alle $x\in A$
als eine offene und zusammenh"angende Teilmenge des Zahlenstrahles 
konstruiert, sodass
\begin{equation}
M_{\Phi}=]0,\Delta_{\Phi}[
\end{equation}
ist. Denn, auch wenn $\Phi$ im Zustand $x\in A$ nicht sensitiv ist, ist
$$M_{\Phi}(x)=\emptyset=]0,-\infty[=]0,\Delta_{\Phi}[$$
eine offene und zusammenh"angende Teilmenge des Zahlenstrahles.
Eigentlich, so f"allt uns hierbei auf,
bezeichnen die Indizierungen der
Ausdr"ucke $M_{\Phi}(x)$ und $\Delta_{\Phi}$
insofern zu kurz greifend, als aus ihnen nicht hervorgeht,
auf welche beiden Normierungen $||\mbox{id}||_{X}$ und $||\mbox{id}||_{Y}$
Bezug genommen wird.
\newline
Wir sehen ferner, dass wir diesen Sensitivit"atsbegriff
einer Wellenfunktion $\Phi$ 
analog f"ur eine 
Wellenfunktion $\Xi\in Q^{P\times\mathbb{R}}$ verfassen 
und damit unmittelbar verallgemeinern  
k"onnen, wenn $(Q,d_{Q})$ und $(P,d_{P})$ metrische R"aume sind.
Wir brauchen dazu lediglich
die norminduzierten Distanzfunktionen
$$||\mathbf{P}_{1}-\mathbf{P}_{2}||_{X}\quad \mbox{bzw.}\quad ||\mathbf{P}_{1}-\mathbf{P}_{2}||_{Y}$$
durch Metriken $d_{Q}$ auf der Menge $Q$ bzw. $d_{P}$ auf der Menge $P$
zu ersetzen.\footnote{F"ur jede nat"urliche Zahl $n\in\mathbb{N}$ und f"ur alle $j\in\{1,2,\dots n\}$
sei 
f"ur jedes $n$-Tupel $a$ mit
$\mathbf{P}_{j}a$ dessen $j$-te Komponente bezeichnet; und f"ur jede Menge $A$
von lauter $n$-Tupeln bezeichnet $\mathbf{P}_{j}A$
die Menge aller deren $j$-ter Komponenten.} 
Uns ist ausserdem klar, dass wir jene norminduzierten Distanzfunktionen auch durch quasimetrische
Distanzfunktionen oder durch
Distanzfunktionen im weitesten Sinn ersetzen k"onnen.  
Sodann liegt nichts n"aher, als folgende Verallgemeinerung zu formulieren:
\newline
\newline
{\bf Definition 1.2:\newline Sensitivit"at bez"uglich zweier nicht-negativer Funktionen}
\index{Sensitivit\"at bez\"uglich zweier nicht-negativer Funktionen}
\newline
{\em F"ur jede Wellenfunktion $\Phi\in B^{A\times\mathbb{R}}$ und f"ur
jedes Paar $(\delta_{A},\delta_{B})$
zweier
auf den Mengen $A$ bzw. $B$ 
erkl"arter Funktionen $\delta_{A}\in [0,\infty[^{A\times A}$ und $\delta_{B}\in [0,\infty[^{B\times B}$
sei
\begin{equation}\label{aalzb}
\begin{array}{c}
M_{\Phi}^{\delta_{A},\delta_{B}}(x):=
\Bigl\{\delta\in\mathbb{R}^{+}:\ \forall\ \varepsilon\in\mathbb{R}^{+}
\ \exists\ y\in \{z\in A:\delta_{A}(y,x)<\varepsilon\}\\ \forall t_{\star}\in\mathbb{R}^{+}\ \exists\ 
t\in ]-\infty,t_{\star}[\cup ]t_{\star},\infty[:\\
\delta_{B}(\Phi(x,t),\Phi(y,t))\ >\ \delta\Bigr\}
\end{array}
\end{equation}
und 
\begin{equation}\label{aalzbb}
\begin{array}{c}
\Delta_{\Phi}^{\delta_{A},\delta_{B}}:A\to\mathbb{R}^{+}\cup\{-\infty\}\ ,\\
\quad\\
x\mapsto\Delta_{\Phi}^{\delta_{A},\delta_{B}}(x):=\sup M_{\Phi}^{\delta_{A},\delta_{B}}(x)
\end{array}
\end{equation}
das verallgemeinerte Aufl"osungsfeld der Wellenfunktion\index{Aufl\"osungsfeld einer Wellenfunktion} $\Phi$ bez"uglich der 
Funktionen $\delta_{A}$ und $\delta_{B}$.\index{verallgemeinertes Aufl\"osungsfeld einer Wellenfunktion}
Auch 
$M_{\Phi}^{\delta_{A},\delta_{B}}(x)$ ist f"ur alle $x\in A$
eine offene und zusammenh"angende Teilmenge des Zahlenstrahles,
sodass
die mengenwertige Abbildung
$$M_{\Phi}^{\delta_{A},\delta_{B}}=]0,\Delta_{\Phi}^{\delta_{A},\delta_{B}}[$$
existiert. Jeden Zustand $x\in A$ nennen wir genau dann einen sensitiven
Zustand der Wellenfunktion $\Phi$ bez"uglich der 
Funktionen $\delta_{A}$ und $\delta_{B}$,
wenn $M_{\Phi}^{\delta_{A},\delta_{B}}(x)$ nicht leer ist. Jede 
Wellenfunktion nennen wir genau dann bez"uglich $\delta_{A}$ und $\delta_{B}$ sensitiv, wenn sie 
bez"uglich dieser Distanzfunktionen sensitive Zust"ande hat.}
\newline
\newline
Der antagonistische Kompagnon der metrischen Sensitivit"at bzw.
der Sensitivit"at bez"uglich zweier nicht-negativer Funktionen
ist die metrische Komanenz bzw. die Komanenz bez"uglich zweier nicht-negativer Funktionen, die wie folgt
definiert sei:
\newline
\newline
{\bf Definition 1.3:\newline Komanenz bez"uglich zweier nicht-negativer 
Funktionen}\index{Komanenz bez\"uglich zweier nicht-negativer Funktionen}\newline
{\em F"ur jede Wellenfunktion $\Phi\in B^{A\times\mathbb{R}}$ und f"ur
jedes Paar $(\delta_{A},\delta_{B})$ zweier
auf den Mengen $A$ und $B$ 
erkl"arter Funktionen $\delta_{A}\in [0,\infty[^{A\times A}$ bzw. $\delta_{B}\in [0,\infty[^{B\times B}$
existiert
die Abbildung
\begin{equation}\label{aoilzbb}
\begin{array}{c}
{\rm B}_{\Phi}^{\delta_{A},\delta_{B}}:\mathbb{R}^{+}\times\mathbb{R}^{+}\to\mathbb{R}^{+}\cup\{-\infty\}\ ,\\
(\delta,t)\mapsto {\rm B}_{\Phi}(\delta,t)\ ,\\
\mbox{f"ur}\\
\sup\{q\in\mathbb{R}^{+}:\delta_{A}(x,y)\leq q\ \land\ 
\vartheta\in[-t,t]\ \land\ (x,y)\in A\times A\\
\Rightarrow\\
\delta_{B}(\Phi(x,\vartheta),\Phi(y,\vartheta))\leq\delta \}\\
=:{\rm B}_{\Phi}^{\delta_{A},\delta_{B}}(\delta,t)\ ,
\end{array}
\end{equation}
exakt welche wir als die Komanenzfunktion der Wellenfunktion $\Phi$
bez"uglich der Funktionen $\delta_{A}$ und $\delta_{B}$
bezeichnen.\index{Komanenzfunktion bez\"uglich zweier nicht-negativer Funktionen}
Genau dann, wenn 
\begin{equation}\label{aoilzb}
{\rm B}_{\Phi}^{\delta_{A},\delta_{B}}>-\infty
\end{equation}
positiv ist, nennen wir die 
Wellenfunktion $\Phi$ bez"uglich der Funktionen $\delta_{A}$ und $\delta_{B}$
komanent.}
\newline
\newline
Wir vereinbaren, dass 
f"ur jede Wellenfunktion $\Phi\in B^{A\times\mathbb{R}}$,
falls $\delta_{A}=\delta_{B}=\delta$ ist,
$${\rm B}_{\Phi}^{\delta_{A},\delta_{B}}:={\rm B}_{\Phi}^{\delta}$$
sei. Von metrischer Komanenz einer 
Wellenfunktion sprechen wir exakt dann, wenn diese Wellenfunktion bez"uglich zweier Metriken komanent 
ist.\index{metrische Komanenz}
Wir weisen auf den entscheidenden Sachverhalt hin, dass 
die Komanenz einer 
Wellenfunktion bez"uglich zweier nicht-negativer Funktionen nicht deren
Sensitivit"at bez"uglich derselben ausschliesst.
Denn, wenn $\Phi\in B^{A\times\mathbb{R}}$ bez"uglich der Funktionen $\delta_{A}$ und $\delta_{B}$
komanent ist,
gilt die Aussage (\ref{aoilzb}), die genau dann wahr ist, wenn f"ur alle 
$(\delta,t)\in \mathbb{R}^{+}\times\mathbb{R}^{+}$ die Menge
\begin{displaymath}
\begin{array}{c}
\{q\in\mathbb{R}^{+}:\delta_{A}(x,y)\leq q\ \land\ 
\vartheta\in[-t,t]\ \land\ (x,y)\in A\times A\Rightarrow\\
\delta_{B}(\Phi(x,\vartheta),\Phi(y,\vartheta))\leq\delta \}\not=\emptyset
\end{array}
\end{displaymath}
ist. Auch, wenn das Supremum aller Aufl"osungsfeldwerte 
$$\sup\Delta_{\Phi}^{\delta_{A},\delta_{B}}(A)>0$$
positiv ist und es daher Elemente $x_{\star}\in A$ geben muss,
f"ur welche die Sensitivit"at bez"uglich der Funktionen $\delta_{A}$ und $\delta_{B}$
gegeben ist, ist damit lediglich behauptet, dass die 
Menge 
$M_{\Phi}^{\delta_{A},\delta_{B}}(x_{\star})\not=\emptyset$ nicht leer ist, sodass es 
positive Zahlen $\delta(x_{\star})$ gibt, f"ur welche die Aussage 
\begin{displaymath}
\begin{array}{c}
\forall\ \varepsilon\in\mathbb{R}^{+}
\ \exists\ y\in \{z\in A:\delta_{A}(y,x_{\star})<\varepsilon\}\\
\forall t_{\star}\in\mathbb{R}^{+}\ \exists\ 
t\in ]-\infty,t_{\star}[\cup ]t_{\star},\infty[:
\delta_{B}(\Phi(x_{\star},t),\Phi(y,t))\ >\ \delta(x_{\star})
\end{array}
\end{displaymath}
wahr ist: Wenn wir nun $\delta(x_{\star})/2$ nehmen und eine beliebige positive reelle Zahl $t\in \mathbb{R}^{+}$, so
ist dennoch keineswegs impliziert, dass die Menge
\begin{displaymath}
\begin{array}{c}
\{q\in\mathbb{R}^{+}:\delta_{A}(x,y)\leq q\ \land\ 
\vartheta\in[-t,t]\ \land\ (x,y)\in A\times A\Rightarrow\\
\delta_{B}(\Phi(x,\vartheta),\Phi(y,\vartheta))\leq\delta(x_{\star}) \}
\end{array}
\end{displaymath}
leer ist und damit die Komanenz der 
Wellenfunktion
$\Phi$ bez"uglich der Funktionen $\delta_{A}$ und $\delta_{B}$
negiert. 
Die Aussage 
\begin{equation}\label{aooilzb}
\inf\{{\rm B}_{\Phi}^{\delta_{A},\delta_{B}}(q,t):t\in\mathbb{R}^{+}\}>0
\end{equation} allerdings
widerspricht der Komanenz der 
Wellenfunktion
$\Phi$ bez"uglich der Funktionen $\delta_{A}$ und $\delta_{B}$ 
f"ur alle 
$$q\geq\sup\Delta_{\Phi}^{\delta_{A},\delta_{B}}(A)$$  
sehr wohl. Die Aussage (\ref{aooilzb}) ist aber 
f"ur keine positive Zahl $q$ impliziert, wenn 
lediglich 
die Komanenz der 
Wellenfunktion
$\Phi$ bez"uglich der Funktionen $\delta_{A}$ und $\delta_{B}$ 
in Form der Reichhaltigkeitsaussage (\ref{aoilzb}) behauptet
ist.
\newline
Weshalb stufen wir diesen Sachverhalt der Nicht-Komplementarit"at von Komanenz und Sensitivit"at
als einen entscheidenden ein, dass 
die Komanenz einer 
Wellenfunktion bez"uglich zweier nicht-negativer Funktionen nicht deren
Sensitivit"at bez"uglich derselben ausschliesst?
\newline
Im ersten Teil der 
der Konzepte der abstrakten Ergodentheorie \cite{rabe}
zeigten wir auf eine h"ochst einfache Weise den insensitven topologischen  
Ergodensatz $(2.2)^{1}$,\index{insensitver topologischer Ergodensatz} welcher behauptet, dass 
das Mengensystem 
$$[[\Phi]]_{\mathbf{T}(\mathbf{P}_{2}\Phi)}=
\{\mathbf{cl}_{\mathbf{T}(\mathbf{P}_{2}\Phi)}(\tau):\tau\in[\Phi]\}$$
abgeschlossener Trajektorienh"ullen 
den Zustandsraum $\mathbf{P}_{2}\Phi$ jeder topologisierten Flussfunktion 
$(\Phi,\mathbf{T}(\mathbf{P}_{2}\Phi))$ partioniere, 
falls die Fl"usse $\Phi^{t}$
f"ur alle $t\in\mathbb{R}$ bez"uglich der 
Zustandsraumtopologie $\mathbf{T}(\mathbf{P}_{2}\Phi)$ stetig sind.\newline
Der insensitve topologische Ergodensatz $(2.2)^{1}$ ist darin 
viel umfassender als der Satz von der Existenz der Zimmer $(2.1.2)^{0}$\index{Satz von der Existenz der Zimmer} der
Abhandlung "uber den elementaren Quasiergodensatz \cite{raab},
dass der insensitve topologische 
Ergodensatz 
auf dem Generalit"atsniveau der allgemeinen Topologie 
verfasst ist, w"ahrend der 
Satz von der Existenz der Zimmer
glatte Trajektorien im endlichdimensionalen und reellen Zustandsraum zum 
Gegenstand hat.
Den elementaren Quasiergodensatz zeigten wir indess etwas
aufwendiger auf der Grundlage der Immanenz und der Komanenz 
trajektorieller Partitionen im endlichdimensionalen und reellen Zustandsraum.
In der 
Abhandlung "uber den elementaren Quasiergodensatz 
mussten wir 
Immanenz und Komanenz 
als Eigenschaften jeweiliger 
trajektorieller Partitionen
zun"acht erst herausstellen und 
einf"uhren.
Wenn wir den zentralen Teil der Argumentation betrachten, die zum Satz von der 
Existenz der Zimmer $(2.1.2)^{0}$ leitet, denjenigen Teil der Argumentation, der
bereits
von Immanenz und der Komanenz ausgeht, so
ist diese Argumentation nicht wesentlich aufw"andiger, als es 
der Beweis des insensitven
topologischen Ergodensatzes $(2.2)^{1}$ ist.
\newline
Die zum
Satz von der 
Existenz der Zimmer leitende 
Argumentation l"asst sich offensichtlich f"ur metrische Zustandsr"aume analog 
finden, weil sich auch f"ur metrische Zustandsr"aume die 
analoge Immanenz und Komanenz
verfassen l"asst, worauf wir in der Abhandlung "uber den elementaren Quasiergodensatz
hinwiesen.
Die
Voraussetzung des Satzes von der 
Existenz der Zimmer ist es, dass eine jeweilige
trajektorielle Partition des Zustandsraumes
immanent und komanent ist.
Dieser Voraussetzung
gen"ugt jede  
trajektorielle Partition $[\Phi]$, die durch eine topologisierte Flussfunktion
$(\Phi,\mathbf{T}(d))$ erzeugt ist, deren 
Zustandsraumtopologie $\mathbf{T}(d)$ durch eine Metrik $d$ 
induziert ist, wenn die Fl"usse $\Phi^{t}$
f"ur alle $t\in\mathbb{R}$ stetig sind bez"uglich der metrischen Zustandsraumtopologie
$\mathbf{T}(d)$.\newline
Es ist aber nicht umgekehrt so, dass jede trajektorielle Partition $[\Psi]$ eines 
metrischen 
Zustandsraumes $(\mathbf{P}_{2}\Psi,d)$, welche 
immanent und komanent ist,
notwendigweise so beschaffen sein muss, dass 
all deren 
Fl"usse $\Psi^{t}$
f"ur $t\in\mathbb{R}$ bez"uglich der metrischen Zustandsraumtopologie stetig sein m"ussen, welche 
durch die Metrik 
$d$ induziert ist:
Innerhalb der metrischen Zustandsraumtopologie besagt die 
Partitionsbehauptung des Satzes von der Existenz der Zimmer $(2.1.2)^{0}$  
deutlich mehr als die 
Partitionsbehauptung des insensitven topologischen
Ergodensatzes $(2.2)^{1}$, wie wir im Laufe dieser Abhandlung 
erfahren werden.
Und zwar besagt die 
Partitionsbehauptung des Satzes von der Existenz der Zimmer $(2.1.2)^{0}$ 
insofern mehr als der insensitve topologische
Ergodensatz $(2.2)^{1}$, als die Voraussetzung des insensitven
topologischen Ergodensatzes 
das Auftreten
mengenweiser Sensititvit"at ausschliesst.
Im Fall metrischer Zustandsr"aume impliziert diese im dritten
Kapitel vorgestellte allgemeine Form der Sensititvit"at
die metrische Sensitivit"at. Metrisch sensitive
Flussfunktionen sind kein Gegenstand des insensitven
topologischen Ergodensatzes. Sie k"onnen aber 
Gegenstand des Satzes von der 
Existenz der Zimmer sein, so dieser 
f"ur metrische Zustandsr"aume formuliert wird.
Und zwar k"onnen metrisch sensitive
Flussfunktionen gerade deshalb 
Gegenstand des Satzes von der 
Existenz der Zimmer sein,
weil metrisch sensitive
Flussfunktionen, wie wir darlegten, 
sehr wohl 
komanent sein k"onnen.
Das ist der Grund, weshalb wir die herausgestellte 
Vereinbarkeit metrischer Komanenz mit der 
metrischen Sensititvit"at als entscheidend einstufen.
\newline\newline
{\bf 1.1.3 Metrische Sensitivit"at als intermittente Negation der\newline\quad\quad\quad\quad\quad\quad Diskontinuit"at}\newline\newline
F"ur den Fall, dass die Distanzfunktionen $\delta_{A}=\delta_{B}$ Metriken 
sind und $\Phi$ eine Flussfunktion ist, erhalten wir den herk"ommlichen Begriff der
metrischen Sensitivit"at einer Flussfunktion.\index{metrische Sensitivit\"at}
In dieser Definition lassen wir aber
Distanzfunktionen im weitesten Sinn zu, unter welchen 
wir exakt alle auf kartesischen Quadraten von Mengen $\mathbf{P}_{1}\mathbf{P}_{1}d$ definierten 
Funktionen $d$ verstehen wollen, deren Werte nicht 
negative reelle Zahlen oder $\infty$ sind\index{Distanzfunktion im weitesten Sinn}
und f"ur welche die Diagonalminimalit"at\index{Diagonalminimalit\"at von Distanzfunktionen}
\begin{equation}\label{aathbb}
d(x,x)\leq d(x,y)
\end{equation}
f"ur alle $x,y\in \mathbf{P}_{1}\mathbf{P}_{1}d$ gilt.
Dabei ist uns aber bewusst,
dass der 
f"ur Distanzfunktionen im weitesten Sinn verfasste
Sensitivit"atsbegriff   
entsprechend
befremdliche
Eigenschaften hat. 
Mit dem herk"ommlichen Begriff der
metrischen Sensitivit"at
verbindet
den befremdlichen 
Sensitivit"atsbegriff f"ur Distanzfunktionen im weitesten Sinn
aber die sukzessive Steigerbarkeit des Anspruches 
an jeweilige Distanzfunktionen bis hin zu der 
Gegebenheit, dass dieselben Metriken sind. In dem Fall, dass 
eine Wellenfunktion $\Phi$ bez"uglich einer Metrik ihrer Definitionsmenge $d_{\mathbf{P}_{1}\Phi}$ und einer
Metrik $d_{\mathbf{P}_{2}\Phi}$ ihrer 
Wertemenge 
sensitiv ist, ist sie unstetig bez"uglich diesen Metriken. Sensitivit"at ist demnach 
eine spezielle Form der Diskontinuit"at, die f"ur Wellenfunktionen verfasst ist.  
\newline 
Die Spezialisierung besteht dabei nicht allein in der Unbefristetheit der Diskontinuit"atsform, welche
Sensitivit"at ist. Wenn f"ur alle $t\in\mathbb{R}$ mit $\Phi^{t}:=\Phi(\mbox{id},t)$ 
der durch die Wellenfunktion $\Phi$ 
festgelegte Pseudofluss $\Phi^{t}$ ist und 
es wahr ist, dass eine beidseitig unbeschr"ankte Folge reeller Zahlen $\{t_{j}\}_{j\in\mathbb{N}}$ existiert, f"ur die
f"ur alle $j\in\mathbb{N}$ zutrifft, dass die Abbildung 
$\Phi^{t_{j}}$ unstetig ist, so ist $\Phi$ damit noch nicht notwendigerweise bez"uglich $d_{\mathbf{P}_{1}\Phi}$
und $d_{\mathbf{P}_{2}\Phi}$
sensitiv. Genau dann, wenn f"ur eine Wellenfunktion $\Phi$
eine beidseitig unbeschr"ankte Folge reeller Zahlen $\{t_{j}\}_{j\in\mathbb{N}}$
existiert, f"ur die $\Phi^{t_{j}}$ unstetig ist, nennen wir sie und die Wellenfunktion $\Phi$
unbefristet unstetig. Es ist hierbei klar, wie die unbefristete Diskontinuit"at 
einer Wellenfunktion
im Bezug auf jeweilige allgemeine Distanzfunktionen analog zu 
verallgemeinern ist.\index{unbefristete Diskontinuit\"at einer Wellenfunktion}
\newline
Die unbefristete Unstetigkeit der Wellenfunktion $\Phi$ impliziert noch nicht, dass 
sie sensitiv ist. 
Denn es kann ja sein, dass es eine Folge $\{Q_{j}\}_{j\in\mathbb{N}}$ von Teilmengen der 
Definitionsmenge $\mathbf{P}_{1}\Phi$ gibt, f"ur die
jeder unstetige Pseudofluss $\Phi^{t_{j}}$
auf dem Bereich $\mathbf{P}_{1}\Phi\setminus {\rm Q}_{j}$ stetig ist, wobei
jede Pseudotrajektorie $\tau$ der Kollektivierung $[\Phi]=\{\Phi(x,\mathbb{R}):x\in\mathbf{P}_{1}\Phi\}$ 
so verl"auft, dass die M"achtigkeit
$$\mathbf{card}(\{k\in\mathbb{N}:{\rm Q}_{k}\cap\tau\not=\emptyset\})\ \in\ \mathbb{N}$$
endlich ist. Die Unbefristetheit der Pseudoflussfamilie $\{\Phi^{t}\}_{t\in\mathbb{R}}$ 
impliziert also noch nicht die noch spezifischere Diskontinuit"atsbehauptung,
welche die Sensitivit"at einer Wellenfunktion aufstellt. Wir pr"agen die Parole von der 
intermittenten Negation,\index{intermittente Negation} f"ur welche die
positive Strukturbehauptung der Sensitivit"at innerhalb der negierten
Kontinuit"at als ein klassisches Beispiel gelte,
wohingegen der Begriff der intermittenten Unstetigkeit  
die Konzeption intermittenter Negation am deutlichsten illustriere:
\newline\newline 
{\bf Definition 1.4: Intermittente Unstetigkeit\index{intermittente Unstetigkeit}}\newline
{\em Seien $(X^{1},T^{1})$ und $(X^{2},T^{2})$
zwei topologische R"aume und sei $f\in (X^{2})^{X^{1}}$ eine Abbildung
des ersten derselben auf den zweiten.
Wir nennen $f$ genau dann im Punkt $x\in X^{1}$ intermittent unstetig, wenn die Abbildung $f$
in $x$ unstetig ist und wenn dabei
in jeder Umgebung ${\rm U}(x)\in T^{1}_{\{x\}}$ ein Punkt $y\in {\rm U}(x)$ und eine in
${\rm U}(x)$ enthaltene
Umgebung ${\rm U}(x,y)\in T^{1}_{\{y\}}$ des Punktes $y\in {\rm U}(x)$
existieren, auf der die Restriktion $f|{\rm U}(x,y)$ stetig ist.\newline 
Genau dann, wenn $f$ im Punkt $x$ intermittent unstetig ist, sagen
wir auch, dass $f$ in $x$ in erster Ordnung intermittent unstetig sei, falls
es eine Umgebung ${\rm U}(x)\in T^{1}_{\{x\}}$ gibt, in der
kein Punkt $y\in {\rm U}(x)\setminus\{x\}$ ist, in dem $f$ 
intermittent unstetig ist.\index{intermittente Unstetigkeit erster Ordnung}
\newline 
F"ur jede nat"urliche Zahl $n>1$
nennen wir $f$ hingegen genau dann im Punkt $x\in X^{1}$ in $n$-ter Ordnung intermittent unstetig,
wenn die Abbildung $f$
in $x$ unstetig ist und dabei
in jeder Umgebung ${\rm U}(x)\in T^{1}_{\{x\}}$ ein Punkt $y\in {\rm U}(x)\setminus\{x\}$
existiert, in dem $f$ in $n-1$-ter Ordnung intermittent unstetig ist.
\index{intermittente Unstetigkeit $n-1$-ter Ordnung}}
\newline 
\newline
Offenbar kann
ein topologische Raum
$(X^{1},T^{1})$, auf dem eine intermittent 
unstetige
Funktion definiert ist, auch eine endliche Topologie haben.
Nehmen wir eine endliche, mehr als ein Element umfassende Partition 
$P$ der Menge $X^{1}$, sodass das Mengensystem $P^{\cup}$
aller Vereinigungen "uber Teilmengen dieser Partition $P$ eine 
Topologie ist. Sind $Q_{1},Q_{2}\in P$ zwei verschiedene Teilmengen der Menge $X^{1}$,
so ist das Mengensystem $(P\setminus\{Q_{1}\})\cup\{Q_{1}\cup Q_{2}\}$ keine
Partition, wohl aber
$$T_{1}:=((P\setminus\{Q_{1}\})\cup\{Q_{1}\cup Q_{2}\})^{\cup}$$
eine Topologie auf $X^{1}$; wobei wir darauf hinweisen, dass
keine Partition die leere Menge als Element hat. Auf der Umgebung
$Q_{1}\cup Q_{2}\in T_{1}$ eines Punktes $x\in Q_{1}$
kann es f"ur denselben leicht eingerichtet werden, dass
$x$ ein Punkt ist, in dem eine Abbildung $g$ intermittent unstetig ist, die 
auf $Q_{2}$ stetig ist.
\newline
Die 
Intermittenz der Unstetigkeit gibt es, die 
Intermittenz der Stetigkeit gibt es nicht:
Wenn wir das analoge Pendant zu der 
intermittenten Unstetigkeit, das wir, wenn es existierte, 
intermittente Stetigkeit nennen sollten,
zu verfassen 
versuchen, indem wir in der Definition 1.1b
Stetigkeit und Unstetigkeit austauschen, stellen
wir fest, dass dies nicht geht. 
Der Leser malt sich zu seiner Orientierung auch "uber 
diesen negativen Sachverhalt
leicht eine im Ursprung bez"uglich der nat"urlichen Topologie $\mathbf{T}(1)$ 
des Zahlenstrahls
intermittent 
unstetige und eine in zweiter Ordnung intermittent unstetige reellwertige 
Funktion des Zahlenstrahls aus.\newline\newline 
{\bf 1.1.4 Von Metriken zu Distanzfunktionen}\newline\newline
Die Verallgemeinerungsschritte von Metriken zu Pseudometriken\index{Pseudometrik} oder von Metriken zu
Quasimetriken\index{Quasimetrik} oder hin zu Semimetriken\index{Semimetrik} sind 
Verallgemeinerungen, die sich bereits als durchaus
erforschenswerte Generalisierungen herausgestellt und etabliert haben; deren
Bezeichnungsweise befindet sich augenscheinlich in einer labilen Verfestigungsphase.
Als Pseudometriken gelten hierbei mit weitgehender Gel"aufigkeit 
Distanzfunktionen, die bis auf die Identitivit"at einer Metrik alle Eigenschaften einer Metrik haben,
w"ahrend Distanzfunktionen, die Metriken bis auf die Symmetrie einer Metrik sind,
meistens als Quasimetriken bezeichnet werden. 
Als Semimetriken werden oft Distanzfunktionen bezeichnet, die zwar
einzig die Dreiecksungleichung verletzen;
diese Bezeichnungsweise die Dreiecksungleichung nicht erf"ullender 
Distanzfunktionen als
Semimetriken
beinhaltet aber, dass eine grosse Vielfalt von Distanzfunktionen, die
auf sehr verschieden wirkende Weisen der 
Dreiecksungleichung nicht gen"ugen allesamt  
Semimetriken genannt werden. 
Wir verweisen zu diesem Thema auf Pervins Lehrbuch \cite{perv}
und auf
die Artikel
H.C. Reichels und W.Rupperts \cite{disa} und \cite{disb}, die ein wenig
von dem Geist aufbrechenden Interesses an Verallgemeinerungen des Begriffes der Metrik vermitteln.
Distanzfunktionen sind heute vielfach 
als "Aquivalente von Programmierungskonzeptionen
in verschiedensten Disziplinen
in
Gebrauch; beispielsweise, um
mit Hilfe von Distanzfunktionen Korreliertheiten zu quantifizieren und die jeweiligen Korrelationswerte
zu Klassifizierungen oder zur automatisierten Musterkennung heranzuziehen. Es gibt daher 
vielerlei neue Literatur zu Distanzfunktionen.
\newline
Mittels des verallgemeinerten Aufl"osungsfeldes\index{verallgemeinertes Aufl\"osungsfeld einer Wellenfunktion} der Wellenfunktion $\Phi$ bez"uglich der 
Funktionen $\delta_{A}$ und $\delta_{B}$
ist die in der Definition 1.2 punktweise definierte Sensitivit"at der Wellenfunktion $\Phi$ bez"uglich der 
Funktionen $\delta_{A}$ und $\delta_{B}$ quantitativ bewertet.
Die Funktionen $\delta_{A}$ und $\delta_{B}$ generieren dabei
metrische Topologien auf die bekannte Weise, falls sie 
Metriken auf $A$ bzw. $B$ sind; wobei $\delta_{A}$ und $\delta_{B}$ nicht nur in dem Fall, dass sie
Metriken sind,  
Topologien generieren.
Es gibt Semimetriken, die
auf die gleiche Weise wie Metriken generalisiert metrische Topologien
festlegen, die wir hier noch nicht vorstellen wollen. 
Nichtsdestotrotz:
Nur dann, wenn
jene Generalisierungen der Metriken dies leisten, dass
alle Topologien als
generalisierte metrische Topologien
dargestellt werden k"onnen,
umfasst der gerade vorgestellte punktweise Sensitivit"atsbegriff
der Wellenfunktion $\Phi$ bez"uglich zweier 
jeweiliger
Funktionen $\delta_{A}$ und $\delta_{B}$
die volle topologische Allgemeinheit.
Nur dann ist in jenem gerade vorgestellten Sensitivit"atsbegriff
bereits
die allgemeine topologische Verfassung des  
Sensitivit"atsbegriffes gefunden.
Andernfalls
aber
k"onnen wir diesen Sensitivit"atsbegriff 
nicht einfach analog und ohne einen Zwischenschritt
f"ur eine 
Flussfunktion binnen irgendeines topologischen Raumes verallgemeinern. 
\newline
Sodass wir eingestehen m"ussen, dass wir die letztgestellte Frage negieren m"ussen:
Nein, der Begriff
der Sensitivit"at eines Attraktors ist uns in allgemeiner
topologischer Generalit"at noch nicht zuhanden. Wir k"onnen aber immerhin den 
im Fall, dass $\delta_{A}$ und $\delta_{B}$ Metriken sind,
durch das verallgemeinerte Aufl"osungsfeld 
\begin{equation}\label{aalzbbb}
\Delta_{\Xi}^{d_{X}}:=\Delta_{\Xi}^{d_{X},d_{X}}
\end{equation}
gegebenen 
Sensitivit"atsbegriff
formulieren, der f"ur eine 
Flussfunktion $\Xi$ binnen eines metrischen Raumes $(X,d_{X})$
verfasst ist. Exakt 
diesen Sensitivit"atsbegriff nennen wir den metrischen Sensitivit"atsbegriff.  
Dieser metrische Sensitivit"atsbegriff\index{metrischer Sensitivit\"atsbegriff} ist
dabei deshalb eine topologische Kovariante,\index{topologische Kovarianz}
weil jeder Hom"oomorphismus eines metrischen Raumes
denselben auf einen metrischen Raum abbildet.
\newline
Wir haben an dieser Stelle vielleicht auch erste Ideen, wie wir diesen 
nunmehr
zuhandenen
metrischen Sensitivit"atsbegriff
zu dem Sensitivit"atsbegriff f"ur eine 
Flussfunktion binnen eines uniformen topologischen Raumes verallgemeinern k"onnen.
Aber, ob wir einen Sensitivit"atsbegriff f"ur eine 
Flussfunktion binnen {\em irgendeines} topologischen Raumes finden, ist fraglich.
Wenn wir den Sensitivit"atsbegriff von der punktweisen Sensitivit"at herkommend
mit Hilfe von Aufl"osungsfeldern 
alleine auf der Grundlage einer gegebenen Topologie des Zustandsraumes
verfassen wollen, so scheint
dies 
aussichtslos,
allenfalls mit erheblichen Umst"andlichkeiten verbunden zu sein.
\section{"Ubertragbarkeit des metrischen Sensitivit"atsbegriffes und 
Grenzen deren Praktizierbarkeit}
Was solls? Jede Menge, jeder Zustandsraum ist metrisierbar: Sollte es uns nicht m"oglich 
sein, einen Sensitivit"atsbegriff zu verfassen, der sich 
f"ur eine beliebige Flussfunktion $\xi$  
direkt auf eine beliebige
"Uberdeckung $\mathcal{A}\subset 2^{\mathbf{P}_{2}\xi}$
des Zustandsraumes $\mathbf{P}_{2}\xi=\bigcup\mathcal{A}$ oder wenigstens auf 
eine beliebige Topologisierung $\mathbf{T}(\mathbf{P}_{2}\xi)$ des Zustandsraumes bezieht,
dann haben wir einfach zuerst 
eine Metrik $d$ auf demselben anzulegen, auf der dann, durch dieselbe vermittelt, 
die metrische Sensitivit"at gem"ass der Definition 
1.2 als Begriff zur Verf"ugung steht.\newline
Und was soll hier "uberhaupt Direktheit 
des Bezuges auf eine Zustandsraumstrukturierung
des Sensitivit"atsbegriffes heissen?
Indirektheit 
dieses Bezuges des Sensitivit"atsbegriffes
heisst hier zwar lediglich, dass die Reihenfolge 
des Vorgehens bei der Anlage einer jeweiligen
Modellierung eines dynamischen Systemes eingehalten ist.
Indirektheit heisst zwar nur,
dass zu dem vorgebenen Paar
$$(\xi,\mathcal{A})$$ 
erst eine Metrisierung festgelegt $d$ werden muss, ehe dann die metrische
Sensitivit"at f"ur das jeweilige dynamischen System
verfasst ist, das dann das Tripel
$$(\xi,\mathcal{A},d)$$
objektiviert;
sodass anschliessend die auf $d$ 
bezogene Sensitivit"at der Flussfunktion $\xi$
befunden werden kann oder aber negiert. \newline
Selbst wenn dies der einzige Unterschied zwischen einem 
verallgemeinerten Sensitivit"atsbegriff einerseits, der sich 
direkt auf eine vorgegebene Zustandsraumstrukturierung
zu beziehen vermag und einem 
applizierten metrischen Sensitivit"atsbegriff andererseits w"are,
der den Umweg "uber eine
vorher erst anzulegende Metrisierung gehen muss: So w"are allein darin schon
der
Mehrarbeitsaufwand angezeigt, eine Metrisierung erst anzulegen.
Dieses unmittelbar einleuchtende, ergonomische Argument zeichnet sich indess
eher dadurch aus, dass es leicht vorgebracht werden kann,
als dadurch, dass das Licht dieses leicht einleuchtenden Argumentes
den haupts"achlichen Zweck der 
Generalisierung des Sensitivit"atsbegriffes erhellte.
\newline
Und jenes Argument ist fragw"urdig.
Nehmen wir an, es sei tats"achlich diejenige Form der Sensitivit"at  
der Flussfunktion $\xi$,
welche auf das jeweilige Ergebnis der Anwendung eines Metrisierungsalgorithmus bezogen ist,
das 
die Metrik $d$ sei,
die f"ur ein von aussen herangetragenes Interesse relevant ist:
Es ist in diesem Fall allemal vorher n"otig, diesen Metrisierungsalgorithmus, der $d$ bestimmt,
auf den Zustandsraum anzuwenden, gleich ob wir f"ur eine vorgelegte 
Zustandsraumstrukturierung $\mathcal{A}$ einen
Sensitivit"atsbegriff haben, der nicht von Interesse ist.
In diesem Fall
trifft jenes plausible ergonomische Argument also nicht.\newline
Es kann aber genausogut auch sein, dass 
pauschal
ein nicht interessenspezifischer Metrisierungsalgorithmus eine Metrik $d$
und eine Form der Sensitivit"at produziert, wo
eigentlich 
ein anderer Sensitivit"atsbegriff 
als der auf $d$ basierende Sensitivit"atsbegriff
der relevante und der den Interessen angemessene w"are. Es ist m"oglich,
dass der auf dem Metrisierungsresultat $d$ basierende Sensitivit"atsbegriff als Instrument gegeben ist, derweil
ein Sensitivit"atsbegriff nat"urlich w"are, der
auf die vorgegebene Zustandsraumstrukturierung $\mathcal{A}$ eingeht.
Das Mengensystem $\mathcal{A}$ steht hier, vermutlich ein wenig deutungsentzogen 
abstrahiert,
etwa f"ur die
Detektorik, die ein Beobachter zur Verf"ugung hat, dem ein 
dynamisches System an die Hand gegeben werden soll:
Die Topologie des 
Zahlenstrahles abstrahiert zusammen mit der 
jeweiligen Strukturierung $\mathcal{A}$ 
die beschreibungsrelevante Detektorik eines 
Entwicklungsvorganges.\newline
Nehmen wir an, uns sei eine Flussfunktion $\Phi$ vorgelegt, deren Zustandsraum
durch eine nicht metrisierbare Topologie
$\mathbf{T}(\mathbf{P}_{2}\Phi)$ strukturiert sei, wobei
die Flussfunktion $\Phi$ bez"uglich dieser Topologie stetig sei, sodass
der im dritten Kapitel des ersten Teiles der Konzepte der abstrakten Ergodentheorie \cite{rabe}
besprochene Kommutator f"ur $\Phi^{t}:=\Phi(\mbox{id},t)$ und f"ur alle $t\in\mathbb{R}$
$$[\mathbf{cl}_{\mathbf{T}(\mathbf{P}_{2}\Phi)}, \Phi^{t}]=\emptyset$$
ist und der insensitve topologische Ergodensatz das gleiche besagt, wie der
verallgemeinerte insensitve Ergodensatz $(3.3)^{1}$: N"amlich, dass
das Mengensystem bez"uglich $\mathbf{T}(\mathbf{P}_{2}\Phi)$
abgeschlossener Trajektorien $\tau\in[\Phi]:=\{\Phi(x,\mathbb{R}):x\in\mathbf{P}_{2}\Phi\}$
$$\{\mathbf{cl}_{\mathbf{T}(\mathbf{P}_{2}\Phi)}(\tau):\tau\in[\Phi]\}=[[\Phi]]_{\mathbf{T}(\mathbf{P}_{2}\Phi)}\in\mathbf{part}(\mathbf{P}_{2}\Phi)$$
den Zustandsraum $\mathbf{P}_{2}\Phi$ partioniert.
Nehmen wir ferner an, es sei uns m"oglich gewesen, den 
Sensitivit"atsbegriff in allgemeiner topologischer Generalit"at auszusprechen. Anschliessend
sei es uns gelungen, mit Hilfe dieses topologischen Sensitivit"atsbegriffes
Aussagen zu formulieren, die f"ur jede 
Flussfunktion 
gelten, die bez"uglich einer jeweiligen Topologie stetig ist.
Diese Aussagen st"unden uns nicht zur Verf"ugung, wenn wir uns darin ersch"opften, 
allein daran zu gehen, auf dem Zustandsraum 
$\mathbf{P}_{2}\Phi$ eine Metrik $d_{\mathbf{P}_{2}\Phi}$ anzulegen, deren
Beziehung zu der nicht metrisierbaren Topologie $\mathbf{T}(\mathbf{P}_{2}\Phi)$
wir dann erst noch zu untersuchen h"atten.\newline
Hier mag der Physiker oder der Numeriker einwenden,
dass
dies doch praxisfremde Verstiegenheiten seien, die
Phasenr"aume in Betracht ziehen, deren Topologie nicht metrisierbar ist.  
Wir teilen diese Auffassung nicht, widersprechen ihr aber an dieser Stelle
auch nicht. Stattdessen weisen wir aber darauf hin, dass es uns der verallgemeinerte
insensitve Ergodensatz $(3.3)^{1}$ zweierlei erlaubt:\newline\newline
{\bf 1.} Der verallgemeinerte insensitve Ergodensatz erlaubt uns,
im Hinblick auf ergodentheoretische Fragen
frei mit beliebigen Netzen zu arbeiten, die durch die Strukturierung
$\mathcal{A}$ darstellbar sind. Die 
Kl"arung der Frage, ob dabei 
im jeweiligen Einzelfall 
die Voraussetzung des verallgemeinerten insensitven Ergodensatzes\index{verallgemeinerter insensitver Ergodensatz}
vorliegen, ist dabei
durch den Explizierungssatz\index{Explizierungssatz} $(3.8)^{1}$
vereinfacht.\newline\newline
{\bf 2.} Ausgerechnet die Numerik ist durch die Konfrontation mit dem Fehlen 
von 
Strukturvoraussetzungen charakterisiert. Dieses Fehlen ist das "Aquivalent 
der
Allgemeinheitsstufe, in welcher
der verallgemeinerte insensitve Ergodensatz verfasst ist: Es ist der verallgemeinerte insensitve Ergodensatz,
der deshalb f"ur die ergodentheoretischen Belange der Numerik relevant ist. Die
Numerik interessiert
sich ausserdem sehr f"ur Sensitivit"atsfragen, die auf derselben
Allgemeinheitsstufe gestellt sind, auf der wir den 
Sensitivit"atsbegriff verfassen wollen.\newline\newline
Im Hinblick auf die Selbstbestimmung\footnote{Selbstverst"andlich hat die Numerik in der Behandlung
kontinuierlicher dynamischer Systeme ein grosses Anwendungsgebiet und daher 
eine Interesse an kontinuierlichen dynamischen System --
als ihren jeweiligen Anwendungsobjekten.} der Numerik
allerdings interessieren den Numeriker  
vorzugsweise diskrete dynamische Systeme, w"ahrend wir
hier Flussfunktionen behandeln und diskrete dynamische Systeme nur sekund"ar untersuchen.
Gerade die durch Flussfunktionen dargestellten kontinuierlichen dynamischen Systeme
erm"oglichen es, dass die nat"urliche Vertrautheit des Menschen mit dem Zahlenstrahl auf beliebige
Zustandsr"aume "ubertragen werden kann. Dies ist ein wichtiger Gesichtspunkt des ideengeschichtlichen
Erfolges des Ansatzes des Determinismus, der vermutlich sogar 
in der
evolutionshistorischen Ideengeschichte der Biologie wurzelt. Diese
evolutionshistorische Wurzel des Determinismus reicht dabei aber "uber die kulturhistorische Ideengeschichte ebenso hinaus, 
wie die Frage unsere thematischen 
Grenzen "uberschreitet, 
was es sei, was die Perspektive des Determinismus so erfolgreich macht.\newline
Jene "Ubertragbarkeit des Zahlenstrahles
auf einen beliebigen Zustandsraum k"onnen wir aber hier f"ur uns nutzen:
Zu jeder Teilmenge $a\subset\mathbb{R}$ des Zahlenstrahles 
und zu jeder Menge von Zust"anden 
$Z\subset\mathbf{P}_{2}\xi$
gibt es
f"ur jede Flussfunktion $\xi$
die Menge
\begin{equation}
a\times_{\xi} Z:=\bigcup_{z\in Z}\xi(a,z)\subset \mathbf{P}_{2}\xi\ ,
\end{equation}
sodass f"ur jedes Mengensystem des Zahlenstrahles 
$$\mathcal{R}\subset 2^{\mathbb{R}}$$
das Zustandsraum-Mengensystem
\begin{equation}
\begin{array}{c}
\mathcal{R}\times_{\xi} Z=
\{a\times_{\xi} Z:a\in\mathcal{R} \}\\
=\{\xi(z,{\rm R}):(z,{\rm R})\in Z\times\mathcal{R}\}\subset2^{\mathbf{P}_{2}\xi} 
\end{array}
\end{equation}
bestimmt
ist. 
Unsere Notation unterscheidet also, so, wie es "ublich ist, die Abbildung $\mathbf{P}_{1}\times_{\xi}\mathbf{P}_{2}$, deren 
Definitionsmenge 
$$\mathbf{P}_{1}(\mathbf{P}_{1}\times_{\xi}\mathbf{P}_{2})=2^{\mathbb{R}}\times 2^{\mathbf{P}_{2}\xi}$$
ist, nicht von der 
auf dem 
kartesischen Produkt der 
Potenzmenge der Potenzmenge des Zahlenstrahles mit der Potenzmenge des 
jeweiligen Zustandsraumes $\mathbf{P}_{2}\xi$ definierten
Abbildung, welche die mengenweise Version der Abbildung $\mathbf{P}_{1}\times_{\xi}\mathbf{P}_{2}$
ist.\footnote{Diese gel"aufige Praktik sorgt manchmal allerdings auch f"ur Konfusion:  
F"ur jede reelle Zahl $t$ und f"ur jedes
Mengensystem
$\mathcal{Z}$ des Zustandsraumes ist
beispielsweise
$$\{\{t\}\}\times_{\xi}\mathcal{Z}=\{\{\{t\}\}\times_{\xi}{\rm L}:{\rm L}\in \mathcal{Z}\}=\{\xi({\rm L},t):{\rm L}\in \mathcal{Z}\}=
\xi^{t}\mathcal{Z} $$
$$=\{t\}\times_{\xi}\mathcal{Z}$$
das entwicklungsverschobene Mengensystem
des Zustandsraumes.
Wenn $\mathcal{Z}$ hierbei eine Zustandsraum"uberdeckung ist, ist
die entwicklungsverschobene Zustandsraum"uberdeckung keineswegs
notwendigerweise 
wieder eine Zustandsraum"uberdeckung.}
Zu jeder Menge von Zust"anden 
$Z\subset\mathbf{P}_{2}\xi$ ist damit
die 
Abbildung
$$\mbox{id}\times_{\xi}Z:2^{\mathbb{R}}\to 2^{\mathbf{P}_{2}\xi}$$
festgelegt. 
Auf der gesamten Potenzmenge $2^{\mathbf{P}_{2}\xi}$ existiert
deren Inversion 
$(\mbox{id}\times_{\xi}Z)^{-1}$ dabei einzig in dem Fall, dass erstens 
der Zustandsraum mit einer einzigen Trajektorie $\xi(\alpha,\mathbb{R})=\mathbf{P}_{2}\xi$
identisch ist und dass zweitens $Z$ einelementig ist. 
Die "Ubertragungen 
F"ur beliebige 
Mengen von Zust"anden $Z$ transferieren
die Abbildungen $\mbox{id}\times_{\xi}Z$
jeweilige Zahlenstrahl-Mengensysteme
$\mathcal{R}\subset 2^{\mathbb{R}}$ 
auf den Zustandsraum in Gestalt der
Zustandsraum-Mengensysteme
$\mathcal{R}\times_{\xi} Z$. F"ur beliebige 
Mengen von Zust"anden $Z$ ist 
diese Vermittlung durch jeweilige Zahlenstrahl-Mengensysteme objektivierter
Strukturierungen des Zahlenstrahles 
vergleichsweise wenig transparent.
Anders ist es hingegen, wenn wir statt der beliebigen Zustandsmenge $Z$ eine Initialbasis der Flussfunktion $\xi$
nehmen:
Es sei f"ur jede Flussfunktion $\xi$
das Mengensystem
\begin{equation}
\begin{array}{c}
\mathbf{\alpha}(\xi):=\Bigl\{X\in 2^{\mathbf{P}_{2}\xi}:\bigcup[\xi]_{X}=2^{\mathbf{P}_{2}\xi}\ \land\\ 
x\in X\Rightarrow \bigcup[\xi]_{X\setminus\{x\}}\not=2^{\mathbf{P}_{2}\xi}\Bigr\}
\end{array}
\end{equation}
aller Zustandsmengen $X$ festgelegt, die so beschaffen sind, dass sie exakt einen Zustand jeder
Trajektorie $\tau\in [\xi]$ als Element enthalten. Wobei wir, so wie im ersten Teil der 
Konzepte der abstraken Ergodentheorie, f"ur jedes Mengensystem $\mathcal{Q}$ und f"ur jede Menge ${\rm A}$ mit
$\mathcal{Q}_{{\rm A}}$ die ${\rm A}$-Auswahl\index{Auswahl} aus $\mathcal{Q}$, d.h.,
das Mengensystem $\{a\in \mathcal{Q}:a\cap{\rm A}\not=\emptyset\}\subset\mathcal{Q}$, notieren.
Exakt jedes Element des Mengensystemes $\mathbf{\alpha}(\xi)$ nennen wir 
eine Initialbasis der Flussfunktion $\xi$.\index{Initialbasis einer Flussfunktion}
Ist $X\in \mathbf{\alpha}(\xi)$ eine Initialbasis, so ist die inverse Abbildung 
$(\mbox{id}\times_{\xi}X)^{-1}$ nur auf jeder Teilmenge 
$$Y\subset \mathcal{R}\times_{\xi}X$$ definiert und 
das Zustandsraum-Mengensystem
$\mathcal{R}\times_{\xi}X$ ist die maximale Definitionsmenge der 
Inversion $(\mbox{id}\times_{\xi}X)^{-1}$.
Es
gibt zu jedem Element $\tilde{a}$ des
Mengensystemes $\mathcal{R}\times_{\xi}X$ des Zustandsraumes genau 
eine Menge reeller Zahlen 
$$(\mbox{id}\times_{\xi}X)^{-1}(\tilde{a})\in \mathcal{R}\ ,$$
sodass
\begin{equation}
(\mbox{id}\times_{\xi}X)^{-1}(\mathcal{R}\times_{\xi}X)=\mathcal{R}
\end{equation}
die Wertemenge der 
auf der maximale Definitionsmenge $\mathcal{R}\times_{\xi}X$ definierten
Inversion $(\mbox{id}\times_{\xi}X)^{-1}$ ist. 
Das Produkt $$\mathcal{R}\times_{\xi}X\subset 2^{\mathbf{P}_{2}\xi}$$
ist f"ur jede   
Initialbasis $X\in \mathbf{\alpha}(\xi)$ der Flussfunktion $\xi$ und f"ur jedes 
Mengensystem $\mathcal{R}$ des Zahlenstrahles $\mathbb{R}$
gleichsam eine interferenzfreie Vereinigung "uber Multiplizierungen der auf jeweilige Trajektorien kopierten Beschaffenheit 
des Mengensystemes $\mathcal{R}$ des Zahlenstrahles:
Die bin"aren Schnittverh"altnisse $\{a\cap b:a,b\in\mathcal{R}\}$, die bin"aren Vereinigungsverh"altnisse 
$\{a\cup b:a,b\in\mathcal{R} \}$ sowie die bin"aren Substraktionsverh"altnisse
$\{a\setminus b:a,b\in\mathcal{R}\}$ werden ebenso auf alle Trajektorien "ubertragen, wie
die entsprechenden endlichen, abz"ahlbaren oder beliebigen
Schnitt-, Vereinigungs-, Substraktionsverh"altnisse. 
Pr"azisieren wir diese Aussage: F"ur alle 
Teilmengen des Zustandsraumes, die Elemente 
$\tilde{a},\tilde{b}\in \mathcal{R}\times_{\xi}X$ sind, gilt
\begin{equation}\label{sioma}
\begin{array}{c}
(\mbox{id}\times_{\xi}X)^{-1}(\tilde{a})\subset(\mbox{id}\times_{\xi}X)^{-1}(\tilde{b})\Leftrightarrow
\tilde{a}\subset\tilde{b}\ ,\\
(\mbox{id}\times_{\xi}X)^{-1}(\tilde{a})\setminus(\mbox{id}\times_{\xi}X)^{-1}(\tilde{b})\in\mathcal{R} \Leftrightarrow
\tilde{a}\setminus\tilde{b}\in \mathcal{R}\times_{\xi}X
\end{array}
\end{equation}
und f"ur alle $\mathcal{Q}\subset\mathcal{R}\times_{\xi}X$ gilt
\begin{equation}\label{siomb}
\begin{array}{c}
\bigcup\{(\mbox{id}\times_{\xi}X)^{-1}(q):q\in \mathcal{Q}\}\in\mathcal{R} \Leftrightarrow
\bigcup\{q\in \mathcal{Q}\}\in \mathcal{R}\times_{\xi}X\ ,\\
\bigcap\{(\mbox{id}\times_{\xi}X)^{-1}(q):q\in \mathcal{Q}\}\in\mathcal{R} \Leftrightarrow
\bigcap\{q\in \mathcal{Q}\}\in \mathcal{R}\times_{\xi}X\ .
\end{array}
\end{equation}
Wir k"onnen also sagen, dass die Abbildung
$(\mbox{id}\times_{\xi}X)^{-1}$ f"ur jede Initialbasis $X\in \mathbf{\alpha}(\xi)$
ein Isomorphismus ist, der das Zustandsraum-Mengensystem
$\mathcal{R}\times_{\xi}X$ auf das Zahlenstrahl-Mengensystem
$\mathcal{R}$ abbildet und dessen Isomorphie darin besteht, dass er
die Inklusionsrelation $\subset$, die Darstellbarkeit einer Menge als Subtraktion
sowie die Darstellbarkeit einer Menge als Schnitt oder als Vereinigung "uber Teilmengen
von $\mathcal{R}\times_{\xi}X$
invariant l"asst. 
Ist
$X\in \mathbf{\alpha}(\xi)$ eine Initialbasis, so ist daher $\mathcal{R}\times_{\xi}X$ genau dann eine 
Partition bzw. Sigmaalgebra bzw. Topologie des 
Zustandsraumes, wenn $\mathcal{R}$ eine Partition bzw. Sigmaalgebra bzw. Topologie des Zahlenstrahles ist;
was sich fortsetzen liesse. Die Isomorphismen
$(\mbox{id}\times_{\xi}X)^{-1}$ f"ur jeweilige Initialbasen $X$
lassen sich zu folgender differenzierten Version derselben abwandeln, wobei 
sie im allgemeinen die beschriebene Isomorphie verlieren:
F"ur jedes Zahlenstrahl-Mengensystem $\mathcal{R}$  
und f"ur jede
Initialbasis $X\in \mathbf{\alpha}(\xi)$ und f"ur jedes Mengensystem $\hat{X}\in 2^{X}$ sei
\begin{equation}
\mathcal{R}\times^{\xi}\hat{X}:=
\bigcup\{\mathcal{R} \times_{\xi[{\rm X}]}{\rm X}:{\rm X}\in\hat{X}\}\ ,
\end{equation}
wobei die Restriktionen 
der Flussfunktion $\xi$
auf die von den jeweiligen Mengen ${\rm X}\in\hat{X}$ getroffenen
Unterzustandsr"aume
$$\xi[{\rm X}]:=\xi|(\bigcup\{\xi(x,\mathbb{R})\in {\rm X}\}\times\mathbb{R})$$
Flussfunktionen sind. Genau dann, wenn
das Zahlenstrahl-Mengensystem $\mathcal{R}$ den Zahlenstrahl
und 
das Mengensystem $\hat{X}\in 2^{X}$
die Initialbasis $X$ "uberdeckt,
"uberdeckt deren beider Produkt $\mathcal{R}\times^{\xi}\hat{X}$ den 
Zustandsraum $\mathbf{P}_{2}\xi$. 
F"ur das Zustandsraum-Mengensystem $\mathcal{R}\times^{\xi}\hat{X}$ finden sich 
nur in Ausnahmen
Entsprechungen zu 
den soeben aufgez"ahlten Isomorphien (\ref{sioma}) und (\ref{siomb}):
Falls das Mengensystem $\hat{X}$ 
von Mengen der Initialbasis $X$
beispielsweise paarweise disjunkt ist, 
gilt allerdings die Implikation (\ref{sioma}) f"ur alle $\tilde{a},\tilde{b}\in \mathcal{R}\times^{\xi}\hat{X}$
und die Implikation
(\ref{siomb}) f"ur alle $\mathcal{Q}\subset\mathcal{R}\times^{\xi}\hat{X}$.
\newline
Falls $X$ eine Initialbasis bzw. $\hat{X}$ ein paarweise disjunktes und nicht leeres Mengensystem
einer Initialbasis
ist, 
geht aus der Isomorphie (\ref{siomb}) zwar hervor,
dass 
die Abbildung $\mbox{id}\times_{\xi}X$ bzw. $\mbox{id}\times^{\xi}\hat{X}$
metrische Topologien des Zahlenstrahles $\mathbf{T}(d_{1})$
f"ur eine jeweilige Metrik $d_{1}$ des Zahlenstrahles  
auf 
Topologien $\mathbf{T}(d_{1})\times_{\xi}X$ 
bzw. $\mathbf{T}(d_{1})\times^{\xi}\hat{X}$
des Zustandsraumes transferriert.
Es geht daraus aber nicht hervor, dass 
$\mathbf{T}(d_{1})\times_{\xi}X$ bzw. $\mathbf{T}(d_{1})\times^{\xi}\hat{X}$ eine metrisierbare Topologie ist. Im Gegenteil, wir erkennen leicht, dass
die Topologien
$\mathbf{T}(d_{1})\times_{\xi}X$ bzw. $\mathbf{T}(d_{1})\times^{\xi}\hat{X}$ gerade Beispiele f"ur nicht metrisierbare Topologien sind.
Ausser in dem trivialen Fall, in dem 
$\mathbf{card}(X)=1$ bzw. $\mathbf{card}(\bigcup\hat{X})=1$
ist: Denn es gilt f"ur jeden Zustand $z\in\mathbf{P}_{2}\xi$
\begin{equation}
\begin{array}{c}
\mathbf{card}\Bigl(\bigcap(\mathbf{T}(d_{1})\times_{\xi}X)_{\{z\}}\Bigr)=\mathbf{card}(X)\ ,\\
\mathbf{card}\Bigl(\bigcap(\mathbf{T}(d_{1})\times^{\xi}\hat{X})_{\{z\}}\Bigr)=\mathbf{card}(\bigcup\hat{X})\ .
\end{array}
\end{equation}
Nichtsdestotrotz finden wir zu jedem Paar 
$(d_{1},d_{X})$, das aus einer
Metrik des Zahlenstrahles $d_{1}$ und einer Metrik $d_{X}$ einer jeweiligen Initialbasis $X\in \mathbf{\alpha}(\xi)$
zusammengesetzt ist,
eine Metrisierung des Zustandsraumes auf die folgende Weise:\newline
Dass es zu jedem Zustand $z\in\mathbf{P}_{2}\xi$
genau ein Initialbasiselement $X^{\xi,\top}z\in X$ gibt, f"ur das
\begin{equation}\label{siomy}
z\in:\xi(X^{\xi,\top}z,\mathbb{R})
\end{equation}
gilt, wenn $X$ eine Initialbasis ist, erm"oglicht es uns, 
f"ur jede Flussfunktion $\xi$ und f"ur jede deren Initialbasen $X\in \mathbf{\alpha}(\xi)$
den Paraprojektor in die jeweilige Initialbasis $X$ als die Abbildung
\begin{equation}\label{siomz}
\begin{array}{c}
X^{\xi,\top}:\mathbf{P}_{2}\xi\to X\ ,\\
z\mapsto X^{\xi,\top}z
\end{array}
\end{equation}
festzulegen,\index{Paraprojektor in eine Initialbasis}
deren jeweiliger 
Wert durch (\ref{siomy}) eindeutig 
festgelegt ist.
Ist 
$d_{X}$ eine Metrik der Initialbasis $X$, so induziert $d_{X}$ die metrische 
Topologie $\mathbf{T}(d_{X})$, f"ur die f"ur jedes Zahlenstrahl-Mengensystem $\mathcal{R}$
\begin{equation}
\mathcal{R}\times^{\xi}\mathbf{T}(d_{X})
\end{equation}
genau dann eine Topologie ist, wenn $\mathcal{R}$ eine Topologie des Zahlenstrahles ist.
Ferner 
ist dann f"ur jede Norm $\eta$ der reellen Ebene $\mathbb{R}^{2}$ die Abbildung
\begin{displaymath}
\mathbf{d}_{\xi}(\eta,d_{1},d_{X}):\mathbf{P}_{2}\xi\times\mathbf{P}_{2}\xi\to[0,\infty]
\end{displaymath}
mit den jeweiligen Werten
\begin{equation}
\begin{array}{c}
\mathbf{d}_{\xi}(\eta,d_{1},d_{X})(q,p)\\
:=\eta\Bigl(d_{X}(X^{\xi,\top}q,X^{\xi,\top}p),\ d_{1}(\xi(X^{\xi,\top}q,\mbox{id})^{-1}(q),\
\xi(X^{\xi,\top}p,\mbox{id})^{-1}(p))\Bigl)
\end{array}
\end{equation}
f"ur alle $(q,p)\in \mathbf{P}_{2}\xi\times\mathbf{P}_{2}\xi$
eine Metrik des Zustandsraumes, da
$\eta$ eine Norm des $\mathbb{R}^{2}$ ist: $\mathbf{d}_{\xi}(\eta,d_{1},d_{X})$ ist in dem 
Sinn identifizierend, dass $\mathbf{d}_{\xi}(\eta,d_{1},d_{X})(q,p)$
genau dann null ist, wenn $q=p$ ist. $\mathbf{d}_{\xi}(\eta,d_{1},d_{X})$ ist
offensichtlich symmetrisch.
Die Dreiecksungleichung f"ur die Norm $\eta$
bedingt, dass die Dreiecksungleichung f"ur die Metrik $\mathbf{d}_{\xi}(\eta,d_{1},d_{X})$
zutrifft.
Offensichtlich
ist die von der Metrik $\mathbf{d}_{\xi}(\eta,d_{1},d_{X})$ des Zustandsraumes
induzierte metrische Zustandsraumtopologie
f"ur alle jeweiligen Normen $\eta$ des $\mathbb{R}^{2}$
\begin{equation}
\mathbf{T}(\mathbf{d}_{\xi}(\eta,d_{1},d_{X}))=\mathbf{T}(d_{1})\times^{\xi}\mathbf{T}(d_{X})
\end{equation}
gleich und insoweit konstruktionsunabh"angig. Wir haben damit die 
Anleitung dazu, den Zustandsraum einer vorgelegten Flussfunktion $\xi$ auf eine
derselben adaptierte Weise zu metrisieren. Die Konkretisierung eines dementsprechend
adaptierten Metrisierungsverfahren fusste dabei allerdings auf der
Metrisierung einer jeweiligen Intialbasis. Dass dabei noch die Frage offen ist, wie 
eine zweckm"assige Metrisierung der jeweiligen Intialbasis gew"ahlt sein soll,
stellt nicht den Kern der Schwierigkeit dar, den Zustandsraum einer jeweiligen Flussfunktion $\xi$ auf 
die beschriebene Weise zu metrisieren. 
Jene noch offene Frage ist
ausserdem
von Interessenvorgaben abh"angig, die thematisch ausserhalb der Untersuchungen
dieser Abhandlung liegen. 
Der Problemzentrum, eine 
Zustandsraummetrisierung mit Hilfe der Metrisierung jeweiliger Intialbasen
zu vollziehen,
besteht n"amlich darin,
"uberhaupt Intialbasen zu bestimmen.
Die Vorstellung ist vermutlich konditioniert,
dass es m"oglich sei, gleichsam einem Gradienten 
in einem unbestimmten Sinn
zu folgen, der in einem Zustand
von einer Trajektorie wegzeigt, um dabei sukzessive auf
die Zust"ande einer Intialbasis zu treffen. 
Jene Vorstellung leitet indess nur in den F"allen
zu einer erfolgreichen Bestimmung einer Intialbasis,
in denen keine Sensitivit"at der Flussfunktion vorliegt.  
Stattdessen sind wir mit dem folgenden Sachverhalt konfrontiert:
Innerhalb eines chaotischen Attraktors beispielsweise ist es bekanntlich
prinzipiell unm"oglich, Intialbasen als solche 
mit numerischen Methoden nachzuweisen.\newline 
Dies gelte als das Schlusswort eines Pl"adoyers f"ur die Anerkennung der Bem"uh-
ung, einen m"oglichst allgemeinen
Sensitivit"atsbegriff zu verfassen. Nun
wollen wir angesichts der Suche nach einem allgemeinen Sensitivit"atsbegriff
den mit der Schwierigkeit 
des trajektoriellen Metriktransfers
thematisch verwandten
Behelf des Begriffes intrinsischer Sensitivit"at\index{intrinsische Sensitivit\"at} vorstellen.
\chapter{Ein Behelf und ein Ausblick}
\section{Fixierte und unfixierte intrinsische Sensitivit"at}
{\small Im Begriff der fixierten intrinsischen 
Sensitivit"at ber"uhren sich die Ph"anomene der Sensitivit"at im Sinne einer 
gewissen
Diskontinuit"at -- n"amlich im Sinn der Grenzwertaussage
(\ref{oweheva}) --  und der Nichtinvertierbarkeit kleinster Unterflussfunktionen. Blosse
Nichtinvertierbarkeit erscheint bei Flussfunktionen im weiteren Sinn
als deren intrinsische Sensitivit"at.\newline
Diese ph"anomenelle Kontamination erw"unschen wir naturgem"ass zun"achst nicht,
wollen wir doch das Ph"anomen der Sensitivit"at f"ur sich alleine erfassen.
Der Kontakt von Nichtinvertierbarkeit und der
Sensitivit"at im Sinne einer Diskontinuit"at kann bei Flussfunktionen
nicht auftreten. Verfassten wir die intrinsische 
Sensitivit"at nur f"ur Flussfunktionen, so vermieden wir zwar
jene ph"anomenelle Kontamination. Wir schl"ossen dann aber auch gleich 
die wichtige Klasse der Flussfunktionen im weiteren Sinn, die 
keine Flussfunktionen sind, aus dem Objektbereich des Begriffes
der intrinsischen Sensitivit"at aus:
\newline
Von jeder Flussfunktion $\theta$ verlangen wir, dass die Inversion 
$\theta(x,\mbox{id})^{-1}$ f"ur alle Zust"ande $x\in\mathbf{P}_{2}\theta$
existiert. Dies schliesst aus, dass Zust"ande $y\in\mathbf{P}_{2}\theta$ existieren,
durch welche 
Zyklen verlaufen. Es gibt die speziellen Flussfunktionen $\xi$ im weiteren Sinn, die wir
nach der Festlegung 1.1
Entwicklungen nennen, welche 
zwar nicht
die Invertierbarkeitsforderung an eine Flussfunktion erf"ullen; 
f"ur diese speziellen Flussfunktionen $\xi$ im weiteren Sinn existiert
aber der jeweilige inverse Fluss
$(\xi^{t})^{-1}=\xi(\mbox{id},t)^{-1}$ f"ur alle $t\in\mathbb{R}$.
Bei den Entwicklungen $\xi$
kann es vorkommen, dass durch einen ihrer Zust"ande $y$
ein Zyklus verl"auft, sodass
es eine Periode gibt, eine kleinste positive Zahl $T(y)\in\mathbb{R}^{+}$, f"ur welche
f"ur alle $t\in\mathbb{R}$ die Gleichung 
$\xi(y,t+\mathbb{Z}T(y))=\{\xi(y,t)\}$ erf"ullt ist.
Ferner sind diskrete Phasenfl"usse leicht in Flussfunktionen im weiteren Sinn einbettbar, 
die aber keine Flussfunktionen sind,
wie
wir im Abschnitt 1.1
im Zusammenhang mit der begrifflichen Abgrenzung seltsamer Attraktoren gegen sensitive Attraktoren  
erl"auterten; in dem Passus,
nachdem wir die Entwicklungen\index{Flussfunktion im weiteren Sinn} einf"uhrten.
Dort betrachteten wir die zu einem jeweiligen Autobolismus $f$
gem"ass
(\ref{seltd}) formulierbare, zu $f$ "aquivalente Flussfunktion im weiteren Sinn $\underline{f}$. Sie nimmt
in Spr"ungen die Werte der Kompositionen der Menge $\{f^{k}:k\in\mathbb{Z}\}$
an. Dabei ist die Funktion
$\underline{f}(x,\mbox{id})$ f"ur alle Zust"ande $x\in \mathbf{P}_{2}f$
exakt auf den Intervallen 
$[j,j+1[$ f"ur $j\in \mathbb{Z}$ konstant und
die maximalen Invertierbarkeitsbereiche der Funktion $\underline{f}(x,\mbox{id})$
sind Teilmengen der Elemente der Partition (\ref{brum}).\newline
Wir sehen: Wir kommen wir kaum darum herum, den Begriff
der intrisischen Sensitivit"at auch f"ur Flussfunktionen im weiteren Sinn
verfassen zu wollen, welche relevant sind;
einmal sind sie, "ahnlich wie $\underline{f}$,
diskrete Phasenfl"usse einbettende Flussfunktionen im weiteren Sinn; 
oder andermal Flussfunktionen im weiteren Sinn,
f"ur die Zyklen und Fixpunkte
auftreten.}
\newline\newline
Die grunds"atzliche Idee der intrinsischen Sensitivit"at ist es, 
die Entwicklung 
$$\xi(x,]a,b[+\{t\})$$
kleiner Intervalle $]a,b[\subset\mathbb{R}$ der phasischen 
Parametrisierung einer jeweilgen Trajektorie $\xi(x,\mathbb{R}) \in [\xi]$ f"ur 
variierende 
$t\in \mathbb{R}$ vom Zustandsraum 
einer Flussfunktion im weiteren Sinn $\xi$
zur"uckzu"ubertragen
auf den Zahlenstrahl und die Menge reeller Zahlen
$$\xi(x,\mbox{id})^{-1}(\xi(x,]a,b[+\{t\}))\subset \mathbb{R}$$
f"ur verschiedene $t\in \mathbb{R}$ zu betrachten.
Auf diese Weise ist es nicht n"otig, sich auf den 
Zustandsraum und eine Strukturierung desselben zu beziehen.
Der Begriff intrinsischer Sensitivit"at begn"ugt sich damit, wenigstens 
eine allgemeine Sensitivit"at f"ur 
einzelne Trajektorien, f"ur
jeweilige kleinste Unterflussfunktionen\footnote{Jede
Flussfunktion im weiteren Sinn
$\xi\subset(\mathbf{P}_{2}\xi\times\mathbb{R})\times\mathbf{P}_{2}\xi$
ist eine Menge von Paaren, sodass die Funktion $\xi(x,\mbox{id})$ f"ur alle $x\in\mathbf{P}_{2}\xi$ die Inklusion
$$\xi(x,\mbox{id})\subset(\{x\}\times\mathbb{R})\times\xi(x,\mathbb{R}))$$
erf"ullt
und insbeondere die Inklusion
$$\xi(x,\mbox{id})=\xi|((\{x\}\times\mathbb{R})\times\xi(x,\mathbb{R})))\subset\xi\ .$$
Daher folgen wir der gel"aufigen Bezeichnungssystematik, wenn wir 
$\xi(x,\mbox{id})$ als eine Unterflussfunktion\index{Unterflussfunktion} von $\xi$, 
der Flussfunktion im weiteren Sinn, bezeichnen.
Obwohl die Definitionsmenge dieser 
Unterflussfunktion
$$\mathbf{P}_{1}\xi(x,\mbox{id})=\mathbb{R}=\mathbf{P}_{2}\mathbf{P}_{1}\xi$$
nur der zweite Definitionsmengenfaktor von $\xi$ ist.
Jede Funktion $\phi$ ist genau dann eine Unterflussfunktion einer Flussfunktion $\xi$ im weiteren Sinn, wenn $\phi$
eine Restriktion der Flussfunktion $\xi$ ist und dabei 
der Zahlenstrahl auf eine Weise
die Definitonsmenge $\mathbf{P}_{1}\phi=\mathbb{R}$ 
oder der zweite Definitonsmengenfaktor 
$\mathbf{P}_{2}\mathbf{P}_{1}\phi=\mathbb{R}$ ist, dass
die Wertemenge $\mathbf{P}_{2}\phi$ ein Unterzustandsraum des Zustandsraumes $\mathbf{P}_{2}\xi$ ist.
$\xi(x,\mbox{id})$ ist in dem Sinn eine kleinste Unterflussfunktion, dass
die Restriktion 
$\xi(x,\mbox{id})|[0,1]$ der Unterflussfunktion $\xi(x,\mbox{id})$, 
wie
jede echte Restriktion der Unterflussfunktion $\xi(x,\mbox{id})$,
keine Flussfunktion im weiteren Sinn ist. Wenn $\phi$ eine Unterflussfunktion einer Flussfunktion $\xi$ im weiteren Sinn
ist, so ist $\phi$ genau dann eine kleinste Unterflussfunktion von $\xi$, wenn
$\phi$ die Definitonsmenge $\mathbf{P}_{1}\phi=\mathbb{R}$ hat.\newline
Ferner sehen wir, dass $\xi$ ist genau dann eine Flussfunktion im weiteren Sinn ist, die keine Flussfunktion ist,
wenn sie eine kleinste Unterflussfunktion hat, die
zwar eine Flussfunktion im weiteren Sinn, jedoch keine Flussfunktion ist.}
zu verfassen,
indem er auf dieselben  
den metrischen Sensitivit"atsbegriff "ubertr"agt, der sich auf die nat"urliche Metrik des Zahlenstrahles bezieht.
Der Begriff 
intrinsischer Sensitivit"at
will erfassen, was das auf eine jeweilige Trajektorie transferierte Pendant zu der 
metrischen  
Sensitivit"at einer Flussfunktion $\xi_{[1]}$ im weiteren Sinn 
ist, deren Zustandsraum der Zahlenstrahl
$$\mathbf{P}_{2}\xi_{[1]}=\mathbb{R}$$
ist. Nehmen wir die nat"urliche Zahlenstrahlmetrik $|\mathbf{P}_{1}-\mathbf{P}_{2}|$,
so ist $\xi_{[1]}$ bez"uglich derselben genau dann 
sensitiv, wenn es einen Zustand $x\in\mathbb{R}$ gibt und eine Zahl $\delta(x)\in ]0,\infty[$, 
f"ur die es f"ur alle 
$\varepsilon\in ]0,\infty[$ eine betragsm"assig beliebig grosse Zahl $t\in\mathbb{R}$ gibt, f"ur die
$$|x-y|<\varepsilon\ \land\ |\xi_{[1]}(x,t)-\xi_{[1]}(y,t)|\geq\delta(x) $$
gilt. 
Zu jeder Flussfunktion im weiteren Sinn $\xi$ und zu jeder Bijektion $Q$ des 
Zustandsraumes $\mathbf{P}_{2}\xi$
der ersteren 
auf eine beliebige Wertemenge $\mathbf{P}_{2}Q$ ist die Komposition
\begin{equation}\label{somx}
\begin{array}{c}
Q\star\xi:\mathbf{P}_{2}Q\times\mathbb{R}\to \mathbf{P}_{2}Q\\
(\alpha,t)\mapsto Q(\xi(Q^{-1}(\alpha),t))
\end{array}
\end{equation}
gegeben, die wieder
eine Flussfunktion im weiteren Sinn ist.
Von jeder Flussfunktion im weiteren Sinn $\Phi$
verlangen wir, dass f"ur alle Zust"ande $a,b\in \mathbf{P}_{2}\Phi$ die "Aquivalenz
$$b\in \Phi(a,\mathbb{R})\Leftrightarrow\Phi(b,\mathbb{R})=\Phi(a,\mathbb{R})$$
gilt, wenn wir verlangen, dass die Kollektivierung 
$[\Phi]=\{\Phi(a,\mathbb{R}):a\in \mathbf{P}_{2}\Phi\}$ eine Partition 
des Zustandsraumes $\mathbf{P}_{2}\Phi$ ist, was nicht m"oglich w"are, wenn
jene "Aquivalenz nicht g"alte. Dieselbe ist also zu der Partitivit"atsforderung an eine 
Flussfunktion "aquivalent. Ist $\xi$
speziell eine Flussfunktion im weiteren Sinn $\phi_{[1]}$,  
deren Zustandsraum der Zahlenstrahl ist, und $q$ eine Bijektion des Zahlenstrahles $\mathbb{R}$
auf eine beliebige Wertemenge $\mathbf{P}_{2}q$ derselben,
ist die Komposition $q\star\phi_{[1]}$ gerade das, was wir 
unter einer eindimensionalen Flussfunktion verstehen wollen.\index{eindimensionale Flussfunktion}
Zyklen entsprechen daher Kompaktifizierungen eindimensionaler Flussfunktionen.
Exakt jede Flussfunktion im weiteren Sinn, deren 
Zustandsraum der Zahlenstrahl ist, nennen wir eine Zahlenstrahlflussfunktion, sodass 
jede eindimensionale Flussfunktion
als eine
Komposition $q\star\phi_{[1]}$ einer Bijektion $q$ des Zahlenstrahles $\mathbb{R}$
mit einer
Zahlenstrahlflussfunktion $\phi_{[1]}$\index{Zahlenstrahlflussfunktion} auffassbar ist.
F"ur jede Flussfunktion $\xi$ und f"ur jedes deren
Kollektivelemente $\tau\in[\xi]$,  
dessen Kardinalit"at die des Kontinuums ist,  
gibt es eine Bijektion 
$q_{\tau}$ des Zahlenstrahles auf dieses Kollektivelement $\tau\in[\xi]$.
$\tau$ ist f"ur jede Zahlenstrahlflussfunktion $\psi$
der Zustandsraum der 
eindimensionalen Flussfunktion $q_{\tau}\star\psi$.
Jede 
eindimensionale Flussfunktion $q\star\phi_{[1]}$
ist genau dann eine Flussfunktion bzw. Entwicklung, wenn die Zahlenstrahlflussfunktion
$\phi_{[1]}$ eine Flussfunktion bzw. Entwicklung ist, wobei dann sowohl die Funktion
$\phi_{[1]}(y,\mbox{id})$ als auch deren Inversion
$\phi_{[1]}(y,\mbox{id})^{-1}$ f"ur jede reelle Zahl $y\in\mathbb{R}$
eine Bijektion des Zahlenstrahles auf denselben ist. Offensichtlich
ist auch, dass jede Zahlenstrahlflussfunktion $\phi_{[1]}$, die
eine finit-additive Entwicklung\index{finit-additive Entwicklung}
gem"ass der Festlegung 1.1 ist, bei der Komposition mit einer beliebigen
Bijektion $g$ die Eigenschaft 
der finiten Additivit"at in dem Sinn
nicht verlieren kann, dass $g\star\phi_{[1]}$ wieder eine
finit-additive Entwicklung ist. \newline
Nicht nur
die Eigenschaft einer Zahlenstrahlflussfunktion, eine Flussfunktion in weiteren Sinn bzw.
eine Flussfunktion bzw. eine Entwicklung bzw. finit-additiv zu sein,
ist f"ur jede Bijektion $q$ des Zahlenstrahles eine 
Invariante der
auf der Menge aller Zahlenstrahlflussfunktionen definierten 
Abbildung 
$q\star\mbox{id}$: F"ur 
jede Flussfunktion 
$\xi$ im weiteren Sinn ist f"ur jede Bijektion $Q$
des 
Zustandsraumes $\mathbf{P}_{2}\xi$
die Abbildung
$q\star\mbox{id}$ eine, die 
die aufgez"ahlten Eigenschaften erh"alt.
\newline
Wir erwarten daher, dass 
die Gegebenheit der Komposition $q\star\phi_{[1]}$ und die  
Kenntnis der Flussfunktion $\phi_{[1]}$ die Bijektion $q$ 
und ihrer Inverser $q^{-1}$ festlegt.
Da f"ur alle $x,t\in\mathbb{R}$
$$q\star\phi_{[1]}(q(x),t)=q(\phi_{[1]}(x,t))$$
ist, gilt
also auch 
$$q\star\phi_{[1]}\Bigl(q(x),\phi_{[1]}(x,\mbox{id})^{-1}(\vartheta)\Bigr)=q(\vartheta)$$
f"ur alle $x,\vartheta\in\mathbb{R}$.
Daher ist f"ur alle $x\in\mathbb{R}$ das 
Bildelement $q(x)$ durch die Gleichung 
$$q\star\phi_{[1]}\Bigl(q(x),\phi_{[1]}(x,\mbox{id})^{-1}(x)\Bigr)=q(x)$$
eindeutig festgelegt
und f"ur alle $\vartheta\in\mathbb{R}$
ist daher beispielsweise f"ur die Spezifizierung $x=0$
$$q\star\phi_{[1]}\Bigl(q(0),\phi_{[1]}(0,\mbox{id})^{-1}(\vartheta)\Bigr)$$
$$=q\star\phi_{[1]}(q(0),\mbox{id})\circ\phi_{[1]}(0,\mbox{id})^{-1}(\vartheta)=q(\vartheta)\ ,$$
also
$$\phi_{[1]}(0,\mbox{id})\circ \Bigl(q\star\phi_{[1]}(q(0),\mbox{id})\Bigr)^{-1}$$
$$=\phi_{[1]}\Bigl(0,\ (q\star\phi_{[1]}(q(0),\mbox{id}))^{-1}\Bigr)=q^{-1}\ .$$
Es ist offenbar f"ur alle $x\in\mathbb{R}$
\begin{equation}\label{somxa}
q^{-1}=\phi_{[1]} \Bigl(x,\ (q\star\phi_{[1]}(q(x),\mbox{id}))^{-1}\Bigr)
\end{equation}
die gesuchte Inversion der Bijektion $q$, welche deren Bild des Zahlenstrahles
$q(\mathbb{R})$ auf denselben abbildet. Und analog ist
\begin{equation}\label{saomxa}
Q^{-1}=\xi \Bigl(x,\ (Q\star\xi(Q(z),\mbox{id}))^{-1}\Bigr)
\end{equation}
die Inversion der Bijektion $Q$, welche das Bild des Zustandsraumes $Q(\mathbf{P}_{2}\xi)$
auf $\mathbf{P}_{2}\xi$ bijeziert. 
\newline
Wir fragen uns, 
ob und gegebenenfalls wie an der 
Komposition $q\star\psi$
in Erscheinung tritt, dass die Zahlenstrahlflussfunktion
$\psi$ gegebenenfalls sensitiv ist.
Exakt im Fall der Sensitivit"at von $\psi$ gibt es
einen Zustand $\alpha\in q(\mathbb{R})$ und eine Zahl $\delta\in ]0,\infty[$, 
f"ur die es f"ur alle 
$\varepsilon\in ]0,\infty[$ eine betragsm"assig beliebig grosse Zahl $t\in\mathbb{R}$ gibt, f"ur die
es einen Zustand $\beta\in q(\mathbb{R})$ gibt, f"ur den 
$|q^{-1}(\alpha)-q^{-1}(\beta)|<\varepsilon$ und dabei
\begin{displaymath}
|\psi(q^{-1}(\alpha),t)-\psi(q^{-1}(\beta),t)|\geq\delta
\end{displaymath}
ist, sodass 
\begin{equation}\label{somxb}
\begin{array}{c}
|q^{-1}(\alpha)-q^{-1}(\beta)|<\varepsilon\\ 
\land\\ 
|q^{-1}(q\star\psi(\alpha,t))-q^{-1}(q\star\psi(\beta,t))|\geq\delta
\end{array}
\end{equation}
gilt. 
F"ur jede Metrik des Zahlenstrahles 
$d_{\mathbb{R}}$ und f"ur jede Bijektion $q^{-1}$ auf denselben ist die 
auf dem Quadrat $\mathbb{R}\times\mathbb{R}$ definierte Abbildung $d_{\mathbb{R}}[q^{-1}]$ 
gegeben, die
mit der Setzung
\begin{equation}
d_{\mathbb{R}}[q^{-1}](x,y):=d_{\mathbb{R}}(q^{-1}(x),q^{-1}(y))
\end{equation}
f"ur alle $(x,y)\in \mathbb{R}\times\mathbb{R}$
festgelegt ist. Diese Abbildung
$d_{\mathbb{R}}[q^{-1}]$ ist 
offenbar eine Metrik. Genau dann, wenn 
$\psi$ eine sensitive Zahlenstrahlflussfunktion 
bez"uglich der Metrik $d_{\mathbb{R}}$
ist,
ist
$q\star\psi$ eine bez"uglich der Metrik $d_{\mathbb{R}}[q^{-1}]$
sensitive eindimensionale Flussfunktion.
Der Sachverhalt, dass $q\star\psi$ bez"uglich der Metrik $d_{\mathbb{R}}[q^{-1}]$
sensitiv 
ist, ist dabei allerdings eine redundante 
Kopie der Aussage, dass 
$\psi$  
bez"uglich der Metrik $d_{\mathbb{R}}$ 
sensitiv ist.\newline
Die verschwommene Frage,
wie an der eindimensionalen Flussfunktion
$q\star\psi$
in Erscheinung tritt, dass
$\psi$ eine sensitive Zahlenstrahlflussfunktion ist,
spitzt sich zu der Vorfrage zu, ob deren Sensitivit"at in einem noch unbestimmten Sinn
eine Invariante gegen"uber der Abbildung $q\star\mbox{id}$ 
ist:
Falls $g$ speziell ein Autobolismus ist, hier also eine Bijektion des Zahlenstrahles auf denselben,
k"onnen wir diese Vorfrage konkretisieren zu der Frage, ob die Sensitivit"at von
Zahlenstrahlflussfunktionen eine Invariante der Abbildung
$g\star\mbox{id}$ ist. Oder wir k"onnen fragen, f"ur welche Autobolismen
des Zahlenstrahles $g$ die Funktion 
$g\star\mbox{id}$ die Sensitivit"at der Zahlenstrahlflussfunktionen erh"alt.
Denn die Sensitivit"at der jeweiligen Bildelemente $g\star\psi$ ist uns
dann ein Begriff, weil in dem Fall, dass $g$ ein 
Autobolismus des Zahlenstrahles ist, $g\star\psi$ eine Zahlenstrahlflussfunktion
ist. \newline
Die klare Frage, f"ur welche Autobolismen
des Zahlenstrahles $g$ die Funktion 
$g\star\mbox{id}$ die Sensitivit"at der Zahlenstrahlflussfunktionen erh"alt,
ist dabei ebenso klar, wie sie gestellt ist, als eine solche einzustufen, die an dieser Stelle in unsere
Thematik eindringt: Wir wollen hier die Abstrahierbarkeit und Modifizierbarkeit der Sensitivit"at 
auskundschaften
und die Grundz"uge 
abstrakter Sensitivit"at 
bestimmen. In unserem Themenbereich liegt es,
zu versuchen, die gestellte Vorfrage dadurch zu einer 
Frage wandeln, dass wir einen Sensitivit"atsbegriff
f"ur beliebige eindimensionale Flussfunktionen
$q\star\psi$
formulieren. Wie sich zeigen wird, gibt es den 
Begriff der intrinsischen 
Sensitivit"at, der die Beschaffenheit 
des Zustandsraumes einer eindimensionalen Flussfunktion gleichsam 
aussen vorl"asst und der dadurch f"ur beliebige eindimensionale Flussfunktionen
$q\star\psi$ gegeben ist; also auch f"ur Zahlenstrahlflussfunktionen $g\star\psi$,
wenn $g$ ein
Autobolismus
des Zahlenstrahles ist. Daher 
ergibt sich 
die Frage, in welchem Verh"altnis die intrinsische Sensitivit"at
zur Sensitivit"at von Zahlenstrahlflussfunktionen steht, mit welcher
wir deren auf die nat"urliche Metrik bezogene 
metrische Sensitivit"at meinen. Denn 
die intrinsische Sensitivit"at erweist sich als
kein "Aquivalent der Sensitivit"at von Zahlenstrahlflussfunktionen.
All diesen Fragen 
gehen wir in dieser Arbeit nicht nach.
Wir wollen den Begriff der intrinsischen 
Sensitivit"at dennoch vorstellen, weil er eine Modifizierung der Sensitivit"at 
darstellt:
Es bezeichne 
\begin{equation}
Q\Delta P:=(Q\setminus P)\cup(P\setminus Q) 
\end{equation}
die symmetrische Differenz zweier Mengen\index{symmetrische Differenz zweier Mengen}
$Q$ und $P$. F"ur jede Wellenfunktion $\psi\in Y^{Y\times \mathbb{R}}$ sei
f"ur jedes deren Definitionsmengenelemente $(z,t)\in \mathbf{P}_{1}\psi$
die beiden mengenwertigen Abbildungen
\begin{equation}\label{sisiaomx}
\begin{array}{c}
\lbrack z,t\rbrack^{\psi}\ \mbox{bzw.}\ \lbrack z,t\rbrack_{\psi}:2^{\mathbb{R}}\to 2^{\mathbb{R}}\ ,\\
{\rm X}\mapsto\lbrack z,t\rbrack^{\psi}({\rm X})\ \mbox{bzw.}\ \lbrack z,t\rbrack_{\psi}({\rm X})
\end{array}
\end{equation}
durch die Festlegung ihrer jeweiligen Werte
\begin{equation}
\begin{array}{c}
\lbrack z,t\rbrack^{\psi}({\rm X}):=\psi(z,\mbox{id})^{-1}\circ\psi(\mbox{id},t)\circ\psi(z,\mbox{id})({\rm X})\ ,\\ 
\lbrack z,t\rbrack_{\psi}({\rm X}):=\lbrack z,t\rbrack^{\psi}({\rm X})\Delta({\rm X}\{t\})
\end{array}
\end{equation}
f"ur alle ${\rm X}\subset \mathbb{R}$ und $(z,t)\in \mathbf{P}_{1}\psi$ bestimmt.
Die Funktion
\begin{equation}\label{sisibomx}
\begin{array}{c}
\psi(z,\mbox{id})^{-1}\psi(\mbox{id},t)\circ\psi(z,\mbox{id})\ \Delta\ (\mbox{id}+\{t\})\\
=\lbrack z,t\rbrack_{\psi}
\end{array}
\end{equation}
kommt 
der durch ihre Schreibweise angedeuteten
Suggestion, die
Beschaffenheit eines Kommutators zu haben, nicht allein insofern nur eingeschr"ankt nach,
dass 
das Wertemengenelement $z\in \mathbf{P}_{2}\psi$ 
und die reelle Zahl $t$ in
verschiedenen Strukturvoraussetzungen 
stehen,\footnote{W"ahrend die Addition auf dem Zahlenstrahl $\mathbb{R}$
definiert ist, ist auf der Wertemenge $z\in \mathbf{P}_{2}\psi$
keine Addition gegeben. Selbst wenn die Wertemenge $\mathbf{P}_{2}\psi$ in einem linearen 
Raum liegt und die Addition seiner Elemente erkl"art ist,
besteht aber weiterhin die Asymmetrie zwischen Wertemengenelementen 
und Elementen des zweiten kartesischen Definitionsmengenfaktors, der $\mathbb{R}$ ist.
Ist $\theta\in Y^{Y\times X}$
eine auf dem kartesischen Produkt zweier Teilmengen zweier
linearer R"aume definierte Funktion, so ist diejenige 
Abwandlung des Ausdruckes $\lbrack z,t\rbrack_{\psi}$ der Gleichung
(\ref{sisibomx}) zu $[a,b]_{\theta}$, bei welcher $\psi$ durch $\theta$ 
und $z$ bzw. $t$ durch $a$ bzw. $b$
ersetzt ist,
erst dann m"oglich, wenn $X=Y$ ist. In diesem Fall existiert dann zwar
$$[b,a]_{\theta}=\theta(\mbox{id},a)\circ\theta(b,\mbox{id})\Delta\theta(b,\mbox{id})\circ(\mbox{id}+\{a\})=
[a,b]_{\theta}$$ 
gleichberechtigt
neben $[a,b]_{\theta}$. Dann ist allerdings die Koexistenz
von $[a,b]_{\theta}$ und $[b,a]_{\theta}$ als solche
in blosser identischer Redundanz nutzlos.}
weswegen wir 
f"ur jede Wellenfunktion $\psi$ exakt
die Abbildung $\lbrack z,t\rbrack_{\psi}$ als den phasischen Pseudokommutator
der Wellenfunktion $\psi$ an der Stelle $(z,t)\in \mathbf{P}_{1}\psi$
bezeichnen.\index{phasischer Pseudokommutator einer Wellenfunktion}
Die jeweiligen Werte des phasischen Pseudokommutators $\lbrack z,t\rbrack_{\psi}({\rm X})$ 
der Wellenfunktion $\psi$
unterscheiden sich zum Zweck seiner Normierung von den
Werten $\lbrack z,t\rbrack^{\psi}({\rm X})$ lediglich dadurch, dass die letzteren die 
symmetrische Differenz der Werte des phasischen Pseudokommutators $\lbrack z,t\rbrack_{\psi}({\rm X})$ und der um $t$ translierten
Teilmenge ${\rm X}+\{t\}\subset\mathbb{R}$ sind.
Weil in die translierten Teilmengen ${\rm X}+\{t\}$ keine Information "uber die jeweilige Wellenfunktion $\psi$ eingeht,
ist die Aufgabe, die jeweilige Wellenfunktion $\psi$ zu beschreiben, schon 
mit der Bestimmung der Werte der Abbildung $\lbrack z,t\rbrack^{\psi}$ getan.
Deshalb nennen wir exakt $\lbrack z,t\rbrack^{\psi}$
f"ur jede Wellenfunktion $\psi$
den\index{aktiver Pseudokommutator einer Wellenfunktion} 
aktiven Pseudokommutator der Wellenfunktion $\psi$ an der Stelle $(z,t)\in \mathbf{P}_{1}\psi$.
Offenbar gilt f"ur den phasischen Pseudokommutator
der Wellenfunktion $\psi$ die "Aquivalenz
\begin{equation}\label{sisicomx}
\begin{array}{c}
\Bigl(\forall\ {\rm X}\subset \mathbb{R},(z,t)\in \mathbf{P}_{1}\psi\quad\  \lbrack z,t\rbrack_{\psi}({\rm X})=\emptyset\Bigr)\\
\Leftrightarrow\\
\Bigl(\forall\ z\in \mathbf{P}_{2}\psi, a,b\in\mathbb{R}\quad \ \psi(\mbox{id},b)(\psi(z,a))=\psi(z,a+b)\Bigr)\ ,
\end{array}
\end{equation}
welche behauptet, dass der phasische Pseudokommutator
an jeder Stelle der Definitionsmenge $(z,t)\in \mathbf{P}_{1}\psi$ 
gem"ass der Identit"at 
$$\lbrack z,t\rbrack_{\psi}=\emptyset$$
genau dann verschwindet, wenn $\psi$ eine 
bis auf deren Zentriertheit\index{Zentriertheit}
additive Wellenfunktion im Sinne der 
rechten Seite dieser "Aquivalenz ist. Wir wollen unter einer 
additiven Wellenfunktion\index{additive Wellenfunktion} exakt eine solche Wellenfunktion $\alpha$
verstehen,
die
so beschaffen ist, dass 
f"ur alle $a,b\in \mathbb{R}$ und f"ur alle $z\in\mathbf{P}_{2}\alpha$
\begin{equation}\label{siomw}
\begin{array}{c}
\alpha(z,0)=z,\\  
\alpha(z,a+b)=\alpha(\alpha(z,a), b)=\alpha(\alpha(z,b), a)
\end{array}
\end{equation}
gilt. Uns interessieren aber insbesondere 
Flussfunktionen und gegenw"artig liegt unser Augenmerk nicht auf denjenigen
Flussfunktionen, welche additiv sind. Im Gegenteil, wir interessieren uns f"ur diejenigen 
Flussfunktionen $(\xi,\mathcal{A})$ im weiteren Sinn, die
innerhalb eines einzelnen Kollektivelementes $\tau\in[\xi]$ Formen
sensitiven Verhaltens zeigen, welche wir begrifflich erfassen wollen. Und zwar ohne, dass
wir uns dabei auf die Strukturierungen $\mathcal{A}$ der 
jeweiligen Flussfunktionen $(\xi,\mathcal{A})$ im weiteren Sinn beziehen.
Es ist aber
der gem"ass (\ref{sisiaomx})-(\ref{sisibomx}) eingef"uhrte aktive
Pseudokommutator, der uns weiterhilft. Mit der
Festlegung des Mengensystemes
\begin{equation}\label{asimove}
\mathcal{B}(1):=\Bigl\{]x-\varepsilon,x+\varepsilon[:x\in\mathbb{R},\varepsilon\in\mathbb{R}^{+}\Bigr\}
\end{equation}
aller offenen Kugeln des Zahlenstrahles
und mit der Festlegung, dass $\mbox{diam}({\rm Q})$ f"ur jede Teilmenge ${\rm Q}\subset\mathbb{R}$ des Zahlenstrahles
deren Durchmesser
\begin{equation}\label{siomv}
\mbox{diam}({\rm Q})=\sup\{|a-b|:a,b\in {\rm Q}\}
\end{equation}
bezeichnen soll,
verhilft uns der aktive Pseudokommutator
zu der folgenden  
\newline\newline
{\bf Definition 2.1: Intrinsische Sensitivit"at einer Flussfunktion}\index{intrinsische Sensitivit\"at}\newline 
{\em  F"ur jede Flussfunktion $\xi$ im weiteren Sinn nennen wir jede deren
kleinster Unterflussfunktionen} $\xi(x,\mbox{id})$ {\em f"ur $x\in \mathbf{P}_{2}\xi$
genau dann eine unfixiert
intrinsisch sensitive Unterflussfunktion\index{unfixiert intrinsisch sensitive Unterflussfunktion} dieser Flussfunktion $\xi$, wenn 
es einen Zustand $z\in \xi(x,\mathbb{R})$ 
gibt, f"ur den
es eine Folge reeller Zahlen $\{t_{j}\}_{j\in\mathbb{N}}$ und eine positive reelle Zahl $q$ gibt, f"ur die}
\begin{equation}\label{sisidomx}
\lim_{j\to\infty}
\sup\Bigl\{\mbox{diam}(Q):Q\in\mathcal{B}(1)\ \land\ Q\subset[z,t_{j}]^{\xi}(]-q,q[)\Bigr\}=0\\
\end{equation}
{\em ist.
Nur dann, wenn dabei die Folge reeller Zahlen $\{t_{j}\}_{j\in\mathbb{N}}$ H"aufungspunkte $t^{\star}$ hat,
bezeichnen wir die unfixiert
intrinsisch sensitive Unterflussfunktion} $\xi(x,\mbox{id})$ {\em als
eine fixiert intrinsisch sensitive Trajektorie.
Und zwar bezeichnen wir} $\xi(x,\mbox{id})$ {\em in diesem Fall als eine im Definitionsmengenelement 
$(z,t^{\star})$ fixiert intrinsisch sensitive Unterflussfunktion. Exakt jede Flussfunktion, die eine 
unfixiert bzw. eine fixiert
intrinsisch sensitive Unterflussfunktion hat, bezeichnen wir als 
eine unfixiert bzw. eine fixiert intrinsisch sensitive Flussfunktion.\index{intrinsisch sensitive Flussfunktion}   
\index{fixiert intrinsisch sensitive Flussfunktion}}
\newline\newline
Jede additive Flussfunktion hat demnach keine intrinsisch sensitiven Unterflussfunktionen,
w"ahrend eine additive Flussfunktion
im weiteren Sinn intrinsisch dann sensitive Unterflussfunktionen
hat, wenn sie kleinste, Unterflussfunktionen hat, die nicht invertierbar sind.
\newline
Ein jeweiliges Kollektivelement $\tau\in[\xi]$ ist f"ur sich genommen eine Punktmenge des
Zustandsraumes. Das Konzept der Sensitivit"at ist indess so angelegt, dass es nicht innerhalb des 
von der Parametrisierung der Kollektivelemente losgel"osten statischen Bildes besteht.
Sensitivit"at ist
wesentlich im Verh"altnis zu der Parametrisierung der Kollektivelemente konzipiert. Wir k"onnen
dennoch von der intrinsisch sensitiven Trajektorie $\tau\in[\xi]$ einer Flussfunktion im weiteren Sinn $\xi$
reden und damit den Sachverhalt meinen, dass eine derjenigen Unterflussfunktionen
$\xi(x,\mbox{id})$ intrinsisch sensitiv sei, deren
Zustandsraum die intrinsisch sensitive Trajektorie $\tau\in[\xi]$ ist.
Denn zustandsraumgleiche kleinste Unterflussfunktionen
haben das gleiche intrinsische Sensitivit"atsverhalten.\index{intrinsisch sensitive Trajektorie}
\newline
Wie so oft, wird aber auch hier der Bezug auf die Konstanten eines jeweiligen Kontextes innerhalb desselben  
nicht mehr ausdr"ucklich erw"ahnt: W"ahrend der Betrachtung einer bestimmten Flussfunktion ist
dann die Sprechweise, einfach von einer intrinsisch sensitiven Trajektorie zu reden, zu deuten
als die Bezeichnung einer intrinsisch sensitiven Unterflussfunktion der kontextuell
bestimmten Flussfunktion. Die "Okonomie dringt hier noch weiter vor, da es keine intrinsisch sensitive Unterflussfunktion
ausserhalb der Beziehung zu einer Flussfunktion im weiteren Sinn oder zu einem Phasenfluss gibt.
Es ist daher m"oglich, isoliert von einer intrinsisch sensitiven Trajektorie zu reden; in dem Wissen,
dass der Bezug zu einer Flussfunktion stillschweigend angenommen, indess lediglich nicht benannt ist.
\newline
Wir beschreiben das einfachste Szenario intrinsischer Sensitivit"at,
eine im Definitionsmengenelement $(z,t^{\star})\in\theta(x,\mathbb{R})\times\mathbb{R}$ fixiert 
intrinsisch sensitive Unterflussfunktion $\theta(x,\mbox{id})$ der Flussfunktion $\theta$,
welche als eine Flussfunktion auf der Trajektorie $\theta(x,\mathbb{R})$ invertierbar sein soll:
Es gibt eine
positive reelle Zahl $q$, f"ur die
das im Zahlenstrahl liegende Urbild
\begin{equation}
\begin{array}{c}
\theta(z,\mbox{id})^{-1}(\theta^{t^{\star}})^{-1}\theta(z,]-\varepsilon,\varepsilon[)\\
\supset ]-q,q[\not\supset\\
\theta(z,\mbox{id})^{-1}(\theta^{t^{\star}})^{-1}\theta(z,]-\varepsilon,\varepsilon[)
\end{array}
\end{equation}
f"ur alle $\varepsilon\in\mathbb{R}^{+}$ nicht in der Menge 
$]-q,q[$ konzentriert werden kann.
F"ur die reelle Funktion
$$\theta(z,\mbox{id})^{-1}(\theta^{t^{\star}})^{-1}\theta(z,\mbox{id}):\mathbb{R}\to\mathbb{R}$$
zeigt sich also bez"uglich der nat"urlichen Topologie des Zahlenstrahles 
im Ursprung sensitives Verhalten. Wie dieses 
sensitive Verhalten auf der Trajektorie $\theta(x,\mathbb{R})$ erscheint, h"angt
von der Strukturierung $\mathcal{A}$ des Zustandsraumes $\mathbf{P}_{2}\theta$
ab. Im Szenario nicht fixierter, jedoch unfixierter intrinsischer Sensitivit"at 
der Unterflussfunktion $\theta(x,\mbox{id})$, 
die invertierbar ist,
wissen wir hingegen lediglich, dass es eine h"aufungspunktfreie, mithin unbeschr"ankte Folge 
reeller Zahlen $\{t_{j}\}_{j\in\mathbb{N}}$ gibt, f"ur die
es zu jedem noch so kleinen $\varepsilon\in\mathbb{R}^{+}$ eine nat"urliche Zahl
$j(\varepsilon)$ gibt, f"ur welche f"ur alle Indizes $j>j(\varepsilon)$
\begin{equation}
\begin{array}{c}
\theta(z,\mbox{id})^{-1}(\theta^{t_{j}})^{-1}\theta(z,]-\varepsilon,\varepsilon[)\\
\supset ]-q,q[\not\supset\\
\theta(z,\mbox{id})^{-1}(\theta^{t_{j}})^{-1}\theta(z,]-\varepsilon,\varepsilon[)
\end{array}
\end{equation}
gilt.
F"ur keines der Folgenglieder der Folge reeller Funktionen
$$\{f_{j}\}_{j\in\mathbb{N}}=\Bigl\{\theta(z,\mbox{id})^{-1}(\theta^{t^{\star}})^{-1}\theta(z,\mbox{id})\Bigr\}_{j\in\mathbb{N}}$$
zeigt sich im Ursprung sensitives Verhalten. F"ur einen hinreichend hohen 
Index $j$ dehnt die Abbildung durch die Funktion $f_{j}$ jedoch jede noch so kleine Ursprungsumgebung 
"uber das Intervall $]-q,q[$ hin aus. Im Fall nicht fixierter, jedoch unfixierter intrinsischer Sensitivit"at
der Unterflussfunktion $\theta(x,\mbox{id})$ findet sich also zu jeder positiven Zahl $\varepsilon$
ein Paar zweier Zust"ande $\alpha,\beta\in\theta(x,\mathbb{R})$
und eine reelle Zahl $t\in \mathbb{R}$, f"ur welche die Ungleichung 
$$\mbox{diam}\Bigl(\theta(\beta,\mbox{id})^{-1}(\theta(\alpha,]t-\varepsilon,t+\varepsilon[))\Bigr)>2q$$
eingehalten ist. Wir 
k"urzen f"ur jede Flussfunktion $\xi$ im weiteren Sinn und f"ur alle deren Zust"ande $\alpha,\beta\in\mathbf{P}_{2}\xi$
und f"ur alle reelle Zahlen $t\in \mathbb{R}$
mit 
\begin{equation}
\xi[\beta,\alpha,t](\varepsilon):=\xi(\beta,\mbox{id})^{-1}(\xi(\alpha,]t-\varepsilon,t+\varepsilon[))
\end{equation}
ab und 
erkennen die folgende Phrasierung der Definition 2.1: Genau dann, wenn
die Unterflussfunktion $\xi(x,\mbox{id})$
einer Flussfunktion $\xi$ im weiteren Sinn fixiert
intrinsisch sensitiv ist,
gibt
es $\alpha,\beta\in\tau$
und eine reelle Zahl $t\in \mathbb{R}$,
sodass 
die 
Grenzwertaussage wahr ist, die das Nichtverschwinden des Limes superior 
\begin{equation}\label{oweheva}
\lim_{\varepsilon\to 0}\sup\{|a-b|:a,b\in\xi[\beta,\alpha,t](\varepsilon)\}>0
\end{equation}
behauptet, wie wir mit Blick auf 
die Festlegung 
(\ref{siomv}) sehen.
Ist $\xi(x,\mbox{id})$ hingegen unfixiert
intrinsisch sensitiv, so gilt die Grenzwertaussage,
dass eine reelle Zahl $t\in \mathbb{R}$ existiert, f"ur welche
\begin{equation}\label{owehevb}
\lim_{\varepsilon\to 0}\inf_{\alpha,\beta\in\tau}\sup\{|a-b|:a,b\in\xi[\beta,\alpha,t](\varepsilon)\}>0
\end{equation}
gilt.
Die Negationen sowohl der fixierten als auch der unfixierten Sensitivit"at sind die 
den Ungleichungen (\ref{oweheva}) und (\ref{owehevb})
entsprechenden homogenen Grenzwertaussagen, bei denen die jeweiligen linken Seiten mit null identifiziert sind.
Diese Negationen sind also ausdr"ucklich Kontinuit"atsbehauptungen, wie 
uns die Phrasierung der Definition 2.1 er"offnet, die 
uns ausserdem weiterleitet zu der folgenden
\newline\newline
{\bf Bemerkung 2.2: Bijektionsinvarianz der intrinsischen Sensitivit"at}\index{Bijektionsinvarianz der intrinsischen Sensitivit\"at}\newline
{\em Es sei $Q$
eine Bijektion des Zustandsraumes $\mathbf{P}_{2}\xi$
einer Flussfunktion im weiteren Sinn
auf eine beliebige Wertemenge $\mathbf{P}_{2}Q$ und} 
$Q\star\mbox{id}$ {\em der Operator, der auf der Menge aller Flussfunktionen im weiteren Sinn
definiert ist, deren Zustandsraum die Definitionsmenge $\mathbf{P}_{1}Q=\mathbf{P}_{2}\xi$ ist und dessen jeweilige Werte
gem"ass} (\ref{somx}) {\em festgelegt sind. Es gilt: Genau dann, wenn} $\xi(x,\mbox{id})$
{\em eine fixiert intrinsisch sensitive Unterflussfunktion bzw. eine 
unfixiert intrinsisch sensitive Unterflussfunktion
ist, ist}
\begin{equation}\label{owthevb}
Q\star\xi(Q(x),\mbox{id})
\end{equation}
{\em eine fixiert intrinsisch sensitive Unterflussfunktion bzw. eine 
unfixiert intrinsisch sensitive Unterflussfunktion.}
\newline\newline
{\bf Beweis:}\newline
F"ur alle Zust"ande $\alpha,\beta\in\mathbf{P}_{2}\xi$ der Flussfunktion
im weiteren Sinn $\xi$
und f"ur alle reelle Zahlen $t\in \mathbb{R}$ ist die f"ur das 
intrinsische Sensitivit"atsverhalten der transformierten Unterflussfunktion
$Q\star\xi(Q(x),\mbox{id})$
massgebliche, 
auf dem positiven Zahlenstrahl definierte
Funktion $Q\star\xi[Q(\beta),Q(\alpha),t]$ von der Bijektion $Q$ unabh"angig, da
f"ur alle $\varepsilon\in\mathbb{R}^{+}$
\begin{displaymath}
\begin{array}{c}
Q\star\xi[Q(\beta),Q(\alpha),t](\varepsilon)\\
=(Q(\xi(\beta,\mbox{id}))^{-1}(Q(\xi(\alpha,]t-\varepsilon,t+\varepsilon[)))\\
=\xi(\beta,\mbox{id})^{-1}(\xi(\alpha,]t-\varepsilon,t+\varepsilon[))
\end{array}
\end{displaymath}
von der Bijektion $Q$ unabh"angig ist.\newline 
{\bf q.e.d.}
\newline\newline 
Wir beschreiben eine fixiert intrinsisch sensitive Unterflussfunktion $\theta(x,\mbox{id})$ einer Flussfunktion 
$\theta$
in der modifizierten Perspektive dieser Paraphrase der Definition 2.1:
Auf $\theta(x,\mathbb{R})$ existiert ein Zustand
$\alpha\in\theta(x,\mathbb{R})$, zu dem es eine feste positive Konstante $\delta(\alpha,t)$
und eine bestimmte Phase $t\in\mathbb{R}$ gibt.
Unter den Entwicklungen dieses Zustandes $\alpha$ aus der Menge $\theta(\alpha,]t-\varepsilon,t+\varepsilon[)$ sind
dabei in jeder 
zeitlichen Umgebung $]t-\varepsilon,t+\varepsilon[$ 
solche zeitlich separablen Entwicklungen
$\theta(\beta,t_{1})\in \theta(\alpha,]t-\varepsilon,t+\varepsilon[)$ und 
$\theta(\beta,t_{2})\in \theta(\alpha,]t-\varepsilon,t+\varepsilon[)$ eines 
Zustandes $\beta$, f"ur die
$$|t_{1}-t_{2}|>\delta(\alpha,t)$$
gilt. Jedes nicht trivale Zimmer\index{nicht trivales Zimmer} $\chi\in[[\Phi]]\setminus[\Phi]$ einer $\mathcal{C}^{1}$-Flussfunktion $\Phi$
binnen
eines kompakten und endlichdimensionalen reellen Zustandsraumes 
beispielsweise
enth"alt 
keine intrinsisch sensitiven Trajektorien.
$\chi$ ist kompakt, sodass das Minimum
$$\min\{||\partial_{2}\Phi(x,0)||:x\in\chi\}\geq 0$$
existiert. W"are dasselbe null, so enthielte das Zimmer $\chi$ einen Fixpunkt,
denn die Unterflussfunktion $\Phi(x,\mbox{id})^{-1}$
existiert f"ur alle $x\in \mathbf{P}_{2}\Phi$, weil
$\Phi$ eine Flussfunktion sein soll.
Da demnach
$$\min\{||\partial_{2}\Phi(x,0)||:x\in\chi\}> 0$$
ist, ist f"ur alle Trajektorien
$\tau\in[\Phi]\cap\chi$ und f"ur alle
$\alpha,\beta\in\tau$
$$\mbox{diam}\Bigl(\Phi(\beta,\mbox{id})^{-1}(\Phi(\alpha,]t-\varepsilon,t+\varepsilon[))\Bigr)$$
$$>\min\{||\partial_{2}\Phi(x,0)||:x\in\chi\}\ .$$
Die intrinsische Sensitivit"at einer Flussfunktion ist hier also als die Negation einer
Form deren Stetigkeit formuliert, wie uns deutlich wird, wenn wir die Verneinung
der intrinsischen Sensitivit"at des 
Kollektivelementes $\theta(x,\mathbb{R})\in[\theta]$
Flussfunktion $\theta$ als den Sachverhalt phrasieren, dass f"ur alle Zust"ande
$\alpha,\beta\in \xi(x,\mathbb{R})$ und f"ur alle $t\in\mathbb{R}$ die 
Grenzwertaussage (\ref{oweheva}) 
gilt. Was wir im Bild dieser Paraphrase der Definition 2.1 auch deutlich sehen, ist,
dass die blosse Nichtinvertierbarkeit 
einer Unterflussfunktion $\xi(x,\mbox{id})$
einer Flussfunktion $\xi$ im weiteren Sinn
allemal die fixierte intrinsische 
Sensitivit"at der Unterflussfunktion $\xi(x,\mbox{id})$ impliziert;
nicht nur dann,
wenn es ein offenes Intervall $J\subset \mathbb{R}$ gibt, auf dem
$\xi(x,\mbox{id})$ konstant ist.\newline 
Dadurch, dass wir in der Definition 2.1 Flussfunktion im weiteren Sinn zulassen,
nehmen wir in Kauf, dass die ledigliche Nichtinvertierbarkeit 
einer Unterflussfunktion als deren fixierte intrinsische 
Sensitivit"at erscheint. F"ur den Begriff der fixierten intrinsischen 
Sensitivit"at ergibt sich damit die eingangs erw"ahnte Kontamination
des Ph"anomens der Sensitivit"at im Sinne der Diskontinuit"at gem"ass
(\ref{oweheva}) und der Nichtinvertierbarkeit.
\newline
Wir betonen hierbei, dass der Formulierung der 
intrinsischen Sensitivit"at losgel"ost von der Kontinuit"at einer 
jeweiligen Flussfunktion im Bezug auf eine 
jeweilige Zustandsraumstrukturierung angelegt ist.
In der Losl"osung von der Kontinuit"at externer Referenz
auf eine Zustandsraumstrukturierung, 
die der Begriff der 
intrinsischen Sensitivit"at hat,
besteht die ihn benennende Intrinsik.

\section{An der Generalisierbarkeitsgrenze des metrischen Sensitivit"atsbegriffes
d"ammert die Idee der Homogenit"at.}
Sollte 
es m"oglich sein, den Sensitivit"atsbegriff alleine auf der Grundlage einer gegebenen Topologie des Zustandsraumes 
auf eine elegante Weise abstrakter zu definieren,
so m"usste die Sensitivit"at vermutlich auf die Weise anders als in der Definition 1.2 formulierbar sein,
dass jener alternative Zugang zum Sensitivit"atsbegriff die Sensitivit"at nicht mit Hilfe von Aufl"osungsfeldern
als globale
punktweise Sensitivit"at fasst.
Wollen wir hierbei nicht von Anfang an punktweise Sensitivit"at
formulieren, so m"ussen wir dann vermutlich von der integralen Gegebenheit eines jeweiligen sensitiven 
Attraktors ausgehen und damit von einer 
Form der Homogenit"at desselben.
Speziell dieser vielleicht allzu naheliegend anmutende Gedanke, 
dass eine Form jener Homogenit"at der sensitiven Attraktoren einer Flussfunktion
als die nicht singul"are Konstanz des Aufl"osungsfeldes derselben greifbar
sein k"onnte, erwiese sich aber als naiv, wenn 
sich
sensitive Attraktoren angeben liessen, deren Aufl"osungsfeld nicht konstant ist.\newline
Der Einwand, dass die Annahme der nicht singul"aren Konstanz
des Aufl"osungsfeldes auf sensitiven Attraktoren
doch durch die einfache Ab"anderung der Metrik $d$
des Zustandsraumes einer jeweiligen Flussfunktion $\Psi$, deren 
Aufl"osungsfeld $\Delta_{\Psi}^{d}$ auf einem sensitiven Attraktor $\chi\in[[\Psi]]_{\mathbb{T}(d)}$
konstant ist, widerlegt werde, macht 
die Idee der attraktorweisen Konstanz des Aufl"osungsfeldes
auf den ersten Blick zunichte:
Wenn nicht f"ur jede Flussfunktion $\Psi$ und f"ur jeden der sensitiven Attraktoren $\chi\in[[\Psi]]_{\mathbb{T}(d)}\setminus[\Psi]$
derselben und f"ur
alle Metriken $d$ und $\tilde{d}$ des Zustandsraumes $\mathbf{P}_{2}\Psi$ die Invarianz der
Urbildmengen des positiven Zahlenstrahles $\mathbb{R}^{+}$ jeweiliger Aufl"osungsfelder
\begin{equation}\label{alizba}
(\Delta_{\Psi}^{d})^{-1}(\mathbb{R}^{+})=(\Delta_{\Psi}^{\tilde{d}})^{-1}(\mathbb{R}^{+})
\end{equation}
gilt, dann ist der Gedanke der attraktorweisen Konstanz des Aufl"osungsfeldes falsch.
\newline
Wenn hingegen diese "ubersichtliche Gegebenheit der attraktorweisen Konstanz des Aufl"osungsfeldes
f"ur jede 
bez"uglich $d$ stetige
Flussfunktion $\Psi$ binnen des Kompaktums $\mathbf{P}_{2}\Psi$ des 
kompakten
metrischen Raumes
$(\mathbf{P}_{2}\Psi,d)$
gelten w"urde, so w"are gerade das Mengensystem
\begin{equation}\label{alizbb}
(\Delta_{\Psi}^{d})^{-1}(\mathbb{R}^{+})=[[\Psi]]_{\mathbb{T}(d)}\setminus[\Psi]
\end{equation}
die Menge aller sensitiven Attraktoren der Flussfunktion $\Psi$ und
\begin{equation}\label{alizbb}
(\Delta_{\Psi}^{d})^{-1}(\{-\infty\})=\bigcup([[\Psi]]_{\mathbb{T}(d)}\cap[\Psi])
\end{equation}
die Menge aller Zust"ande des Zustandsraumes $\mathbf{P}_{2}\Psi$, die nicht sensitiv sind. Das Mengensystem
$$(\Delta_{\Psi}^{d})^{-1}(\mathbb{R}^{+})\cup\{(\Delta_{\Psi}^{d})^{-1}(\{-\infty\})\}\in\mathbf{part}(\mathbf{P}_{2}\Psi)$$
w"are eine Partition des Zustandsraumes $\mathbf{P}_{2}\Psi$.
Dabei wandelte sich also die darstellerische Funktion des Aufl"osungsfeldes $\Delta_{\Psi}^{d}$.
Es beschriebe dann nicht mehr eine Differenzierung zwischen Zust"anden innerhalb eines 
sensitiven Attraktors, sondern zwischen Zust"anden verschiedener Attraktoren und 
die 
gem"ass $(1.12)^{1}$ im ersten Teil \cite{rabe} konstruierte
Niveaulinientopologie $\mathbf{T}_{N}(\Delta_{\Psi}^{d})$ des Aufl"osungsfeldes $\Delta_{\Psi}^{d}$
\begin{equation}\label{aulinia}
\mathbf{T}_{N}(\Delta_{\Psi}^{d})=([[\Psi]]_{\mathbb{T}(d)}\setminus[\Psi])^{\cup}
\end{equation}
fiele mit der selbstdualen Topologie $([[\Psi]]_{\mathbb{T}(d)}\setminus[\Psi])^{\cup}$
der sensitiven Attraktoren in eins.
\newline 
Wie wollen wir aber jenen Einwand entkr"aften, dass die   
attraktorweise Konstanz des Aufl"osungsfeldes verworfen werde durch die 
von der Vorgabe einer jeweiligen Flussfunktion unabh"angige
Abwandlungsfreiheit 
der Metrisierung des Zustandsraumes;
und wie sollen wir uns stattdessen das Zustandekommen der Invarianz der
Urbildmengen des positiven Zahlenstrahles $\mathbb{R}^{+}$ jeweiliger Aufl"osungsfelder
gegen"uber der Wahl der Metrisierung des Zustandsraumes gem"ass (\ref{alizba}) vorstellen?
\newline
Stellen wir jenen Einwand vor den Hintergund der
Vermutung maximaler Sensitivit"at, welche in der
Darstellung des elementaren Quasiergodensatzes \cite{raab}
im Abschnitt 2.3 ge"aussert wurde:
\newline
\newline
\glqq {\em Dar"uber hinaus k"onnte 
der Aufl"osungsfeldwert jedes Zimmers\index{Zimmer} $\chi$ als dessen jeweilige Aufl"osungskonstante 
$\Delta(\chi)$
in dem Sinn existieren, dass 
das 
Aufl"osungsfeld\index{Aufl\"osungsfeld einer Wellenfunktion} 
der auf das jeweilige Zimmer $\chi$ restringierten Wellenfunktion $\Psi|\chi\times\mathbb{R}$
konstant ist: Es nimmt den singul"aren Wert $-\infty$ f"ur jedes triviale Zimmer an, keiner
dessen Zust"ande
sensitiv ist.
F"ur 
jedes nicht triviale Zimmer $\chi$ hingegen k"onnte als jeweilige 
Aufl"osungskonstante\index{Aufl\"osungskonstante eines Zimmers}
f"ur alle $x\in\chi$
die Zahl
\begin{equation}
\Delta(\chi):=\Delta_{\Psi|\chi\times\mathbb{R}}(x)=||\chi||\in \mathbb{R}^{+}\ ,
\end{equation} 
vorliegen und einfach mit 
dem Durchmesser des betrachteten nicht trivialen Zimmers $\chi$ "ubereinstimmen.}\grqq\newline\newline
Wir redeten dabei von Zimmern, die in endlichdimensionalen und reellen Zustandsr"aumen liegen.
Diese Zimmer sind diejenigen Abschl"usse $\mathbf{cl}(\tau)$ der 
Kollektivierungselemente $\tau\in[\Psi]$ f"ur
reelle Flussfunktionen $\Psi$ im weiteren Sinn, welche Mannigfaltigkeiten 
sind. Genau dann, wenn $\mathbf{cl}(\tau)\setminus\tau\not=\emptyset$ ist,
ist $\mathbf{cl}(\tau)$ ein nicht triviales Zimmer.
Nehmen wir den Fall, dass es Kollektivierungselemente $\tau\in[\Psi]$
vorliegen, die stetige Bilder offener oder halboffener Intervalle sind. 
Es ist dann zwar $\mathbf{cl}(\tau)\setminus\tau\not=\emptyset$. 
$\mathbf{cl}(\tau)$ ist dann aber keine Mannigfaltigkeit.
\newline
Die
Verallgemeinerung 
der 
Eigenschaft nicht trivialer reeller Zimmer, die im Hinblick auf die
Sensitivit"at relevant ist,\index{nicht triviale Zimmer} ist die Eigenschaft 
derjeniger
Zustandsraumteilmengen, die
von Poincar$\acute{e}$schen Entwicklungen\index{Poincar$\acute{e}$sche Entwicklung} erzeugt werden: 
\newline
\newline
{\bf Definition 2.3:\newline Poincar$\mathbf{\acute{e}}$sche und global Poincar$\mathbf{\acute{e}}$sche Entwicklungen}\newline
{\em Als eine Poincar$\acute{e}$sche Entwicklung bezeichnen wir exakt jede 
strukturierte Entwicklung $(\xi,\mathcal{A})$, f"ur welche die Aussage} 
\begin{equation}\label{thisbec}
\begin{array}{c}
\forall\ z\in \mathbf{P}_{2}\xi,{\rm L}\in\mathcal{A}_{\xi(z,\mathbb{R})}, t_{\star}\in\mathbb{R}\ \exists\ t\in ]-\infty,t_{\star}[\cup]t_{\star},\infty[\ :\\
\xi(z,t)\in {\rm L}
\end{array}
\end{equation}
{\em wahr ist. Jede Poincar$\acute{e}$sche Entwicklung $(\xi,\mathcal{A})$ nennen wir 
genau dann global Poincar$\acute{e}$sch, wenn "uber} (\ref{thisbec}) {\em hinaus f"ur jede Lokalisierung ${\rm L}\in\mathcal{A}$ die Aussage} 
\begin{equation}\label{thisbic}
\begin{array}{c}
\forall t_{\star}\in\mathbb{R}\ \exists\ t\in ]-\infty,t_{\star}[\cup]t_{\star},\infty[\ :\\
\xi(z,t)\in {\rm L}
\end{array}
\end{equation}
{\em gilt.}\index{global Poincar$\acute{e}$sche Entwicklung}
\newline
\newline
Weshalb wir Poincar$\acute{e}$sche Entwicklungen nach dem Pionier Henri Poincar$\acute{e}$\index{Poincar$\mathbf{\acute{e}}$, Henri}
benennen, ist evident. Es war Henri Poincar$\acute{e}$, der Ludwig Boltzmanns Versuch,
Irreversibilit"at aus der statistischen Mechanik abzuleiten, den ber"uhmt gewordenen Wiederkehreinwand\index{Wiederkehreinwand}
entgegenhielt.
Wenn wir von Poincar$\acute{e}$schen Entwicklungen sprechen,
ist es hierbei allerdings kaum m"oglich, von Entwicklungen nur in dem 
in der Definition 2.3
objektivierten Sinn zu reden; in dem Sinn, dass
das Wort \glqq Entwicklung\grqq aussschliesslich die
Objekte 
der Klasse der Entwicklungen meint, welche die Klasse der Flussfunktionen
umfasst und dabei von der Klasse der Flussfunktionen im weiteren Sinn umfasst ist.
Wir sagen z.B. auch, dass 
jede Entwicklung eines Zustandes $z\in \mathbf{P}_{2}\xi$ einer Poincar$\acute{e}$schen Entwicklung $(\xi,\mathcal{A})$
gem"ass (\ref{thisbec})
so beschaffen ist, dass sie jede Lokalisierung ${\rm L}\in\mathcal{A}$ immer wieder 
trifft, falls die Entwicklung des Zustandes $z\in \mathbf{P}_{2}\xi$ einmal in die 
Lokalisierung ${\rm L}\in\mathcal{A}$ f"uhrte; was genau dann der Fall ist, wenn 
${\rm L}$ ein Element der Auswahl $\mathcal{A}_{\xi(z,\mathbb{R})}$ ist. F"ur jede global Poincar$\acute{e}$sche 
Entwicklung $(\xi,\mathcal{A})$ ist jede Entwicklung eines Zustandes $z\in \mathbf{P}_{2}\xi$ gem"ass (\ref{thisbic}) 
so beschaffen ist, dass sie sogar jede Lokalisierung ${\rm L}\in\mathcal{A}$ immer wieder 
trifft. Und damit hat
eine jeweilige Zustandsentwicklung
einer global Poincar$\acute{e}$schen 
Entwicklung mit einer Zustandsentwicklung innerhalb eines reellen Zimmers deren Charakteristikum gemein.
\newline
Reelle Zimmer\index{nicht triviale Zimmer} $\chi\in[[\Psi]]$ sind
durch die Restriktion der jeweiligen reellen Flussfunktion $\Psi|(\chi\times \mathbb{R})$ erzeugt, welche
der Entwicklung $$(\Psi|(\chi\times \mathbb{R}),\mathbf{T}(\mathbf{P}_{2}\Psi)\cap\chi)$$
entspricht, wobei $\mathbf{T}(\mathbf{P}_{2}\Psi)$ die nat"urliche Zustandsraumtopologie des endlichdimensionalen Zustandsraumes
benenne. Gleich, ob $\chi\in[\Psi]$ trivial oder $\chi\not\in[\Psi]$
nicht trivial ist:
In jedem Fall ist $(\Psi|(\chi\times \mathbb{R}),\mathbf{T}(\mathbf{P}_{2}\Psi)\cap\chi)$ eine 
global
Poincar$\acute{e}$sche Entwicklung.
Der Unterschied zwischen 
nicht trivialen reellen Zimmern $\chi\in[[\Psi]]\setminus[\Psi]$ und trivialen reellen Zimmern $\chi\in[\Psi]$ zeigt sich 
darin, dass f"ur $\chi\in[[\Psi]]\setminus[\Psi]$ die Aussage
\begin{equation}
\begin{array}{c}
\forall\ {\rm U}\in\mathbf{T}(\mathbf{P}_{2}\Psi)\cap\chi, z\in \chi\\
\mathbf{card}\Bigl(\Bigl\{{\rm X}\subset\mbox{con}(\mathbb{R}):\bigcup{\rm X}=\xi(z,\mbox{id})^{-1}({\rm U})\Bigr\}\Bigr)=\aleph_{0}
\end{array}
\end{equation}
wahr ist, die falls $\chi\in[\Psi]$ trivial ist, falsch ist.
Wobei $\aleph_{0}=\mathbf{card}(\mathbb{N})$ die Kardinalit"at der abz"ahlbar unendlichen Mengen sei und
$\mbox{con}(\mathbb{R})$ das
Mengensystem aller bez"uglich der nat"urlichen Zahlenstrahltopologie einfach zusammenh"angenden und
mehrelementigen
Teilmengen des Zahlenstrahles.
$\mbox{con}(\mathbb{R})$ ist demnach die Menge aller
mehrelementigen 
Intervalle des Zahlenstrahles. Die Intervalle in $\mbox{con}(\mathbb{R})$ k"onnen beidseitig offen bzw. geschlossen, halbseitig offen bzw. geschlossen
sein, nicht aber punktf"ormig oder leer. $\mbox{con}(\mathbb{R})$ ist damit die Menge aller 
nicht punktf"ormigen
Konvexa des Zahlenstrahles. 
Wir bemerken hierzu: F"ur alle durch eine Vereinigung 
$\bigcup{\rm X}={\rm Y}$
"uber eine Teilmenge ${\rm X}\subset\mbox{con}(\mathbb{R})$
darstellbaren Teilmengen ${\rm Y}\subset\mathbb{R}$ ist die 
M"achtigkeit $\mathbf{card}({\rm X})<\aleph_{1}$ unterhalb
der Kardinalit"at des Zahlenstrahles $\aleph_{1}$.
\newline
Im Fall einer reellen Flussfunktion $\Psi$
w"are der Sachverhalt 
der attraktorweisen Konstanz des Aufl"osungsfeldes
mit der
Vermutung maximaler Sensitivit"at\index{Vermutung maximaler Sensitivit\"at}
in Einklang zu bringen. Die maximale Sensitivit"at
w"are ein Spezialfall der attraktorweisen Konstanz des Aufl"osungsfeldes.
Dabei w"are die maximale Sensitivit"at aber nicht nur 
als Spezialfall mit der 
attraktorweisen Konstanz des Aufl"osungsfeldes
in Einklang bringbar, sondern 
die maximale Sensitivit"at zeigte uns dabei auch,
wie wir uns das Zustandekommen der Invarianz der
Urbildmengen des positiven Zahlenstrahles $\mathbb{R}^{+}$ jeweiliger Aufl"osungsfelder
gegen"uber der Wahl der Metrisierung des Zustandsraumes zu denken
haben. Sogleich sehen wir nun, dass
der mit der 
Abwandlungsfreiheit 
der Metrisierung des Zustandsraumes argumentierende
Einwand die Idee der attraktorweisen Konstanz des Aufl"osungsfeldes
nicht zunichte macht: \newline 
Wenn die vermutete maximale Sensitivit"at
gilt,
dann "andert sich mit der Abwandlung der Metrik des Zustandsraumes $d$ zu $\tilde{d}$
im allgemeinen 
zwar der Wert des Durchmessers eines jeweiligen Attraktors. 
Er ist aber
f"ur alle Zust"ande $z$ desselben auch nach jeder Abwandlung der Metrik
gleich und zwar der Wert $\Delta_{\Psi}^{\tilde{d}}(z)$ des Aufl"osungsfeldes
an der Stelle $z$. Vor der Leinwand der Vorstellung und Vermutung maximaler Sensitivit"at
tritt die Nicht-Stichhaltigkeit des mit der Abwandlungsfreiheit der Zustandsraummetrisierung argumentierenden Einwandes
gegen die nicht singul"are Konstanz des Aufl"osungsfeldes  
hervor.
\newline 
Die 
maximale Sensitivit"at endlichdimensionaler und reeller Zimmer, die Konstanz deren Aufl"osungsfeldes, das
den maximal m"oglichen Wert annimmt,
ist der gleichermassen markante wie kompakte Sachverhalt, 
der eine Form der
Homogenit"at 
sensitiver Attraktoren
anzeigt, die wir hier in ihrer F"ulle innerhalb des thematischen Fokus 
dieses Traktates nicht
beschreiben.\footnote{Die diversen
Auspr"agungen der Homogenit"at oder der topologischen Indifferenz
sensitiver Attraktoren
bleiben allerdings in diesem Traktat in ihrer Denkbarkeit nach noch unbestimmt.
Die topologische Indifferenz sensitiver Attraktoren stellen wir hier
ebenso wie die probabilistische Beschaffenheit eines sensitiven Attraktors
nicht dar: Beide 
Themen haben n"amlich eine Vielzahl von Gesichtspunkten. Sie bed"urfen
eigener Abhandlungen.}
Wir w"urden die Vermutung maximaler Sensitivit"at freilich gerne
beweisen.\newline
Gelte uns dieser Beweis als Ziel, gar als das Hauptanliegen dieses Traktates?
Nein, das Hauptanliegen ist dieser Beweis ausdr"ucklich nicht! Wir haben uns vorgenommen, allgemein zu konzipieren. 
Und zwar die Sensitivit"at kontinuierlicher dynamische Systeme:
Wir objektivieren
allgemeine Sichtweisen, indem wir Begriffe verfassen, d.h., Objektklassen formulieren,
die es uns erlauben, die Sensitivit"at kontinuierlicher dynamischer Systeme 
auf der Grundlage beliebiger Zustandsraumstrukturierungen  
zu beschreiben. Dabei stossen wir auch auf Bemerkungen und S"atze.
Das Ziel, die vermutete
maximale Sensitivit"at zu beweisen, sei uns aber eines der Gestirne, an denen wir uns dabei orientieren.
In dieser Arbeit erreichen wir dieses Ziel noch nicht.
\chapter{Generalisierte Sensitivit"at}
{\small Nachdem sich im ersten Kapitel unsere Vorstellung 
in gleichsam einer perspektivischen Systole 
sehr auf die Spannung 
einliess, die mit der Ungewissheit verbunden ist, 
ob es "uberhaupt einen Sensitivit"atsbegriff gibt, der "uber die Allgemeinheit verallgemeinerter metrischer 
Sensitivit"at hinaus verfassbar ist,  
l"osen wir diese Spannung nun auf.
Im Zuge der anschliessenden Diastole geht 
die Perspektive auf, in der die allgemeine Weite des 
generalisierten Sensitivit"atsbegriffes sichtbar wird.
Wir f"uhren im insofern diastolischen Auftakt dieses Kapitels sogleich die  
fixierte Sensitivit"at f"ur beliebige strukturierte Flussfunktionen ein, 
die eine Verallgemeinerung des Sensitivit"atsbegriffes 
"uber den Allgemeinheitsgrad der allgemeinen Topologie hinaus ist.
Wir bemerken den Pluralismus der Sensitivit"at, welcher dem Pluralismus der Formen 
des Chaos entspricht und finden, dass wir uns daran orientieren k"onnen, dass Chaos\index{Chaos}
eine intermittente Negation ist.\index{intermittente Negation}\newline
Eine Reihe
verallgemeinerter Sensitivit"atsbegriffe, die uns alle die Verallgemeinerung des Sensitivit"atsbegriffes 
"uber den Allgemeinheitsgrad der allgemeinen Topologie hinaus an die Hand
geben, f"ugt sich in die Hierarchie\index{Hierarchie der Sensitivit\"at}
der Sensitivit"at, die wir im dritten Abschnitt vorstellen. Diese Sensitivit"atsbegriffe
sind mengenweise und keine zustandsweise verfassten Sensitivit"atsbegriffe,\index{mengenweise Sensitivit\"at} wie es die 
metrische Sensitivit"at\index{metrische Sensitivit\"at} ist.\index{zustandsweise Sensitivit\"at}
Diese Hierarchie ist offenbar das zentrale Konzept, 
das den Pluralismus der Sensitivit"at\index{Pluralismus der Sensitivit\"at} zusammenh"alt, der
ph"anomenell ist, wo dies Konzept fundamental ist.
\newline 
Dann,
im vierten  
Abschnitt dieses Kapitels pr"asentieren wir
die Ultrakolokalisation\index{Ultrakolokalisation} in Gestalt einer 
zu jeder jeweiligen strukturierten Flussfunktion $\Lambda=(\xi, \mathcal{A})$ im weiteren Sinn gegebenen
Relation 
$$\sim_{\Lambda}\ \subset\ \mathbf{P}_{2}\xi\times \mathbf{P}_{2}\xi\ ,$$ die
zwischen zwei jeweiligen Zust"anden des Zustandsraumes $\mathbf{P}_{2}\xi$ besteht. Die Ultrakolokalisation 
vermittelt zwar
eine spezielle Form der Sensitivit"at, jawohl. Indess, diese ist eine spezielle Form der 
verallgemeinerten Sensitivit"at. Denn die Ultrakolokalisation ist dabei f"ur 
beliebige Paare $\Lambda=(\xi,\mathcal{A})$ formuliert, deren erste Komponente eine Flussfunktion im weiteren Sinn 
ist und deren zweite Komponente ein 
den Zustandsraum $\mathbf{P}_{2}\xi$
"uberdeckendes Mengensystem $\mathcal{A}\subset 2^{\mathbf{P}_{2}\xi}$ ist. Daher ist 
die von der Ultrakolokalisation 
abgeleitete Sensitivit"at
eine in h"ochster Allgemeinheit verfasste Spielart der verallgemeinerten Sensitivit"at.
Die von der Ultrakolokalisation 
abgeleitete Sensitivit"at
ist eine Form 
zustandsweiser Sensitivit"at.
"Uber die Abschw"achung der Ultrakolokalisation kommen wir zur Kolokalisation,\index{Kolokalisation} welche es erm"oglicht,
die verallgemeinerte zustandsweise Sensitivit"at zu verfassen.
F"ur Poincar$\acute{e}$sche Entwicklungen\index{Poincar$\acute{e}$sche Entwicklung} gelingt es uns schliesslich, diese 
verallgemeinerte zustandsweise Sensitivit"at mit der mengenweisen Sensitivit"at zu identifizieren.}
\section{Die fixierte Sensitivit"at strukturierter\\ Flussfunktionen}
Unsere Jagd auf Ignoranz setzt als die Jagd auf das 
eigene Dunkel an. 
Wir stellen uns der Unklarheit, um mit der damit verbundenen
Exposition die Grundlage desjenigen Kontrastes zu schaffen, der es erm"oglicht, grosse Deutlichkeit
herauszuarbeiten.
Wir stellen uns also, sofern es dessen "uberhaupt bedarf, dumm,
indem wir ebenso fragen, um uns 
Sachverhalte mit gr"osstm"oglicher Klarheit
vor Augen zu f"uhren.\newline
Aber haben wir uns nun eigentlich dumm gestellt und aus
der Vorstellungsklamm heraus, dass der Begriff der Sensitivit"at lediglich mit dem 
Begriff der Metrik kogeneralisierbar sein k"onnte,
dumm gefragt?
L"asst sich die Frage, ob der Begriff der Sensitivit"at lediglich mit dem 
Begriff der Metrik kogeneralisierbar sei, nicht leicht 
konstruktiv 
verneinen, wenn
wir bemerken, dass wir die Generalisierung des Sensitivit"atsbegriffes im wesentlichen schon 
im ersten Teil der Konzepte der abstrakten Ergodentheorie verfassten und
untersuchten? An welcher Stelle?
Wo wir die mit der Cantor-Stetigkeit die
Negation der Cantor-Stetigkeit\index{Cantor-Stetigkeit} eines jeweiligen
Phasenflusses $\{\xi^{t}\}_{t\in\mathbb{R}}$ bez"uglich einer 
Zustandsraum"uberdeckung $\mathcal{A}$ fassen, die wir
im dritten Kapitel des ersten Teiles der Konzepte der abstrakten 
Ergodentheorie formulierten. F"ur jedes Mengensystem $\mathcal{Z}$ bezeichne
\begin{equation}\label{blzb}
\mathcal{Z}^{\not 2}:=\{({\rm A},{\rm B})\in\mathcal{Z}\times\mathcal{Z}:{\rm A}\cap{\rm B}=\emptyset\}
\end{equation}
die Menge aller Paare disjunkter Mengen des Mengensystemes $\mathcal{Z}$, sodass
$$\mathbf{P}_{1}\mathcal{Z}^{\not 2}=\mathbf{P}_{2}\mathcal{Z}^{\not 2}$$
das Mengensystem aller Mengen ${\rm A}$ des Mengensystemes $\mathcal{Z}$ ist,
zu denen es Mengen ${\rm B}$ desselben Mengensystemes $\mathcal{Z}$ gibt, welche
zu ${\rm A}$ disjunkt sind. $\mathbf{P}_{1}\mathcal{Z}^{\not 2}$ ist nicht mit einem
in $\mathcal{Z}$ enthaltenen Mengensystem lauter paarweise disjunkter Mengen zu verwechseln,
ferner differiert
$$\mathbf{P}_{1}\mathcal{Z}^{\not 2}\times\mathbf{P}_{2}\mathcal{Z}^{\not 2}
\not=\mathcal{Z}^{\not 2}\ .$$
"Uber jedes Mengensystem $\mathcal{Z}$
sagen wir dabei genau dann, dass es
die verallgemeinerte Hausdorffsche Trennungseigenschaft\index{verallgemeinerte Hausdorffsche Trennungseigenschaft}
habe, wenn f"ur alle $$(x,y)\in\bigcup\mathcal{Z}\times\bigcup\mathcal{Z}$$
die "Aquivalenz 
\begin{equation}\label{bluzb}
\begin{array}{c}
\exists({\rm X},{\rm Y})\in \mathcal{Z}^{\not 2}:\ x\in{\rm X}\land y\in{\rm Y}\\
\Leftrightarrow x\not=y
\end{array}
\end{equation}
wahr ist.
Jede Distanzfunktion 
\begin{displaymath}
\begin{array}{c}
d^{\circ}:\mathbf{P}_{1}\mathcal{Z}^{\not 2}\times\mathbf{P}_{2}\mathcal{Z}^{\not 2}\to
[0,\infty[\ ,\\
({\rm A},{\rm B})\mapsto d^{\circ}({\rm A},{\rm B})
\end{array}
\end{displaymath}
ist hierbei genau dann eine auf dem Mengensystem $\mathbf{P}_{1}\mathcal{Z}^{\not 2}$
definierte Metrik,\index{Metrik} wenn
f"ur alle ${\rm A},{\rm B},{\rm C}\in \mathbf{P}_{1}\mathcal{Z}^{\not 2}$
\begin{displaymath}
\begin{array}{c}
d^{\circ}({\rm A},{\rm B})=0\Leftrightarrow{\rm A}={\rm B}\ ,\\
d^{\circ}({\rm A},{\rm B})=d^{\circ}({\rm B},{\rm A})\ ,\\
d^{\circ}({\rm A},{\rm B})+d^{\circ}({\rm B},{\rm C})\geq d^{\circ}({\rm A},{\rm C})
\end{array}
\end{displaymath}
gilt.
\newline
\newline
{\bf Definition 3.1: Fixierte Sensitivit"at\index{fixierte Sensitivit\"at}}\newline  
{\em Wir nennen jede strukturierte Flussfunktion im weiteren Sinn $(\xi,\mathcal{A})$ genau dann im Zustand 
$\alpha\in\mathbf{P}_{2}\xi$ fixiert
sensitiv, wenn\index{fixierte Sensitivit\"at in einem Zustand}
die Aussage 
\begin{equation}\label{blzba}
\begin{array}{c}
\exists\ ({\rm X},{\rm Y})\in\mathcal{A}^{\not 2}:\ 
\forall\ {\rm A}\in\mathcal{A}_{\{\alpha\}}, t_{\star}\in\mathbb{R}\\ 
\exists\ (x,y)\in{\rm A}\times{\rm A} ,\ t\in ]-\infty,t_{\star}[\cup]t_{\star},\infty[\ :\\
\xi(x,t)\in{\rm X}\ \land\ \xi(y,t)\in{\rm Y}
\end{array}
\end{equation}
wahr ist. 
Genau dann, wenn es einen Zustand $x\in\mathbf{P}_{2}\xi$ gibt, in dem $(\xi,\mathcal{A})$ 
fixiert sensitiv ist,
nennen wir die strukturierte Flussfunktion im weiteren Sinn $(\xi,\mathcal{A})$ fixiert sensitiv.}
\newline
\newline
Das Verh"altnis Cantor-stetiger Flussfunktionen im weiteren Sinn $(\xi,\mathcal{A})$
zur fixierten Sensitivit"at
ist so einfach und doch Orientierung gebend, wie der Beweis der folgenden 
\newline
\newline
{\bf Bemerkung 3.2:}\newline
{\em Jede
Cantor-stetige Flussfunktion im weiteren Sinn $(\xi,\mathcal{A})$ ist nirgends
von fixiert
sensitiv.}
\newline
\newline
{\bf Beweis:}\newline
Wenn f"ur den 
Zustand 
$\alpha\in\mathbf{P}_{2}\xi$ die Aussage
(\ref{blzba}) gilt, trifft
es zu, dass
\begin{equation}\label{balizbb}
\begin{array}{c}
\exists\ {\rm B}\in\mathcal{A}:\ 
\forall\ {\rm A}\in\mathcal{A}_{\{\alpha\}}, t_{\star}\in\mathbb{R}\ \exists\ t\in ]-\infty,t_{\star}[\cup]t_{\star},\infty[\ :\\
\xi({\rm A},t)\not\subset{\rm B}
\end{array}
\end{equation}
wahr ist. Mithin trifft
es erst recht zu, dass f"ur den 
Zustand $\alpha\in\mathbf{P}_{2}\xi$ die Gegebenheit 
\begin{equation}\label{blizbb}
\begin{array}{c}
\exists\ {\rm B}\in\mathcal{A}:\ 
\forall\ {\rm A}\in\mathcal{A}_{\{\alpha\}}, \exists\ t\in\mathbb{R} \ :\\
\xi({\rm A},t)\not\subset{\rm B}
\end{array}
\end{equation}
vorliegt. Die Aussage, dass es 
einen Zustand 
$\alpha\in\mathbf{P}_{2}\xi$ gibt, f"ur den die Aussage (\ref{blizbb}) erf"ullt ist,
ist hierbei gerade die Negation der Cantor-Stetigkeit\index{Cantor-Stetigkeit} 
desjenigen 
Phasenflusses $\{\xi^{t}\}_{t\in\mathbb{R}}$, dessen Existenz in jener Aussage 
(\ref{blizbb}) behauptet ist. Denn, wenn die 
strukturierte Flussfunktion im weiteren Sinn $(\xi,\mathcal{A})$
Cantor-stetig ist,
gilt f"ur jedes Mitglied $\xi^{t}=\xi(\mbox{id},t)$
des von der 
strukturierten Flussfunktion $(\xi,\mathcal{A})$
festgelegten Phasenflusses $\{\xi(\mbox{id},t)\}_{t\in\mathbb{R}}$ 
\begin{equation}\label{kalizbb} 
\begin{array}{c}
\forall\ ({\rm B},t)\in\mathcal{A}\times\mathbb{R}\ \exists\ {\rm A}\in\mathcal{A}\ :\\
\xi({\rm A},t)\subset{\rm B}\ .
\end{array}
\end{equation}
{\bf q.e.d.}
\newline\newline
Wir bezeichnen hierbei jede strukturierte Flussfunktion im weiteren Sinn $(\xi,\mathcal{A})$ 
genau dann als im Zustand 
$\alpha\in\mathbf{P}_{2}\xi$
von erster Art protosensitiv, wenn\index{Protosensitivit\"at erster Art}
wenn f"ur $\alpha$ und $(\xi,\mathcal{A})$
die Aussage
(\ref{blizbb}) wahr ist.
Wir nennen $(\xi,\mathcal{A})$
hingegen genau dann im Zustand 
$\alpha\in\mathbf{P}_{2}\xi$
von zweiter Art protosensitiv, wenn\index{verallgemeinerte Protosensitivit\"at zweiter Art}
\begin{equation}\label{ealizbb}
\begin{array}{c}
\forall\ {\rm A}\in\mathcal{A}_{\{\alpha\}}, t_{\star}\in\mathbb{R}\ \exists\ 
t\in ]-\infty,t_{\star}[\cup]t_{\star},\infty[,\ {\rm B}\in\mathcal{A}\ :\\
\xi({\rm A},t)\not\subset{\rm B}
\end{array}
\end{equation}
gilt. Die Protosensitivit"at zweiter Art ist also f"ur jede strukturierte Flussfunktion im weiteren Sinn $(\xi,\mathcal{A})$
genau dann gegeben, wenn dieselbe die
gem"ass (\ref{ealizbb}) formulierte Eigenschaft hat, welche die lediglich in der Quantorensequenzierung abgewandelte
Eigenschaft ist,
welche $(\xi,\mathcal{A})$ h"atte, wenn 
$(\xi,\mathcal{A})$ von erster Art protosensitiv w"are.
Daher impliziert die Protosensitivit"at erster Art die Protosensitivit"at zweiter Art.
Die Protosensitivit"at erster Art ist eine Vorform fixierter Sensitivit"at. Insofern,
als die Negation der Cantor-Stetigkeit die Protosensitivit"at erster Art impliziert,
ist die Negation der Cantor-Stetigkeit wiederum eine 
Vorform der Protosensitivit"at erster Art.
Die Protosensitivit"at zweiter Art scheint ein zun"achst wenig zweckm"assig verfasstes 
logisches
Pendant der Protosensitivit"at erster Art sein,
gilt doch f"ur jede Flussfunktion im weiteren Sinn $(\xi,\mathcal{A})$, deren
zweite Komponente $\mathcal{A}$ 
die gem"ass (\ref{bluzb}) durch blosse Zitierung verallgemeinerte Hausdorffsche Trennungseigenschaft hat,
dass $(\xi,\mathcal{A})$ von zweiter Art
protosensitiv ist. Doch damit noch nicht genug. 
Stellen wir die Quantorensequenzierung bei der
fixierten Sensitivit"at in der selben Weise um, wie
wir sie beim Schritt von der Protosensitivit"at erster Art zur Protosensitivit"at zweiter Art
umstellen, so ergibt
sich das in der Quantorensequenzierung analog
abgewandelte Pendant der 
fixierten Sensitivit"at.
Dieses analoge Pendant der 
fixierten Sensitivit"at
definiert in dem Fall, dass das strukturierende 
Mengensystem $\mathcal{A}$ 
die Hausdorffsche Trennungseigenschaft hat,
den Determinismus der Flussfunktionen im weiteren Sinn $(\xi,\mathcal{A})$: 
Exakt ein Paar $(\psi,\mathcal{A})$ dessen erste Komponente eine Wellenfunktion ist,
deren Wertemenge $\mathbf{P}_{2}\psi$ mit deren Definitionsmengenfaktor $\mathbf{P}_{2}\mathbf{P}_{1}\psi$ 
"ubereinstimmt,
und deren zweite Komponente eine "Uberdeckung $\mathcal{A}\subset 2^{\mathbf{P}_{2}\psi}$
deren Wertmenge $\mathbf{P}_{2}\psi=\bigcup\mathcal{A}$ ist, nennen wir eine 
strukturierte Wellenfunktion.\index{strukturierte Wellenfunktion}
Jede strukturierte Wellenfunktion $(\psi,\mathcal{A})$, deren zweite 
Komponente $\mathcal{A}$ die Hausdorffsche Trennungseigenschaft hat,
ist genau dann eine Flussfunktionen im weiteren Sinn,
wenn 
die Aussage
\begin{equation}\label{bllzba}
\begin{array}{c}
\forall\ {\rm A}\in\mathcal{A}_{\{\alpha\}}, t_{\star}\in\mathbb{R}\\ 
\exists\ (x,y)\in{\rm A}\times{\rm A} ,\ t\in ]-\infty,t_{\star}[\cup]t_{\star},\infty[,\ ({\rm X},{\rm Y})\in\mathcal{A}^{\not 2}\ :\\
\psi(x,t)\in{\rm X}\ \land\ \xi(y,t)\in{\rm Y}
\end{array}
\end{equation}
wahr ist. Diese Aussage ist eine 
Quantorenumsequenzierung
der Aussage
(\ref{blzba}) in der Definition 3.1 der fixierten Sensitivit"at.
Die Quantoren sind
so umgestellt, wie bei der 
Vertauschung der Quantorensequenzierung, welche die
Protosensitivit"at erster Art
in die
Protosensitivit"at zweiter Art
wandelt.
All ihrer fraglichen Zweckm"assigkeit zum Trotz:
Die 
Protosensitivit"at zweiter Art erweist 
sich im Unterschied zu 
ihrem Pendant,
der Vorform fixierter Sensitivit"at, welche die
Protosensitivit"at erster Art ist,
als die 
Vorform gerade derjeniger metrisierter Verallgemeinerung metrischer 
Sensitivit"at, welche dem Permanenzprinzip folgt: Denn es ist die im Rahmen
Hausdorffscher Trennbarkeit m"ogliche, (\ref{bllzba}) gem"asse
Formulierung des Determinismus der Flussfunktionen, deren nachfolgende Versch"arfung die 
mit dem 
Permanenzprinzip konform gehende 
Verallgemeinerung der Sensitivit"at bez"uglich einer Distanzfunktion
ist.
\newline
Ist jedes Mitglied $\xi^{t}=\xi(\mbox{id},t)$
des von der 
strukturierten Flussfunktion $(\xi,\mathcal{A})$
festgelegten Phasenflusses $\{\xi(\mbox{id},t)\}_{t\in\mathbb{R}}$ ein Autobolismus, d.h.,
invertierbar, und ist jenes Mitglied $\xi^{t}$ dabei
bez"uglich 
der Zustandsraum"uberdeckung $\mathcal{A}$ Cantor-stetig,\index{Cantor-Stetigkeit} so
gilt nach der Bemerkung $(3.10)^{1}$ des ersten Teiles der Konzepte der abstrakten 
Ergodentheorie
f"ur alle reellen Zahlen $t\in\mathbb{R}$ die Kommutatorrelation\index{Kommutatorrelation} 
\begin{equation}\label{calizbb}
[\mathbf{cl}_{\mathcal{A}},\xi^{t}]=\emptyset\ .
\end{equation}
Denn die G"ultigkeit der Kommutatorrelationen (\ref{calizbb})
beinhaltet nach dem Explizierungssatz $(3.8)^{1}$ des ersten Teiles, dass alle Autobolismen $\xi^{t}$ sowohl
Cantor-stetig als auch konvers Cantor-stetig sind, 
wobei die konverse Cantor-Stetigkeit eines Autobolismus das 
"Aquivalent der Cantor-Stetigkeit seiner Inversion ist.
Jede 
Flussfunktion im weiteren Sinn $(\xi,\mathcal{A})$, welche f"ur alle $t\in\mathbb{R}$ der Kommutatorrelation (\ref{calizbb})
gen"ugt, ist erst recht total insensitiv.
Insensitivit"at ist also, wie wir an (\ref{calizbb}) verdeutlicht sehen, die
Konsequenz dynamischer Symmetrie. 
\section{Ein Aufsatz:\\ Die intermittente Konzeption des Chaos}
Genau dann, wenn Sensitivit"at neben der sogenannten topologischen Transitivit"at\index{topologische Transitivit\"at}
vorliegt, spricht die Majorit"at von Chaos.
Die sogenannte topologische Transitivit"at ist dabei
im Bezug auf eine topologische Zustandsraumstrukturierung 
die Gegebenheit dessen, was das von uns im zweiten Kapitel des ersten Teiles der Konzepte der abstrakten 
Ergodentheorie eingehend untersuchte Koh"arenzkriterium\index{Koh\"arenzkriterium}
von jedem verallgemeinerten Attraktor verlangt. Exakt dann, wenn Sensitivit"at {\em innerhalb eines Attraktors} auftritt,
zeigt sich also dies Ph"anomen, welches gel"aufigerweise als 
Chaos\index{Chaos} bezeichnet wird und das insofern exakt in jenen
sogenannten chaotischen Attraktoren vorkommt.\index{chaotischer Attraktor}
Die Majorit"atsauffassung von Chaos als der Sensitivit"at innerhalb eines Attraktors
zersplittert aber dadurch, dass es diverse Auffassungen davon gibt, was als Sensitivit"at
gelten soll; welche aber allemal als ein Begriff vorliegt, 
welcher seine Eigenheit gegen"uber 
der konkurrierenden einfachen Negation der Kontinuit"at
vertritt. \newline
Nehmen wir die Streuung des Sensitivit"atsbegriffes an, wenngleich
deren Nachweis eine statistische Erhebung "uber die Fachliteratur 
erforderte, die unseren Rahmen sprengt, so kommt der Verdacht 
und das Unbehagen auf, dass diese Streuung des Sensitivit"atsbegriffes
durch eine sachliche und belangliche Divergenz verschiedener Sensitivit"atsbegriffe
bedingt ist; dass es einen Sensitivit"atspluralismus tats"achlich gibt,
der die verschiedenen Formen des Chaos definiert.
Diese sachliche und belangliche Divergenz 
der Sensitivit"atsbegriffe
ist aber bislang selber in keiner "Ubersicht dargelegt und abstrahiert.
Die Schwierigkeit, den Sensitivit"atsbegriff zu verallgemeinern,
schreckt dabei naturgem"ass davon ab, es anzugehen, jene "Ubersicht
verfassen zu wollen.
Die diversen Sensitivit"atsbegriffe unterscheidet zwar meistens eine positive Strukturbehauptung 
von der einfachen Kontinuit"atsnegation, deren jeweilige 
positive Strukturbehauptung f"allt aber verschieden aus. Letztere ist dabei
durchaus, jeweils praktikabilt"atszweckm"assig, an jeweilige
Aufgabestelllungen angepasst. \newline
Unabh"angig von der Divergenz des Sensitivit"atsbegriffes geht 
die Majorit"atsauffassung von Chaos
innerhalb der Theorie dynamischer 
Systeme aber bei genauerem Hinsehen "uber die volkst"umliche Vorlage der Idee des 
Chaos 
hinweg; und zwar auf so resolute Weise, dass diese
Hinwegsetzung "uber 
die volkst"umliche Vorlage 
sich diagnostizieren l"asst, obwohl diese Vorlage naturgem"ass nur
unscharf gefasst ist. 
\newline
Der 
Begriff des Chaos als der Sensitivit"at innerhalb eines Attraktors 
wird dem 
namengebenden Chaos der unscharfen, mythologischen Charakterisierung
des
Altgriechentumes  
als Strukturlosigkeit auf den ersten Blick dadurch gerecht, als wir
dynamische Symmetrien benennen k"onnen, die 
dies begrifflich formulierte Chaos negieren. Hier ist indess
einwendbar, dass der volkst"umliche
Gedanke 
des Chaos,
nicht etwa einzelne Symmetrien ausschliesst, sondern per conceptionem absolut sein will und
ein undenkliches Utopos des Fehlens jeglicher  
Struktur behaupten m"ochte. Sodass
es damit nicht getan ist, dass wir einzelne Symmetrien aufz"ahlen k"onnen, 
welche jeweils alle nicht vorliegen,
wenn
dies absolute Chaos gegeben w"are, das jedwede
Symmetrie ausschliesst.\newline 
Wiederum
gegen diesen Einwand l"asst sich
vorbringen, dass
dies absolut gemeinte Chaos die Strukturlosigkeit schlechthin bezeichnete; 
sodass es nicht m"oglich sein d"urfte, Chaos als entgrenzte Strukturlosigkeit innerhalb einer formalen 
Wissenschaft zu objektivieren. Dies allerdings stellt sich allemal heraus: Mythoskonform modelliertes
Chaos ist eine einfache Negation -- darstellbar nur innerhalb eines gesetzten Rahmens, der also nicht negiert ist, sondern der
positiv vorausgesetzt ist.
Innerhalb dieses Rahmens aber
ist mythoskonform modelliertes Chaos insofern eine einfache Negation, als
es nicht etwa
eine intermittente Negation einer Struktur behauptet, bei der innerhalb 
dieser Negation
zugleich schon
eine andere positive
Strukturbehauptung mitaufgestellt ist.\newline  Und darin zeigt sich nun die Emanzipation des 
Begriffes des Chaos als der Sensitivit"at innerhalb eines Attraktors
von dem Mythos des Chaos, dass 
die Sensitivit"at von vorneherein die Rolle eines Positivums
spielt; wobei die metrische Sensitivit"at beispielsweise, wie wir agitierten, eine 
intermittente Negation Kontinuit"at formuliert. Darin tritt 
ein Zug des wissenschaftlichen Begriffes des Chaos eklatant hervor, der
denselben aus der schlichten Negativit"at\index{intermittente Negation} hebt, die ein mythoskonformer Begriff des Chaos
h"atte.
\newline 
Warum wir diesen Sachverhalt so versessen betonen, dass mythoskonform modelliertes Chaos eine 
einfache Negation innerhalb eines positiven Rahmens w"are?
Dieser Sachverhalt, dass Sensitivit"at eine 
intermittente Negation\index{intermittente Negation} der Kontinuit"at ist es, welcher 
die divergierenden der Sensitivit"atsbegriffe zusammenh"alt und welcher die Grundlage 
eines Konzeptes einer "Ubersicht "uber die  
Sensitivit"atsbegriffe bildet.
Denken wir an unsere unscharfe Ahnung der Homogenit"at oder der topologischen Indifferenz
sensitiver Attraktoren, die zugleich eine Hoffnung formuliert, weil 
die Homogenit"at sensitiver Attraktoren eine 
Gegebenheit w"are, auf deren Grundlage wir gerne arbeiten w"urden:
Homogenit"atsformen decken sich nicht mit der volkst"umlichen Idee des Chaos. Homogenit"atsformen "uberschneiden  
sich vielmehr mit dem ideellen Archetyp des Nichts. Wobei
wir in unserer Wissenschaft vielerlei Formen
der 
Homogenit"at kennen, welche oft
Gegebenheiten benennen, die Symmetrien entsprechen.
(Diese Aussage will lediglich 
die Korrelation von Symmetrie und Homogenit"at
einer "ausserlichen Statistik der Namengebung 
in der Mathematik behaupten.)
Symmetrien lassen sich als das Verschwinden 
eines Kommutators formulieren und 
oft werden Nullstellenszenarien als 
die jeweiligen 
homogenen 
Szenarien bezeichnet. Chaos ist in der Tat keine nirwanische Form der Homogenit"at, geben wir sowohl dem
Volksmund als der theoretischen Essayistik anderer Disziplinen recht, nichtsdestotrotz 
freuen wir uns dar"uber, wenn sich f"ur manche Formen des Chaos noch unbestimmte Formen der Homogenit"at ergeben,
weil diese Bearbeitungsm"oglichkeiten verheissen. \newline    
Vor dem verschwommenen Hintergrund des Auseinandergehens von
geahnter Homogenit"at und der volkst"umlichen Idee des Chaos
f"allt uns der eigent"umliche Zug der fixierten Sensitivit"at auf, dass 
deren Vorform der Protosensitivit"at erster Art
nicht einfach die 
Aufl"osung beliebig kleiner Lokalisierungen behauptet. Stattdessen verlangt fixierte Sensitivit"at, dass diese 
Aufl"osung w"ahrend des Verlaufes jeweiliger Entwicklungen 
in beiden Entwicklungsrichtungen
immer wieder auftritt. Dies unterscheidet die Gegebenheit 
der Protosensitivit"at erster Art von wiederum deren 
Vorform, die exakt dann gegeben ist,
wenn es einen Zustand $\alpha\in\mathbf{P}_{2}\xi$ gibt, f"ur den (\ref{blizbb}) gilt.
Letztere
ist das blosse
Fehlen der Cantor-Stetigkeit.
Die Aussage, dass es 
einen Zustand 
$x\in\mathbf{P}_{2}\xi$ gebe, f"ur den (\ref{balizbb}) gilt, differiert von der lediglichen
Negation der Cantor-Stetigkeit, welche besagt, dass (\ref{kalizbb})
gilt.\newline 
Durch die betonte Tatsache, dass der wissenschaftliche Begriff des  
Chaos bereits als keine positiv eingerahmte, einfache Negation festgelegt und verfestigt ist, 
sind wir
nicht davon abgehalten, 
den allgemeinen Begriff des Chaos
auf dem Begriff der 
fixierten Sensitivit"at zu gr"unden, im Gegenteil.
Wir k"onnten somit die Suche nach dem 
allgemeinen Begriff der Sensitivit"at voreilig abbrechen,
welcher derjenige Begriff sein soll, der den allgemeinen Begriff des Chaos
determiniert. Denn nachdem ja der allgemeine Begriff des Attraktors 
im ersten Teil der Konzepte der abstrakten Ergodentheorie bereits
formuliert und untersucht ist, indiziert der 
Sensitivit"atsbegriff den Begriff des Chaos: Der Rahmen der nicht einfachen Negation,
die Chaos ist, ist der jeweilige Attraktor, in dem Chaos als die jeweilige nicht einfache Negation
der Kontinuit"at, welche die jeweilige Sensitivit"at ist,
auftritt. 
Gr"unden wir einen Begriff des Chaos 
auf dem Begriff der fixierten Sensitivit"at, so ist dies dementsprechend fixiertes
Chaos.\newline
Indess, wir merken, dass es einen 
Pluralismus der Sensitivit"at gibt, der den l"angst in der Untersuchung befindlichen 
Pluralismus der Formen des Chaos
generiert. Wir suchen also mehr als nur einen einzigen Begriff der verallgemeinerten Sensitivit"at.
Wir suchen diejenigen Begriffe der verallgemeinerten Sensitivit"at, welche die allgemeinen Begriffe des Chaos
bestimmen. Sch"on w"are es dabei, wenn es uns gl"uckte, dennoch zu {\em einer} 
begrifflichen Dominante von h"ochstem Generalisierungsniveau zu gelangen,
und
den Pluralismus der Sensitivit"at in einer "ubersichtlichen Hierarchie anzuordnen.
\newline
Wir haben f"urwahr viele W"unsche. Fehlt noch was?
Ja, wir h"atten es zudem gerne, wenn 
sich 
die Hierarchie der Sensitivit"at
besser
an die Absolutheit f"ugte, welche
der Mythos des Chaos 
meint, eine Absolutheit,
die keine
einzelnen, durch Chaos ausgeschlossenen Symmetrien aufz"ahlt.
Die gesuchte Hierarchie der Sensitivit"at\index{Hierarchie der Sensitivit\"at} sollte sich
besser, als beispielsweise
der recht einfache, auf dem Begriff der fixierten Sensitivit"at gegr"undete Begriff des fixierten Chaos
an die besagte Absolutheit f"ugen.
Dies ist denkbar
in der Form, dass
die Sensitivit"atsbegriffe 
{\em eine} Schl"usselsymmetrie
auf intermittente Weise negieren. Denn dann deuten wir den Mythos des Chaos
als ein Exemplar der rethorischen Standardfigur des Mythos,
welche die abstrakte Spitze unter der Impression ozeanischer Ph"anomenf"ulle lallt.
Der besagten mythogischen Absolutheit
werden wir dadurch noch
am ehesten gerecht,
wenn wir eine Schl"usselsymmetrie als diejenige ausmachen,
exakt deren Negation die allgemeine Sensitivit"at ist, die den allgemeinen Begriff des Chaos
bestimmt. 
"Ahnlich wie bei der paradigmatischen, intermittenten Unstetigkeit gem"ass der Definition 1.1b
liege 
dabei das, was den Pluralismus der Sensitivit"at aufgespannt und zusammenh"alt
in der Intermittenz der Negation der Schl"usselsymmetrie.
Der Kandidat f"ur die intermittent zu negierende Schl"usselsymmetrie steht dabei schon fest.
\section{Die Hierarchie der Sensitivit"at}
Denn, wie uns die Augen aufgehen,\index{Hierarchie der Sensitivit\"at}
halten wir den allgemeinen Begriff der 
Sensitivit"at seit geraumer Zeit in der Hand.
\newline\newline
Exakt ein Paar $(\phi,\mathcal{Z})$ eines Autobolismus $\phi$ einer Menge und einer
"Uberdeckung derselben durch ein in deren Potenzmenge enthaltenes Mengensystem $\mathcal{Z}$ 
nennen wir einen strukturierten Autobolismus.\index{strukturierter Autobolismus} F"ur jeden
strukturierten Autobolismus $(\phi,\mathcal{Z})$
ist die durch
\begin{equation}\label{thisbaa}
[\mathbf{cl}_{\mathcal{Z}},\phi]:=\mathbf{cl}_{\mathcal{Z}}\phi\Delta\phi\mathbf{cl}_{\mathcal{Z}}
\end{equation}
auf der Potenzmenge der Wertemenge $\mathbf{P}_{2}\phi$ definierte mengenwertige 
Abbildung gegeben. Dabei notiere $\Delta$ die symmetrische Differenz zweier Mengen und
$\mathbf{cl}_{\mathcal{Z}}$ sei der H"ullenoperator des Mengensystemes $\mathcal{Z}$, den wir
im dritten Kapitel des ersten Teiles der Konzepte der abstrakten Ergodentheorie
einf"uhrten. Dementsprechend ist
f"ur jede
strukturierte Wellenfunktion\index{strukturierte Wellenfunktion} 
$\Psi=(\psi,\mathcal{Z})$
jeweils 
f"ur alle $t\in\mathbb{R}$
und f"ur alle Autobolismen $\psi^{t}:=\psi(\mbox{id},t)$
durch
\begin{equation}\label{thisbea}
[\mathbf{cl}_{\mathcal{Z}},\psi^{t}]:=\mathbf{cl}_{\mathcal{Z}}\psi^{t}\Delta\psi^{t}\mathbf{cl}_{\mathcal{Z}}
\end{equation}
eine Abbildung
auf der Potenzmenge der Wertemenge $\mathbf{P}_{2}\psi$ definiert. Wobei wir daran erinnern, dass f"ur
strukturierte Wellenfunktionen deren Wertemenge mit deren zweitem Definitionsmengenfaktor
"ubereinstimmt. Von strukturierten Flussfunktionen im weiteren Sinn 
unterscheidet sie, dass die Kollektivierung $[\gamma]$ einer strukturierten Wellenfunktion
$\Gamma=(\gamma,\mathcal{G})$ nicht notwendigerweise deren 
Wertemenge $\mathbf{P}_{2}\gamma$ partioniert. 
Genau dann, wenn es eine Menge ${\rm X}\subset\mathbf{P}_{2}\phi$
gibt, f"ur die $[\mathbf{cl}_{\mathcal{Z}},\phi]({\rm X})\not=\emptyset$ ist,
nennen wir den strukturierten Autobolismus $(\phi,\mathcal{Z})$ mengenweise sensitiv.\index{mengenweise Sensitivt\"at strukturierter Autobolismen} In diesem Fall gibt es nach dem 
Explizierungssatz\index{Explizierungssatz} $(3.8)^{1}$ des ersten Teiles 
der Konzepte der abstrakten Ergodentheorie \cite{rabe}
die Menge $Q\in \mathcal{Z}$, deren Bild
$$\hat{Q}\in\{\phi^{-1}Q,\phi Q\}$$ 
im Mengensystem nicht $\mathcal{Z}$-haltiger Mengen
\begin{equation}\label{athisbea}
\mathbf{un}(\mathcal{Z}):=\Bigl\{{\rm X}\subset\bigcup\mathcal{Z}:{\rm L}\in
\mathcal{Z}\setminus\{\emptyset\}\Rightarrow{\rm L}\not\subset{\rm X}\Bigr\} 
\end{equation}
ist, das wir bereits vom ersten Teil her kennen. Es sei 
$$\hat{Q}\in: \mathbf{un}(\mathcal{Z})\cap\{\phi^{q}Q\}$$
und $q\in\{-1,1\}$:
Genau dann, wenn es ein solches nicht $\mathcal{Z}$-haltiges 
Bild $\hat{Q}$ einer Lokalisierung $Q\in \mathcal{Z}$ 
gibt und "uberdies einen Punkt $x\in \hat{Q}$ dieses Lokalisierungsbildes
sowie
eine Lokalisierung ${\rm L}\in \mathcal{Z}_{\{x\}}$  
dieses Punktes $x\in \hat{Q}$, die so beschaffen ist, dass die Restriktion $\phi^{-q}|{\rm L}$ 
Cantor-stetig ist,
nennen wir den strukturierten Autobolismus $(\phi,\mathcal{Z})$
f"ur das Lokalisierungsbild $\hat{Q}$ von erstem Grad  
mengenweise sensitiv\index{mengenweise Sensitivit\"at ersten Grades} und notieren genau dann
\begin{equation}\label{bthisbea}
\mbox{s-ord}_{(\phi,\mathcal{Z})}(Q)\geq 1\ .
\end{equation}
Als f"ur das Lokalisierungsbild $\hat{Q}$ von zweitem Grad mengenweise sensitiv bezeichnen 
wir den strukturierten Autobolismus $(\phi,\mathcal{Z})$ hingegen 
nur dann,
wenn $(\phi,\mathcal{Z})$ f"ur $\hat{Q}$ von erstem Grad
mengenweise sensitiv ist, sodass  
es einen Punkt $x\in \hat{Q}$ und eine Lokalisierung ${\rm L}\in \mathcal{Z}_{\{x\}}$ gibt, 
auf der die Bijektion $\phi^{-q}|{\rm L}$ Cantor-stetig ist.
Damit dann aber $(\phi,\mathcal{Z})$ f"ur $\hat{Q}$ als von zweitem Grad mengenweise sensitiv gilt,
sei sowohl notwendig als auch hinreichend, dass ferner ein Punkt $x\in \hat{Q}$ und ${\rm L}\in \mathcal{Z}_{\{x\}}$ und ein Punkt
$$x_{1}\in \phi^{-q}{\rm L}$$
existieren, 
f"ur den eine Lokalisierung ${\rm L}_{1}$ existiert, welche die Bedingung 
\begin{equation}\label{cthis}
\begin{array}{c}
{\rm L}_{1}\in \mathcal{Z}_{\{x_{1}\}}\ \land\\ 
{\rm L}_{1}\subset\phi^{-q}{\rm L}\ \land\\ 
\phi^{q}|{\rm L}_{1}\not\in\mathcal{C}_{+}(\mathcal{Z}\cap{\rm L}_{1})
\end{array}
\end{equation}
erf"ullt, 
auf der die Restriktion $\phi^{q}|{\rm L}_{1}$ nicht
Cantor-stetig ist. 
Denn f"ur jedes Mengensystem $\mathcal{A}$ bezeichnet $\mathcal{C}_{+}(\mathcal{A})$ 
die Menge aller 
Cantor-stetigen Abbildungen der Menge $\bigcup \mathcal{A}$ in dieselbe.
Dies ist dadurch erm"oglicht, dass die Abbildung einer Lokalisierung durch eine Cantor-stetige Funktion ein
lokalisierungshaltiges Bild hat, sodass
$$\phi^{-q}{\rm L}\not\in\mathbf{un}(\mathcal{Z})$$ 
gilt.
Exakt dann, wenn die Elemente 
$$x\in \hat{Q},\ {\rm L}\in \mathcal{Z}_{\{x\}},\ x_{1}\in \phi^{-q}{\rm L}$$
existieren, f"ur welche die Aussage (\ref{cthis}) wahr ist,
schreiben wir
\begin{equation}\label{cthisbea}
\mbox{s-ord}_{(\phi,\mathcal{Z})}(Q)\geq 2\ .
\end{equation}
F"ur jede nat"urliche Zahl $n>1$ bezeichnen wir den strukturierten Autobolismus $(\phi,\mathcal{Z})$
als f"ur das Lokalisierungsbild $\hat{Q}$ von $n$-tem Grad mengenweise sensitiv, wenn
es f"ur alle $j\in\{1,2\dots n-1\}$
\begin{equation}\label{cthisbea}
\begin{array}{c}
x_{0}\in \hat{Q},\ {\rm L}_{0}\in \mathcal{Z}_{\{x_{0}\}}\ ,\\
x_{j}\in \phi^{q^{j}}{\rm L}_{j-1} ,\ {\rm L}_{j}\in \mathcal{Z}_{\{x_{j}\}}
\end{array}
\end{equation}
gibt, wobei f"ur alle $j\in\{1,2\dots n-1\}$ 
\begin{equation}
{\rm L}_{j}\subset\phi^{-q}{\rm L}_{j-1}\ \Leftrightarrow\ q^{j}\not=q
\end{equation}
und
die Aussage
\begin{equation}\label{dthisbea}
\begin{array}{c}
\phi^{q}|{\rm L}_{j}\in\mathcal{C}_{+}(\mathcal{Z}\cap{\rm L}_{j})\ \Leftrightarrow\ q^{j}=q\\
\land\\
\phi^{q}|{\rm L}_{j}\not\in\mathcal{C}_{+}(\mathcal{Z}\cap{\rm L}_{j})\ \Leftrightarrow\ q^{j}\not=q
\end{array}
\end{equation}
wahr ist. 
Exakt in diesem Fall notieren wir
\begin{equation}\label{cthisbea}
\mbox{s-ord}_{(\phi,\mathcal{Z})}(Q)\geq n\ .
\end{equation}
Als die Sensitivit"atsordnung des strukturierten Autobolismus $(\phi,\mathcal{Z})$ f"ur die
Menge $\hat{Q}$ bezeichnen wir dabei 
das Supremum 
\begin{equation}
\mbox{s-ord}_{(\phi,\mathcal{Z})}(Q):=\sup\{k\in\mathbb{N}:\mbox{s-ord}_{(\phi,\mathcal{Z})}(Q)\geq k\}\ ,
\end{equation}
exakt welches wir als die Sensitivit"atsordnung des strukturierten Autobolismus $(\phi,\mathcal{Z})$
f"ur die Lokalisierung $Q\in \mathcal{Z}$ bezeichnen.
Dabei verstricken wir uns keineswegs in eine implizite Definition, wie vielleicht der
oberfl"achliche Blick meint.\index{Sensitivit\"atsordnung f"ur eine Lokalisierung}
Exakt das globale Supremum 
\begin{equation}\label{ethisbea}
\mbox{s-ord}(\phi,\mathcal{Z}):=\sup\{\mbox{s-ord}_{(\phi,\mathcal{Z})}(Q):Q\in \mathcal{Z}\}
\end{equation}
schliesslich nennen wir 
f"ur jeden strukturierten Autobolismus $(\phi,\mathcal{Z})$
die Sensitivit"atsordnung desselben.\index{Sensitivit\"atsordnung eines strukturierten Autobolismus}
Wir beobachten, dass wir die Inversionen $\phi^{-1}$ w"ahrend dieser Konstruktion
immer nur auf Mengen anwandten. Wir k"onnen daher den Begriff der 
Sensitivit"atsordnung 
der Pseudofl"usse $\gamma^{t}$
einer strukturierten Wellenfunktion 
$\Gamma=(\gamma,\mathcal{G})$
ganz analog zu dem
Begriff der Sensitivit"atsordnung eines 
strukturierten Autobolismus aufbauen und die Konstruktionen (\ref{athisbea})-(\ref{ethisbea})
entsprechend auch f"ur jeden
Pseudofluss $\gamma^{t}$ einer
strukturierten Wellenfunktion $\Gamma=(\gamma,\mathcal{G})$
verfassen.\newline
Damit steht uns f"ur alle nat"urlichen Zahlen $n\in\mathbb{N}$ und 
f"ur jeden Pseudofluss $\gamma^{t}$ einer strukturierten Wellenfunktion $\Gamma=(\gamma,\mathcal{G})$ 
der Begriff dessen 
mengenweiser Sensitivit"at $n$-ten Grades
f"ur jeweilige Lokalisierungsbilder $(\gamma^{t})^{-1}G$ oder 
$\gamma^{t} G$ f"ur jede Lokalisierung $G\in \mathcal{G}$ zur Verf"ugung.
Das gleiche gilt f"ur alle $t\in\mathbb{R}$ sowohl f"ur den Term der 
Sensitivit"atsordnung $\mbox{s-ord}_{(\gamma^{t},\mathcal{G})}(G)$ f"ur eine jeweilige Lokalisierung
$G\in \mathcal{G}$ als auch f"ur die 
Sensitivit"atsordnung $\mbox{s-ord}(\gamma^{t},\mathcal{G})$ eines jeweiligen Pseudoflusses.
\newline  
\newline 
{\bf Definition 3.3: Mengenweise Sensitivit"at und deren Ordnung.\index{mengenweise Sensitivit\"at}} \newline
{\em Es sei $\Psi=(\psi,\mathcal{Z})$ eine strukturierte Wellenfunktion.\index{strukturierte Wellenfunktion}
Genau dann, wenn es eine reelle Zahl $t\in \mathbb{R}$ und eine Menge ${\rm X}\subset\mathbf{P}_{2}\psi$ gibt,
f"ur die
\begin{equation}\label{thisbe}
[\mathbf{cl}_{\mathcal{Z}},\psi^{t}]({\rm X})\not=\emptyset
\end{equation}
ist,
nennen wir die strukturierte Wellenfunktion $\Psi$
mengenweise sensitiv. Dabei nennen wir ${\rm X}$ eine sensitive Menge der 
strukturierten Wellenfunktion $\Psi$ und 
stufen die strukturierte Wellenfunktion $\Psi$ als von nulltem Grad mengenweise sensitiv
\index{sensitive Menge einer strukturierten Wellenfunktion} ein.\index{strukturierte Wellenfunktion nullten Grades}\newline 
Von der unbefristeten mengenweisen Sensitivit"at reden wir 
genau dann, wenn es eine beidseitig unbeschr"ankte Folge reeller Zahlen $\{t_{j}\}_{j\in\mathbb{N}}$
gibt, f"ur die f"ur alle $j\in\mathbb{N}$ eine sensitive Menge ${\rm X}_{j}\subset\mathbf{P}_{2}\psi$ existiert, f"ur die
\begin{equation}\label{thisbe}
[\mathbf{cl}_{\mathcal{Z}},\psi^{t_{j}}]({\rm X}_{j})\not=\emptyset
\end{equation}
ist.\index{mengenweise Sensitivt\"at strukturierter Wellenfunktionen}\index{unbefristete mengenweise Sensitivt\"at strukturierter Wellenfunktionen}
Wir sagen, dass }
\begin{equation}\label{thisbe}
\mbox{s-ord}(\Psi):=\sup_{t\in\mathbb{R}}\mbox{s-ord}(\psi^{t},\mathcal{Z})
\end{equation}
{\em die Sensitivit"atsordnung der strukturierten\index{Sensitivit\"atsordnung einer strukturierten Wellenfunktion}
Wellenfunktion $\Psi$ ist.}
\newline
\newline
Um nun aber die Bindung der mengenweisen Sensitivit"at nullten Grades an jeweilige Lokalisierungen
ausdr"ucken zu k"onnen, legen wir fest, dass
wir $\Psi$ genau dann als 
in jedem Bild $(\psi^{t})^{-1}Z\in\mathbf{un}(\mathcal{Z})$ oder 
$\psi^{t} Z\in\mathbf{un}(\mathcal{Z})$ von nulltem Grad mengenweise sensitiv
bezeichnen, wenn $Z\in \mathcal{Z}$
eine Lokalisierung und $t\in\mathbb{R}$ ist. Die Lokalisierung $Z\in \mathcal{Z}$ nennen wir exakt in diesem Fall eine
sensitive Lokalisierung von nulltem Grad.\index{sensitive Lokalisierung von nulltem Grad} 
Und konsequenterweise bezeichnen wir sie f"ur jede nicht negative ganze Zahl $n$ genau dann als eine 
sensitive Lokalisierung von $n$-tem Grad, wenn
$$\mbox{s-ord}_{\Psi}(Z)\geq n$$
gilt und nennen sie genau dann eine sensitive Lokalisierung von $n$-ter Ordnung, wenn
$$\mbox{s-ord}_{\Psi}(Z)= n$$ 
ist.\index{sensitive Lokalisierung von $n$-tem Grad}\index{sensitive Lokalisierung $n$-ter Ordnung} 
Betrachten wir beispielsweise den Fall einer nat"urlich topologisierten Entwicklung 
$(\xi,\mathbf{T}(n)\cap\mathbf{P}_{2}\xi)$,
so sehen wir, dass f"ur jede 
sensitive Lokalisierung ${\rm U}\in \mathbf{T}(n)\cap\mathbf{P}_{2}\xi$ eine Zahl $t\in\mathbb{R}$
existiert, f"ur welche das Lokalisierungsbild 
$\xi^{t}{\rm U}\in \mathbf{un}(\mathbf{T}(n)\cap\mathbf{P}_{2}\xi)$ keine
offene Kugel enth"alt: Das Lokalisierungsbild 
$\xi^{t}{\rm U}$ hat damit, grob aber eindringlich gesagt, seine Dimension verringert.
Es kann dabei ein Fraktal\index{Fraktal} sein.\newline
Genau dann, wenn f"ur $n\in\mathbb{N}_{0}$ eine sensitive Lokalisierung von $n$-tem Grad
in einem Attraktor $a\in\mbox{{\bf @}}(\xi,\mathcal{A})$ 
einer strukturierten Flussfunktion $(\xi,\mathcal{A})$ im weiteren Sinn enthalten ist,
nennen wir den Attraktor $a$ einen chaotischen Attraktor von $n$-tem Grad.
Wir beachten: Wir legen damit nicht fest, dass
jeder chaotischer Attraktor von nulltem Grad oder von irgendeinem Grad 
als ein chaotischer Attraktor gelte.
Wir sagen auch nicht, dass ein Attraktor $a$ von $n$-ter Ordnung sei, wenn
er eine sensitive Lokalisierung $n$-ter Ordnung enth"alt.
Die chaotische Ordnung eines Attraktors $a$ sei als die Sensitivit"atsordnung 
\begin{equation}\label{aprthis}
\mbox{s-ord}((\xi|(a\times\mathbb{R}),\mathcal{A}\cap a))
\end{equation}
der Restriktion $(\xi|(a\times\mathbb{R}),\mathcal{A}\cap a)$ festgelegt.\index{chaotische Ordnung eines Attraktors}
Ein chaotischer Attraktor $n$-ter Ordnung sei exakt ein solcher, dessen chaotische Ordnung $n$ ist.
Was als ein chaotischer Attraktor gelte, werden wir erst sp"ater, nach der Einf"uhrung der
verallgemeinerten zustandsweisen 
Sensitivit"at in der Definition 3.5, bestimmen.\index{chaotischer Attraktor}
\index{chaotischer Attraktor von 
$n$-tem Grad}\index{chaotischer Attraktor $n$-ter Ordnung}
\newline
Anders als die metrische Sensitivit"at, die, wie 
wir ausf"uhrten,
eine sehr spezielle Form negierter Stetigkeit ist, ist die
mengenweise Sensitivit"at -- oder die mengenweise Sensitivit"at nullten Grades
oder die mengenweise Sensitivit"at nullter Sensitivit"atsordnung --  einer strukturierten Entwicklung
eine einfache Negation einer Kontinuit"atsform, n"amlich die einfache Negation der Cantor-Stetigkeit.
Statt einer strukturierten Wellenfunktion\index{Cantor-Stetigkeit}
zu diagnotizieren, dass sie mengenweise sensitiv sei, k"onnen wir auch sagen, dass sie nicht
kommutativ Cantor-stetig\index{kommutative Cantor-Stetigkeit} sei.
Dass die kommutative Cantor-Stetigkeit durch eine strukturierte Entwicklung $(\xi,\mathcal{A})$ verletzt ist,
ist die notwendige Voraussetzung daf"ur, dass bei 
einer strukturierten Entwicklung $(\xi,\mathcal{A})$
Chaos 
auftritt.
Die notwendige Voraussetzung f"ur Chaos
verletzt damit ausdr"ucklich 
die Voraussetzung des 
verallgemeinerten insensitven Ergodensatzes $(3.3)^{1}$ des ersten Teiles,\index{verallgemeinerter insensitver Ergodensatz}
den die Implikation
\begin{equation}\label{bprthis}
\begin{array}{c}
[\mathbf{cl}_{\mathcal{A}},\xi^{\mathbb{R}}]=\{\emptyset\}
\Rightarrow[[\xi]]_{\mathcal{A}}\in\mathbf{part}(\mathbf{P}_{2}\xi)
\end{array}
\end{equation}
darstellt. Das Mengensystem $[[\xi]]_{\mathcal{A}}$ ist nicht notwendigerweise eine Partition des 
Zustandsraumes, wenn die strukturierte Entwicklung $(\xi,\mathcal{A})$ mengenweise Sensitivit"at
zeigt.
Es kann aber
trotz auftretender mengenweiser Sensitivit"at
dennoch vorkommen, dass keine verschiedenen Vorzimmer $\chi_{1},\chi_{2}\in [[\xi]]_{\mathcal{A}}$
schneiden, weil die Implikation (\ref{bprthis}) offensichtlich nicht umkehrbar ist:
Haben wir n"amlich eine feste Zustandsraumstrukturierung $\mathcal{A}$ und eine feste 
Kollektivierung $[\xi]$, ist mit denselben bereits das Mengensystem
$[[\xi]]_{\mathcal{A}}$ bestimmt, die 
Wahl der Flussfunktionen $\xi$ im weiteren Sinn der Faser kollektivierungsgleicher 
Flussfunktionen im weiteren Sinn
\begin{equation}\label{dprthis}
[\mbox{id}]^{-1}(\{[\xi]\})=
\Bigl\{\psi\in (\bigcup[\xi])^{(\bigcup[\xi])\times\mathbb{R}}:[\psi]=[\xi]\Bigr\} 
\end{equation}
ist indess noch frei. 
Wenn wir uns beispielsweise eine Partition $[\vartheta]$ einer offenen Kreisscheibe 
der Ebene $\mathbb{R}^{2}$ in glatte und geschlossene 
Trajektorien der Ebene und in einen einzigen Fixpunkt $\{\omega\}$ vorstellen, so sehen wir leicht ein, dass 
wir das Mengensystem $[\vartheta]$ auf verschiedene Weise 
durch Funktionen $\vartheta_{0},\vartheta_{1}\in [\mbox{id}]^{-1}(\{[\vartheta]\})$ parametrisieren k"onnen.
Der Leser verifiziert,
dass es m"oglich ist, zu erreichen, dass
\begin{displaymath}
[\mathbf{cl},\vartheta_{0}^{\mathbb{R}}]=\{\emptyset\}\ \land
[\mathbf{cl},\vartheta_{1}^{\mathbb{R}}]\not=\{\emptyset\}
\end{displaymath}
ist. Es gibt demnach Kollektivierungen $[\xi]$, deren mengenweise Sensitivit"at 
f"ur eine jeweilige Zustandsraumstrukturierung $\mathcal{A}$
in dem Sinn
wegtransformierbar ist, dass eine Flussfunktion im weiteren Sinn $\xi_{0}\in [\mbox{id}]^{-1}(\{[\xi]\})$ 
existiert, f"ur die
$[\mathbf{cl}_{\mathcal{A}},\xi_{0}^{\mathbb{R}}]=\{\emptyset\}$ ist, w"ahrend 
es zugleich andere Flussfunktionen $\xi_{1}\in [\mbox{id}]^{-1}(\{[\xi]\})$
gibt, f"ur die $[\mathbf{cl}_{\mathcal{A}},\xi_{1}^{\mathbb{R}}]\not=\{\emptyset\}$
ist. Es liegt nahe, f"ur jede Partition $\mathcal{P}$ und f"ur jede 
"Uberdeckung $\mathcal{A}\subset 2^{\bigcup\mathcal{P}}$ mit 
\begin{equation}\label{eprthis}
\bigcup\mathcal{P}=\mathcal{A} 
\end{equation}
exakt die Menge
von Funktionen 
\begin{equation}
\mathbf{can}(\mathcal{P},\mathcal{A}):=
\Bigl\{\xi_{0}\in [\mbox{id}]^{-1}(\{\mathcal{P}\}):
[\mathbf{cl}_{\mathcal{A}},\vartheta_{0}^{\mathbb{R}}]=\{\emptyset\}\Bigr\}
\end{equation}
als die Menge Cantor-geeichter Flussfunktionen\index{Cantor-geeichte Flussfunktion}
des Paares $(\mathcal{P},\mathcal{A})$ einzuf"uhren. Die Faser 
\begin{equation}
[\mbox{id}]^{-1}(\{\mathcal{P}\})=\Bigl\{\psi\in (\bigcup[\xi])^{(\bigcup[\xi])\times\mathbb{R}}:[\psi]=\mathcal{P}\Bigr\} 
\end{equation}
von Flussfunktionen im weiteren Sinn, deren Kollektivierung die Partition $\mathcal{P}$ ist,
ist genau dann leer, wenn die Partition $\mathcal{P}$ ein Element hat, dessen Kardinalit"at
"uber $\aleph_{1}$ liegt. Wenn die Faser $[\mbox{id}]^{-1}(\{[\xi]\})$ nicht leer ist, kann 
die Menge Cantor-geeichter Flussfunktionen
$\mathbf{can}([\xi],\mathcal{A})$ durchaus leer sein. 
Dann gibt es keine Cantor-geeichte Flussfunktion $\tilde{\xi}$, die 
$[\xi]$ in dem Sinn erzeugt, dass $[\xi]=[\tilde{\xi}]$ ist. Exakt 
eine solche Partition $\mathcal{P}$, die als eine
Kollektivierung $[\xi]=\mathcal{P}$ einer Flussfunktion im weiteren Sinn
$\xi\in[\mbox{id}]^{-1}(\{\mathcal{P}\})$ darstellbar ist, die jedoch
f"ur eine "Uberdeckung $\mathcal{A}\subset 2^{\bigcup\mathcal{P}}$ mit (\ref{eprthis})
nicht Cantor-geeicht werden kann, nennen wir eine bez"uglich der 
Strukturierung $\mathcal{A}$
Cantor-unstetige Kollektivierung.\index{Cantor-unstetige Kollektivierung}
Wir stossen auf die folgende Frage: Gilt die logische Vervollst"andigung 
des verallgemeinerten insensitven Ergodensatzes in Form der f"ur alle strukturierten Flussfunktionen
$(\xi,\mathcal{A})$ g"ultigen
"Aquivalenz
\begin{equation}\label{fprthis}
\mathbf{can}([\xi],\mathcal{A})\not=\emptyset\Leftrightarrow[[\xi]]_{\mathcal{A}}\in\mathbf{part}(\mathbf{P}_{2}\xi)\ ?
\end{equation}
Nein. Es erweist sich als keine einfache Aufgabe, eine 
bez"uglich einer
nat"urlichen 
Topologie Cantor-unstetige Kollektivierung explizit zu benennen.
Indess, der Satz von der Sensitivit"at nicht trivialer Zimmer $(2.2.2)^{0}$ der
Abhandlung des elementaren Quasiergodensatzes \cite{raab} sagt uns, dass die 
nicht trivialen Zimmer, die Gegenstand dieses Satzes sind, so beschaffen sind,
dass sie bez"uglich der
nat"urlichen 
Topologie ihres Zustandsraumes von Cantor-unstetigen trajektoriellen
Partitionen erzeugt sind. F"ur sie 
darf es allemal
keine bez"uglich der nat"urlichen 
Topologie Cantor-geeichten Flussfunktionen im weiteren Sinn geben, deren Kollektiverung sie sind,
sodass die "Aquivalenz (\ref{fprthis}) nicht gilt.
\newline
Mit dem Sachverhalt, dass 
die mengenweise Sensitivit"at ausdr"ucklich die Voraussetzung des 
verallgemeinerten insensitven Ergodensatzes verletzt,
verlagert sich die Relevanz von den Zimmern\index{Zimmer}
einer jeweiligen strukturierten Flussfunktion im weiteren Sinn $(\xi,\mathcal{A})$ 
auf deren Vorzimmer, welche exakt die Mengen des Mengensystemes $[[\xi]]_{\mathcal{A}}$ 
sind. Jedes Vorzimmer $\chi\in [[\xi]]_{\mathcal{A}}$ ist genau dann ein Zimmer,
wenn die Auswahl
\begin{equation}\label{cprthis}
([[\xi]]_{\mathcal{A}})_{\chi}=\{\chi\}
\end{equation}
ist. Damit verlagert sich gleichzeitig die Richtung unseres Interesses
in der Form, dass wir 
angehalten sind,
nach den Bedingungen zu fragen, unter denen die Umkehrbarkeit
der Implikation (\ref{bprthis}) vorliegt: Wenn die Implikationsrichtung 
der Implikation (\ref{bprthis}) umkehrbar ist, schliesst mengenweise
Sensitivit"at -- und damit schon Chaos nullten Grades -- die partitive 
Architektur der Zimmer aus. Wir interessieren uns daher 
f"ur die Szenarien, in denen die Umkehrbarkeit
der Implikation (\ref{bprthis}) nicht gegeben ist, sodass
Chaos nullten Grades nicht 
ausgeschlossen ist, das sich aber dennoch in die 
durch die partitive 
Architektur der Zimmer gesetzte 
Ordnung f"ugt. Und insbesondere solche F"alle reizten, in denen 
extreme Formen der Sensitivit"at und des Chaos dennoch in der
partitiven Zimmerordnung bleiben, so es solche  
Szenarien geben sollte.   
\newline
Inwieweit "aussert sich die schlichte Negativit"at der mengenweisen Sensitivit"at
darin, dass die mengenweise Sensitivit"at
unspezifisch ist und
dass die mengenweise Sensitivit"at nicht nur insofern 
allgemein ist, dass sie f"ur beliebige 
strukturierte Wellenfunktionen konstituiert ist?
Die mengenweise Sensitivit"at behauptet im konkreten Fall eine in folgendem Sinn zun"achst
unbestimmte 
Sensitivit"at:\newline
Sei $\Xi=(\xi,\mathbf{T}(d))$ eine strukturierte Entwicklung, deren Strukturierung 
eine metrische Topologie ihres Zustandsraumes ist,
auf dem die Metrik $d$ definiert ist, womit auch 
das gem"ass der Definition 1.2 verfasste Aufl"osungsfeld\index{Aufl\"osungsfeld einer Wellenfunktion} $\Delta_{\xi}^{d}$ der Entwicklung $\xi$
bez"uglich 
dieser Metrik $d$ bestimmt ist.
Es gilt die Implikation
\begin{equation}\label{thisbeb}
\begin{array}{c}
z\in\bigcup(\Delta_{\xi}^{d})^{-1}(\mathbb{R}^{+})\\
\Rightarrow\\
\exists\ {\rm U}\in\mathbf{T}(d)_{\{z\}},t\in\mathbb{R} :\ 
[\mathbf{cl}_{\mathcal{Z}},\xi^{t}]({\rm U})\not=\emptyset\ .
\end{array}
\end{equation}
Denn, wie wir wissen,
schliesst metrische Sensitivit"at die metrische Stetigkeit aus. 
Metrische Stetigkeit ist aber notwendig
f"ur die Cantor-Stetigkeit bez"uglich der metrischen Topologie $\mathbf{T}(d)$. 
Nach dem Explizierungssatz 
\index{Explizierungssatz} $(3.8)^{1}$ des ersten Teiles ist die kommutative Cantor-Stetigkeit einer  
Entwicklung zu deren Cantor-Stetigkeit "aquivalent, welche zu der metrischen Stetigkeit "aquivalent
ist. Wir erwarten, dass die Implikation (\ref{thisbeb})
dabei im allgemeinen nicht umkehrbar ist.\newline
Dass die Implikation (\ref{thisbeb}) im allgemeinen nicht umkehrbar ist, wiese uns darauf hin,
dass wir mit der Verfassung der mengenweisen Sensitivit"at in der Definition 3.3 
das Permanenzprinzip\index{Permanenzprinzip} zwar nicht verletzen;
dass wir aber 
in der Definition 3.3
den Begriff der metrischen Sensitivit"at und der
Sensitivit"at bez"uglich zweier Distanzfunktionen
auf eine unidentische Weise ausdehnen:
Der Begriff der
mengenweisen Sensitivit"at ist nicht mit dem 
Begriff deren metrischer Sensitivit"at identisch, wenn wir die mengenweise Sensitivit"at
auf die F"alle 
strukturierter Entwicklungen $\Xi=(\xi,\mathbf{T}(d))$ eingeschr"anken, die
"Aquivalente 
metrisierter Entwicklungen sind. F"ur $(\xi,\mathbf{T}(d))$ kann mengenweise Sensitivit"at im Bezug
auf die metrische Topologie $\mathbf{T}(d)$ vorliegen,
obwohl $(\xi,d)$ nicht metrisch sensitiv im Sinn der Definition 1.2 ist.
Wir wollen folgenden Sachverhalt unmissverst"andlich herausstellen und notieren hierzu
Inklusion und Schnitt f"ur Klassen wie f"ur Mengen:\newline
Es sei 
f"ur alle $n\in\mathbb{N}\cup\{0,\infty\}$ mit
$\mbox{SF}^{n}$ die Klasse strukturierter Flussfunktionen im weiteren Sinn bezeichnet, 
deren Sensitivit"atsordnung mindestens $n$ ist.
Es sei $\mbox{MF}$ die Klasse metrisierter 
Flussfunktionen im weiteren Sinn. Letztere sind alle von der Form $(\xi,d)$, sodass
$$(\mathbf{P}_{1}\oplus\mathbf{T})\mbox{MF}\subset\mbox{SF}^{0}$$
die Klasse derjeniger strukturierter Flussfunktionen im weiteren Sinn ist, die von der 
Form $(\xi,\mathbf{T}(d))$ sind und die 
zu metrisierten Flussfunktionen $(\xi,d)$
"aquivalent sind,
wenn dabei $\mathbf{T}$ als der Operator gelesen wird, der 
jeder Metrik $d$ die von ihr induzierte Topologie zuweist.
Die Klasse strukturierter Flussfunktionen im weiteren Sinn mit mengenweiser Sensitivit"at, die zu
metrisierten Flussfunktionen "aquivalent sind, ist also
$$\mbox{SF}^{0}\cap(\mathbf{P}_{1}\oplus\mathbf{T})\mbox{MF}\ .$$
Ferner sei $\mbox{SMF}$ die Klasse metrisierter Flussfunktionen im weiteren Sinn, die
metrisch sensitiv sind, und 
$\mbox{PE}$ bezeichne die Klasse Poincar$\acute{e}$scher Entwicklungen.\index{Poincar$\acute{e}$sche Entwicklung}
Dann gilt der folgende 
\newline 
\newline 
{\bf Eingrenzungssatz\index{Eingrenzungssatz} 3.4:\newline Die eingrenzende Permanenz der mengenweisen Sensitivit"at}\index{Eingrenzende Permanenz der mengenweisen Sensitivit\"at}\newline 
{\em Der Begriff der mengenweisen Sensitivit"at
strukturierter Flussfunktionen im weiteren Sinn 
setzt den Begriff metrischer Sensitivit"at metrisierter Flussfunktionen im weiteren Sinn
auf der Klasse Poincar$\acute{e}$scher Entwicklungen eingrenzend fort: Es gelten die Inklusionen}
\begin{equation}\label{bebthis}
\begin{array}{c}
\Bigl(\mbox{SF}^{1}\cap(\mathbf{P}_{1}\oplus\mathbf{T})\mbox{MF}\Bigr)\cap\mbox{PE}\\
\subset\Bigl((\mathbf{P}_{1}\oplus\mathbf{T})\mbox{SMF}\Bigr)\cap\mbox{PE}\\
\subset\Bigl(\mbox{SF}^{0}\cap(\mathbf{P}_{1}\oplus\mathbf{T})\mbox{MF}\Bigr)\cap\mbox{PE}\ .
\end{array}
\end{equation}
{\bf Beweis:}\newline 
Wir brauchen lediglich die 
Inklusion 
$$\Bigl(\mbox{SF}^{1}\cap(\mathbf{P}_{1}\oplus\mathbf{T})\mbox{MF}\Bigr)\cap\mbox{PE}\subset \Bigl((\mathbf{P}_{1}\oplus\mathbf{T})\mbox{SMF}\Bigr)\cap\mbox{PE}$$
zu verifzieren, denn
die Klasse strukturierter
Entwicklungen der Form $\Xi=(\xi,\mathbf{T}(d))$, f"ur die
$(\xi,d)$ eine metrisierte Entwicklung ist, die bez"uglich $d$ metrisch
sensitiv ist und die zur Klasse $\mbox{SMF}$ metrisch sensitiver metrisierter Entwicklungen z"ahlt,
ist gem"ass (\ref{thisbeb}) mengenweise sensitiv.
\newline 
Im Vergleich zu der Sensitivit"at von Wellenfunktionen bez"uglich zweier Distanzfunktionen
gem"ass der Definition 1.2 besteht vor allem der 
Unterschied der mengenweisen Sensitivit"at 
gem"ass der Definition 3.3 zu der Sensitivit"at gem"ass der Definition 1.2, dass die mengenweise Sensitivit"at
nicht zustandsweise verfasst ist. F"ur eine jeweilige sensitive Menge ${\rm X}\subset\bigcap\mathcal{A}$
einer strukturierten Flussfunktion $(\psi,\mathcal{A})$ im weiteren Sinn zeigt sich
mengenweise Sensitivit"at genau dann, wenn 
$$[\mathbf{cl}_{\mathcal{A}},\psi^{\mathbb{R}}]({\rm X})\not=\{\emptyset\}$$
gilt. Falls $\psi$ eine Entwicklung  
ist, heisst das
nach dem Explizierungssatz\index{Explizierungssatz} $(3.8)^{1}$, 
dass es eine Menge ${\rm A}^{\star}\in \mathcal{A}$ und eine Zahl $t^{\star}\in\mathbb{R}$ gibt, f"ur die
$\psi^{t^{\star}}{\rm A}^{\star}\in\mathbf{un}(\mathcal{A})$ eine Menge ist, in der keine Lokalisierung aus $\mathcal{A}$
Platz findet. F"ur alle 
Lokalisierungen ${\rm L}\in\mathcal{A}_{\{z\}}$ f"ur
Zust"ande $z\in\psi^{t^{\star}}{\rm A}^{\star}$
gilt dann daher, dass
$$(\psi^{t^{\star}})^{-1}({\rm L})\not\subset{\rm A}^{\star}$$
ist. Ist $(\psi,d)$ eine metrisierte Entwicklung und $\mathcal{A}=\mathbf{T}(d)$ zu der
durch $d$ induzierten Topologie spezifiziert, so ist f"ur jeden Zustand
$z\in\psi^{t^{\star}}{\rm A}^{\star}$ wahr, dass 
das Bild $(\psi^{t^{\star}})^{-1}{\rm U}$
jeder seiner Umgebungen ${\rm U}\in\mathbf{T}(d)_{\{z\}}$
die offene Menge ${\rm A}^{\star}\in \mathbf{T}(d)$ schneidet, weil $(\psi^{t^{\star}})^{-1}(z)\in{\rm A}^{\star}$ 
und $z\in {\rm U}$ ist. Zugleich ist es aber 
ausgeschlossen, dass ${\rm U}$ in $(\psi^{t^{\star}})^{-1}{\rm A}^{\star}$ enthalten ist. 
W"ahlen wir den Zustand $z\in\psi^{t^{\star}}{\rm A}^{\star}$ so, dass f"ur eine positive reelle Zahl $\delta$
\begin{equation}\label{cebthis}
\inf\{d(z,x):x\in \partial{\rm A}^{\star}\}>\delta
\end{equation}
ist, so ist
also die Aussage
\begin{equation}\label{aebthis}
\begin{array}{c}
\exists\ \delta\in\mathbb{R}^{+}:\ \forall\ {\rm U}\in\mathbf{T}(d)_{\{z\}}\ \exists\ y\in {\rm U}:\\
d((\psi^{t^{\star}})^{-1}(z),(\psi^{t^{\star}})^{-1}(y))>\delta
\end{array}
\end{equation}
wahr. Wenn $(\psi,\mathbf{T}(d))$ eine strukturierte Entwicklung ist, die 
von erster Ordnung sensitiv ist, so gibt es gem"ass der Definition 3.3 nicht nur eine Menge ${\rm A}^{\star}\in \mathbf{T}(d)$ 
und eine Zahl $t^{\star}\in\mathbb{R}$, f"ur die
$\psi^{t^{\star}}{\rm A}^{\star}\in\mathbf{un}(\mathbf{T}(d))$ eine Menge ist, in der keine Umgebung der Topologie $\mathbf{T}(d)$
Platz findet. Es gibt im Fall der Sensitivt"at erster Ordnung eine offene Menge ${\rm A}^{\star}\in \mathbf{T}(d)$ 
und eine Zahl $t^{\star}\in\mathbb{R}$, f"ur welche die Aussage
$$\psi^{t^{\star}}{\rm A}^{\star}\in\mathbf{un}(\mathbf{T}(d))\ \land$$
$$\exists\ z\in \psi^{t^{\star}}{\rm A}^{\star},\ {\rm U}\in\mathbf{T}(d)_{\{z\}}\ :$$
$$(\psi^{t^{\star}})^{-1}|{\rm U}\in\mathcal{C}_{+}(\mathbf{T}(d)\cap{\rm U})$$
wahr ist. Daher finden wir einen Zustand $z\in\psi^{t^{\star}}{\rm A}^{\star}$, der so beschaffen ist, 
dass f"ur eine positive reelle Zahl $\delta$ die
Ungleichung
(\ref{cebthis}) eingehalten ist, sodass (\ref{aebthis}) gilt, wobei "uberdies eine Umgebung
dieses Zustandes $z$ existiert, auf der die
Restriktion $(\psi^{t^{\star}})^{-1}|{\rm U}$ stetig ist. Denn
die Cantor-Stetigkeit ist nach dem Explizierungssatz\index{Explizierungssatz} $(3.8)^{1}$ eine identische Fortsetzung der topologischen Stetigkeit. 
Betrachten wir diese Umgebung ${\rm U}\in\mathbf{T}(d)_{\{z\}}$ und einen Punkt $y$
derselben, f"ur den die Aussage (\ref{aebthis}) zutrifft, weil 
$$(\psi^{t^{\star}})^{-1}(y)\not\in {\rm A}^{\star}$$
ist:
Falls $\psi$ Poincar$\acute{e}$sch ist,\index{Poincar$\acute{e}$sche Entwicklung} gibt es 
f"ur alle Umgebungen ${\rm U}_{y}\in \mathbf{T}(d)_{\{y\}}$, die in ${\rm U}\supset {\rm U}_{y}$ enthalten sind,
eine beidseitig 
unbeschr"ankte Folge reeller Zahlen $\{t_{j}(y,{\rm U}_{y})\}_{j\in\mathbb{N}}$, f"ur die f"ur alle $j\in\mathbb{N}$ 
$$\psi(y,t_{j}(y,{\rm U}_{y}))\in{\rm U}_{y}$$ 
ist. Weil $(\psi^{t^{\star}})^{-1}$ auf der Umgebung ${\rm U}$ stetig ist, in der
$y$ ist und in der ${\rm U}_{y}\subset{\rm U}$ liegt, gibt es f"ur alle positiven Zahlen $\varepsilon\in \mathbb{R}^{+}$
eine Umgebung ${\rm U}_{y}(\varepsilon)\in \mathbf{T}(d)_{\{y\}}$, f"ur die alle Zust"ande $\hat{y}\in {\rm U}_{y}(\varepsilon)$
der Ungleichung
$$d((\psi^{t^{\star}})^{-1}(\hat{y}),(\psi^{t^{\star}})^{-1}(y))<\varepsilon$$
gen"ugen. W"ahlen wir $\varepsilon=\delta/2$,
so gilt daher f"ur alle Zust"ande $\hat{y}\in {\rm U}_{y}(\delta/2)$
$$d((\psi^{t^{\star}})^{-1}(z),(\psi^{t^{\star}})^{-1}(\hat{y}))>\delta/2\ .$$
Wenn $\psi$ Poincar$\acute{e}$sch ist, gibt es auch eine beidseitig 
unbeschr"ankte Folge reeller Zahlen 
\begin{equation}\label{debthis}
\{t_{j}(y,\delta)\}_{j\in\mathbb{N}}:=\{t_{j}(y,{\rm U}_{y}(\delta/2))\}_{j\in\mathbb{N}}\ ,
\end{equation}
f"ur die f"ur alle $j\in\mathbb{N}$ 
$$\psi(y,t_{j}(y,{\rm U}_{y}(\delta/2)))\in{\rm U}_{y}(\delta/2)$$
ist,
sodass f"ur die beidseitig unbeschr"ankte Folge (\ref{debthis}) die 
Ungleichung
$$d((\psi^{t^{\star}})^{-1}(z),(\psi^{t^{\star}})^{-1}(\psi(y,t_{j}(y,\delta)))>\delta/2$$
erf"ullt wird. Es
gibt also den Zustand 
$z\in\psi^{t^{\star}}{\rm A}^{\star}$, ferner eine Umgebung ${\rm U}\in \mathbf{T}(d)_{\{z\}}$, auf der $\psi^{t^{\star}}$ stetig ist;
und ausserdem ist in dieser Umgebung ${\rm U}$ eine Umgebung $\hat{{\rm U}}\in \mathbf{T}(d)_{\{z\}}$ enthalten,
in der 
Zust"ande 
$y$ existieren, f"ur die es eine Folge $\{\tilde{t}_{j}\}_{j\in\mathbb{N}}$ gibt,
f"ur die f"ur alle $j\in\mathbb{N}$ 
$$d((\psi^{t^{\star}})^{-1}(z),(\psi^{t^{\star}})^{-1}(\psi(y,\tilde{t}_{j}))>\delta/2$$
gilt; sodass es wegen der Stetigkeit der Restriktion $(\psi^{t^{\star}})^{-1}|U$ 
eine positive Zahl $\hat{\delta}$ gibt, f"ur die f"ur alle $j\in\mathbb{N}$
$$d(z,\psi(y,\tilde{t}_{j}))>\hat{\delta}$$
gilt, wobei $\{\tilde{t}_{j}\}_{j\in\mathbb{N}}$ beidseitg unbeschr"ankt ist.
\newline
{\bf q.e.d.}
\newline
\newline
Wir beachten, dass
der Eingrenzungssatz\index{Eingrenzungssatz} 3.4 
beschr"ankt auf die Klasse Poincar$\acute{e}$scher Entwicklungen $\mbox{PE}$ gilt, 
welche die Klasse global Poincar$\acute{e}$scher Entwicklungen echt umfasst.\index{global Poincar$\acute{e}$sche Entwicklung}
\index{Poincar$\acute{e}$sche Entwicklung}
F"ur den eineindeutigen universellen Operator $(\mathbf{P}_{1}\oplus\mathbf{T})$ finden sich triviale Beispiele daf"ur,
dass die Verneinung
$$\neg\Bigl(\mbox{SF}^{1}\subset(\mathbf{P}_{1}\oplus\mathbf{T})\mbox{SMF}\subset\mbox{SF}^{0}\Bigr)$$
wahr ist. 
Weil die mengenweise Sensitivit"at verschiedener Ordnungen keine
identische Fortsetzung metrischer Sensitivit"at ist,  
sind wir durchaus noch an anderen Formen 
verallgemeinerter Sensitivit"at interessiert, die
auf 
der Grundlage beliebiger Zustandsraumstrukturierungen verfassbar sind;
insbesondere an solchen Formen 
verallgemeinerter Sensitivit"at, 
die dem Permanenzprinzip auf die Weise folgen, dass sich 
die metrische Sensitivit"at als eine 
restringierte 
Form jener verallgemeinertern Sensitivit"at ausdr"ucken l"asst. Gleichermassen pr"azisiert wie objektiviert gesprochen:
Es sei $\mbox{X}=\mbox{X}^{+}\cup\mbox{X}^{-}$ die Klasse exakt der Objekte, die
jene gesuchte Form verallgemeinerter Sensitivit"at haben k"onnen und sie entweder haben oder aber nicht haben. 
Exakt 
diejenigen Objekte der Klasse $\mbox{X}$, die jene gesuchte Form verallgemeinerter Sensitivit"at haben, seien die
Objekte der Klasse $\mbox{X}^{+}$, w"ahrend exakt diejenigen Objekte der Klasse $\mbox{X}$, 
welchen jene verallgemeinerte Sensitivit"at fehlt, die
Objekte der Klasse $\mbox{X}^{-}$ seien.
Es gebe eine Unterklasse $\mbox{MX}\subset\mbox{X}$, welche die eineindeutige Konstruktion $\phi$ auf die Klasse
metrisierter  Flussfunktionen im weiteren Sinn
$$\mbox{MF}=\phi\mbox{MX}$$ abbilde. F"ur den eineindeutigen universellen Operator $\phi$ gelte, dass
$$\phi(\mbox{MX}\cap\mbox{X}^{+})=\mbox{SMF}$$
gerade die Klasse metrisch sensitiver metrisierter Flussfunktionen im weiteren Sinn sei.
\newline
Um eine solche permanenzprinziptreue Fortsetzung 
metrischer Sensitivit"at zu anzulegen, welche in Gestalt der beschriebenen Klasse $\mbox{X}$ 
objektiviert ist,
sprachen wir von Protosensitivit"at erster und zweiter Art.
Darum stellten wir der Protosensitivit"at erster Art
das triviale Pendant der Protosensitivit"at zweiter Art gegen"uber:
Die Gegebenheit
fixierter Sensitivit"at
einer Flussfunktion im weiteren Sinn $(\xi,\mathbf{T}(d))$, 
deren Zustandsraumstrukturierung $\mathbf{T}(d)$ eine durch die Metrik 
$d$ induzierte Topologie ist, "uberbestimmt offenbar sogar eine ausgezeichnete 
Form der metrischen Sensitivit"at bez"uglich $d$, welche die 
Funktion $\xi$ in diesem Fall hat. Die 
logische Lokalisierung innerhalb der Protosensitivit"at erster Art, welche 
die fixierte Sensitivit"at ist, folgt also dem
Permanenzprinzip nicht auf die Weise, die durch die Klasse $\mbox{X}$ objektiviert ist,
die offenbar auch zu keiner der Klassen $\mbox{SF}^{n}$ f"ur $n\in\mathbb{N}\cup\{0,\infty\}$
"aquivalent ist.
Daher soll uns sowohl die durch die Klassen $\mbox{SF}^{n}$ f"ur $n\in\mathbb{N}\cup\{0,\infty\}$
objektivierte Hierarchie der Sensitivit"at\index{Hierarchie der Sensitivit\"at} als auch der 
Begriff
fixierter Sensitivit"at als die nicht exakt
treffende Verallgemeinerung der metrischen Sensitivit"at noch nicht genug sein.\index{Permanenzprinzip}
F"ur jedes Paar 
$(\xi,\mathbf{T}(d))$ liegt pauschal 
Protosensitivit"at der zweiten Art vor und sogar
das entsprechende Pendant fixierter Sensitivit"at gem"ass (\ref{bllzba}), das den
Determinismus formuliert. Was
ganz im Gegensatz zu der mit der fixierten Sensitivit"at verbundenen "Uberbestimmung steht.
Wir hoffen,
dass
an dem Beispiel 
der durch die metrische Topologie $\mathbf{T}(d)$ strukturierten 
Flussfunktion im weiteren Sinn $(\xi,\mathbf{T}(d))$ sehr deutlich wird, dass
die folgende Definition der quantifizierten Sensitivit"at einer 
strukturierten Wellenfunktion $(\psi,\mathcal{Z})$ die gesuchte 
geradlinige Fortsetzung des Begriffes der punktweisen Sensitivit"at
einer Wellenfunktion bez"uglich zweier Distanzfunktionen gem"ass der
Definition 1.2 ist. 
\newline
Wenn wir auf dem zu der zweiten Komponente $\mathcal{Z}$
einer strukturierten Wellenfunktion $(\psi,\mathcal{Z})$
gem"ass (\ref{bluzb})
festgelegten Mengenssystem $\mathbf{P}_{1}\mathcal{Z}^{\not 2}$
eine Distanzfunktion $d^{\circ}$
anlegen, erhalten wir die folgende Verallgemeinerung, die mit dem
Permanenzprinzip harmoniert:
\newline
\newline{\bf Definition 3.5: Quantifizierte Sensitivit"at\index{quantifizierte Sensitivit\"at}} \newline
{\em Es sei das Tripel $(\psi,\mathcal{Z},d^{\circ})$ aus einer 
Wellenfunktion $\psi$, einer "Uberdeckung deren Wertemenge $\mathcal{Z}$ und 
einer Distanzfunktion $d^{\circ}$ gem"ass} (\ref{bluzb}) {\em zusammengesetzt.
Wir nennen die Wellenfunktion $\psi$ genau dann in einem Punkt $\alpha\in\mathbf{P}_{2}\xi=\bigcup\mathcal{Z}$
durch die Distanzfunktion $d^{\circ}$
quantifiziert sensitiv, wenn\index{quantifizierte Sensitivit\"at in einem Punkt}
die 
Aussage
\begin{equation}\label{bllzba}
\begin{array}{c}
\exists\ \delta(\alpha)\in ]0,\infty[:\ \forall\ {\rm A}\in\mathcal{A}_{\{\alpha\}}, t_{\star}\in\mathbb{R}\\ 
\exists\ (x,y)\in{\rm A}\times{\rm A} ,\ t\in ]-\infty,t_{\star}[\cup]t_{\star},\infty[,\ ({\rm X},{\rm Y})\in\mathcal{A}^{\not 2}\ :\\
\quad\\
\psi(x,t)\in{\rm X}\ \land\ \xi(y,t)\in{\rm Y}\ \land\\
d^{\circ}({\rm X},{\rm Y})\geq\delta(\alpha)
\end{array}
\end{equation}
wahr ist.}
\newline
\newline
Angewandt
auf unser Beispiel 
der durch die metrische Topologie $\mathbf{T}(d)$ strukturierten 
Flussfunktion im weiteren Sinn $(\xi,\mathbf{T}(d))$
zeigt sich: 
Genau dann, wenn ein Punkt $x\in\bigcup\mathbf{T}(d)$
durch die von der Metrik $d$ 
abgeleitete Hausdorff-Distanz $d^{[1]}$ quantifiziert sensitiv ist, ist er
bez"uglich $d$ metrisch sensitiv im Sinn der
Definition 1.2.\newline
Der Begriff metrischer Sensitivit"at liegt nicht zwischen den 
metrischen Spezifzierungen der allgemeinen Begriffe der 
quantifizierten
Sensitivit"at und der fixierten Sensitivit"at,
sondern 
die quantifizierte
Sensitivit"at setzt
den Begriff metrischer Sensitivit"at fort, w"ahrend die fixierte Sensitivit"at
von demselben
abweicht.
Das Begriffspaar der fixierten Sensitivit"at einerseits und  
der quantifizierten Sensitivit"at
andererseits
generalisiert
demnach den Begriff metrischer Sensitivit"at auf eine im flachen Sinn
dialektische Weise:
W"ahrend der Begriff der fixierten Sensitivit"at der Flussfunktion im weiteren Sinn $(\xi,\mathcal{A})$
totale Allgemeinheit
erreicht, 
dabei aber aus der Spur des Permanenzprinzipes l"auft,
folgt der Begriff der quantifizierten Sensitivit"at des
Tripels $(\psi,\mathcal{Z},d^{\circ})$ 
zwar
dem Permanenzprinzip. Er kommt aber nicht ohne 
das Hinzuziehen der jeweiligen Distanzfunktion
$d^{\circ}$ aus. Wir bezeichnen die beobachtete Dialektik 
als die Dialektik fixierter und quantifizierter Sensitivit"at.\index{Dialektik fixierter und quantifizierter Sensitivit\"at}
\newline
Dies n"otige Additiv der Distanzfunktion
$d^{\circ}$
macht dabei den 
Begriff der quantifizierten Sensitivit"at an sich keineswegs weniger
allgemein. In der Wahl strukturierender Mengensysteme $\mathcal{Z}$ sind wir frei.
In der theoretischen Praxis der Untersuchung der Frage, ob unter bestimmten Voraussetzungen 
quantifizierte Sensitivit"at auftritt, m"ussen wir 
allerdings auf die jeweiligen Wahlen
der Distanzfunktionen $d^{\circ}$ eingehen.\newline
Haben wir, wenn wir auf $\mathcal{Z}$ beispielsweise 
eine Metrik $d^{\circ}$ anlegen wollen, nicht erst die Metrisierungsaufgabe
abzuarbeiten, die wir im letzten Abschnitt beschrieben?
Wir brauchen, um eine metrische
quantifizierte Sensitivit"at anwenden zu k"onnen,
keine Intialbasis\index{Intialbasis}
zu bestimmen. Und deren Bestimmung ist es, welche bei dem Transfer 
von Zahlenstrahlmetriken auf 
Zustandsr"aume
sensitiver Flussfunktionen vor
numerische Ungangbarkeiten stellt. Insoweit 
k"onnen wir mit dem Begriff der 
quantifizierten Sensitivit"at als 
Verallgemeinerung zufrieden sein. 
\newline
Die Dialektik fixierter und quantifizierter Sensitivit"at birgt aber noch einen
beunruhigenden weiteren Gesichtspunkt.
Wir sagten bereits, dass die 
es die Detektorik sei, welche als die jeweiligen strukturierenden Mengensysteme
abstrahiert ist. Der Begriff der quantifizierte Sensitivit"at wird insofern nicht handgreifbar,
als die 
quantifizierten Sensitivit"at
keine Aussage macht, wo im Zustandsraum die jeweilige  
Befindlichkeit
$$\psi(x,t)\in{\rm X}\ \land\ \xi(y,t)\in{\rm Y}$$
der Aussage (\ref{bllzba}) konstatiert wird.
Es wird in (\ref{bllzba}) lediglich die Existenz 
der Paare $({\rm X},{\rm Y})\in\mathcal{A}^{\not 2}$
behauptet, deren Komponenten irgendwo im Zustandsraum liegen und 
welche
dabei hinreichende Distanz voneinander haben. 
Der Begriff der quantifizierten Sensitivit"at ist darin 
gegens"atzlich zum Begriff
der fixierten Sensitivit"at, bei dem zu einem jeweiligen Zustand $\alpha$
ein Paar $({\rm X}(\alpha),{\rm Y}(\alpha))\in\mathcal{A}^{\not 2}$ fixiert wird;
das ist der Grund, weshalb wir die fixierte Sensitivit"at so nennen.
Es mag 
an dieser Stelle das Motiv der Divergenz zwischen der konstruktiven Mathematik 
und der nicht konstruktiven Mathematik anklingen, durcheinandert"onend:
Wenn wir uns darauf besinnen, dass es die 
Detektorik ist, welche die Abstraktion strukturierender Mengensysteme 
gleichermassen objektiviert, wie die Abstraktion
die detektorische Datenerhebung
auch verh"ullt, dann bedeutet diese konzeptionelle
Differenz zwischen 
quantifizierter und fixierter Sensitivit"at dies:
Die fixierte Sensitivit"at modelliert, dass die 
Detektoren f"ur den jeweiligen Nachweis der fixierten 
sensitiven Dynamik des Zustandes $\alpha$ 
aufgebaut stehen, w"ahrend
die Detektoren f"ur den jeweiligen Nachweis der 
quantifizierten sensitiven Dynamik des Zustandes $\alpha$ auf g"ottliche Weise
allgegenw"artig sind. Wir sind nicht und keine Ger"atschaft ist Gott.
Deshalb nehmen wir auch den gleichsam mobilit"atsbehinderten
Begriff der fixierten Sensitivit"at ernst. 
Untersuchen wir,
f"ur welche Gegebenheiten sich 
fixiert sensitives Verhalten zeigt,
so brauchen wir an keine Distanzfunktion $d^{\circ}$ zu denken.

\section{Ultrakolokalisation und Kolokalisation als nachvollziehbare Sensitivit"at}
Wir passen unsere Redeweise diesem Allgemeinheitsgrad insofern schrittweise
an, als wir uns erst die Gegebenheiten f"ur Flussfunktionen auf dem 
Generalit"atsniveau der allgemeinen Topologie vor Augen f"uhren,
um dann einen Interpretationshintergrund anloger, noch allgemeinerer 
Konstruktionen vergegenw"artigt zu haben, wenn
wir anschliessend diese noch allgemeineren Konstruktionen
vorstellen.\newline  
Zun"achst
legen wir fest, dass
wir jedes Paar $(\xi, \mathbf{T}(\mathbf{P}_{2}\xi))$ in Analogie zu der Sprechweise 
von metrisierten Flussfunktionen 
genau dann als eine 
topologisierte Flussfunktion\index{topologisierte Flussfunktion} bezeichnen, wenn
$\xi$ eine Flussfunktion ist, deren Zustandsraum $\mathbf{P}_{2}\xi$ mit der Topologie 
$\mathbf{T}(\mathbf{P}_{2}\xi)$
versehen ist. Wir betrachten w"ahrend dieses 
Abschnittes
Flussfunktionen, die nicht nur
insofern normiert seien, dass die Zentriertheit
$\xi(\mbox{id},0)=\mbox{id}$ gegeben ist. Wir betrachten 
Flussfunktionen, deren Definitionsmenge auf 
das kartesische Produkt $\mathbf{P}_{1}\xi=\mathbf{P}_{2}\xi\times\overline{\mathbb{R}}$
erweitert ist, wobei, wie "ublich,
\begin{equation}
\overline{\mathbb{R}}:=\mathbb{R}\cup\{-\infty,\infty\}
\end{equation}
die Tr"agermenge des kompaktifizierten  
Zahlenstrahles\index{kompaktifizierter Zahlenstrahl} $(\overline{\mathbb{R}},\overline{\mathbf{T}(1)})$ bezeichne, dessen Topologie 
$\overline{\mathbf{T}(1)}\supsetneqq\mathbf{T}(1)$ und dessen algebraische Struktur beispielsweise\index{kompaktifizierte Zahlenstrahltopologie}
in \cite{elst} kurz erl"autert ist. Wir nennen hierbei $\overline{\mathbb{R}}$ einfach den kompaktifizierten Zahlenstrahl, ohne,
dass wir damit gegen verfestigte Redeweisen verstossen.
Es sei dabei 
\begin{equation}\label{tante}
\xi(\mbox{id},0)=\xi(\mbox{id},\pm\infty)=\mbox{id}\ .
\end{equation}
Dies erm"oglicht es uns, jeder der in dieser Weise normierten Flussfunktionen $\xi$ leicht deren
zeitreziproke Flussfunktion $\overline{\xi}$ zuzuordnen,\index{zeitreziproke Flussfunktion} die
auf der jeweils gleichen Definitionsmenge wie $\xi$
\begin{equation}
\mathbf{P}_{1}\overline{\xi}=\mathbf{P}_{1}\xi
\end{equation}
festgelegt sein soll; und zwar dadurch, dass
f"ur alle $t\in\overline{\mathbb{R}}$ 
\begin{equation}
\overline{\xi}(\mbox{id},t):=\left\{
\begin{array}{cll}
\xi(\mbox{id},1/t)&\ f.\ &t\in\mathbb{R}\setminus\{0\} \\
\mbox{id}&\ f.\ &t\in\{-\infty,0,\infty\}
\end{array}\right .
\end{equation}
sei. 
Es gilt offenbar 
\begin{equation}
[\xi]=[\overline{\xi}]
\end{equation}
und damit auch, dass die Mengensysteme der
Vorzimmer  
\begin{equation}
[[\xi]]_{\mathcal{A}}=[[\overline{\xi}]]_{\mathcal{A}}
\end{equation}
des Paares $(\xi,\mathcal{A})$
und des Paares
$(\overline{\xi},\mathcal{A})$ f"ur alle 
jeweiligen Zustandsraum"uberdeckungen $\mathcal{A}$
"ubereinstimmen. 
Damit haben wir in der Abbildung $\overline{\mbox{id}}$,
deren jeweilige Werte zu $\xi$ 
gerade $\overline{\xi}$ seien,
eine universelle Involution auf der Klasse normierter 
Flussfunktionen konstruiert. Der Zweck der 
Normierung (\ref{tante}) von Flussfunktionen ist es, die   
Zuordnung jeweils zugeh"origer zeitreziproker Flussfunktionen
zu erleichtern.
Der Zweck zeitreziproker Flussfunktionen $\overline{\xi}$
wiederum ist es, es
im Bedarfsfall zu erleichtern, die asymptotisch 
langfristigen Entwicklungen zu beschreiben.
\newline
Wir nennen dabei die topologisierte Flussfunktion $(\xi, \mathbf{T}(\mathbf{P}_{2}\xi))$
gem"ass der gel"aufigen Bezeichnungsweise  
genau dann stetig, wenn sie stetig bez"uglich 
der Topologie $\mathbf{T}(\mathbf{P}_{1}\xi)$
ihrer Definitionsmenge und
der 
Zustandsraumtopologie $\mathbf{T}(\mathbf{P}_{2}\xi)$ ist, welche
als die zweite Komponente
dieser topologisierten Flussfunktion $(\xi, \mathbf{T}(\mathbf{P}_{2}\xi))$
angegeben ist. 
Die Topologie $\mathbf{T}(\mathbf{P}_{1}\xi)$ der Definitionsmenge
ist hierbei die Produkttopologie der Zustandsraumtopologie $\mathbf{T}(\mathbf{P}_{2}\xi)$
und der kompaktifizierten Zahlenstrahltopologie $\overline{\mathbf{T}(1)}\supset\mathbf{T}(1)$, welche
die 
nat"urliche Zahlenstrahltopologie $\mathbf{T}(1)$ echt enth"alt.\index{kompaktifizierte Zahlenstrahltopologie} Die Produkttopologie
$\mathbf{T}(\mathbf{P}_{1}\xi)$
schreiben wir von der gel"aufigen Notation abweichend. 
Wir wollen 
nicht mit Hilfe
des selben Zeichens \glqq$\times$\grqq, mit dem wir das kartesische Produkt
darstellen, Produkttopologien notieren, weil diese gel"aufige Schreibweise 
irref"uhrend sein kann und es ja
durchaus auch vorkommt, dass wir das kartesische Produkt zweier Topologien bezeichnen wollen und nicht deren
Produkttopologie.
Daher schreiben wir die Produkttopologie
einfach hintereinander: Es ist 
$$\mathbf{T}(\mathbf{P}_{1}\xi)=\mathbf{T}(\mathbf{P}_{2}\xi)\overline{\mathbf{T}(1)}\supsetneqq \mathbf{T}(\mathbf{P}_{2}\xi)\mathbf{T}(1)\ .$$ 
Im Hinblick auf die Stetigkeit der 
topologisierten Flussfunktion $(\xi, \mathbf{T}(\mathbf{P}_{2}\xi))$, f"ur die (\ref{tante}) gilt, spielt die
Kompaktifizierung des Zahlenstrahles keine Rolle. 
Die Stetigkeit der 
topologisierten Flussfunktion $(\xi, \mathbf{T}(\mathbf{P}_{2}\xi))$
heisst also, dass erstens die Stetigkeit in der ersten Ver"anderlichen 
bez"uglich der Zustandsraumtopologie $\mathbf{T}(\mathbf{P}_{2}\xi)$
in Form der G"ultigkeit
der Aussage 
\begin{equation}\label{ydnula}
\begin{array}{c}
\forall\ {\rm U}\in \mathbf{T}(\mathbf{P}_{2}\xi)_{\{\xi(x,t)\}}\ \exists\ {\rm V}\in 
\mathbf{T}(\mathbf{P}_{2}\xi)_{\{x\}}:\\\
\bigcup\xi({\rm V},t)\subset{\rm U} 
\end{array}
\end{equation}
f"ur alle Definitionsmengenelemente $(x,t)\in \mathbf{P}_{1}\xi$ gegeben ist. Dabei notiert 
f"ur alle Zust"ande $y\in \mathbf{P}_{2}\xi$ die Auswahl 
$\mathbf{T}(\mathbf{P}_{2}\xi)_{\{y\}}\not\ni\emptyset$ das Umgebungssystem des Zustandes
$y$.\index{Umgebungssystem eines Zustandes} Wollen wir diese Stetigkeit in der ersten Ver"anderlichen
einzeln
benennen, so sprechen wir 
von der zust"andlichen Stetigkeit der 
topologisierten Flussfunktion $(\xi, \mathbf{T}(\mathbf{P}_{2}\xi))$. 
Dieselbe ist also genau dann gegeben, wenn (\ref{ydnula}) 
gilt.\index{zust\"andliche Stetigkeit einer topologisierten Flussfunktion}
Und zweitens besagt die Stetigkeit der topologisierten Flussfunktion
im gel"aufigen Sinn, dass
auch die Stetigkeit in der zweiten Ver"anderlichen 
bez"uglich $\mathbf{T}(\mathbf{P}_{2}\xi)$
in Form der G"ultigkeit
der Aussage 
\begin{equation}\label{ydnulb}
\begin{array}{c}
\forall\ {\rm U}\in \mathbf{T}(\mathbf{P}_{2}\xi)_{\{\xi(x,t)\}}\ \exists\ 
\delta\in\mathbb{R}^{+}\ :\\
\xi(x,]t-\delta,t+\delta[)\subset{\rm U}
\end{array}
\end{equation}
f"ur alle $(x,t)\in \mathbf{P}_{1}\xi$
vorliegt, falls
\begin{equation}\label{tnez}
\xi(\mbox{id},0)=\mbox{id}
\end{equation}
ist. 
Wenn wir auch diese Stetigkeit in der zweiten Ver"anderlichen
f"ur sich genommen
bezeichnen wollen, so sei von der phasischen Stetigkeit der 
topologisierten Flussfunktion $(\xi, \mathbf{T}(\mathbf{P}_{2}\xi))$ die Rede.
Wenn wir von \index{phasische Stetigkeit einer topologisierten Flussfunktion}
phasischer Stetigkeit sprechen,
gehen wir demnach von der Zentriertheit\index{Zentriertheit} gem"ass (\ref{tnez}) aus.
Phasische Stetigkeit liegt nur vor, wenn (\ref{tnez}) gilt.
Die Stetigkeit
einer topologisierten Flussfunktion umfasst also zum einen 
deren zust"andliche Stetigkeit, die 
innerhalb des gem"ass (\ref{ydnula}) bestimmten Selbstbezuges der
Zustandsraumtopologie $\mathbf{T}(\mathbf{P}_{2}\xi)$
bleibt.\index{Stetigkeit einer topologisierten Flussfunktion}
Zum anderen 
umfasst die Stetigkeit
einer topologisierten Flussfunktion auch deren 
phasische Stetigkeit, welche
sich auf die nat"urliche Zahlenstrahltopologie $\mathbf{T}(1)$ und 
die angegebene Zustandsraumtopologie $\mathbf{T}(\mathbf{P}_{2}\xi)$ bezieht.
Zu der phasischen Stetigkeit 
einer topologisierten Flussfunktion kann "uberdies deren 
konverse 
phasische Stetigkeit\index{konvers phasische Stetigkeit einer 
topologisierten Flussfunktion} gegeben sein, die
genau dann vorliege, wenn die 
Zentriertheit gem"ass (\ref{tnez}) gegeben ist und
die
Aussage 
\begin{equation}\label{ydnunulb}
\begin{array}{c}
\forall\ \delta\in\mathbb{R}^{+}\ \exists{\rm U}\in \mathbf{T}(\mathbf{P}_{2}\xi)_{\{\xi(x,t)\}}\ :\\
{\rm U}\cap\xi(x,\mathbb{R})\subset\xi(x,]t-\delta,t+\delta[)
\end{array}
\end{equation}
f"ur alle $(x,t)\in \mathbf{P}_{1}\xi$
gilt.
Genau dann, wenn eine topologisierte Flussfunktion
sowohl phasische Stetigkeit als auch konverse 
phasische Stetigkeit besitzt, nennen wir sie perfekt 
phasisch stetig.\index{perfekte phasische Stetigkeit einer topologisierten Flussfunktion}
Die Stetigkeit der topologisierten Flussfunktion $(\xi, \mathbf{T}(\mathbf{P}_{2}\xi))$ 
ist von der 
auf die Zustandsraumtopologie
$\mathbf{T}(\mathbf{P}_{2}\xi)$ bezogenen
Stetigkeit ihrer jeweiligen Fl"usse $\xi^{t}=\xi(\mbox{id},t)$ f"ur alle $t\in\overline{\mathbb{R}}$
begrifflich unterschieden, welche die
G"ultigkeit
der Aussage 
\begin{equation}\label{ydnulc}
\begin{array}{c}
\forall\ x\in\mathbf{P}_{2}\xi, {\rm U}\in \mathbf{T}(\mathbf{P}_{2}\xi)_{\{\xi(x,t)\}}\ \exists\ 
{\rm V}\in \mathbf{T}(\mathbf{P}_{2}\xi)_{\{x\}}:\\
\xi^{t}{\rm V}\subset{\rm U}
\end{array}
\end{equation}
f"ur alle $t\in \overline{\mathbb{R}}$ behauptet. Offensichtlich 
existieren Beispiele daf"ur, dass
ein Phasenfluss $\{\xi^{t}\}_{t\in\overline{\mathbb{R}}}=\{\xi(\mbox{id},t)\}_{t\in\overline{\mathbb{R}}}$
lauter bez"uglich $\mathbf{T}(\mathbf{P}_{2}\xi)$ stetige Fl"usse $\xi^{t}$ hat
und dass dieser Phasenfluss $\{\xi^{t}\}_{t\in\overline{\mathbb{R}}}$ dabei
eine topologisierte Flussfunktion $(\xi, \mathbf{T}(\mathbf{P}_{2}\xi))$ festlegt,
die nicht stetig ist, wenn wegen der 
Diskontinuit"at der
flusswertigen Funktion $\xi^{\mbox{id}}$
des Zahlenstrahles $\overline{\mathbb{R}}$ die Aussage (\ref{ydnulb}) nicht gilt.
Wenn $b$ eine Bijektion des 
Zahlenstrahles auf denselben ist, die kein 
Hom"oomorphismus ist,
und wenn $\{\xi^{t}\}_{t\in\overline{\mathbb{R}}}$
ein bez"uglich $\mathbf{T}(\mathbf{P}_{2}\xi)$ stetiger Phasenfluss ist, so definiert
\begin{equation}\label{endula} 
\{\xi^{b(t)}\}_{t\in\overline{\mathbb{R}}}=:\{\xi_{b}(\mbox{id},t)\}_{t\in\overline{\mathbb{R}}}
\end{equation}
die Flussfunktion $\xi_{b}$ die bez"uglich $\mathbf{T}(\mathbf{P}_{2}\xi)$ im allgemeinen nicht stetig ist:
$\xi_{b}$ ist im Fall, dass $b$ kein Hom"oomorphismus ist, genau dann stetig, wenn
$b$ nur die Indizes $t$ jeweiliger Fl"usse $\xi^{t}$, nicht aber die Fl"usse $\xi^{t}= \xi^{b(t)}$ vertauscht.
\newline
Die Stetigkeit des Phasenflusses $\{\xi^{t}\}_{t\in\overline{\mathbb{R}}}$ bez"uglich einer jeweiligen 
Zustandsraumtopologie $\mathbf{T}(\mathbf{P}_{2}\xi)$
haben wir bereits im ersten Teil der abstrakten Ergodentheorie 
zu der Cantor-Stetigkeit\index{Cantor-Stetigkeit} des Phasenflusses $\{\xi^{t}\}_{t\in\overline{\mathbb{R}}}$  
bez"uglich einer jeweiligen
Zustandsraum"uberdeckung $\mathcal{A}$ verallgemeinert:
Ist 
$\mathcal{A}$ eine Zustandsraum"uberdeckung, d.h., ein
Mengensystem
$\mathcal{A}\in 2^{\mathbf{P}_{2}\xi}$,
dessen Vereinigung  
$\bigcup\mathcal{A}=\mathbf{P}_{2}\xi$
der Zustandsraum ist, so gilt der Phasenfluss $\{\xi^{t}\}_{t\in\overline{\mathbb{R}}}$ genau dann
als bez"uglich der Zustandsraum"uberdeckung $\mathcal{A}$ kommutativ
Cantor-stetig, wenn der 
im ersten Teil der abstrakten Ergodentheorie eingef"uhrte
Kommutator
$[\mathbf{cl}_{\mathcal{A}},\xi^{t}]$ f"ur alle $t\in\overline{\mathbb{R}}$
leer ist.
Da die Generalisierbarkeit "uber den Verallgemeinertheitsgrad der allgemeinen Topologie hinaus
die Formulierung des als ein verallgemeinerter Quasiergodensatz auffassbaren 
Satzes $(3.3)^{1}$ "uber verallgemeinerte
Zimmer\index{verallgemeinerter insensitver Ergodensatz} erm"oglicht,
erscheint 
uns die "uber die Allgemeinheitsstufe der allgemeinen Topologie hinausgehende
Verallgemeinerung
verheissungsvoll. Wir erweitern  
schliesslich die Sprechweise 
von metrisierten und topologisierten Flussfunktionen dahingehend,
dass wir jedes Paar $(\xi,\mathcal{A})$  
genau dann als eine 
strukturierte Flussfunktion\index{strukturierte Flussfunktion} bezeichnen, wenn
$\xi$ eine Flussfunktion ist
und $\mathcal{A}$ eine deren Zustandsraum"uberdeckungen.
Spannen wir uns nun nicht l"anger auf die Folter!
Pr"asentieren wir endlich den Begriff der Kolokalisation, den wir am Anfang des Kapitels ank"undigten
als eine nachvollziehbare Spezialform der mengenweisen Sensitivit"at
einer strukturierten Flussfunktion im weiteren Sinn; welche 
zwar eine Spezialform der mengenweisen Sensitivit"at,\index{mengenweise Sensitivit\"at}
nichtsdestotrotz eine solche Spezialform ist, die auf dem gleichen 
Generalit"atsniveau vorliegt, auf dem wir die
mengenweise Sensitivit"at verfassten. Der Weg zum Begriff der Kolokalisation
f"uhrt "uber wiederum deren Spezialform der Ultrakolokalisation, die wir ebenfalls in dem
Allgemeinheitsgrad beliebiger Zustandsraumstrukturierungen formulieren:
\newline
F"ur jede strukturierte Flussfunktion $\Lambda=(\xi,\mathcal{A})$
im weiteren Sinn
und 
f"ur jedes Paar zweier Zust"ande 
$(z,y)\in \mathbf{P}_{2}\xi\times \mathbf{P}_{2}\xi$ 
bezeichnen wir $z$ genau dann als
durch die strukturierte Flussfunktion im weiteren Sinn\index{strukturierte Flussfunktion im weiteren Sinn} $\Lambda$ mit $y$ 
ultrakolokalisierend, wenn es 
f"ur das Zustandspaar\index{ultrakolokalisierendes Zustandspaar} $(z,y)$
einen Zustand $\omega(z,y)$ gibt, f"ur den
\begin{equation}\label{ydnulf}
\begin{array}{c}
\forall\ {\rm U}\in\mathcal{A}_{\{\omega(z,y)\}},\ t_{\star}\in\mathbb{R}\ \exists\ (t_{-}({\rm U}),t_{+}({\rm U}))\in]-\infty,t_{\star} [\times ]t_{\star},\infty[:\\
\{\xi(z,t_{\pm}({\rm U})),\xi(y,t_{\pm}({\rm U}))\}\subset {\rm U}\\
\end{array}
\end{equation}
gilt. Exakt diesen Sachverhalt notieren wir in der Form
\begin{equation}\label{ygdnulf}
z,y\ \sim>_{\Lambda}\ \omega(z,y)\ .
\end{equation}
Den 
Sachverhalt (\ref{ygdnulf})
fassen wir mit Hilfe der Sprechweise in Worte, dass
genau dann, wenn $z,y\ \sim>_{\Lambda}\ \omega(z,y)$ gilt,
$z$ und $y$ in $\omega(z,y)$ durch $\Lambda$ ultrakolokalisieren.
Die Schreibweise 
\begin{equation}\label{yhdnulf}
\begin{array}{c}
z\sim_{\Lambda} y\\
:\Leftrightarrow\\
\exists\ \omega(z,y)\in \mathbf{P}_{2}\xi\times \mathbf{P}_{2}\xi:\ z,y\ \sim>_{\Lambda}\ \omega(z,y)
\end{array}
\end{equation}
dr"uckt demnach den Sachverhalt aus, dass $z$
durch die strukturierte Flussfunktion\index{strukturierte Flussfunktion} $\Lambda$ zu $y$ 
ultrakolokalisiert ist. Dabei lassen wir die Erw"ahnung der strukturierten Flussfunktion fallen,
wenn sie eine Konstante des Kontextes ist.
Es sei
\begin{equation}\label{yfdnulf}
\mbox{coloc}(\Lambda):=
\Bigl\{\omega\in \mathbf{P}_{2}\xi:\exists (z,y)\in \mathbf{P}_{2}\xi^{2}:z,y\ \sim>_{\Lambda}\ \omega\Bigr\}
\end{equation}
und 
\begin{equation}\label{yfdyulf}
\begin{array}{c}
\mbox{coloc}^{0}(\Lambda)\\
:=\Bigl\{\omega\in \mathbf{P}_{2}\xi:\exists (z,y)\in \mathbf{P}_{2}\xi\times 
(\mathbf{P}_{2}\xi\setminus\{z\})
:z,y\ \sim>_{\Lambda}\ \omega\Bigr\}\ .
\end{array}
\end{equation}
Exakt die Zust"ande 
der Menge
\begin{equation}\label{yftisdlf}
\begin{array}{c}
\mbox{coloc}^{+}(\Lambda):=\Bigl\{\omega\in \mbox{coloc}(\Lambda):
\exists (z,y)\in \mathbf{P}_{2}\xi^{2}:z,y\ \sim>_{\Lambda}\ \omega\ \land\\ 
\exists\ {\rm L}(z),{\rm L}(y)\in\mathcal{A}:\  z\not\in{\rm L}(y)\ni y\ \land\ y\not\in{\rm L}(z)\ni z\Bigr\}
\end{array}
\end{equation}
bezeichnen wir hierbei als die ultrasensitven Zust"ande der strukturierten 
Flussfunktion $\Lambda$ im weiteren Sinn.\index{ultrasensitver Zustand einer strukturierten Flussfunktion im weiteren Sinn}
Dass wir die ultrasensitven Zust"ande der Menge $\mbox{coloc}^{+}(\Lambda)$ gerade so benennen, liegt nahe.
Denn, wenn $z$ ein sensitiver Zustand einer metrisierten
Flussfunktion $(\xi,d)$ im Sinne der Definition 1.2 ist,
ist $z$
offensichtlich keineswegs notwendigerweise ein ultrasensitver Zustand der Flussfunktion 
$\Xi=(\xi,\mathbf{T}(d))$, welche diejenige topologisierte Flussfunktion $\Xi=(\xi,\mathbf{T}(d))$ sei,
die dadurch bestimmt ist,
dass $\mathbf{T}(d)$ die von $d$ induzierte metrische Topologie ist.
Wohingegen die Inklusion
\begin{equation}\label{ylos}
\mbox{coloc}^{+}(\Xi)\subset(\Delta_{\Psi}^{d})^{-1}(\mathbb{R}^{+})\ ,
\end{equation}
evident ist,
in der $\Delta_{\Psi}^{d}$ gem"ass der Festlegung (\ref{aalzbbb}) nach der Definition 1.2 das 
Aufl"osungsfeld der Flussfunktion $\Psi$ bez"uglich einer Distanzfunktion $d$ notiert, die (\ref{aathbb}) erf"ullt.\index{Aufl\"osungsfeld einer Wellenfunktion}
Alle ultrasensitven Zust"ande der topologisch strukturierten 
Flussfunktion $\Xi=(\xi,\mathbf{T}(d))$ sind sensitive Zust"ande der metrisierten
Flussfunktion $(\xi,d)$ im Sinn der Definition 1.2. Dabei sind, freiz"ugig gesprochen, ultrasensitve Zust"ande 
einer jeweiligen, mit der metrischen Topologie $\mathbf{T}(d)$ topologisierten Flussfunktion $\Xi=(\xi,\mathbf{T}(d))$
Rarit"aten. Denn, in der zeitlichen Interpretation 
stellt sich das Szenario ultrasensitver Zust"ande folgendermassen dar:
Da $\mathbf{T}(d)$ dem Hausdorffschen Trennungsaxiom gen"ugt, ist ein 
ultrasensitver Zustand
des Paares $\Xi=(\xi,\mathbf{T}(d))$ ein solcher Zustand, dem zwei separable Zust"ande verschiedener 
Trajektorien gleichzeitig immer wieder beliebig nahe kommen.
\newline
Genau dann, wenn ein
ultrasensitver Zustand $z\in\mbox{coloc}^{+}((\xi,\mathcal{A}))\cap a$
in einem Attraktor $a\in\mbox{{\bf @}}(\xi,\mathcal{A})$ 
einer strukturierten Flussfunktion $(\xi,\mathcal{A})$ im weiteren Sinn ist,
nennen wir den Attraktor $a$ einen ultrachaotischen Attraktor.\index{ultrachaotischer Attraktor}
Als Ultrachaos bezeichnen wir den Sachverhalt, dass in einem Attraktor ein
ultrasensitver Zustand existiert.\index{Ultrachaos}
\newline
Es mag die Ultrasensitivit"at 
einer mit einer metrischen Topologie strukturierten 
Flussfunktion
zwar eine Rarit"at sein. Indess, die Ultrasensitivit"at ist 
insofern nachvollziehbar, als sie uns den kinematischen Mechanismus zeigt, wie Sensitivit"at 
im Fall der Ultrasensitivit"at 
zustande 
kommt. Die Ultrakolokalisation betrachtet das Zustandekommen der Ultrasensitivit"at
aus der Perspektive, die der Entwicklungsrichtung der Lokalisiserungsaufl"osung entgegengesetzt blickt.
\newline
Ist nun aber die Ultrasensitivit"at eines Zustandes 
einer solchen strukturierten Flussfunktion $(\xi,\mathcal{A})$, deren Strukturierung 
$\mathcal{A}$ so beschaffen ist, dass es Zust"ande $\omega\in\bigcup\mathcal{A}$ gibt, f"ur 
welche
die echte Inklusion
\begin{equation}\label{yftislf}
\bigcap\mathcal{A}_{\{\omega\}}\supsetneqq\{\omega\}
\end{equation}
gilt, deutlich 
weniger dramatisch?
Selbst in dem Fall, dass
eine
topologisierte Flussfunktion $(\xi,\mathbf{T}(\mathbf{P}_{2}\xi))$ vorliegt,
kann $(\xi,\mathbf{T}(\mathbf{P}_{2}\xi))$ 
auf h"ochst unspektakul"are Weise
Zust"ande in 
$\mbox{coloc}^{0}((\xi,\mathbf{T}(\mathbf{P}_{2}\xi)))$
haben:
Es gibt ja Topologien $\mathbf{T}(\mathbf{P}_{2}\xi)$, f"ur die keineswegs notwendigerweise 
f"ur alle Zust"ande $x\in \mathbf{P}_{2}\xi$ die punktweise Lokalisierbarkeit in Form der Gleichung
\begin{equation}\label{ydnulx}
\bigcap\mathbf{T}(\mathbf{P}_{2}\xi)_{\{x\}}=\{x\}
\end{equation}
gegeben ist. 
Exakt jede ${\rm T}_{0}$-Topologie\index{${\rm T}_{0}$-Topologie} ist eine Topologie,
f"ur die f"ur keinen der Punkte des zu ihr "aquivalenten topologischen Raumes 
punktweise Lokalisierbarkeit gem"ass 
der Identit"at (\ref{ydnulx}) gegeben ist. Diese 
Identit"at formuliert
das Hausdorffsche Trennungsaxiom\index{Hausdorffsches Trennungsaxiom} f"ur den einzelnen Punkt.
Demgegen"uber ist exakt jede ${\rm T}_{1}$-Topologie\index{${\rm T}_{1}$-Topologie}
eine Topologie von der Art, dass Punkte des zu ihr "aquivalenten topologischen Raumes existieren,
in denen die punktweise Lokalisierbarkeit gem"ass (\ref{ydnulx}) verletzt ist.
Die ${\rm T}_{2}$-Topologien\index{${\rm T}_{2}$-Topologie}
sind gerade diejenigen Topologien $\mathbf{T}$, f"ur welche die punktweise Lokalisierbarkeit in dem Sinn
global gegeben ist, dass f"ur jeden ihrer Punkte $x$ die 
Gleichung $\{x\}=\bigcap\mathbf{T}_{\{x\}}$ gilt. Die 
globale Version der punktweisen Lokalisierbarkeit formuliert
also das Hausdorffsche Trennungsaxiom.
Die Gegebenheiten und Einschr"ankungen, die ohne die G"ultigkeit des   
Hausdorffschen Trennungsaxiomes vorliegen, hat S.Willard in \cite{will}
untersucht. Wenn wir allgemeine Mengensysteme statt Topologien betrachten, so fallen auch
die in \cite{will} dargestellten Gegebenheiten weg, weswegen die 
Gegebenheit globaler punktweiser Lokalisierbarkeit
f"ur die Betrachtung strukturierter Flussfunktionen $\Lambda=(\xi, \mathcal{A})$
eine um so gr"ossere Rolle spielt.\newline Zum Beispiel formulieren wir mehrere der S"atze
"uber Flussfunktionen $\Lambda=(\xi, \mathcal{A})$ 
f"ur den Fall, dass 
ein gewisser zentraler Teil $(\mathbf{P}_{2}\xi)^{\circ}\subset \mathbf{P}_{2}\xi$ des Zustandsraumes 
so beschaffen ist, dass
dessen durch das Spurmengensystem $\mathcal{A}\cap (\mathbf{P}_{2}\xi)^{\circ}$ 
gegebene Strukturierung
ein 
deterministisches Mengensystem ist. 
Im ersten Kapitel der Konzepte der abstrakten Ergodentheorie \cite{rabe} untersuchten 
wir z.B. die invarianten gemeinsamen Topologien $\mathbf{T}(\Xi)=\mathbf{T}(\Xi)^{c}$ 
jeweiliger 
Autobolismenmengen $\Xi$, die nach dem Korollar $1.2^{1}$ in dem, dass sie
selbstduale Topologien sind, ${\rm T}_{0}$-Topologien sind. 
Nach dem Satz $1.1^{1}$ des ersten Teiles der allgemeinen Ergodentheorie
sind alle 
nicht mit Potenzmengen identischen,
selbstdualen Topologien $\mathbf{T}$ solche Topologien, f"ur die nicht f"ur alle ihre Punkte $x$ 
punktweise Lokalisierbarkeit 
in Form der Identit"at $\bigcap\mathbf{T}_{\{x\}}=\{x\}$
gegeben ist.\newline
Die punktweise Lokalisierbarkeit 
der Zust"ande des Zustandsraumes 
einer jeweiligen  
strukturierten Flussfunktion $\Lambda=(\xi, \mathcal{A})$
modelliert die deterministische Aufl"osbarkeit\index{deterministische Aufl\"osbarkeit}  
der durch die Flussfunktion $\xi$ dargestellten
Vorg"ange w"ahrend des
kontinuierlichen Ablaufes der Zeit, die durch jeweilige Zustandsentwicklungen $\xi(z,\mbox{id})$ 
eines jeweiligen Zustandes $z$ des Zustandsraumes $\bigcup\mathcal{A}$
objektiviert sind.
Die Zustandsraum"uberdeckung $\mathcal{A}$ stellt hierbei dessen vorgegebene  
Aufl"osbarkeit dar, weshalb wir die zweite Komponente einer 
strukturierten Flussfunktion $\Lambda=(\xi, \mathcal{A})$, das Mengensystem $\mathcal{A}$,
als deren Strukturierung bezeichnen.\index{Strukturierung einer strukturierten Flussfunktion}
Die Modellierung der deterministischen Aufl"osbarkeit
gerade durch Mengensysteme mit punktweiser Lokalisierung bildet demnach den Hintergrund, vor dem
wir 
jedes Mengensystem $\mathcal{Z}$
genau dann als deterministisches Mengensystem\index{deterministisches Mengensystem} bezeichnen, wenn
f"ur es die Implikation 
\begin{equation}\label{ydnalx}
z\in\bigcup\mathcal{Z} \Rightarrow \bigcap\mathcal{Z}_{\{z\}}=\{z\}
\end{equation}
wahr ist.
Sei $\omega$
ein Zustand des Zustandsraumes einer 
strukturierten Flussfunktion $\Lambda=(\xi, \mathcal{A})$, f"ur den
die
Menge $\bigcap\mathcal{A}_{\{\omega\}}\not=\{\omega\}$ "uber $\{\omega\}$ hinaus
ausdehnt ist. Zu $\omega$ gebe es zwei Zust"ande $x,y\in \mathbf{P}_{2}\xi$,
deren Entwicklungen so verlaufen, dass f"ur jede als Zeit interpretierbare Zahl $t_{\star}\in\mathbb{R}^{+}$
eine Zahl $t_{\star}<t_{+}$ und eine Zahl $t_{-}<-t_{\star}$ gibt, 
f"ur die
$\xi(x,t_{\pm})$ und $\xi(y,t_{\pm})$ 
gleichzeitig in $\mathcal{A}_{\{\omega\}}$ sind.
Dann ist der Zustand $\omega$ genau dann ein 
ultrasensitver Zustand der Flussfunktion $\Lambda=(\xi, \mathcal{A})$, wenn 
$x$ und $y$ bez"uglich der Strukturierung $\mathcal{A}$ separierbar sind. Nehmen wir eine bez"uglich 
$\mathcal{A}$ kommutativ Cantor-stetige Flussfunktion $(\xi, \mathcal{A})$, f"ur welche
das Mengensystem nicht trivialer verallgemeinerter Zimmer
$$[[\xi]]_{\mathcal{A}}\setminus[\xi]$$ nicht leer ist, sodass
die Restriktion 
$$\xi|(\bigcup[[\xi]]_{\mathcal{A}}\setminus[\xi])\times\overline{\mathbb{R}}$$
eine 
auf bez"uglich dem Spurmengensystem
$\mathcal{A}\cap(\bigcup[[\xi]]_{\mathcal{A}}\setminus[\xi])$ kommutativ Cantor-stetige Flussfunktion 
ist, so ist 
nach dem Satz $3.3^{1}$ des ersten Teiles der abstrakten Ergodentheorie
\index{verallgemeinerter insensitver Ergodensatz} 
das Paar
$$\Gamma=\Bigl(\xi|(\bigcup[[\xi]]_{\mathcal{A}}\setminus[\xi])\times\overline{\mathbb{R}},[[\xi]]_{\mathcal{A}}\setminus[\xi]\Bigr)$$
eine strukturierte Flussfunktion, wobei auf h"ochst triviale Weise
\begin{equation}\label{trivialia}
\bigcup([[\xi]]_{\mathcal{A}}\setminus[\xi])=\mbox{coloc}^{0}(\Gamma)
\end{equation}
Die Menge aller ultrasensitven Zust"ande des strukturierten
Flussfunktion $\Gamma$ hingegen ist leer, denn 
es gilt f"ur alle 
Flussfunktionen $(\xi,\mathcal{P})$ im weiteren Sinn,
deren Strukturierung $\mathcal{P}\in\mathbf{part}(\mathbf{P}_{2}\xi)$ den Zustandsraum partioniert
\begin{equation}
\mbox{coloc}^{+}((\xi,\mathcal{P}))=\emptyset\ .
\end{equation}
Wohingegen 
keineswegs f"ur 
jede Flussfunktionen $(\xi,\mathcal{A})$ im weiteren Sinn, f"ur welche die 
Hausdorffsche Trennbarkeit in dem Sinn nirgends gegeben ist, dass die Implikation
$$z\in\mathbf{P}_{2}\xi\Rightarrow\mathcal{A}_{\{z\}}\setminus\{z\}\not=\emptyset$$
zutrifft, gelten muss, dass $\mbox{coloc}^{+}((\xi,\mathcal{A}))$ leer ist.
Die Ultrasensitivt"at kann durch Strukturierungsvergr"oberungen beliebig abgeschw"acht
werden, weswegen wir anstreben, eine Abwandlung der Ultrakolokalisation zu finden,
welche die verallgemeinerte zustandsweise Sensitivt"at darstellt.\newline
F"ur jedes Mengensystem $\mathcal{Z}$
sei f"ur jede Menge ${\rm X}$ 
\begin{equation}
\mathcal{Z}_{({\rm X})}:=\{{\rm Z}\cup{\rm X}:{\rm Z}\in\mathcal{Z}_{{\rm X}}\}\cup
\{{\rm Z}\not\in\mathcal{Z}_{{\rm X}}\}
\end{equation}
die Vergr"oberung des Mengensystemes $\mathcal{Z}$ zu dem Mengensystem $\mathcal{Z}_{({\rm X})}$,
exakt welches wir die Vergr"oberung des Mengensystemes $\mathcal{Z}$ auf die Menge ${\rm X}$
nennen.\index{Vergr\"oberung eines Mengensystemes auf eine Menge}
Die echte Inklusion (\ref{yftislf}) ist demnach f"ur alle Vergr"oberungen $\mathcal{A}_{({\rm A})}$
der Strukturierung 
$\mathcal{A}$ auf ${\rm A}\in \mathcal{A}_{\{\omega\}}$ erf"ullt,
wenn $\{\omega\}\subsetneqq{\rm A}$ gilt. 
Genau dann, wenn ${\rm X}\subset\bigcup\mathcal{Z}$ ist,
ist die Vereinigung "uber die Vergr"oberung des Mengensystemes $\mathcal{Z}$ zu dem Mengensystem $\mathcal{Z}_{({\rm X})}$
die Invariante 
$$\bigcup\mathcal{Z}_{({\rm X})}=\bigcup\mathcal{Z}\ ,$$
weshalb wir exakt jede Vergr"oberung eines Mengensystemes, welche die
Vereinigung "uber es invariant l"asst, als eine vereinigungsinvariante 
Vergr"oberung des urspr"unglichen Mengensystemes bezeichnen. Es gilt 
offensichtlich die Implikation 
$$z\in\bigcup\mathcal{Z}\ \land\ \bigcap\mathcal{Z}_{\{z\}}={\rm H}(z)\Rightarrow$$
$$\mathcal{Z}_{({\rm H}(z))}=\mathcal{Z}\ ,$$
welche allerdings nicht umkehrbar ist. Damit haben wir auf der Klasse der Paare 
$(\mathcal{Z},{\rm X})$
von Mengensystemen $\mathcal{Z}$ und Mengen ${\rm X}$ den universellen Vergr"oberungsoperator\index{Vergr\"oberungsoperator} 
$$\mathbf{P}_{1\quad(\mathbf{P}_{2})}$$
festglegt, den wir folgendermassen auf die Klasse der Flussfunktionen im weiteren Sinn
"ubertragen: Es sei f"ur jede strukturierte Flussfunktion $\Lambda=(\xi, \mathcal{A})$
im weiteren Sinn und f"ur jede Menge ${\rm X}$
\begin{equation}\label{tislf}
\Lambda_{({\rm X})}:=(\xi, \mathcal{A}_{({\rm X})})
\end{equation}
die Vergr"oberung der Flussfunktion $\Lambda$ auf ${\rm X}$.\index{Vergr\"oberung einer Flussfunktion im weiteren Sinn auf eine Menge}
Mit Hilfe der Vergr"oberung jeweiliger Flussfunktionen im weiteren Sinn
k"onnen wir die Ultrakolokalisation zur Kolokalisation abschw"achen und finden 
damit endlich die verallgemeinerte Form der zustandsweisen Sensitivit"at\index{zustandsweise Sensitivit\"at}  
einer strukturierten Flussfunktion im weiteren Sinn und damit den allgemeinen
Begriff des chaotischen Attraktors:
\newline\newline
{\bf Definition 3.6: Kolokalisation und zustandsweise Sensitivit"at}\newline
{\em F"ur jede strukturierte Flussfunktion $\Lambda=(\xi, \mathcal{A})$
im weiteren Sinn
bezeichnen wir mit} $\mbox{coloc}_{-}(\Lambda)$
{\em die Menge
aller ihrer sensitiven Zust"ande, wobei
jeder Zustand $\omega\in\mathbf{P}_{2}\xi$
exakt in folgendem Fall als ein sensitiver Zustand\index{sensitiver Zustand einer Flussfunktion 
im weiteren Sinn} der Flussfunktion im weiteren Sinn $\Lambda$ gelte:
F"ur eine Folge von Lokalisierungen 
$$\{{\rm L}_{j}\}_{j\in\mathbb{N}}\in\mathcal{A}_{\{\omega\}}^{\mathbb{N}}\ ,$$
f"ur die 
\begin{equation}
\bigcap_{j\in\mathbb{N}}{\rm L}_{j}=\bigcap\mathcal{A}_{\{\omega\}}
\end{equation}
gilt, existiert eine Folge   
$\{(z_{j},y_{j})\}_{j\in\mathbb{N}}$ von Zustandspaaren,
f"ur welche f"ur jeden Index $j\in\mathbb{N}$ die jeweilige nebenbedingte vergr"oberte Ultrakolokalisation}
\begin{equation}\label{thiseva}
\begin{array}{c}
z_{j},y_{j}\ \sim>_{\Lambda_{({\rm L}_{j})}}\omega\\
\land\\
\exists\ t\in\mathbb{R},\ {\rm L}\in\mathcal{A}_{\{\xi^{t}(\omega)\}}:\ \forall\ j\in\mathbb{N}\
\xi^{t}(z_{j}),\xi^{t}(y_{j})\not\in {\rm L}
\end{array}
\end{equation}
{\em gegeben ist. Exakt jede derartige Folge $\{(z_{j},y_{j})\}_{j\in\mathbb{N}}$ von 
Zustandspaaren bezeichnen wir als im Zustand 
$\omega$ kolokalisierend.\index{kolokalisierende Zustandspaarfolge}}   
\newline
\newline
Es f"ugt sich in unsere Benennungssystematik, nun exakt jeden Attraktor
$a\in\mbox{{\bf @}}(\xi,\mathcal{A})$ 
einer strukturierten Flussfunktion im weiteren Sinn $(\xi,\mathcal{A})$,
welcher einen sensitiven Zustand $z\in\mbox{coloc}_{-}((\xi|a\times\mathbb{R},\mathcal{A}\cap a))$ hat
als einen chaotischen Attraktor zu bezeichnen.\index{chaotischer Attraktor}
Wobei wir mit Nachdruck darauf hinweisen, dass im allgemeinen die Differenz
\begin{equation}\label{thisava}
\mbox{coloc}_{-}((\xi,\mathcal{A}))\cap a\not=
\mbox{coloc}_{-}((\xi|a\times\mathbb{R},\mathcal{A}\cap a))\\
\subset \mbox{coloc}_{-}((\xi,\mathcal{A}))\cap a
\end{equation}
zu beachten ist.
Wir betonen hierbei ausserdem, dass wir es sein kann, dass eine Zustandspaarfolge
$\{(z_{j},y_{j})\}_{j\in\mathbb{N}}$ im Zustand 
$\omega$ kolokalisiert, deren Folgengliedkomponenten $z_{j}$ oder $y_{j}$
auf der Trajektorie $\xi(\omega,\mathbb{R})$ liegen, die durch $\omega$ verl"auft.
Ausserdem haben wir nicht ausgeschlossen, dass auch eine Zustandspaarfolge jeweils komponentengleicher Zustandspaare
$$\{(z_{j},y_{j})\}_{j\in\mathbb{N}}=\{(z_{j},z_{j})\}_{j\in\mathbb{N}}$$
als in einem Zustand kolokalisierend $\omega$ gelten kann. Falls 
$\{(z_{j},z_{j})\}_{j\in\mathbb{N}}$ im Zustand $\omega$ kollokalisiert, sagen wir einfach, dass
die Zustandsfolge $\{z_{j}\}_{j\in\mathbb{N}}$ in
$\omega$ kolokalisiere.
Diese M"oglichkeit komponentengleicher Zustandspaare ist f"ur eine separiert kolokalisierende Zustandspaarfolge 
$\{(\hat{z}_{j},\hat{y}_{j})\}_{j\in\mathbb{N}}$ ausgeschlossen:
Genau dann, wenn 
f"ur eine Zustandspaarfolge $\{(\hat{z}_{j},\hat{y}_{j})\}_{j\in\mathbb{N}}$ statt (\ref{thiseva})
f"ur alle $j\in\mathbb{N}$ die vergr"oberte Ultrakolokalisation
\begin{equation}
\begin{array}{c}
\hat{z}_{j},\hat{y}_{j}\ \sim>_{\Lambda_{({\rm L}_{j})}}\omega\\
\land\\
\exists\ {\rm L}(\hat{z}_{j}),{\rm L}(\hat{y}_{j})\in\mathcal{A}:\  \hat{z}_{j}\not\in{\rm L}(\hat{y}_{j})\ni \hat{y}_{j}\ \land\ \hat{y}_{j}\not\in{\rm L}(\hat{z}_{j})\ni \hat{z}_{j}
\end{array}
\end{equation}
mit der Nebenbedingung der Separiertheit der jeweiligen Komponenten $\hat{z}_{j}$ und $\hat{y}_{j}$
gilt, nennen wir die Zustandspaarfolge $\{(\hat{z}_{j},\hat{y}_{j})\}_{j\in\mathbb{N}}$ im Zustand 
$\omega$ separiert kolokalisierend.\index{separiert kolokalisierende Zustandspaarfolge}
Wir erl"autern die Definition 3.6 ferner in folgenden Punkten:
\newline
\newline
{\bf 1.} Diese Verallgemeinerung des zustandsweisen Sensitivit"atsbegriffes gen"ugt nun dem
Permanenzprinzip. 
Denn die f"ur Wellenfunktionen\footnote{Die Definition 1.2 gemahnt uns daran, dass wir die allgemeine punktweise Sensitivit"at
noch nicht f"ur Wellenfunktionen verfasst haben.
Offenbar ist die symmetrisch angelegte Relation $\sim_{\Lambda}$ der Ultrakolokalisation
noch weiter 
f"ur Wellenfunktionen anstelle der 
Flussfunktionen
generalisierbar.
Die entsprechende Generalsierungsrichtung interessiert uns aber gegenw"artig nicht.} verfasste
Sensitivit"at bez"uglich zweier positiver Funktionen\index{Sensitivit\"at bez\"uglich zweier nicht-negativer Funktionen}
der 
Definition 1.2 ist
auf folgende Weise
als Verallgemeinerung der Sensitivit"at gem"ass dieser Definition 3.6 auffassbar, wenn
der Fall vorliegt, dass erstens
die jeweilige Wellenfunktion
eine Flussfunktion $\Psi$ ist und dabei zweitens die beiden 
Distanzfunktionen $d_{A}$ und $d_{B}$ der Definition 1.2
zwei identische 
Distanzfunktionen $d\in[0,\infty]^{\mathbf{P}_{2}\Psi\times \mathbf{P}_{2}\Psi}$
des Zustandsraumes $\mathbf{P}_{2}\Psi$
sind; welche demnach der Bedingung (\ref{aathbb}) gen"ugen: F"ur die Zustandsraum"uberdeckung durch das Mengensystem der offenen Pseudokugeln\index{Mengensystem der offenen Pseudokugeln}\index{offene Pseudokugel}
\begin{displaymath}
\mathcal{B}(\mathbf{P}_{2}\Psi,d)=\Bigl\{\{x\in \mathbf{P}_{2}\Psi:d(x,y)<r\}:(y,r)\in\mathbf{P}_{2}\Psi\times[0,\infty]\Bigr\}\ ,
\end{displaymath}
zu denen wir auch die leere Menge z"ahlen,
ist
\begin{equation}\label{ydnuly}
\mbox{coloc}_{-}((\Psi,\mathcal{B}(\mathbf{P}_{2}\Psi,d)))=\bigcup(\Delta_{\Psi}^{d})^{-1}(\mathbb{R}^{+})\ ,
\end{equation}
wobei wir mit $\Delta_{\Psi}^{d}$ das 
Aufl"osungsfeld der Flussfunktion $\Psi$ bez"uglich der Distanzfunktion\index{Distanzfunktion} $d$
gem"ass der Festlegung (\ref{aalzbbb}) notieren.\index{Aufl\"osungsfeld einer Wellenfunktion}
Die offenen Pseudokugeln k"onnen einelementig sein. Ist $\{x\}$ eine einelementige offene Pseudokugel, 
durch welche der Zyklus $\Psi(x,\mathbb{R})$ verl"auft,
so ist dieser Zyklus
$\Psi(x,\mathbb{R})$ in der Menge sensitiver Zust"ande 
$\mbox{coloc}_{-}((\Psi,\mathcal{B}(\mathbf{P}_{2}\Psi,d)))$ enthalten.
\newline
{\bf 2.} Die Verallgemeinerung des Sensitivit"atsbegriffes der Definition 3.6 gen"ugt also zwar dem
Permanenzprinzip, nichtsdestotrotz geht eine Differenz
zwischen intuitiven Vorstellungen von der Sensitivit"at
und der gem"ass der Definition 3.6 verallgemeinerten Sensitivit"at
auf. Wie wir sehen:
Schon die 
Elemente der Menge $\bigcup(\Delta_{\Psi}^{d})^{-1}(\mathbb{R}^{+})$, die
wir als die
bez"uglich der Distanzfunktion $d$ sensitiven 
Zust"ande der Flussfunktion $\Psi$  
ansehen,  
zeigen
ein 
Verhalten, das
von dem Verhalten abweicht, das wir intuitiv mit dem Sensitivit"atsph"anomen meinen und 
das wir
von den durch nat"urliche Zustandsraumtopologien
beschriebenen Szenarien 
kennen.
Und zwar besteht diese Abweichung in gerade den jeweiligen Hinsichten, in denen die Distanzfunktion $d$ 
davon abweicht, eine
Metrik zu sein. Dementsprechend ist die
Differenz zwischen 
intuitiven Vorstellungen von der Sensitivit"at und der allgemeinen Sensitivit"at
gem"ass der Definition 3.6 
erst recht zu beachten.
\newline
{\bf 3.}
Nichtsdestotrotz 
ist der gem"ass der Definition 3.6 verfasste Begriff der allgemeinen und zustandsweisen
Sensitivit"at in folgender Hinsicht nachvollziehbar:
Wie die Ultrakolokalisation beobachtet die 
Kolokalisation gleichsam das Zustandekommen der zustandsweisen Sensitivit"at.
Sie schaut dabei
aus derjenigen Perspektive nach vorne, welche aus der Entwicklungsrichtung 
betrachtet, die zu der jeweiligen Entwicklungsrichtung entgegengesetzt ist, 
in der die Sensitivit"at als Lokalisierungsaufl"osung erscheint.
\newline
{\bf 4.} Die Reichhaltigkeit an sensitiven Zust"anden, dass $\mbox{coloc}_{-}(\Lambda)\not=\emptyset$
ist, impliziert f"ur jede strukturierte Flussfunktion $\Lambda=(\xi, \mathcal{A})$ im weiteren Sinn die
mengenweise Sensitivit"at gem"ass der Definition 3.3. Das ist offensichtlich.
Es gilt die Implikation
\begin{equation}\label{ebsiht}
\begin{array}{c}
\mbox{coloc}_{-}(\Lambda)\not=\emptyset\Rightarrow\\
\exists\ t\in\mathbb{R}:\ [\mathbf{cl}_{\mathcal{A}},\xi^{t}]\not=\emptyset
\end{array}
\end{equation} 
f"ur jede strukturierte Flussfunktion $\Lambda=(\xi, \mathcal{A})$ im weiteren Sinn.
\newline\newline
Wir schliessen diesen Abschnitt mit der 
folgenden Bemerkung, welche eine Symmetriebedingung benennt, unter der die Implikation (\ref{ebsiht}) umkehrbar ist:
\newline\newline
{\bf Satz 3.6:\newline Die "Aquivalenz mengenweiser und zustandsweiser Sensitivit"at}\newline
{\em Es sei $\Lambda=(\xi, \mathcal{A})$ eine Poincar$\acute{e}$sche Entwicklung\index{Poincar$\acute{e}$sche Entwicklung}
gem"ass der Definition} 2.3, {\em deren
Strukturierung $\mathcal{A}$ gem"ass} (\ref{ydnalx}) {\em deterministisch\index{deterministisches Mengensystem} ist. 
Genau dann, wenn $\Lambda$ 
mengenweise Sensitivit"at gem"ass der Definition} 3.3 {\em zeigt, hat $\Lambda$ sensitive Zust"ande: Die "Aquivalenz}
\begin{equation}\label{topit}
\begin{array}{c}
\exists\ t\in\mathbb{R}:\ [\mathbf{cl}_{\mathcal{A}},\xi^{t}]\not=\emptyset\\
\Leftrightarrow\\
\mbox{coloc}_{-}(\Lambda)\not=\emptyset
\end{array}
\end{equation}
{\em gilt f"ur alle 
Poincar$\acute{e}$schen Entwicklungen $(\xi, \mathcal{A})$ mit deterministischer Strukturierung.}
\newline\newline
{\bf Beweis:}\newline
Dass die Existenz sensitiver Zust"anden
Sensitivit"at impliziert, ist evident. Im Fall, dass
$\Lambda$ eine Poincar$\acute{e}$sche Entwicklung ist,
f"ur die
es eine reelle Zahl $t\in\mathbb{R}$ gibt, f"ur die
$$[\mathbf{cl}_{\mathcal{A}},\xi^{t}]\not=\emptyset$$
ist, gibt es nach dem Explizierungssatz $(3.8)^{1}$ des ersten Teiles 
der Konzepte der abstrakten Ergodentheorie \cite{rabe}
eine nicht leere Menge 
${\rm A}\in \mathcal{A}$,
f"ur die entweder die Implikation
$$\tilde{{\rm A}}\in \mathcal{A}\Rightarrow\tilde{{\rm A}}\not\subset\xi^{t}{\rm A}$$
oder die Implikation
$$\tilde{{\rm A}}\in \mathcal{A}\Rightarrow\tilde{{\rm A}}\not\subset(\xi^{t})^{-1}{\rm A}$$ 
wahr ist, wobei, da $\xi$ eine Entwicklung ist,
$$(\xi^{t})^{-1}\in\xi^{\mathbb{R}} $$
ist. Ohne Einschr"ankung der Allgemeinheit k"onnen wir also annehmen, dass
die erste der beiden letztgenannten Implikationen der Fall ist. Es gilt 
daher f"ur jeden Zustand $z$ der Menge $\xi^{t}{\rm A}$
die Implikation
\begin{equation}\label{ruelle}
{\rm Z}\in \mathcal{A}_{\{z\}}\Rightarrow(\xi^{t})^{-1}{\rm Z}\not\subset{\rm A}\ni(\xi^{t})^{-1}(z)\ .
\end{equation}
F"ur jede Folge von Lokalisierungen 
$$\{{\rm L}_{j}\}_{j\in\mathbb{N}}\in\mathcal{A}_{\{z\}}^{\mathbb{N}}\ ,$$
f"ur die wegen der Deterministik der Strukturierung $\mathcal{A}$   
$$\bigcap_{j\in\mathbb{N}}{\rm L}_{j}=\bigcap\mathcal{A}_{\{z\}}=\{z\}$$
ist, gilt daher
dass eine Folge  
$\{z_{j}\}_{j\in\mathbb{N}}$ von Zust"anden existiert,
f"ur welche f"ur jeden Index $j\in\mathbb{N}$ 
$$z_{j}\in ((\xi^{t})^{-1}{\rm L}_{j})\setminus{\rm A}\not=\emptyset$$  
ist, weil wegen der Implikation (\ref{ruelle}) die Aussage
$$(\xi^{t})^{-1}{\rm L}_{j}\not\subset{\rm A}$$
wahr ist.
Weil $\Lambda$ eine Poincar$\acute{e}$sche Entwicklung ist, trifft gem"ass der Definition 2.3 jede Zustandsentwicklung
jede Lokalisierung ${\rm L}\in \mathcal{A}$ immer wieder, wenn die Zustandsentwicklung
auf einer  
Trajektorie verl"auft, welche die Lokalisierung ${\rm L}$  
schneidet. Es gilt (\ref{thisbec}).
Daher trifft die Zustandsentwicklung des Zustandes $z_{j}$ 
die Lokalisierung ${\rm L}_{j}$
immer wieder, sodass
die jeweilige vergr"oberte Ultrakolokalisation
$$z_{j},z_{j}\ \sim>_{\Lambda_{({\rm L}_{j})}}z$$
f"ur alle $j\in\mathbb{N}$ gilt.
Hierbei gibt es aber
die Lokalisierung
${\rm A}\in\mathcal{A}$ und da $\xi$ eine Entwicklung ist, existiert die reelle Zahl
$\hat{t}\in\mathbb{R}$, f"ur die
$$(\xi^{t})^{-1}=\xi^{\hat{t}}$$
ist,
sodass f"ur alle $j\in\mathbb{N}$ die Aussage
\begin{displaymath}
\xi^{\hat{t}}(z_{j})\not\in{\rm A}\
\land\ \xi^{\hat{t}}(z)\in{\rm A} 
\end{displaymath}
wahr ist.\newline
{\bf q.e.d.}
\newline
\newline
Wegen der "Aquivalenz mengenweiser und zustandsweiser Sensitivit"at
sind also exakt alle Attraktoren $a\in\mbox{{\bf @}}(\xi,\mathcal{A})$ 
einer Poincar$\acute{e}$schen Entwicklung $(\xi,\mathcal{A})$,
welche einen sensitiven Zustand $z\in\mbox{coloc}_{-}((\xi|a\times\mathbb{R},\mathcal{A}\cap a))$ haben,
chaotische Attraktoren.\index{chaotischer Attraktor}
Wenn $(\xi,\mathcal{A})$ eine Poincar$\acute{e}$sche Entwicklung ist
koinzidieren der Begriff des chaotischen Attraktors
und des 
chaotischen Attraktors von nulltem Grad.\index{chaotischer Attraktor von nulltem Grad}
Dass $(\xi,\mathcal{A})$ Poincar$\acute{e}$sch ist, heisst noch nicht, dass
$(\xi|a\times\mathbb{R},\mathcal{A}\cap a)$ Poincar$\acute{e}$sch ist.
Ist 
$(\tilde{\xi},\tilde{\mathcal{A}})$ keine Poincar$\acute{e}$sche Entwicklung, 
so ist 
wegen der Implikation
(\ref{ebsiht}) 
jeder Attraktor $b$, in dem ein sensitiver Zustand 
$y\in\mbox{coloc}_{-}((\xi|b\times\mathbb{R},\mathcal{A}\cap b))$
ist,
ein chaotischer Attraktor. Indess, er hat die Besonderheit, 
ein 
chaotischer Attraktor zu sein,
der einen sensitiven Zustand
im Sinn der Definition 3.6 hat. 
\newline
Zu beachten ist, dass f"ur die 
Attraktoren $a\in \mbox{{\bf @}}(\xi,\mathcal{A})$ einer beliebigen strukturierten Flussfunktion 
im weiteren Sinn $(\xi,\mathcal{A})$ im allgemeinen
die Differenz (\ref{thisava})
besteht, was wir folgendermassen illustrieren:
Es kann sein, dass ein sensitiver
Zustand in einem Zustandsraum liegt, den lauter geschlossene glatte
Trajektorien partionieren, wie uns dies der Blick in die 
endlichdimensionalen reellen R"aume der Dimenison $n>1$ veranschaulicht. Die 
zu diesem Szenario geh"orende Entwicklung ist dabei bez"uglich der nat"urlichen 
Zustandsraumtopologie
Poincar$\acute{e}$sch.
Der sensitive Zustand liegt dann in einem trivialen Zimmer, das ein
Attraktor ist. Nichtsdestotrotz ist derselbe nicht von nulltem Grad chaotisch, weil
in demselben keine sensitive Lokalisierung liegen kann.\newline
Die
Existenz sensitiver Zust"ande f"ur sich allein genommen
bedeutet noch nicht, dass Chaos vorliegt: Die 
Existenz sensitiver Zust"ande in einem Attraktor
heisst noch nicht, dass Sensitivit"at
innerhalb dieses Attraktors vorliegt, dass n"amlich solche
Sensitivit"at gegeben ist, die sich auf dessen Spurstrukturierung bezieht. 
\newline
Ein
ultrasensitver Zustand kann in einem endlichdimensionalen und reellen Zustandsraum, den lauter geschlossene 
glatte
Trajektorien partionieren, offensichtlich nicht auftreten. Wir wollen uns im folgenden
auf die Ultrasensitivit"at und die Ultrakolokation
konzentrieren.
\newline

\chapter{Ordnungsspuren im Ultrachaos}
{\small
Wir beginnen unsere Suche nach der Ordnung im Ultrachaos, indem 
wir im ersten Abschnitt 
den Gegenspieler der verallgemeinerten Sensitivit"at, die Komanenz, untersuchen, die
eine Form der Kontinuit"at ist. 
Komanenz kann die Sensitivit"at und damit Chaos von vorne herein ausschliessen.
Ein Abwandlung der Komanenz hingegen, die wir als
konverse Komanenz bezeichnen, schliesst 
weder Sensitivit"at noch
Chaos von vorne herein aus und stellt dabei doch 
die Grundlage daf"ur dar, dass sich folgende Ordnungspuren zeigen: 
Im zweiten Abschnitt besch"aftigt uns der Sachverhalt, dass
die Relation $\sim_{\Lambda}$ sich 
als eine "Aquivalenzrelation erweist, 
wenn die Flussfunktion $\Lambda$ bestimmte Kontinuit"atsvoraussetzungen 
erf"ullt und sie -- im wesentlichen -- konvers komanent ist.}
\section{Komanenz als Sensitivit"atsantagonist}
\subsection{Wo siedelt die logische Konstitution der Komanenz?}
Gleich, ob wir dies seltsam oder bemerkenswert nennen wollen:
Zu der Formulierung deutlich weniger trivialer Aussagen "uber die generalisierte Sensitivit"at,
als es beispielsweise die
Aussage (\ref{trivialia}) ist,
leitet uns die Betrachtung gerade des 
Gegenspielers der Sensitivit"at, n"amlich die Betrachtung der allgemeinen 
Komanenz.newline
In der "Uberschrift \glqq Komanenz als Sensitivit"atsantagonist\grqq liegt 
hierbei die Betonung nicht auf dem \glqq als\grqq.
Wir schauen in diesem Unterabschnitt die Komanenz weitgehend losgel"ost von der dynamischen Interpretation an.
Wir interessieren uns in diesem Unterabschnitt
nicht etwa daf"ur, nach welchem kinematischen Mechanismus die Komanenz als ein Gegenspieler erscheint, sondern wir interessieren uns haupts"achlich f"ur die Komanenz, weil sie
als ein Gegen"uber der Sensitivit"at auffassbar ist,
das uns Einsichten "uber die Sensitivit"at
vermittelt.  
Nichtsdestotrotz wollen wir dennoch
erl"autern,
was uns dazu bringt,
die Komanenz als den Sensitivit"atsantagonisten
zu titulieren.
Um zu sagen,
was es ist,
was uns zu der groben und unscharfen Gegeneinandersetzung
von Komanenz und Sensitivit"at veranlasst,
gen"ugt es, zu bestimmen, wie die 
unabstrahierten Formen der 
Sensitivit"at und der Komanenz,
die Sensitivit"at und Komanenz stetiger reeller Flussfunktionen  
einander gegen"uberstehen:
Seien $$X_{1},X_{2},\dots X_{n}$$
$n$ Mengen, $n\in\mathbb{N}$ eine nat"urliche Zahl und es sei f"ur alle Tupel
$$(x_{1},x_{2},\dots x_{n})\in X_{1}\times X_{2}\times \dots X_{n}$$
die Aussage $A(x_{1},x_{2},\dots x_{n})$ formulierbar, sodass 
die alternierenden Sequenzen 
\begin{equation}\label{quanto}
\begin{array}{c}
\exists\ x_{1}\in X_{1}\ \forall\  x_{2}\in X_{2}\ \exists\ x_{3}\in X_{3}\ \forall\ \dots
x_{n}\in X_{n}\\
A(x_{1},x_{2},\dots x_{n})\\
\Leftrightarrow:\\
(X_{1},X_{2},\dots X_{n})^{\exists}A(x_{1},x_{2},\dots x_{n})\\
:=\Bigl((X_{1},X_{2},\dots X_{n})^{\forall}A(x_{1},x_{2},\dots x_{n})\Bigr)^{\star}\\
\mbox{und entsprechend}\\
\forall\ x_{1}\in X_{1}\ \exists\  x_{2}\in X_{2}\ \forall\ x_{3}\in X_{3}\exists\ \dots
x_{n}\in X_{n}\\
A(x_{1},x_{2},\dots x_{n})\\
\Leftrightarrow:\\
(X_{1},X_{2},\dots X_{n})^{\forall}A(x_{1},x_{2},\dots x_{n})\\
:=\Bigl((X_{1},X_{2},\dots X_{n})^{\exists}A(x_{1},x_{2},\dots x_{n})\Bigr)^{\star}
\end{array}
\end{equation}
ebenfalls Aussage darstellen. Die Negation der ersten dieser beiden Aussagen ist demnach 
$$\neg\Bigl((X_{1},X_{2},\dots X_{n})^{\exists}A(x_{1},x_{2},\dots x_{n})\Bigr)=
(X_{1},X_{2},\dots X_{n})^{\forall}\neg A(x_{1},x_{2},\dots x_{n})\ ,$$
was impliziert, dass
$$\neg\Bigl((X_{1},X_{2},\dots X_{n})^{\forall}A(x_{1},x_{2},\dots x_{n})\Bigr)=
(X_{1},X_{2},\dots X_{n})^{\exists}\neg A(x_{1},x_{2},\dots x_{n})$$
die Negation der zweiten der beiden Aussagen (\ref{quanto}) ist.
F"ur jede
nat"urliche Zahl $k\in\{2,3,\dots n-2\}$ sind auch
\begin{equation}\label{quanta}
\begin{array}{c}
\exists\ x_{1}\in X_{1}\ \forall\  x_{2}\in X_{2}\dots (x_{k},x_{k+1})\in X_{k}\times X_{k+1}\dots x_{n}\in X_{n}\\ 
A(x_{1},x_{2},\dots x_{n})\\
=(X_{1},X_{2},\dots,X_{k-1},X_{k}\times X_{k+1},X_{k+2}\dots X_{n})^{\exists}A(x_{1},x_{2},\dots x_{n})\\
=:\mathbf{F}_{k}(X_{1},X_{2},\dots X_{n})^{\exists}A(x_{1},x_{2},\dots x_{n})\\
\mbox{bzw.}\\
\forall\ x_{1}\in X_{1}\ \exists\  x_{2}\in X_{2}\dots (x_{k},x_{k+1})\in X_{k}\times X_{k+1}\dots x_{n}\in X_{n}\\ 
A(x_{1},x_{2},\dots x_{n})\\
=(X_{1},X_{2},\dots,X_{k-1},X_{k}\times X_{k+1},X_{k+2}\dots X_{n})^{\forall}A(x_{1},x_{2},\dots x_{n})\\
=:\mathbf{F}_{k}(X_{1},X_{2},\dots X_{n})^{\forall}A(x_{1},x_{2},\dots x_{n})
\end{array}
\end{equation}
zwei Aussagen.
Es k"urze $a$ die Aussage $(X_{1},X_{2},\dots X_{n})^{\exists}A(x_{1},x_{2},\dots x_{n})$ 
ab.
Sofern es keine nat"urliche Zahl $k\in\{2,3,\dots n-2\}$ gibt, oder sofern $k=1$ ist,
seien auch die Aussagen $\mathbf{F}_{1}a$, $\mathbf{F}_{1}a^{\star}$, $\mathbf{F}_{k}a$, $\mathbf{F}_{k}a^{\star}$
in der angedeuteten Weise festgelegt, die der geneigte Leser erkenne. 
Wir nennen jede Aussage ${\rm X}$
der Form (\ref{quanto}) 
genau dann eine 
Einfaltung der Aussage ${\rm Y}$ an der $k$-ten Stelle, wenn\index{Einfaltung einer Aussage} 
$${\rm X}=\mathbf{F}_{k}{\rm Y}$$
ist. 
Jede Einfaltung kann r"uckg"angig gemacht werden und f"ur alle $k\in\mathbb{N}$ bezeichne 
$\mathbf{F}_{k}^{1}=\mathbf{F}_{k}$ und es sei
$\mathbf{F}_{k}^{-1}$
der Operator, der auf der Klasse derjeniger
Aussagen $A$ definiert ist, die als an der Stelle $k$ 
eingefaltet aufgefasst werden k"onnen; f"ur die also 
eine Aussage
$B$
existiert, f"ur die
$A=\mathbf{F}_{k}B$ gilt; exakt 
$B$ sei dabei $\mathbf{F}_{k}^{-1}A$.
Wenn zu einer Aussage 
${\rm Y}$ eine endliche Menge
$\{k(1),k(2),\dots k(m)\}\subset \mathbb{N}$
und ein Tupel $(e(1),e(2),\dots e(m))\in \{-1,1\}^{m}$
existieren, f"ur die 
$${\rm X}=\neg\ \mathbf{F}_{k(1)}^{e(1)}\mathbf{F}_{k(2)}^{e(2)}\dots\mathbf{F}_{k(m)}^{e(m)}{\rm Y}$$
ist, so bezeichnen wir
die Aussage ${\rm X}$ als eine Schiefnegation der 
Aussage ${\rm Y}$.\index{Schiefnegation einer Aussage}  
\newline
Es sei $\Psi$ eine normierte, stetige, reelle Flussfunktion, deren Zustandsraum $\mathbf{P}_{2}\Psi$ im $\mathbb{R}^{n}$
endlicher Dimension $n\in \mathbb{N}$ liege, dessen offene, euklidische Kugeln des Radius $r$ um jeden
jeweiligen Punkt $x\in \mathbb{R}^{n}$ 
f"ur alle $r\in \mathbb{R}^{+}$
mit $\mathbb{B}_{r}(x)$
bezeichnet seien.
Ein sensitiver Zustand $z$ einer normierten, stetigen, reellen Flussfunktion $\Psi$
liegt genau dann
vor, wenn
\begin{equation}\label{quantb}
\begin{array}{c}
\mathbf{F}_{1}(\mathbb{R}^{+},\mathbb{R}^{+},\mathbb{R}^{+})^{\exists} {\rm S}_{\Psi}(\delta,\varepsilon,t_{\star}) \\
=(\mathbb{R}^{+},\mathbb{R}^{+}\times\mathbb{R}^{+})^{\exists} {\rm S}_{\Psi}(\delta,\varepsilon,t_{\star})
\end{array}
\end{equation}
f"ur die Setzung
$$ \exists\ (x,t)\in\mathbb{B}_{v}(z)\times ]-\infty,w[\cup]w,\infty[:\Psi(x,t)\not\in \mathbb{B}_{u}(z) $$
$$:= {\rm S}_{\Psi}(u,v,w)$$
gilt. Die Aussage, dass 
$z$ ein sensitiver Zustand der Flussfunktion $\Psi$ ist, sei durch 
$z{\rm S}_{\Psi}$ bezeichnet, wohingegen
$z{\rm K}_{\Psi}$ die Aussage bezeichne, dass
die Komanenz der normierten, stetigen, reellen Flussfunktion $\Psi$ im Zustand $z$ gegeben ist.
$z{\rm K}_{\Psi}$ ist genau dann
gegeben, wenn die Aussage
\begin{equation}\label{quantc} 
\begin{array}{c}
(\mathbb{R}^{+}\times\mathbb{R}^{+},\mathbb{R}^{+})^{\forall}\neg\ {\rm S}_{\overline{\Psi}}(\delta,\varepsilon,t_{\star})\\
=\neg \Bigl((\mathbb{R}^{+}\times
\mathbb{R}^{+},\mathbb{R}^{+})^{\exists}{\rm S}_{\overline{\Psi}}(\delta,\varepsilon,t_{\star})\Bigr)\\
\neg\ \mathbf{F}_{2}\mathbf{F}_{1}^{-1}z{\rm S}_{\overline{\Psi}}
\end{array}
\end{equation}
f"ur die zu der Flussfunktion $\Psi$ zeitreziproke Flussfunktion $\overline{\Psi}$
gilt.\index{zeitreziproke Flussfunktion} Denn die Aussage
$$\forall\ \delta,t_{\star}\in\mathbb{R}^{+}\ \exists\varepsilon\in\mathbb{R}^{+} :$$ 
$$(x,t)\in\mathbb{B}_{\varepsilon}(z)\times ]-1/t_{\star},1/t_{\star}[\Rightarrow\Psi(x,t)\in
\mathbb{B}_{\delta}(z)$$
erkennen wir als die Komanenz der Flussfunktion $\Psi$ im Zustand $z$
und
es gilt 
f"ur alle Tripel
$(\delta,\varepsilon,t_{\star})\in \mathbb{R}^{+}\times\mathbb{R}^{+}\times\mathbb{R}^{+}$
die "Aquivalenz
$$\Bigl((x,t)\in\mathbb{B}_{\varepsilon}(z)\times ]-1/t_{\star},1/t_{\star}[\Rightarrow\Psi(x,t)\in 
\mathbb{B}_{\delta}(z)\Bigr)$$
$$\Leftrightarrow\Bigl((x,1/t)\in\mathbb{B}_{\varepsilon}(z)\times ]-\infty,t_{\star}[\cup]t_{\star},\infty[\Rightarrow\Psi(x,1/t)\in 
\mathbb{B}_{\delta}(z)\Bigr)\ ,$$
wobei die rechte Seite dieser "Aquivalenz zu der Negation
$$\Bigl((x,t)\in\mathbb{B}_{\varepsilon}(z)\times ]-\infty,t_{\star}[\cup]t_{\star},\infty[\Rightarrow
\overline{\Psi}(x,t)\in 
\mathbb{B}_{\delta}(z)\Bigr)$$
$$=\neg\ {\rm S}_{\overline{\Psi}}(\delta,\varepsilon,t_{\star})$$
"aquivalent ist. Wir erkennen leicht, dass auch
f"ur alle normierten und metrisierten Flussfunktionen $\Psi=(\psi,d)$ und alle deren 
Zust"ande $z$
die "Aquivalenzen
\begin{equation}\label{quantd}
\begin{array}{c}
z{\rm K}_{\Psi}\Leftrightarrow \neg\ \mathbf{F}_{2}\mathbf{F}_{1}^{-1}z{\rm S}_{\overline{\Psi}}\ ,\\
z{\rm K}_{\overline{\Psi}}\Leftrightarrow \neg\ \mathbf{F}_{2}\mathbf{F}_{1}^{-1}z{\rm S}_{\Psi}\ ,\\
z{\rm S}_{\Psi}\Leftrightarrow  \neg\ \mathbf{F}_{1}\mathbf{F}_{2}^{-1}z{\rm K}_{\overline{\Psi}}\ ,\\
z{\rm S}_{\overline{\Psi}}\Leftrightarrow \neg\ \mathbf{F}_{1}\mathbf{F}_{2}^{-1}z{\rm K}_{\Psi}
\end{array}
\end{equation}
gelten.
Die Sensitivit"at eines Zustandes $z$ einer
normierten und metrisierten Flussfunktion $\Psi=(\psi,d)$
ist die 
Negation der Einfaltung $\mathbf{F}_{1}\mathbf{F}_{2}^{-1}$ der
Komanenz der zu der Flussfunktion $\Psi$
zeitreziproken metrisierten Flussfunktion $\overline{\Psi}=(\overline{\psi},d)$ im Zustand $z$.
Jene Sensitivit"at ist also eine Schiefnegation dieser Komanenz. 
Das ist also die  Weise, wie
die Komanenz als ein \glqq Sensitivit"atsantagonist\grqq in der Logik der 
Begriffskonstitution erscheint.
\subsection{Allgemeine und interne Komanenz}
Zur Vermittlung der
allgemeinen Komanenz
sei
f"ur jeweils jede topologisierte Flussfunktion $\Xi=(\xi, \mathbf{T}(\mathbf{P}_{2}\xi))$ im weiteren Sinn
f"ur alle reellen Zahlen $t\in\mathbb{R}$
und f"ur alle Zust"ande $z\in \mathbf{P}_{2}\xi$
das Mengensystem
\begin{equation}\label{ydnuld}
\begin{array}{c}
\Bigl\{\bigcup\mathbf{P}_{2}\Omega:\Omega\in\mathbf{T}(\mathbf{P}_{2}\xi)_{\xi(z,[0,t]\cup[t,0])}^{[0,t]\cup[t,0]} \Bigr\}\\
\quad\\
=:\mathbf{dl}_{\mathbf{T}(\mathbf{P}_{2}\Xi)}(\xi,z,t)\ \subset\ \mathbf{T}(\mathbf{P}_{2}\xi)
\end{array}
\end{equation}
von Vereinigungen "uber die Wertemengen $\mathbf{P}_{2}\Omega\subset \mathbf{T}(\mathbf{P}_{2}\xi)$ jeweiliger umgebungswertiger Funktionen 
$$\Omega\in\mathbf{T}(\mathbf{P}_{2}\xi)_{\xi(z,[0,t]\cup[t,0])}^{[0,t]\cup[t,0]}$$
herausgestellt.
Da die Vereinigungen 
$$\bigcup\mathbf{P}_{2}\Omega\in\mathbf{T}(\mathbf{P}_{2}\xi)$$ "uber die Wertemengen $\mathbf{P}_{2}\Omega$ als Vereinigungen 
"uber Teilmengen der Topologie $\mathbf{T}(\mathbf{P}_{2}\xi)$ Elemente derselben sind,
ist das Mengensystem $\mathbf{dl}_{\mathbf{T}(\mathbf{P}_{2}\Xi)}(\xi,z,t)$ 
ein Mengensystem, das in der Topologie $\mathbf{T}(\mathbf{P}_{2}\xi)$ enthalten ist.
Wir nennen exakt dieses Konstrukt $\mathbf{dl}_{\mathbf{T}(\mathbf{P}_{2}\Xi)}(\xi,z,t)$
das Mengensystem der 
$t$-gedehnten Umgebungen 
des Zustandes $z\in \mathbf{P}_{2}\xi$
der topologisierten Flussfunktion $\Xi$.\index{gedehnte Umgebungen 
eines Zustandes} Die Konstruktion dieses 
Mengensystemes ist dabei offenbar f"ur beliebige
Zustandsraum"uberdeckungen $\mathcal{A}$ formulierbar. 
Auf der Klasse der Tripel $(\xi,\mathcal{A},z,t)$,
deren erste Komponente $\xi$ eine Flussfunktion ist, deren
zweite Komponente eine Zustandsraum"uberdeckung $\mathcal{A}$ ist und
deren dritte und vierte Komponente ein Element $(z,t)$ der
Definitionsmenge $\mathbf{P}_{1}\xi$ bilden,
gibt es den universellen Operator
$$\mathbf{dl}_{\mathbf{P}_{2}}(\mathbf{P}_{1},\mathbf{P}_{3},\mathbf{P}_{4})\ ,$$ 
dessen jeweilige Werte
\begin{equation}\label{dynula}
\begin{array}{c}
\mathbf{dl}_{\mathcal{A}}(\xi,z,t)\\
\Bigl\{\bigcup\mathbf{P}_{2}\Omega:\Omega\in\mathcal{A}_{\xi(z,[0,t]\cup[t,0])}^{[0,t]\cup[t,0]} \Bigr\}\\
\end{array}
\end{equation}
sind. Wir wissen hierbei aber jeweils nicht, ob das Mengensystem $\mathcal{A}$ 
die Vereinigungsabgeschlossenheit einer Topologie hat, 
welche
gegen"uber
der Vereinigung "uber beliebige ihrer Teilmengen abgeschlossen ist.
Daher wissen wir auch nicht, ob $\mathbf{dl}_{\mathcal{A}}(\xi,z,t)$ in $\mathcal{A}$ liegt.
Gleich, ob die topologisierte Flussfunktion 
$\Xi$ im weiteren Sinn stetig ist, oder nicht, ist zwar das Mengensystem der 
$t$-gedehnten Umgebungen 
des Zustandes $z\in \mathbf{P}_{2}\xi$ paraphrasierbar als 
\begin{equation}\label{dynulb}
\mathbf{dl}_{\mathbf{T}(\mathbf{P}_{2}\Xi)}(\xi,z,t)=\Bigl\{{\rm U}\in\mathbf{T}(\mathbf{P}_{2}\xi):
\xi(z,[0,t]\cup[t,0])\subset {\rm U}\Bigr\}\ .
\end{equation}
Dem gegen"uber gibt
es offensichtlich Zustandsraum"uberdeckungen $\mathcal{A}$, f"ur die das Mengensystem
$$\Bigl\{{\rm U}\in\mathcal{A}:
\xi(z,[0,t]\cup[t,0])\subset {\rm U}\Bigr\}$$
leer ist. In dem Sachverhalt, dass das Mengensystem der 
$t$-gedehnten Umgebungen $\mathbf{dl}_{\mathbf{T}(\mathbf{P}_{2}\Xi)}(\xi,z,t)$ die
Gleichung (\ref{dynulb}) erf"ullt, "aussert sich
die Kontinuit"at derjeniger Mengensysteme, die 
Topologien sind:
Wir haben ja bereits im ersten Teil der abstrakten Ergodentheorie 
gezeigt, dass die Eigenschaften von Mengensystemen,
die Topologien $\mathbf{T}$ und deren Dualisierungen $\mathbf{T}^{c}$ als solche auszeichnen,
Eigenschaften sind, die
als eine Form der
Cantor-Stetigkeit\index{Cantor-Stetigkeit} angesehen werden
k"onnen. 
Diese Eigenschaften der Topologien und deren Dualisierungen
k"onnen wir
n"amlich als die Cantor-Stetigkeit der Vereinigung und des Schnittes der Elemente
der jeweiligen Topologie $\mathbf{T}$ oder deren Dualisierung $\mathbf{T}^{c}$ formulieren, die wir
dabei
als Verk"upfungen auf der jeweiligen Topologie oder deren Dualisierung 
auffassen.
\newline
Nun k"onnen wir leicht festlegen, dass jede strukturierte Flussfunktion $(\xi,\mathcal{A})$ 
genau dann im Zustand $z\in \mathbf{P}_{2}\xi$ als intern komanent gelten soll, wenn f"ur $z$ die Aussage 
\begin{equation}\label{ydnule}
\begin{array}{c}
\forall\ t\in\mathbb{R}^{+},{\rm L}^{+}\in \mathbf{dl}_{\mathcal{A}}(\xi,z,t)\ \exists\ {\rm L}\in \mathcal{A}_{\{z\}}:\\
\bigcup\{\xi(y,[0,t]\cup[t,0]):y\in {\rm L}\}\subset {\rm L}^{+}
\end{array}
\end{equation}
wahr ist. Die interne Komanenz jeder strukturierten Flussfunktion $(\xi,\mathcal{A})$ bestehe schliesslich darin,
dass sie in jedem ihrer Zust"ande\index{interne Komanenz einer strukturierten Flussfunktion}
intern komanent ist.\index{interne Komanenz in einem Zustand}\newline 
Mit
dem Begriff der internen Komanenz
einer 
strukturierten
Flussfunktion 
kommen wir
"uber die Allgemeinheitsstufe der allgemeinen Topologie hinaus, aber nicht viel weiter. Wir
geraten in Verlegenheit, wenn wir 
die Komanenz auf dieser Allgemeinheitsstufe  mit
Hilfe einer 
Komanenzfunktion 
quantifizieren wollen. Diese 
Schwierigkeit,
eine allgemeine Komanenzfunktion zu formulieren,
sollte uns hier eigentlich
ebensowenig wie die 
Er"orterung der internen
Komanenz 
auf dem sehr allgemeinen Generalit"atsniveau beliebiger 
Paare $(\xi,\mathcal{A})$ einer Flussfunktion und einer beliebigen
Zustandsraum"uberdeckung $\mathcal{A}$
in weiterf"uhrende Untersuchungen verwickeln.
\newline
Das Anliegen dieses Abschnittes ist es ja nicht, die allgemeine interne Komanenz
zu untersuchen. Wir
wollten z.B.
weder darlegen, ob, und gegebenfalls, wie die 
Stetigkeit einer topologisierten Flussfunktion $\Xi$ die interne
Komanenz gem"ass (\ref{ydnule}) impliziert. Noch wollten wir untersuchen,
inwieweit sich eine allgemeine Komanenz
in Form eines
funktionellen "Aquivalentes der Komanenzfunktion quantifizieren l"asst;
wobei uns als funktionelles "Aquivalent der Komanenzfunktion
gelte,
was uns, wie diese, eine Argumentationsgundlage anlegt.
Das bescheidene Anliegen dieses Abschnittes ist es, lediglich vorzuf"uhren,
dass es die Relation der Ultrakolokalisation gibt,
die uns eine Generalisierung der Sensitivit"at "uber 
die metrische Sensitivit"at hinaus darstellt,
um anschliessend m"oglichst allgemeine Bedingungen
zu zeigen, unter denen die Ultrakolokalisation
eine "Aquivalenzrelation ist.
\subsection{Konverse Komanenz}
Wir kommen aber nicht darum herum, die allgemeine Komanenz
in einer Form zu verfassen, die bessere Argumentationsgundlagen
schafft, als dies die interne Komanenz zu leisten vermag:
Um uns die Formulierung jenes allgemeinen Komanenzbegriffes 
zu erleichtern,
der uns eine bessere Argumentationsbasis als die interne Komanenz gibt,
legen wir f"ur jedes Mengensystem $\mathcal{A}$
und f"ur die Menge $\mathcal{A}^{\mathbb{N}}$ aller Folgen $\{A_{j}\}_{j\in \mathbb{N}}$ von Mengen $A_{j}$ aus $\mathcal{A}$
fest, dass f"ur alle $A_{\star}\in \mathcal{A}$ und 
f"ur jede 
Folge $\{A_{j}\}_{j\in \mathbb{N}}\in \mathcal{A}^{\mathbb{N}}$ 
\begin{equation}\label{ydnulw}
\begin{array}{c}
\forall\ A^{+}\in\{A\in\mathcal{A}:A\supset A_{\star}\}\setminus\{A_{\star}\} \exists\ j(A^{+})\in\mathbb{N}:\ \forall\ j>j(A^{+})\\
A_{\star}\subset A_{j+1}\subset A_{j}\subset A^{+}\\
\Leftrightarrow:\\
A_{j}-\mathcal{A}\to A_{\star}
\end{array}
\end{equation}
sei. Und davon abgeleitet sei f"ur jede Folge $\{x_{j}\}_{j\in \mathbb{N}}$ von Elementen aus der Vereinigung $\bigcup\mathcal{A}\ni y$
\begin{equation}\label{ydnuww}
\begin{array}{c}
x_{j}-\mathcal{A}\to y\\
:\Leftrightarrow\\
\exists\ \{A_{j}\}_{j\in \mathbb{N}}\in \mathcal{A}^{\mathbb{N}}:\ A_{j}-\mathcal{A}\to \{y\}\ \land\\
(j\in\mathbb{N}\Rightarrow x_{j}\in A_{j})
\end{array}
\end{equation}
eine kurze Schreibweise, die immer in dem Kontext zu stehen hat, dass $x_{j}$ ein Zustand einer mit $j$ indizierten 
Zustandsfolge ist.
Nun heben wir diejenigen
strukturierten Flussfunktionen $(\xi, \mathcal{A})$ heraus, f"ur welche f"ur alle
Folgen $\{x_{j}\}_{j\in \mathbb{N}}$ von Zust"anden $x_{j}\in \mathbf{P}_{2}\xi$ und f"ur alle $y\in\mathbf{P}_{2}\xi$ 
f"ur jede Lokalisierung $A\in \mathcal{A}$
die Implikation
\begin{equation}\label{ydnulv}
x_{j}-\mathcal{A}\to y\Rightarrow \lim_{j\to\infty} {\rm D}_{2}\Bigl(\xi(x_{j},\mbox{id})^{-1}(A),\xi(y,\mbox{id})^{-1}(A)\Bigr)=0
\end{equation}
gilt. Exakt diese besonderen strukturierten Flussfunktionen $(\xi, \mathcal{A})$, f"ur die f"ur alle
Folgen $\{x_{j}\}_{j\in \mathbb{N}}$ von Zust"anden $x_{j}\in \mathbf{P}_{2}\xi$ und f"ur alle $y\in\mathbf{P}_{2}\xi$ 
und f"ur alle Lokalisierungen $A\in \mathcal{A}$ die Implikation (\ref{ydnulv}) wahr ist,
bezeichnen wir als konvers komanente Flussfunktionen $(\xi, \mathcal{A})$.\index{konvers komanente Flussfunktion}
Wie wir zu dieser Bezeichnung kommen, wird sich aus den folgenden Ausf"uhrungen dieses 
Unterabschnittes ergeben.
\newline
Da wir hierbei keinerlei Bestimmungen "uber Distanzfunktionen festlegten, ausser, was deren
Definitionsmenge sein soll und, dass deren Werte Elemente aus $\overline{\mathbb{R}}$ sein sollen, ist
evident, dass die konverse Komanenz im allgemeinen die interne Komanenz nicht impliziert.
L"auft dieser Sachverhalt logischer Beziehungslosigkeit aber nicht gerade gegen unsere Absicht, den
Begriff der
konversen Komanenz einzuf"uhren, um 
mit der konversen Komanenz
eine bessere Argumentationsgundlage
zu generieren?\newline
Wenn wir 
die logische Beziehungslosigkeit zwischen dem Begriff der
konversen Komanenz und dem Begriff der internen Komanenz 
auf den ersten Blick f"ur ein Zeichen daf"ur halten, 
dass die Verfassung des Begriffes der
konversen Komanenz
im Hinblick auf unsere Motivation
kontraindiziert ist, 
dann verkennen wir, dass 
der Vorzug der konversen Komanenz 
gegen"uber der internen Komanenz
in ihrer Regulierbarkeit durch die Wahl einer 
jeweiligen mengenweisen Zahlenstrahldistanz ${\rm D}_{2}$ besteht.
F"ur alle $z\in\mathbf{P}_{2}\xi$ 
und f"ur alle $X\subset\mathbf{P}_{2}\xi$ und f"ur alle $T\subset\mathbb{R}$
ist die in (\ref{ydnuvv}) spezifizierte Menge
$$\xi(z,\mbox{id})^{-1}(X)\subset\mathbb{R}$$
deutbar als
der Aufenthaltszeitraum der durch $\xi$ bestimmten Entwicklung
des Zustandes $z$ in der Lokalisation $X$.
${\rm D}_{2}$ sei dabei eine Distanzfunktion, die auf dem kartesischen Produkt 
\begin{equation}\label{ydnuvv}
\mathbf{P}_{1}{\rm D}_{2}=2^{\overline{\mathbb{R}}}\times 2^{\overline{\mathbb{R}}}=2^{\mathbf{P}_{2}\mathbf{P}_{1}\xi}\times 
2^{\mathbf{P}_{2}\mathbf{P}_{1}\xi}
\end{equation}
des zweiten kartesischen Faktors $\mathbf{P}_{2}\mathbf{P}_{1}\xi=\overline{\mathbb{R}}$ der 
Definitionsmenge $\mathbf{P}_{1}\xi$ der Flussfunktion $\xi$
definiert ist, worauf der Index der Bezeichnung dieser Distanzfunktion ${\rm D}_{2}$ hinweisen soll. 
Die jeweiligen Werte dieser Distanzfunktion seien in dem Intervall $[0,\infty]$.
Die Distanzfunktion ${\rm D}_{2}$ 
k"onnte also beispielsweise
die sich von der nat"urlichen Metrik des Zahlenstrahles $d_{\overline{\mathbb{R}}}:=|\mathbf{P}_{1}-\mathbf{P}_{2}|$ ableitende 
Hausdorffsche Distanzfunktion $d_{\overline{\mathbb{R}}}^{[1]}$ sein,
die auf dem kartesischen Produkt (\ref{ydnuvv})
definiert ist, die aber bekanntlich auf demselben keine Metrik ist. 
$d_{\overline{\mathbb{R}}}^{[1]}$ ist hingegen nicht 
nur eine 
Metrik auf dem Mengensystem aller in dem Zahlenstrahl liegenden Kompakta, sondern 
$d_{\overline{\mathbb{R}}}^{[1]}$ ist auf 
dem Mengensystem aller
bez"uglich der Metrik des Zahlenstrahles $d_{\overline{\mathbb{R}}}$ abgeschlossenen 
Mengen eine Metrik. 
Sie ist aber keine Metrik
auf der nat"urlichen Topologie $\mathbf{T}(1)$ des Zahlenstrahles und erst recht keine solche
auf dem 
Mengensystem
$$\kappa\mathbf{T}(1):=\{Y\subset\mathbb{R}:\exists\ {\rm U}\in \mathbf{T}(1):{\rm U}\subset Y\}\ ,$$
exakt dessen Elemente wir als keuliale Teilmengen\index{keuliale Menge} des Zahlenstrahles bezeichnen, denn $\kappa o\iota\lambda\iota \acute{\alpha}$
heisst auf altgriechisch Bauch.
Demgegen"uber nennen wir exakt das Mengensystem 
$\kappa^{\star}\mathbf{T}(1)$
all derjeniger Teilmengen des Zahlenstrahles, deren offener Kern eine abz"ahlbare Vereinigung offener Intervalle ist,
das Mengensystem
aller total keulialen Teilmengen des Zahlenstrahles. Dieses ist also
\begin{equation}
\kappa^{\star}\mathbf{T}(1)=\Bigl\{Y\subset\mathbb{R}:{\rm U}\in \mathbf{T}(1)\Rightarrow Y\cap{\rm U}\in\kappa\mathbf{T}(1)\cup\{\emptyset\}\Bigr\}\ .   
\end{equation}
Das Mengensystem
\begin{equation}
\kappa_{\star}\mathbf{T}(1):=\{Y\subset\mathbb{R}:\forall\ t_{\star}\in
\mathbb{R}^{+}\ \exists\ {\rm U}\in \mathbf{T}(1):{\rm U}\subset Y\setminus]-t_{\star},t_{\star}[\}\ ,
\end{equation}
der unbefristet keulialen Teilmengen\index{unbefristet keuliale Menge} des Zahlenstrahles
k"onnen wir wegen dessen Wohlgeordnetheit f"ur die Topologie $\mathbf{T}(1)$ formulieren.\footnote{Es gilt
$$\mathbf{T}(1)\subset\kappa^{\star}\mathbf{T}(1)\subset\kappa_{\star}\mathbf{T}(1)\subset\kappa\mathbf{T}(1)\ .$$}
Eine Abstraktion der unbefristet keulialen Teilmengen ist daher wohl schwerer zu finden
als eine Abstraktion der keulialen und der total keulialen Teilmengen des Zahlenstrahles:
F"ur jedes Mengensystem $\mathcal{A}$ sei
das Mengensystem
\begin{equation}
\begin{array}{c}
\kappa\mathcal{A}:=\Bigl\{Y\subset\bigcup\mathcal{A}:\exists\ {\rm A}\in \mathcal{A}:{\rm A}\subset Y\Bigr\}\ ,\\
\mbox{bzw.}\\
\kappa^{\star}\mathcal{A}:=\Bigl\{Y\subset\mathcal{A}:{\rm L}\in\mathcal{A} \Rightarrow Y\cap{\rm L}\in\kappa\mathcal{A}\cup\{\emptyset\}\Bigr\}
\end{array}
\end{equation}
dasjenige, das
als das Mengensystem aller bez"uglich $\mathcal{A}$ keulialer bzw. total keulialer Teilmengen\index{total keuliale Menge}\index{totale Keulialit\"at}
gelte; womit wir die 
universellen Operatoren $\kappa$ 
und $\kappa^{\star}$
auf der Klasse der Mengensysteme festlegen. Dabei 
leitet sich die Benennung keulialer Teilmengen vom griechischen Wort \glqq$\kappa o\iota\lambda\iota\acute{\alpha}$\grqq
ab, das den Bauch benennt.
Bezeichnet $\mu_{1}$ das Lebesguesche Mass des Zahlenstrahles, so ist
das Mengensystem keulialer Teilmengen $\kappa\mathbf{T}(1)$ nicht mit dem 
Mengensystem $\mu^{-1}_{1}(\mathbb{R}^{+})$ identisch.
Die keulialen Teilmengen aus $\kappa\mathbf{T}(1)$ fallen weder mit den offenen, noch mit den
abgeschlossenen Teilmengen, noch mit
denjenigen Teilmengen des Zahlenstrahles in eins, die positives Lebesguesches Mass haben.
Nichtsdestotrotz ist gerade die f"ur keuliale Teilmengen
gegebene Eigenschaft ein immer wieder auftretendes Reichhaltigkeitskriterium.
\newline
Wieso reden wir von der 
\glqq konversen Komanenz\grqq
exakt jener besonderen strukturierten Flussfunktionen $(\xi, \mathcal{A})$, f"ur die f"ur alle
Folgen $\{x_{j}\}_{j\in \mathbb{N}}$ von Zust"anden $x_{j}\in \mathbf{P}_{2}\xi$ und f"ur alle $y\in\mathbf{P}_{2}\xi$ 
und f"ur alle Lokalisierungen $A\in \mathcal{A}$ die Implikation (\ref{ydnulv}) wahr ist?
\index{konvers komanente Flussfunktion}
Diese Benennung liegt nahe, weil die konverse Komanenz auf die folgende 
Weise das konverse Pendant folgender verallgemeinerten Form der Komanenz ist:
Es sei ${\rm D}_{1}$ eine der Distanzfunktion ${\rm D}_{2}$
entsprechende Distanzfunktion, die
auf dem kartesischen Produkt 
\begin{equation}\label{ydnvvv}
\mathbf{P}_{1}{\rm D}_{2}=2^{\mathbf{P}_{2}\xi}\times 2^{\mathbf{P}_{2}\xi}=2^{\mathbf{P}_{1}\mathbf{P}_{1}\xi}\times 2^{\mathbf{P}_{1}\mathbf{P}_{1}\xi}
\end{equation}
des nun ersten kartesischen Faktors $\mathbf{P}_{1}\mathbf{P}_{1}\xi=\mathbf{P}_{2}\xi$ der 
Definitionsmenge $\mathbf{P}_{1}\xi$ der Flussfunktion $\xi$
definiert ist. 
Wir nennen jede strukturierte Flussfunktion $(\xi, \mathcal{A})$ genau dann bez"uglich 
der Distanzfunktion 
${\rm D}_{1}$ 
komanent, wenn f"ur alle Teilmengen 
$T\subset \mathbb{R}$
die Implikation
\begin{equation}\label{zvnulv}
x_{j}-\mathcal{A}\to y\Rightarrow \lim_{j\to\infty} {\rm D}_{1}\Bigl(\xi(x_{j},T),\xi(y,T)\Bigr)=0
\end{equation}
wahr ist und deren Analogie zu der f"ur die konverse Komanenz konstitutive 
Implikation (\ref{ydnulv}) transparent ist.
Denn die Eigenschaft
einer jeweiligen strukturierten Flussfunktion $(\xi, \mathcal{A})$
f"ur ${\rm D}_{1}$ f"ur jede Teilmenge 
$T\subset \mathbb{R}$
die Implikation (\ref{zvnulv}) zu erf"ullen,
ist eine verallgemeinerte Form der Komanenz, die
den Komanenzbegriff umfasst, den wir
im Traktat "uber den elementaren Quasiergodensatz \cite{raab} f"ur
stetige reelle Flussfunktionen $\psi$ verfassten, deren 
Zustandsraum $\mathbf{P}_{2}\psi$ im $\mathbb{R}^{n}$
endlicher Dimension $n\in\mathbb{N}$ liegt:
Wir zeigten in jenem Traktat, dass die stetigen reellen Flussfunktionen $\psi$
mit kompaktem Zustandsraum $\mathbf{P}_{2}\psi$
in dem Sinn komanent sind, dass ihre Komanenzfunktion
$B_{\psi}$ auf der Definitionsmenge $\mathbb{R}^{+}\times \mathbb{R}^{+}$ Werte 
in $\mathbb{R}^{+}$
hat, f"ur welche 
f"ur
alle Quadrupel
$(x,y,\delta,t)\in \mathbf{P}_{2}\psi\times\mathbf{P}_{2}\psi\times\mathbb{R}^{+}\times \mathbb{R}^{+}$
die Implikation
\begin{equation}\label{zvzulv}
\begin{array}{c}
||x-y||<B_{\psi}(\delta,t)\Rightarrow\\
(\vartheta\in[0,t]\Rightarrow ||\psi(x,\vartheta)-\psi(y,\vartheta)||<\delta)
\end{array}
\end{equation}
wahr ist, wobei $||\mbox{id}||$ die euklidische Norm des $\mathbb{R}^{n}$ sei.
Das Bestehen dieses Sachverhaltes ist offensichtlich zu der Gegebenheit "aquivalent, dass
beispielsweise 
f"ur die Distanzfunktion $\tilde{{\rm D}}_{1}$ mit den jeweiligen Werten 
$$\tilde{{\rm D}}_{1}(X,Y):=\int_{X}\inf\{||x-y||: y\in Y\}d\mu(x)$$
f"ur jeweilige Zustandsraumteilmengen $X,Y\subset \mathbf{P}_{2}\psi$
die Implikation 
\begin{equation}\label{znlv}
x_{j}-\mathbf{T}(n)\to y\Rightarrow \lim_{j\to\infty} 
\tilde{{\rm D}}_{1}\Bigl(\psi(x_{j},T),\psi(y,T)\Bigr)=0
\end{equation}
f"ur alle beschr"ankten Teilmengen $T\subset\mathbb{R}$
und f"ur alle Folgen
$\{x_{j}\}_{j\in \mathbb{N}}$ von Zust"anden $x_{j}\in \mathbf{P}_{2}\psi$ und f"ur alle $y\in\mathbf{P}_{2}\psi$
wahr ist. Die ausdr"uckliche Formulierung einer 
im Bezug auf $\psi$ gewichtenden
Distanzfunktion ${\rm D}_{1}^{\psi}$, f"ur die
die zu (\ref{znlv}) analoge Implikation
auch f"ur alle 
nicht beschr"ankten Teilmengen $T\subset\mathbb{R}$  
und f"ur jede Zustandsfolge
$\{x_{j}\}_{j\in \mathbb{N}}$ und jeden Zustand $y$
gilt, ist m"oglich, jedoch umst"andlich. Wir unterlassen es, sie hinzuschreiben.
"Uber die Hausdorffsche Distanzfunktion $d^{[1]}_{\mathbb{R}^{n}}$, die
sich von der euklidischen Metrik $d_{\mathbb{R}^{n}}$ des $\mathbb{R}^{n}$ ableitet,
bemerken wir allerdings dies:
Die bez"uglich dieser Distanzfunktion $d^{[1]}_{\mathbb{R}^{n}}$
komanenten topologisierten Flussfunktionen $(\psi,\mathbf{T}(n))$,
f"ur welche also
die Implikation
\begin{equation}\label{znlv}
x_{j}-\mathbf{T}(n)\to y\Rightarrow \lim_{j\to\infty} 
d^{[1]}_{\mathbb{R}^{n}}\Bigl(\psi(x_{j},T),\psi(y,T)\Bigr)=0
\end{equation}
f"ur alle Teilmengen $T\subset\mathbb{R}$
und f"ur alle Folgen
$\{x_{j}\}_{j\in \mathbb{N}}$ von Zust"anden $x_{j}\in \mathbf{P}_{2}\psi$ und f"ur alle $y\in\mathbf{P}_{2}\psi$
gilt, sind gerade diejenigen 
Flussfunktionen $\psi$, welche Laplace-stetig im Sinne der Definition des 
Anhanges des 
Traktats "uber den elementaren Quasiergodensatz \cite{raab}
sind: 
Die Laplace-stetigen Flussfunktionen 
mit kompaktem Zustandsraum 
erwiesen sich als diejenigen Flussfunktionen,
all deren Trajektorien geschlossen sind und
die insofern nur triviale Zimmer haben. Es ist vor dem Hintergrund des elementaren Quasiergodensatzes 
evident, dass die Laplace-stetigen Flussfunktionen mit kompaktem Zustandsraum exakt diejenigen reellen
Flussfunktionen $\psi$ mit kompaktem Zustandsraum 
$\mathbf{P}_{2}\psi\subset\mathbb{R}^{n}$
sind, 
f"ur welche 
die zugeh"orige topologisierte Flussfunktion $(\psi,\mathbf{T}(n))$
bez"uglich der Hausdorffschen Distanzfunktion $d_{1}^{[1]}$
konvers komanent ist, was die dargelegte Korrespondenz 
der 
Komanenz bez"uglich einer Zustandsraummengen bewertenden Distanz und der auf eine 
entsprechende
Distanzfunktion des Zahlenstrahles bezogenen
konversen 
Komanenz 
illustriert.
\newline
Wir sehen hierbei, dass die Wahl der 
jeweiligen 
Distanzfunktionen die Strenge sowohl der 
Komanenz bez"uglich denselben als auch der
konversen 
Komanenz bez"uglich denselben bestimmt.
Nehmen wir
statt der Hausdorffschen Distanzfunktion $d_{1}^{[1]}$
die Distanzfunktion $g$
mit den Werten
\begin{equation}\label{auznlv}
g(X,Y):=\frac{1}{\sqrt{2\pi}}\int_{\mathbb{R}}|\mathbf{1}_{X}(x)-\mathbf{1}_{Y}(x)| e^{-x^{2}/2}d\mu(x)
\end{equation}
f"ur alle $X,Y\subset \overline{\mathbb{R}}$, wobei
$\mathbf{1}_{Z}$ f"ur alle $Z\subset \overline{\mathbb{R}}$
die Indikatorfunktion der Menge $Z$ bezeichne:
Wenn wir 
der konversen 
Komanenz diese Distanzfunktion $g$
zugrundelegen, so sind es nicht nur 
Laplace-stetige Flussfunktionen, welche die
reellen Flussfunktionen
mit kompaktem Zustandsraum sind,
f"ur welche die topologisierte Flussfunktion $(\psi,\mathbf{T}(n))$
bez"uglich der Distanzfunktion $g$
konvers komanent ist: Es sind
exakt die stetigen reellen Flussfunktionen
mit kompaktem Zustandsraum, die 
bez"uglich der Distanzfunktion $g$
konvers komanent sind.\newline 
Die Gegebenheit konverser Komanenz 
f"ur eine jeweilige Flussfunktion impliziert
je nach der Wahl der mengenweisen Zahlenstrahldistanz
das Vorliegen verschiedener 
Stetigkeitsklassen der jeweiligen
Flussfunktion. Dieser Sachverhalt
illustriert also auch, dass die verschiedenen Stetigkeitsklassen
verschiedenen Graden der Strenge der zugrundegelegten
mengenweisen Zahlenstrahldistanzen 
entsprechen:
Wir beobachten, dass der hohe 
Stetigkeitsgrad der
Laplace-stetigkeit zu der hohen Strenge der Hausdorffschen Distanzfunktion $d^{[1]}_{\mathbb{R}^{n}}$
geh"ort, w"ahrend
der niedrigere Stetigkeitsgrad der einfachen
Stetigkeit der 
niedrigeren Strenge der Distanzfunktion $g$
entspricht. Offenbar gibt es eine Menge solcher
mengenweiser Zahlenstrahldistanzen, die so,
wie die Distanzfunktion $g$,
in einem noch zu fassenden Sinn gem"assigt sind.
Wir legen Charakterisierungen mengenweiser Distanzfunktionen an:
\newline
Es sei $\Theta$ eine Menge, $(\Theta,\mathbf{T}(\Theta))$ ein topologischer Raum, 
${\rm D}$ eine Abbildung, deren Definitionsmenge $\mathbf{P}_{1}{\rm D}=2^{\Theta}\times 2^{\Theta}$ sei und deren
Wertemenge $\mathbf{P}_{2}{\rm D}\subset[0,\infty]$ nicht negativ sei. 
${\rm D}$ ist demnach eine mengenweise Distanzfunktion auf der Menge $\Theta$; wir k"onnen
auch sagen, dass ${\rm D}$ eine mengenweise Distanzfunktion 
auf dem topologischen Raum $(\Theta,\mathbf{T}(\Theta))$ sei.
Exakt, falls
f"ur alle $X,Y\subset \Theta$
\begin{equation}
{\rm D}(\mathbf{cl}_{\mathbf{T}(\Theta)}(X),Y)={\rm D}(X,\mathbf{cl}_{\mathbf{T}(\Theta)}(Y))={\rm D}(X,Y)
\end{equation}
ist, bezeichnen wir ${\rm D}$ als eine stetige mengenweise Distanzfunktion\index{stetige mengenweise Distanzfunktion} 
auf dem topologischen Raum $(\Theta,\mathbf{T}(\Theta))$. 
Unabh"angig davon, ob ${\rm D}$ stetig ist oder nicht, nennen wir ${\rm D}$
genau dann monoton, wenn
f"ur alle $Z\subset \Theta$ die Beschr"ankung 
\begin{equation}
{\rm D}(X\cap Z,Y\cap Z)\leq {\rm D}(X,Y)
\end{equation}
gilt.\index{monotone mengenweise Distanzfunktion} Falls das Paar 
$(\Theta,d)$ ein metrischer Raum ist,
bezeichnen wir ${\rm D}$ genau dann als eine nicht eskalierte 
mengenweise Distanzfunktion 
auf dem metrischen Raum $(\Theta,d)$, wenn
f"ur alle $X,Y\subset \Theta$ der Limes superior
\begin{equation}
\lim_{a\to\infty}\sup\{{\rm D}(X\cap Z,Y\cap Z):d(Z)>a\}=
{\rm D}(X,Y)
\end{equation}
existiert. F"ur jede Menge $Z\subset \Theta$ bezeichnet
$d(Z)$ hierbei deren Durchmessser $\sup\{d(x,y):x,y\in Z\}$.
Ist
$(\Theta,||\mbox{id}||)$ ein normierter Raum, so
nennen wir ${\rm D}$ dementsprechend genau dann eine nicht eskalierte\index{nicht eskalierte mengenweise Distanzfunktion} 
mengenweise Distanzfunktion 
auf dem normierten Raum $(\Theta,||\mbox{id}||)$, wenn
${\rm D}$ eine nicht eskalierte 
mengenweise Distanzfunktion 
auf dem durch $||\mbox{id}||$ induzierten metrischen Raum $(\Theta,||\mathbf{P}_{1}-\mathbf{P}_{2}||)$
ist. Exakt jede monotone und stetige mengenweise Distanzfunktion 
auf dem topologischen Raum $(\Theta,\mathbf{T}(d))$, die
eine nicht eskalierte 
mengenweise Distanzfunktion 
auf dem metrischen Raum $(\Theta,d)$ ist, bezeichnen wir 
eine 
moderate mengenweise Distanzfunktion 
auf dem metrischen Raum $(\Theta,d)$.\index{moderate mengenweise Distanzfunktion}
Die  Distanzfunktion $g$
charakterisieren wir demnach als eine auf dem Zahlenstrahl
moderate mengenweise Distanzfunktion -- oder pedantisch ausgesprochen: $g$ ist
eine auf dem metrischen Raum $(\mathbb{R}, d_{1})$ moderate mengenweise Distanzfunktion,
wobei $d_{1}=|\mathbf{P}_{1}-\mathbf{P}_{2}|$ die euklidische Distanz ist. 
\newline
Die auf eine mengenweise Zustandsraumdistanz bezogene
Komanenz und die auf eine mengenweise Zahlenstrahldistanz bezogene
konverse 
Komanenz stehen 
auf der strukturellen Ebene 
immer in dem konstituellen Verh"altnis zueinander, analoge
Konstruktionen zu sein. 
Im konkreten Fall einer mitsamt einer vorliegenden mengenweisen Zustandsraumdistanz ${\rm D}_{1}$
und einer jeweiligen mengenweisen Zahlenstrahldistanz ${\rm D}_{2}$
gegebenen, strukturierten Flussfunktion $(\xi, \mathcal{A})$
ist aber "uber dieses konstante strukturelle Verh"altnis hinaus
ein jeweiliges logisches Verh"altnis konkretisiert.
Dieses konkrete Verh"altnis ist
von der jeweiligen Wahl der Distanzfunktionen bestimmt. Dabei steht 
die jeweilige Komanenz und die konkretisierte konverse Komanenz neben der
von der Wahl der Distanzfunktionen losgel"osten internen
Komanenz, die alleine durch die Zustandsraum"uberdeckung 
$\mathcal{A}$ festgelegt ist.\newline
Obwohl wir unsere Untersuchungen eng an unsere Zielsetzung binden,
Kriterien zu bestimmen, die uns sagen sollen,
wann die durch die Zustandsraum"uberdeckung
$\mathcal{A}$ festgelegte Ultrakolokalisation  
eine "Aquivalenzrelation ist, ist die nun folgende Bemerkung 4.1 durch die Tatsache gerechtfertigt, dass
wir im Beweis des Satzes 4.2 die Bemerkung 4.1 geltend machen werden.
Die Bemerkung 4.1
thematisiert
das 
logische Verh"altnis, in dem die
konverse Komanenz zu der kommutativen Cantorstetigkeit des Phasenflusses
einer vorliegenden strukturierten Flussfunktion $(\xi, \mathcal{A})$ steht, die durch
die zus"atzliche Gegebenheit
der Stetigkeit der jeweiligen strukturierten Flussfunktion $(\xi, \mathcal{A})$
in der ersten Ver"anderlichen gem"ass (\ref{vorwoolc}) erg"anzt ist.\newline
Wir schicken dies vorweg:
Das bei der Festlegung (\ref{endula}) ausgesprochene 
Argument, weshalb die Stetigkeit eines 
topologisierten 
Phasenflusses
im allgemeinen
nicht die Stetigkeit der durch denselben 
festgelegten topologisierten Flussfunktion impliziert, gilt offensichtlich allgemeiner:  Die kommutative
Cantorstetigkeit des Phasenflusses $\{\xi^{t}\}_{t\in\overline{\mathbb{R}}}$ bez"uglich der
Zustandsraum"uberdeckung 
$\mathcal{A}$
einer vorliegenden strukturierten Flussfunktion $(\xi, \mathcal{A})$ impliziert zwar, dass
f"ur alle $t\in\overline{\mathbb{R}}$ und alle $x\in \mathbf{P}_{2}\xi$
das verallgemeinerte Pendant zu der Stetigkeit der 
Flussfunktion in der ersten Ver"anderlichen gem"ass (\ref{ydnulb}) wahr ist, n"amlich die
Aussage
\begin{equation}\label{vorwoolb}
\forall\ {\rm A}\in \mathcal{A}_{\{\xi(x,t)\}}\ \exists\ 
\bar{{\rm A}}\in \mathcal{A}_{\{x\}}\ :
\xi(\bar{{\rm A}},t)\subset{\rm A}\ .
\end{equation}
Die auf die Zustandsraum"uberdeckung 
$\mathcal{A}$ bezogene kommutative Cantorstetigkeit des Phasenflusses $\{\xi^{t}\}_{t\in\overline{\mathbb{R}}}$
impliziert aber im allgemeinen
nicht die verallgemeinerte Form der Stetigkeit der 
Flussfunktion in der zweiten Ver"anderlichen gem"ass (\ref{ydnulc}), welche
die Aussage 
\begin{equation}\label{vorwoolc}
\begin{array}{c}
\forall\ t\in \overline{\mathbb{R}},\ x\in\mathbf{P}_{2}\xi,\ {\rm A}\in \mathcal{A}_{\{\xi(x,t)\}} \exists\ 
\delta\in\mathbb{R}^{+}:\\
\xi(x,]t-\delta,t+\delta[)\subset{\rm A}
\end{array}
\end{equation}
ist. F"ur jede strukturierte Flussfunktion im weiteren Sinn legen wir daher die zu der Bezeichnungsweise
der Stetigkeit topologisierter Flussfunktionen analoge 
Benennung fest: Jede strukturierte Flussfunktion $(\xi, \mathcal{A})$ im weiteren Sinn
gelte genau dann als stetig, wenn f"ur sie sowohl die zust"andliche 
Stetigkeit\index{zust\"andliche Stetigkeit einer strukturierten Flussfunktion im weiteren Sinn}
in der ersten Ver"anderlichen gem"ass
(\ref{vorwoolb}) als auch die phasische Stetigkeit\index{phasische Stetigkeit einer strukturierten Flussfunktion im weiteren Sinn} gem"ass (\ref{vorwoolc})
gegeben ist, bei welcher zudem die Zentriertheit\index{Zentriertheit} $\xi(\mbox{id},0)=\mbox{id}$ gelte.\index{stetige strukturierte Flussfunktion}
Von der konversen phasischen Stetigkeit des Paares $(\xi,\mathcal{A})$
reden wir genau dann,
wenn f"ur es $\xi(\mbox{id},0)=\mbox{id}$ gilt und die Aussage 
\begin{equation}\label{corwoolv}
\begin{array}{c}
\forall\ t\in \overline{\mathbb{R}},\ x\in\mathbf{P}_{2}\xi,\ \delta\in\mathbb{R}^{+}\
\exists\ 
{\rm A}\in \mathcal{A}_{\{\xi(x,t)\}}:\\
{\rm A}\cap\xi(x,\mathbb{R})\subset\xi(x,]t-\delta,t+\delta[)
\end{array}
\end{equation}
wahr ist. Jede strukturierte Flussfunktion im weiteren Sinn $(\xi,\mathcal{A})$
nennen wir exakt in dem Fall perfekt phasisch stetig, wenn die 
Zentriertheit $\xi(\mbox{id},0)=\mbox{id}$ gegeben ist und wenn
die Aussagen (\ref{vorwoolc}) und (\ref{corwoolv}) beide wahr sind, wenn also $(\xi,\mathcal{A})$\index{perfekt phasische Stetigkeit einer strukturierten Flussfunktion im weiteren Sinn} sowohl phasisch als auch konvers phasisch stetig ist.\index{konvers phasische Stetigkeit einer strukturierten Flussfunktion im weiteren Sinn} Wir wollen an dieser Stelle noch die perfekte Phasenstetigkeit 
einer Flussfunktion $(\xi,\mathcal{A})$ im weiteren Sinn mit deren
intrinsischer Insensitivit"at
vergleichen:
Ist $(\xi,\mathcal{A})$ perfekt phasisch stetig,
gilt erstens deren Phasenstetigkeit in Form
der
Aussage
\begin{displaymath}
\begin{array}{c}
\forall\ t\in\overline{\mathbb{R}},\ \alpha\in\mathbf{P}_{2}\xi,\ {\rm L}^{+}\in \mathcal{A}_{\{\xi(\alpha,t)\}} \exists\ 
\varepsilon\in\mathbb{R}^{+}\\
\xi(\alpha,]t-\varepsilon,t+\varepsilon[)\subset{\rm L}^{+}\ ,
\end{array}
\end{displaymath}
die wegen der Zentriertheit $\xi(\mbox{id},0)=\mbox{id}$ ihre Paraphrase
\begin{displaymath}
\begin{array}{c}
\forall\ t\in\overline{\mathbb{R}},\ \alpha\in\mathbf{P}_{2}\xi,\ {\rm L}^{+}\in \mathcal{A}_{\{\xi(\alpha,t)\}} \exists\ 
\varepsilon\in\mathbb{R}^{+}\\
\xi(\xi(\alpha,t),]-\varepsilon,\varepsilon[)\subset{\rm L}^{+}
\end{array}
\end{displaymath}
impliziert.
Zweitens ist f"ur $(\xi,\mathcal{A})$ konverse Phasenstetigkeit gegeben,
welche
besagt, dass
die Aussage
\begin{displaymath}
\begin{array}{c}
\forall\ t\in \overline{\mathbb{R}},\ \alpha\in\mathbf{P}_{2}\xi,\ \varepsilon\in\mathbb{R}^{+}\ \exists\ 
{\rm L}^{-}\in \mathcal{A}_{\{\xi(\alpha,t)\}}\\
{\rm L}^{-}\cap\xi((\alpha,\mathbb{R})\subset\xi(\alpha,]t-\varepsilon,t+\varepsilon[)
\end{array}
\end{displaymath}
wahr ist,
sodass mit der Phasenstetigkeit zusammen die Aussage
\begin{displaymath}
\begin{array}{c}
\forall\ t\in \mathbb{R},\ \alpha\in\mathbf{P}_{2}\xi,\ \varepsilon\in\mathbb{R}^{+}\ \exists\ 
\delta\in\mathbb{R}^{+}\\
\xi(\xi(\alpha,t),]-\delta,\delta[)\subset\xi(\alpha,]t-\varepsilon,t+\varepsilon[)
\end{array}
\end{displaymath}
resultiert, welche aber von der Aussage 
\begin{displaymath}
\begin{array}{c}
\exists\ \delta\in\mathbb{R}^{+}\ \forall\ t\in \mathbb{R},\ \alpha\in\mathbf{P}_{2}\xi,\ \varepsilon\in\mathbb{R}^{+}\ \\
\xi(\xi(\alpha,t),]-\delta,\delta[)\subset\xi(\alpha,]t-\varepsilon,t+\varepsilon[)
\end{array}
\end{displaymath} 
differiert, die
gerade die intrinsische Insensitivit"at der Flussfunktion $(\xi,\mathcal{A})$ im weiteren Sinn behauptet.
\index{intrinsische Sensitivit\"at}
\newline
\newline
{\bf Bemerkung 4.1:}\newline
{\em Jede phasisch stetige strukturierte Flussfunktion $(\xi, \mathcal{A})$ ist
konvers komanent\index{konverse Komanenz} bez"uglich jeder moderaten mengenweisen Zahlenstrahldistanz ${\rm D}_{2}$.}
\newline
\newline
{\bf Beweis:}\newline
Es sei $y\in\mathbf{P}_{2}\xi$ und $\{x_{j}\}_{j\in \mathbb{N}}$ ein Folge von Zust"anden $x_{j}\in \mathbf{P}_{2}\xi$,
f"ur die
$$x_{j}-\mathcal{A}\to y$$
gelte. Dann gibt es f"ur jede Lokalisierung $A\in \mathcal{A}$ zu jeder 
reellen Zahl  
$$t\in\xi(y,\mbox{id})^{-1}(A)$$
einen Index $j^{\star}(t)\in \mathbb{N}$, f"ur den f"ur alle $j\in \mathbb{N}$ die
Implikation
\begin{equation}\label{rinulv}
j> j^{\star}(t)\Rightarrow \xi(x_{j},t)\in A
\end{equation}
wahr ist. Wegen der phasischen Stetigkeit der 
Flussfunktion in der zweiten Ver"anderlichen gem"ass
(\ref{vorwoolc}) gibt es zu jeder Zahl 
$t\in\xi(y,\mbox{id})^{-1}(A)$ eine 
positive reelle Zahl $\delta(t)$, f"ur welche
die Inklusion
$$\xi(x_{j^{\star}(t)},]t-\delta(t),t+\delta(t)[)\subset A$$
zutrifft.  
F"ur alle $a\in\mathbb{R}^{+}$ 
ist daher
$$\Bigr\{]t-\delta(t),t+\delta(t)[:t\in\mathbf{cl}(\xi(y,\mbox{id})^{-1}(A)\cap [-a,a])\Bigl\}$$
eine "Uberdeckung des Kompaktums
$$\mathbf{cl}(\xi(y,\mbox{id})^{-1}(A)\cap [-a,a]\ .$$
Alle Elemente dieser "Uberdeckung sind offen, sodass es
eine endliche Menge $J(a)\subset \mathbf{cl}(\xi(y,\mbox{id})^{-1}(A)\cap [-a,a]$ gibt, f"ur die
$$\bigcup\Bigr\{]t-\delta(t),t+\delta(t)[:t\in J(a)\Bigl\}\supset\mathbf{cl}(\xi(y,\mbox{id})^{-1}(A)\cap [-a,a]$$
gilt. Also gibt es das Maximum 
$$\max\{j^{\star}(t):t\in J(a)\}$$
all der Indizes $j^{\star}(t)\in \mathbb{N}$, f"ur die jeweils 
die Implikation (\ref{rinulv}) wahr ist.
Da die mengenweise Zahlenstrahldistanz ${\rm D}_{2}$ moderat ist,
gibt es f"ur alle $\varepsilon\in\mathbb{R}^{+}$ 
und f"ur alle $j\in \mathbb{N}$
eine positive Zahl $a(\varepsilon)\in\mathbb{R}^{+}$, f"ur welche die Ungleichungskette
$${\rm D}_{2}\Bigl(\xi(x_{j},\mbox{id})^{-1}(A),\xi(y,\mbox{id})^{-1}(A)\Bigr)-\varepsilon$$
$$\leq{\rm D}_{2}\Bigl(\xi(x_{j},\mbox{id})^{-1}(A)\cap[-a(\varepsilon),
a(\varepsilon)],\ \xi(y,\mbox{id})^{-1}(A)\cap
[-a(\varepsilon),a(\varepsilon)])\Bigr)$$
$$\leq {\rm D}_{2}\Bigl(\xi(x_{j},\mbox{id})^{-1}(A),\ \xi(y,\mbox{id})^{-1}(A)\Bigr)$$
gilt. Wenn 
$k\in\mathbb{N}$ dabei ein Index ist, welcher der Beschr"ankung
$$k\leq\max\{j^{\star}(t):t\in J(a(\varepsilon))\}$$
gen"ugt, so gilt die Inklusion
$$\xi(x_{k},\mbox{id})^{-1}(A)\supset\xi(y,\mbox{id})^{-1}(A)\cap
[-a(\varepsilon),a(\varepsilon)]\ ,$$
sodass
$${\rm D}_{2}\Bigl(\xi(x_{k},\mbox{id})^{-1}(A),\ \xi(y,\mbox{id})^{-1}(A)\Bigr)\ -$$
$${\rm D}_{2}\Bigl(\xi(x_{k},\mbox{id})^{-1}(A)
\cap[-a(\varepsilon),a(\varepsilon)],\ \xi(y,\mbox{id})^{-1}(A)\cap[-a(\varepsilon),a(\varepsilon)]\Bigr)$$
$$<\varepsilon$$
ist.\newline
{\bf q.e.d.}
\section{Ultrakolokalisation als "Aquivalenzrelation}
Wir hoben die 
konvers komanenten Flussfunktionen
hervor, denn es gilt f"ur sie folgender
\newline
\newline
{\bf Satz 4.2: Transitivit"at der Ultrakolokalisation\index{Transitivit\"at der Ultrakolokalisation}}\newline
{\em Es sei $\Lambda=(\xi, \mathcal{A})$ eine strukturierte Flussfunktion und
\begin{equation}\label{ydnulg}
\Lambda^{\star}:=(\xi|\mathbf{P}_{2}\xi^{\star}\times\overline{\mathbb{R}}, \mathcal{A}\cap\mathbf{P}_{2}\xi^{\star})
\end{equation} deren Restriktion 
auf den Reflexivit"atsbereich der Ultrakolokalisation\index{Reflexivit\"atsbereich der Ultrakolokalisation
einer strukturierten Flussfunktion}
\begin{equation}\label{ydnulh}
\mathbf{P}_{2}\xi^{\star}:=\{z\in \mathbf{P}_{2}\xi:z\sim_{\Lambda} z\}
\end{equation}
der strukturierten Flussfunktion $\Lambda=(\xi, \mathbf{T}(\mathbf{P}_{2}\xi))$.
Dann ist
die Ultrakolokalisation $\sim_{\Lambda^{\star}}$ eine
"Aquivalenzrelation, falls f"ur das Spurmengensystem
\begin{equation}\label{ydlalh}
\mathcal{A}^{\star}:=\mathcal{A}\cap \mathbf{P}_{2}\xi^{\star}
\end{equation}
Folgendes gilt:\newline
(a) Das Mengensystem} 
\begin{equation}\label{yhdnulh}
\Bigl\{\xi(z,\mbox{id})^{-1}(A):(z,A)\in \mathbf{P}_{2}\xi^{\star}
\times\mathcal{A}^{\star} \Bigr\}\\
\subset \kappa^{\star}\mathbf{T}(1) 
\end{equation}
{\em  hat nur unbefristet keuliale Mengen\index{keuliale Menge}\index{unbefristet keuliale Menge}
als Elemente 
und die Restriktion $\Lambda^{\star}$
ist
konvers komanent bez"uglich der in} (\ref{auznlv}) {\em festgelegten Distanz $g$.\newline
(b) Die Restriktion $\Lambda^{\star}$ ist
bez"uglich der Zustandsraum"uberdeckung $\mathcal{A}^{\star}$ phasisch stetig.}
\newline
\newline
{\bf Beweis:}\newline
Zuallererst stellt sich die Frage,
ob die Restriktion 
(\ref{ydnulg}) "uberhaupt eine strukturierte Flussfunktion ist,
ob also $\xi|\mathbf{P}_{2}\xi^{\star}\times\overline{\mathbb{R}}$ eine Flussfunktion ist. Das ist
genau dann der Fall, wenn f"ur alle $t\in \overline{\mathbb{R}}$
\begin{equation}\label{ydnuli}
\xi^{t}\mathbf{P}_{2}\xi^{\star}=\mathbf{P}_{2}\xi^{\star}
\end{equation}
invariant und $\xi|\mathbf{P}_{2}\xi^{\star}\times\overline{\mathbb{R}}$
kinematisch abgeschlossen ist. Wenn ein Zustand $x\in \mathbf{P}_{2}\xi\setminus \mathbf{P}_{2}\xi^{\star}$
ausserhalb des Reflexivit"atsbereiches der Ultrakolokalisation der
Flussfunktion $\Lambda$ ist,
dann ist per definitionem f"ur alle $t\in \overline{\mathbb{R}}$ auch der Zustand 
$\xi(x,t)$ kein Element des Reflexivit"atsbereiches $\mathbf{P}_{2}\xi^{\star}$, sodass die Implikation
$$x\in \mathbf{P}_{2}\xi\setminus \mathbf{P}_{2}\xi^{\star}\Rightarrow
\xi(x,\overline{\mathbb{R}})\not\subset\mathbf{P}_{2}\xi^{\star}$$
wahr ist und der Reflexivit"atsbereiches der Ultrakolokalisation
$$\mathbf{P}_{2}\xi^{\star}=\mathbf{P}_{2}\xi\setminus\bigcup\Bigl\{\xi(x,\overline{\mathbb{R}}):
x\in \mathbf{P}_{2}\xi\setminus \mathbf{P}_{2}\xi^{\star}\Bigr\}$$
ist, sodass (\ref{ydnuli}) gilt.
Die Symmetrie der Ultrakolokalisation $\sim_{\Lambda^{\star}}$ ist per constructionem angelegt. Der Kern dieses Beweises
ist es, im Fall (a) die
Transitivit"at der Ultrakolokalisation $\sim_{\Lambda^{\star}}$
zu zeigen.
\newline
Zu (a): Dass $x\sim_{\Lambda^{\star}}z$ gilt,
ergibt sich im Fall der unbefristeten
Keulialit"at gem"ass (\ref{yhdnulh}) aus der
konversen Komanenz der Restriktion $\Lambda^{\star}$:
Wenn f"ur drei Zust"ande $x,y,z\in \mathbf{P}_{2}\xi^{\star}$ die Aussage
\begin{equation}\label{vora}
x\sim_{\Lambda^{\star}}y\quad\land\quad y\sim_{\Lambda^{\star}}z
\end{equation}
wahr ist, gibt 
es Zust"ande
$\omega(x,y),\omega(y,z)\in \mathbf{P}_{2}\xi^{\star}$, f"ur welche 
es zwei 
beidseitig bestimmt divergente
Zahlenfolgen $\{t_{j}\}_{j\in\mathbb{N}}$ und $\{s_{j}\}_{j\in\mathbb{N}}$ gibt, f"ur die
\begin{equation}\label{wolnua}
\begin{array}{c}
\xi(y,t_{j})-\mathcal{A}\to\omega(x,y)\ \land\\
\xi(x,t_{j})-\mathcal{A}\to\omega(x,y)
\end{array}
\end{equation}
und entsprechend
\begin{equation}\label{solnua}
\begin{array}{c}
\xi(y,s_{j})-\mathcal{A}\to\omega(y,z)\ \land\\
\xi(x,s_{j})-\mathcal{A}\to\omega(y,z)
\end{array}
\end{equation}
gilt. Die Aussage (\ref{wolnua}) bestimmt 
wegen der vorausgesetzten konversen Komanenz der Restriktion $\Lambda^{\star}$ bez"uglich der 
Distanzfunktion $g$ f"ur alle ${\rm L}\in \mathcal{A}^{\star}$ die Grenzwerte als
\begin{equation}\label{wolnub}
\begin{array}{c}
\lim_{j\to\infty}g\Bigl(\xi(\xi(y,t_{j}),\mbox{id})^{-1}({\rm L}),\xi(\omega(x,y),\mbox{id})^{-1}({\rm L})\Bigr)=0\ ,\\ 
\lim_{j\to\infty}g\Bigl(\xi(\xi(y,t_{j}),\mbox{id})^{-1}({\rm L}),\xi(\omega(x,y),\mbox{id})^{-1}({\rm L})\Bigr)=0\ ,
\end{array}
\end{equation}
sodass der Grenzwert
\begin{equation}\label{wolnuc}
\lim_{j\to\infty}g\Bigl(\xi(\xi(y,t_{j}),\mbox{id})^{-1}({\rm L}),\xi(\xi(x,t_{j}),\mbox{id})^{-1}({\rm L})\Bigr)=0
\end{equation}
ist.
Es
gilt f"ur die Zahlenfolge $\{s_{j}\}_{j\in\mathbb{N}}$ und 
f"ur alle ${\rm L}\in \mathcal{A}^{\star}$
dementsprechend
\begin{equation}\label{wolnuf}
\lim_{j\to\infty}g\Bigl(\xi(\xi(y,s_{j}),\mbox{id})^{-1}({\rm L}),\xi(\xi(z,s_{j}),\mbox{id})^{-1}({\rm L})\Bigr)=0\ . 
\end{equation}
Dabei gilt f"ur alle $\alpha,\beta\in\mathbb{R}$
und f"ur alle ${\rm V}\in \mathcal{A}^{\star}$ und f"ur jeden Zustand $w$
\begin{equation}\label{wolnug}
\xi(\xi(w,\alpha),\mbox{id})^{-1}({\rm V})-\{\beta-\alpha\}=
\xi(\xi(w,\beta),\mbox{id})^{-1}({\rm V})\ ,
\end{equation}
falls der Phasenfluss $\{\xi^{t}|\mathbf{P}_{2}\xi^{\star}\}_{t\in\overline{\mathbb{R}}}$
abelsch ist, d.h., falls die Gruppe 
$$(\{\xi^{t}|\mathbf{P}_{2}\xi^{\star}:t\in\mathbb{R}\},\circ)$$ eine 
additive Gruppe ist und f"ur alle $\alpha,\beta\in\mathbb{R}$ die Komposition $\xi^{\alpha}\circ\xi^{\beta}= \xi^{\alpha+\beta}$. Allemal gibt es aber f"ur alle $\alpha,\beta\in\mathbb{R}$ eine reelle Zahl $r(\alpha,\beta)$, f"ur die
f"ur alle Lokalisierungen ${\rm V}\in \mathcal{A}^{\star}$ und f"ur jeden Zustand $w$ 
\begin{equation}\label{solnug}
\xi(\xi(w,\alpha),\mbox{id})^{-1}({\rm V})+\{r(\alpha,\beta)\}=
\xi(\xi(w,\beta),\mbox{id})^{-1}({\rm V})
\end{equation}
ist. Es gilt also nicht nur, dass f"ur alle ${\rm L}\in \mathcal{A}^{\star}$
\begin{equation}\label{wolnuh}
\lim_{j\to\infty}g\Bigl(\xi(\xi(y,t_{j}),\mbox{id})^{-1}({\rm L}),\xi(\xi(z,t_{j}),\mbox{id})^{-1}({\rm L})\Bigr)=0\
\end{equation}
ist, sondern 
schliesslich, dass auch f"ur alle ${\rm L}\in \mathcal{A}^{\star}$
\begin{equation}\label{wolnuf}
\lim_{j\to\infty}g\Bigl(\xi(\xi(x,t_{j}),\mbox{id})^{-1}({\rm L}),\xi(\xi(z,t_{j}),\mbox{id})^{-1}({\rm L})\Bigr)=0
\end{equation}
ist.
Wegen
der vorausgesetzten\index{unbefristet keuliale Menge} unbefristeten Keulialit"at
(\ref{yhdnulh}) ist also nicht nur f"ur alle ${\rm Q}\in\mathcal{A}^{\star}_{\{\omega(y,z)\}}$ die Reichhaltigkeitsaussage
\begin{equation}\label{wolnug}
\xi(x,\mbox{id})^{-1}({\rm Q})\cap\xi(z,\mbox{id})^{-1}({\rm Q})\not=\emptyset
\end{equation}
wahr: Die Zahlenfolgen $\{t_{j}\}_{j\in\mathbb{N}}$ und $\{s_{j}\}_{j\in\mathbb{N}}$, f"ur die
(\ref{wolnua}) und (\ref{solnua}) gilt, sind hierbei
beidseitig divergent. 
Es gibt also 
f"ur
alle nicht negativen reellen Zahlen $t_{\star}\in \mathbb{R}^{+}_{0}$ streng monotone und einseitig bestimmt divergente 
Zahlenfolgen $\{q_{j}^{\pm}\}_{j\in\mathbb{N}}$, f"ur welche die Beschr"ankung
$q_{j}^{+}>t_{\star}$ bzw. die Ungleichung $q_{j}^{-}<-t_{\star}$ eingehalten ist,
wobei f"ur alle ${\rm L}\in \mathcal{A}^{\star}$
\begin{equation}\label{wnuh}
\lim_{j\to\infty}g\Bigl(\xi(\xi(x,q_{j}),\mbox{id})^{-1}({\rm L}),\xi(\xi(z,q_{j}),\mbox{id})^{-1}({\rm L})\Bigr)=0
\end{equation}
ist.
"Uber (\ref{wolnug}) hinaus gilt also f"ur alle $t_{\star}\in \mathbb{R}^{+}$ und f"ur alle ${\rm Q}\in\mathcal{A}^{\star}_{\{\omega(y,z)\}}$
wegen
der vorausgesetzten\index{unbefristet keuliale Menge} unbefristeten Keulialit"at die Reichhaltigkeitsaussage
\begin{equation}\label{wolnuh}
\Bigl(\xi(x,\mbox{id})^{-1}({\rm Q})\cap\xi(z,\mbox{id})^{-1}({\rm Q})\Bigr)\setminus[-t_{\star},t_{\star}]\not=\emptyset\ ,
\end{equation}
welche 
die 
Aussage
$$x\sim_{\Lambda^{\star}}z\ .$$
definiert.
Zu (b): 
Wenn die Restriktion $\Lambda^{\star}$ bez"uglich $\mathcal{A}^{\star}$
phasisch stetig ist,
gilt f"ur alle $t\in \mathbb{R}$ die Stetigkeit in der
zweiten Ver"anderlichen in Form der Aussage
\begin{equation}\label{last}
\forall\ x\in\mathbf{P}_{2}\xi, {\rm A}\in \mathcal{A}^{\star}_{\{\xi(x,t)\}}\ \exists\ 
\delta\in\mathbb{R}^{+}:\ 
\xi(x,[t-\delta,t+\delta])\subset{\rm A}\ .
\end{equation}
Nach der Bemerkung 4.1 ist die Restriktion $\Lambda^{\star}$     
konvers komanent bez"uglich der in (\ref{auznlv}) festgelegten Distanz $g$,
die moderat ist. Wegen der Argumentation unter (a) gilt daher die mit der Gleichung
(\ref{wnuh}) 
verbundende Aussage
und wegen (\ref{last})
$x\sim_{\Lambda^{\star}}z$ .\newline
{\bf q.e.d.}
\newline\newline
F"uhren wir auf der Klasse der strukturierten Flussfunktionen im weiteren Sinn
den Operator $\langle\langle\mbox{id}\rangle\rangle$ ein, der jeder jeweiligen
strukturierten Flussfunktion $\Theta=(\vartheta,\theta)$ das Mengensystem
\begin{equation}\label{dramquo}
\langle\langle\Theta\rangle\rangle:=
\Bigl\{\{x\in\mathbf{P}_{2}\vartheta:x\sim_{\Theta}a\}:a\in\mathbf{P}_{2}\vartheta\Bigr\}
\end{equation}
als Wert zuordnet, so stellt sich uns gegebenenfalles
der Sachverhalt, dass $\Theta$
die "Aquivalenzrelation $\sim_{\Theta}$ auf
$\mathbf{P}_{2}\vartheta$ bestimmt,
als die Partitivit"at
$$\langle\langle\Theta\rangle\rangle\in\mathbf{part}(\mathbf{P}_{2}\vartheta)$$
dar. Wie chaotisch ein jeweiliger Attraktor ist, haben wir bislang 
einzig unter dem Gesichtspunkt betrachtet, welche Form der 
Sensitivit"at auf demselben auftritt. Diejenige Sensitivit"atsform, welche
eine andere impliziert, ist eine st"arkere Form der 
Sensitivit"at. Deren st"arkste ist die Ultrasensitivit"at 
$x\sim_{\Lambda}z$ jeweils voneinander 
separabler Zust"ande $x$ und $z$.\newline
Es gibt aber auch den Gesichtspunkt der Ausdehnung der jeweiligen Sensitivit"at,
der bei den Formen mengenweiser Sensitivit"at allerdings 
erst begrifflich erschlossen werden m"usste, 
sofern die Ausdehnungsbewertung f"ur Formen mengenweiser Sensitivit"at "uberhaupt konzipierbar 
ist. Bei den Formen der 
zustandsweisen Sensitivit"at und der Ultrasensitivit"at ist dies m"oglich. Wir haben aber bislang keinerlei
Aussage, die uns eine Auskunft "uber die Beschaffenheit der Bereiche 
$\mbox{coloc}_{-}((\xi|a\times\mathbb{R},\mathcal{A}\cap a))$
innerhalb eines 
Attraktors $a$ 
einer jeweiligen strukturierten Flussfunktion im weiteren Sinn $(\xi,\mathcal{A})$ 
gibt, in denen Zust"ande $z$ liegen, 
die sensitiv sind;
f"ur die es gem"ass
der Definition 3.6 eine Zustandsfolge 
$\{z_{j}\}_{j\in\mathbb{N}}$
gibt, welche in den sensitiven Zust"anden $z$
einer jeweiligen strukturierten Flussfunktion im weiteren Sinn 
kolokalisiert. Bei der Ultrasensitivit"at 
haben wir 
nun immerhin den Sachverhalt, dass konverse 
Komanenz 
deren Bereiche zu Partitionen organisiert. Wie viel 
hilft uns dies im Hinblick auf eine Bewertung 
der Extension des mit der
Ultrasensitivit"at
verbundenen Ultrachaos weiter?
Das Ausmass der trivialen Ultrakolokalisation einer 
strukturierten Flussfunktion $\Theta=(\vartheta,\theta)$ schl"agt sich in der 
Beschaffenheit der Partition $\langle\langle\Theta\rangle\rangle$ folgendermassen nieder:
Genau dann, wenn f"ur 
ein jeweiliges Vorzimmer $\chi\in[[\vartheta]]_{\theta}$ jede Menge
ultrakolokalisierungs"aquivalenter Zust"ande 
$\alpha\in\langle\langle\Theta\rangle\rangle\cap\chi$ einelementig
ist, zeigt sich nirgends im Vorzimmer $\chi$ ultrasensitves 
Verhalten.
Diesen Fall g"anzlich ausgeschlossener Ultrasensitivit"at k"onnen wir daher 
eineindeutig
mit der suggestiven
Redeweise von der totalen
Dissoziation der Ultrakolokalisation\index{totale Dissoziation der Ultrakolokalisation}
auf einem jeweiligen Vorzimmer
verbinden. Die Ultrasensitivit"at innerhalb des Vorzimmers $\chi$ k"onnte aber auf 
die weniger radikale Weise dissoziiert sein, dass beispielsweise die
Implikation
$$\alpha\in\langle\langle\Theta\rangle\rangle\cap\chi 
\Rightarrow \mathbf{card}([\vartheta]_{\alpha})=1$$  
zutrifft und trajektorienweise Dissoziation der Ultrasensitivit"at gegeben ist,
dass keinerlei Ultrasensitivit"at innerhalb einer Trajektorie auftritt.
Oder dass umgekehrt nur innerhalb jeweiliger Trajektorien Ultrasensitivit"at auftritt.
So Zimmer vorliegen, ist gewiss, dass zimmerweise Dissoziation der
Ultrakolokalisation auftritt.
Es gilt dann f"ur jede strukturierte Flussfunktion $\Theta=(\vartheta,\theta)$ und
f"ur deren gem"ass (\ref{dramquo}) festgelegten Wert des 
Operators $\langle\langle\mbox{id}\rangle\rangle$ die Implikation
\begin{equation}
\alpha\in\langle\langle\Theta\rangle\rangle\Rightarrow
\alpha\subset([[\vartheta]]_{\theta})_{\alpha}\ .
\end{equation}
Mehr k"onnen wir im Hinblick auf die extensive 
Bewertung des Ultrachaos nicht unmittelbar sagen. Bei totaler 
Extension
gilt die Aussage
\begin{equation}\label{drama}
\begin{array}{c}
\forall\ t_{\star}\in\mathbb{R}^{+}\exists\ t\in ]-\infty,t_{\star}[\cup]t_{\star},\infty[:\\
\xi(a,t)\in {\rm L}\ \land\ \xi(a,t)\in  {\rm L} 
\end{array}
\end{equation}
f"ur alle Lokalisierungen 
${\rm L}\in\mathcal{A}$ und f"ur alle $a,b\in \mathbf{P}_{2}\xi$ und es
liegt dann f"ur die jeweilige  
strukturierte Flussfunktion im weiteren Sinn
$(\xi, \mathcal{A})$ diejenige Form von Sensitivit"at vor, die 
wir als 
die 
ultimative Sensitivit"at einer strukturierten Flussfunktion im weiteren Sinn $\Lambda=(\xi,\mathcal{A})$
bezeichnen.\index{ultimative Sensitivit\"at einer strukturierten Flussfunktion}
Die ultimative Sensitivit"at kann nur innerhalb eines Attraktors 
auftreten und sie legt das fest, was wir 
dementsprechend als ultimatives Chaos\index{ultimatives Chaos}
bezeichnen.
Denn, wenn (\ref{drama}) gilt, sehen wir diejenige Form von Sensitivit"at vorliegen, die
wir uns "uberhaupt als die dramatischste deren Formen vorstellen k"onnen. 
Dass 
ultimative Sensitivit"at auftritt und
f"ur alle ${\rm L}\in\mathcal{A}$ und $a,b\in \mathbf{P}_{2}\xi$
die 
Aussage 
(\ref{drama})
gilt, dies 
paraphrasiert lediglich die Aussage, dass jedes Vorzimmer $\chi$
des Mengensystemes $[[\xi]]_{\mathcal{A}}$
so beschaffen ist,
dass f"ur alle $x,y,z\in \chi$
$$x,y\sim>_{\Lambda}z$$
gilt, dass also
$$\langle\langle\Lambda\rangle\rangle=\mathbf{P}_{2}\xi$$ 
ist. Ultimatives Chaos kann keinesfalls auf einem 
Zustandsraum auftreten, welcher verschiedene Vorzimmer hat.
Ultimatives Chaos kann aber f"ur einzelne Vorzimmer $\chi$ auftreten
als das ultimative Chaos jeweiliger
Restriktionen $(\xi|(\chi\times\mathbb{R}),\mathcal{A}\cap\chi)$. 
Es ist dabei dann ausgeschlossen, dass ein anderes Vorzimmer 
$\tilde{\chi}\in [[\xi]]_{\mathcal{A}}\setminus\{\chi\}$
existiert, welches $\chi$ schneidet. Ein 
Vorzimmer, auf dem sich 
ultimatives Chaos zeigt, ist demnach ein Zimmer 
des Mengensystemes 
der Vorzimmer $[[\xi]]_{\mathcal{A}}$ der  
strukturierten Flussfunktion im weiteren Sinn $\Lambda=(\xi,\mathcal{A})$.

\newpage
{\Large {\bf Symbolverzeichnis}}\newline\newline
Neben den v"ollig etablierten Bezeichnungen benutzen wir auch in dieser Abhandlung
spezielle, f"ur dieses Traktat eigent"umliche Schreibweisen. 
Diese weniger gel"aufigen Schreibweisen f"uhren wir auch im zweiten Teil der Konzepte der 
abstrakten Ergodentheorie
im fortlaufenden
Text ein und listen sie hier auf.\newline
Wenn wir nun jeweilige idiomatische Notationen erkl"aren, kommt es 
freilich vor, dass diese Erkl"arung eine spezielle Begriffsbildung 
dieser Anbhandlung gebraucht. Die jeweilige spezielle Begriffsbildung findet sich aber im Index vermerkt,
sodass sich im Text die Stelle auffinden l"asst, wo diese 
spezielle Begriffsbildung formuliert wird.
\newline
Es sei $n,m\in\mathbb{N}$, $r\in \mathbb{R}^{+}$, $t\in \mathbb{R}$, $(X,d)$ ein metrischer Raum, $\mathbf{T}$
eine Topologie, 
${\rm E}=\{e\}$
eine einelementige Menge, $j\in\{1,2,\dots n\}$,
$\theta=(\theta_{1},\theta_{2},\dots \theta_{n})$
ein $n$-Tupel, 
$\lambda=(\lambda_{1},\lambda_{2},\dots \lambda_{m})$ ein $m$-Tupel,
$\Psi$ eine Flussfunktion im weiteren Sinn,
$\Phi$ eine metrisierte Flussfunktion, 
$\Lambda$ eine strukturierte Flussfunktion,
$\psi$ eine reelle Flussfunktion, $\xi$ ein Autobolismus, $\Xi$ 
eine Menge von Autobolismen.
$\phi$ sei eine Funktion und es sei ${\rm A}$ und ${\rm B}$ eine Menge, $\mathcal{A}$
hingegen ein
Mengensystem, wobei
$$\mathcal{A}^{c}:=\Bigl\{Y\setminus A:A\in\mathcal{A}\Bigr\}$$
das zu $\mathcal{A}$ komplement"are Mengensystem und $Y=\bigcup\mathcal{A}$ sei.
F"ur eine Folge von Mengen $\{A_{j}\}_{j\in\mathbb{N}}$ des Mengensystemes 
$\mathcal{A}$ und $x\in\bigcup\mathcal{A}$
sei der Sinn der Aussage
$$A_{j}-\mathcal{A}\to x$$
gem"ass (\ref{ydnulw}) festgelegt. Die Notation mit Hilfe von
$$\sim_{\Lambda}\ , \ \sim>_{\Lambda}$$
wird vor der Definition 3.5 eingef"uhrt. 
Ist die mit $\sim_{\Lambda}$ bezeichnete Ultrakolokation eine "Aquivalenzrelation ist, bezeichnet
$\langle\langle\Lambda\rangle\rangle$ deren "Aquivalenzklassen.\newline
Es ist 
\newline
\begin{tabbing}
Wir benutzen zudem folgende \= Etwas \kill
${\rm A}^{\mathbb{N}}$ \> die Menge aller Folgen $\{a_{j}\}_{j\in\mathbb{N}}$,\\ 
\quad \> deren Glieder Elemente der Menge ${\rm A}$ sind,\\
$\mbox{{\bf @}}(\Phi,\mathcal{A})$\> die Menge freier Attraktoren bzgl. $\mathcal{A}$, \\
$\mathcal{A}^{\cup}$ \> das Mengensystem aller Vereinigungen\\ \quad \> "uber Teilmengen von $\mathcal{A}$,\\
$\mathcal{A}_{{\rm A}}$ \> die ${\rm A}$-Auswahl aus $\mathcal{A}$, d.h\\
\quad \> $\{a\in \mathcal{A}:a\cap{\rm A}\not=\emptyset\}\ ,$\\
$\Delta_{\Phi}$ \> das Aufl"osungsfeld von $\Phi$,\\
$\mathbb{B}_{r}(x)$ \> die offene Kugel des $(X,d)$, in dem $x$ ist, deren\\  
\quad \>  Mittelpunkt $x$ ist und die den Radius $r$ hat,\\
$\Psi^{t}$ \> der durch $\Psi$ und $t$ festgelegte Phasenfluss,\\   
\quad \> d.h die Bijektion $\Psi(\mbox{id},t)$\\ 
\quad \>  des Zustandsraumes auf denselben,\\ 
$[\Psi]$  \> die durch $\Psi$ festgelegte trajektorielle Partition,\\ 
\quad \> d.h $[\Psi]:=\{\Psi(x,\mathbb{R}):x\in\mathbf{P}_{2}\Psi\}\ ,$\\ 
$[[\Psi]]_{\mathcal{A}}$ \> die durch $\Psi$ festgelegte Partition \\ 
\quad \>  $[[\Psi]]_{\mathcal{A}}:=\{\mathbf{cl}_{\mathcal{A}}(\tau):\tau\in[\Psi]\}\ ,$\\ 
\quad \> wobei\\
$\mathbf{cl}_{\mathcal{A}}$ \> der H"ullenoperator bzgl. $\mathcal{A}$ ist,\\ 
$:=\bigcap\{A\in\mathcal{A}:A\supset\mbox{id}\}$\> \quad \\
$\downarrow {\rm E} =e$ \> das Element der einelementige Menge ${\rm E}$, \\
$\mathbf{part}({\rm A})$  \> die Menge aller Partitionen der Menge ${\rm A}$,\\
$\mathbf{P}_{j}\Lambda=\lambda_{j}$ , \> \quad\\
$\mathbf{P}_{1}\phi$ \> die Definitionsmenge von $\phi$,\\
$\mathbf{P}_{2}\phi$ \>  die Wertemenge von $\phi$,\\
$\oplus$ \> die Kokatenation,\\ \quad \> sodass $\theta\oplus\lambda=(\theta_{1},\theta_{2},\dots \theta_{n},
\lambda_{1},\lambda_{2},\dots \lambda_{m})$ ist,\\
$\mathbf{T}(n)$ \quad \> die nat"urliche Topologie des $\mathbb{R}^{n}$,\\
$\mathbf{un}(\mathcal{A})$ \quad \> das Mengensystem der Teilmengen von $\bigcup\mathcal{A}$,\\
\quad \> die keine Menge aus $\mathcal{A}$ enthalten. \\
\end{tabbing}

\renewcommand{\indexname}{Begriffsregister}
\printindex

\end{document}